\documentclass[14pt]{elsarticle}

\usepackage{lineno}
\usepackage{hyperref}
\usepackage{amsmath,amssymb}
\usepackage{graphicx}
\usepackage{epstopdf}
\usepackage{subcaption}
\usepackage{textcomp}
\usepackage{mathrsfs}
\usepackage{bbm}
\usepackage{color,soul}
\usepackage[margin=0.9in]{geometry}
\usepackage[ruled]{algorithm2e}
\usepackage{xcolor}
\usepackage{mathtools}
\usepackage{cases} 
\usepackage{multirow} 
\newtheorem{remark}{Remark}
\usepackage[title]{appendix}
\modulolinenumbers[5]
\newcommand{\vect}[1]{\boldsymbol{#1}}
\newcommand{\matr}[1]{\mathbf{#1}}

\graphicspath{ {./graphics/} }
\journal{CMAME}

\newcommand{\cP}{\mathcal{P}}

\newcommand{\vF}{\vect{F}}
\newcommand{\vI}{\vect{I}}
\newcommand{\vJ}{\vect{J}}
\newcommand{\vk}{\vect{k}}

\newcommand{\vP}{\vect{P}}

\newcommand{\vx}{\vect{x}}

\newcommand{\vxi}{\vect{\xi}}

\begin{document}
\begin{frontmatter}
		
\title{Isogeometric analysis with local adaptivity based on a posterior error estimation for elastodynamics}
\author[mymainaddress]{Peng Yu
\corref{mycorrespondingauthor}}
\ead{glpengyu@gmail.com}
\author[secondaddress]{Cosmin Anitescu}
\author[thirdaddress]{Satyendra Tomar}
\author[thirdaddress,mymainaddress]{St\'ephane Pierre Alain Bordas}
\author[mymainaddress]{Pierre Kerfriden	\corref{mycorrespondingauthor}}	
\cortext[mycorrespondingauthor]{Corresponding author}
\ead{pierre.kerfriden@gmail.com}
\address[mymainaddress]{Institute of Mechanics and Advanced Materials, School of Engineering, Cardiff University}
\address[secondaddress]{Institute of Structural Mechanics, Bauhaus Universit\"at Weimar, Germany}
\address[thirdaddress]{Institute of Computational Engineering, University of Luxembourg, Faculty of Sciences Communication and Technology, Luxembourg}

\begin{abstract}
This paper presents a novel methodology of local adaptivity for the frequency-domain analysis of the vibrations of Reissner-Mindlin plates. The adaptive discretization is based on the recently developed Geometry Independent Field approximaTion (GIFT) framework, which may be seen as a generalisation of the Iso-Geometric Analysis (IGA). Within the GIFT framework, we describe the geometry of the structure exactly with NURBS (Non-Uniform Rational B-Splines), whilst independently employing Polynomial splines over Hierarchical T-meshes (PHT)-splines to represent the solution field. The proposed strategy of local adaptivity, wherein a posteriori error estimators are computed based on inexpensive hierarchical $h-$refinement, aims to control the discretisation error within a frequency band. The approach sweeps from lower to higher frequencies, refining the mesh appropriately so that each of the free vibration mode within the targeted frequency band is sufficiently resolved. Through several numerical examples, we show that the GIFT framework is a powerful and versatile tool to perform local adaptivity in structural dynamics. We also show that the proposed adaptive local $h-$refinement scheme allows us to achieve significantly faster convergence rates than when using a uniform $h-$refinement.
\end{abstract}
		
		
\end{frontmatter}
	
\section{Introduction}
Isogeometric analysis (IGA) was proposed in \cite{hughes2005} to integrate Computer Aided Design (CAD) and analysis in Computer Aided Engineering (CAE). Due to the high continuity order of NURBS basic functions \cite{hughes2005, Nguyen201589}, NURBS-based IGA has been successfully used to investigate many problems, and in particular problems related to plate vibrations, including Kirchoff plate \cite{Cottrell20065257, Shojaee2012} and Reissner–Mindlin plate \cite{sobota2017implicit, NME:NME4282}). The results obtained when using IGA are often more accurate than those obtained using the traditional finite element method (FEM). The previously mentioned studies of plate vibrations with IGA are mostly dedicated to homogeneous structures, whereby the vibrations occur globally so that the uniform refinement of NURBS is an adequate method to control the discretisation error. However, when the dynamic solution exhibits local features, due to e.g. sharp geometrical feature and/or varying material properties, the uniform NURBS-based refinement may become inefficient. This is because NURBS basis functions are defined by a tensor product form, which leads to globally structured grid (see Fig.1(a)), which in turns result in computational wastage when trying to capture the local features of interest.
\begin{figure}
	\begin{center}
		\includegraphics[width=0.7\columnwidth]{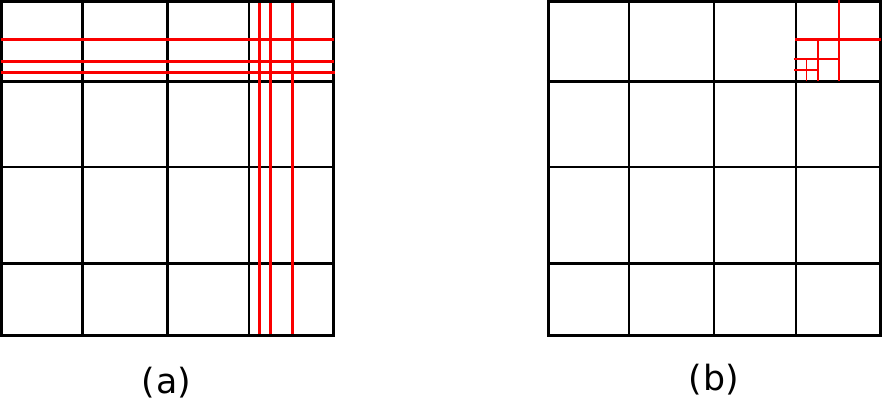}	
		\caption{(a) NURBS global refinement and (b) PHT splines local refinement}
		\label{fig:NURBSvsPHT}
	\end{center}
\end{figure}

To overcome these limitations, splines with local refinement properties such as (truncated) hierarchical B-splines \cite{GiannelliJS-12-THBSplines,giannelli2016thb}, hierarchical NURBS \cite{schillinger2012isogeometric}, locally refined (LR) B-splines \cite{johannessen2014isogeometric}, T-splines \cite{Sederberg2004,Bazilevs2010}, and polynomial/rational splines over hierarchical T-meshes (PHT/RHT)-splines \cite{Deng2008,Nguyen-Thanh2011} were developed. In this study, we choose to apply PHT-splines, as they inherit the main merits of both B-splines and T-splines: basis functions can be represented by Bézier-Bernstein polynomials over a set of Hermite finite elements, and mesh refinement is local and simple (as seen in Fig.\ref{fig:NURBSvsPHT}(b)). In the recent past, PHT-splines have been successfully used to solve static elastic solid problems. The numerical results of \cite{Nguyen-Thanh2011} showed that the adaptive PHT refinement delivers a higher convergence rate than uniform NURBS refinement. However, since PHT-splines are polynomial splines and not rational splines, they are not able to exactly represent the basic geometrical features, e.g. circles, ellipses, that typically arise in engineering design and analysis. This problem may be circumvented by making use of RHT-splines, as proposed in \cite{Nguyen-Thanh2011}. However, RHT-splines, unlike NURBS and T-splines, cannot be seamlessly extracted from existing CAD softwares. Besides, in the context of adaptivity, the updating of the weights during the refinement process require dedicated numerical developments. These difficulties prompted us to look into another direction. 

Inspired by the work proposed in \cite{marussig2015fast,toshniwal2017smooth,anitescu2017recovery,atroshchenko2017weakening}, we will employ the Geometry-Independent Field approximaTion (GIFT) to deal with the aforementioned issue, by allowing the geometry and solution fields to be described using different functional spaces. The GIFT framework was first developed within the context of the boundary element method \cite{marussig2015fast}. Later, Toshniwal et al. \cite{toshniwal2017smooth} established a scheme for unstructured quadrilateral meshes, where the space of geometric design $\mathcal{S}_{\mathbb{D}}$ and the space of solution analysis $\mathcal{S}_{\mathbb{A}}$ were different. The GIFT model, wherein NURBS basis functions are used to describe the geometry without approximation and PHT splines are utilized for analysis, will be used in this paper. This combination is compatible with  state-of-the-art CAD technology (the net of control points is inherited from CAD directly), whilst allowing local mesh refinement to take place in a non-degenerate manner. It is worth noticing that the GIFT scheme may not to satisfy the isogeometric compatibility condition \cite{toshniwal2017smooth}, which requires the solution space to be adequately rich compared to the functional space used to represent the geometry.  This is because PHT splines are polynomials while NURBS are rational functions.
However, the NURBS/PHT combo has been successfully used to develop (adaptive) GIFT schemes  and achieve optimal convergent rates in the context of linear elasticity  \cite{atroshchenko2017weakening,anitescu2017recovery}. It should also be noted that with increasing refinement level of the solution space, the space of PHT-splines will get closer to encompassing the NURBS-based functional space used to represent the geometry.

The main contributions of the paper are twofold. Firstly, we develop a novel methodology of local adaptivity based on GIFT for structural vibration problems. Secondly, we propose a novel frequency-domain adaptation strategy based on a posteriori error estimation and mode sweeping. Closely related to the proposed  adaptivity scheme is that described in \cite{nguyen2014adaptive}, whereby RHT splines are used to obtain a higher convergence rate of free vibration frequencies when compared to that observed when using tensor-product-based NURBS. However, the local refinement of the aforementioned study is driven by a priori error estimation, which does not provide any quantitative measure of accuracy (i.e. it only provides a spatial map of error sources). We aim to develop a comprehensive refinement strategy, which  includes a reliable stopping criterion as provided by a posteriori error estimation (see for instance \cite{stein2007finite,ladeveze2013new,bangerth2010adaptive,gonzalez2014mesh}). More precisely, we will define a hierarchical a posteriori error estimator that makes the best of the PHT-spline local refinement capabilities. More precisely, the accuracy of GIFT solutions will be estimated by computing refined solutions using a finer mesh generated by systematically subdividing every GIFT element. The mesh adaptivity will be performed for every free vibration mode (i.e. frequency and associated mode shape) within a frequency band, sweeping from lower to higher frequencies. The algorithm requires for coarse and fine GIFT estimations of the modes to be put in correspondence in order to be compared. This is not a trivial task. We propose a new algorithm inspired by the Modal Assurance Criterion (MAC) strategy, which is widely used in experimental structural dynamics \cite{allemang2003modal,PASTOR2012543}. We will show that the proposed algorithm is robust with respect to the order of multiplicity of the fine and coarse modes.

The organization of this paper is as follows. In Section \ref{sec:Problem} and Section \ref{sec:discrete}, we formulate the variational form of the free vibrations of Reissner-Mindlin plates, in the frequency domain. We then describe how to discretise such problems using IGA and GIFT. In Section \ref{sec:adpone}, we describe our proposed error estimation strategy, in the context of the adaptation of one, clearly isolated, free vibration mode. The strategy is then extended to the accuracy control of multiple modes in Section \ref{sec:adpfreq}, which includes a technical discussion regarding the control of discretisation errors in the context of free vibration modes of order of multiplicity larger than one. In Section \ref{sec:Numtests}, several numerical examples are presented to evaluate the efficiency of the proposed methodology, and conclusions are drawn in Section \ref{sec:conclusion}.

\section{Problem Statement} \label{sec:Problem}
\begin{figure} 
	\centering
	\def\svgwidth{0.6\columnwidth}
	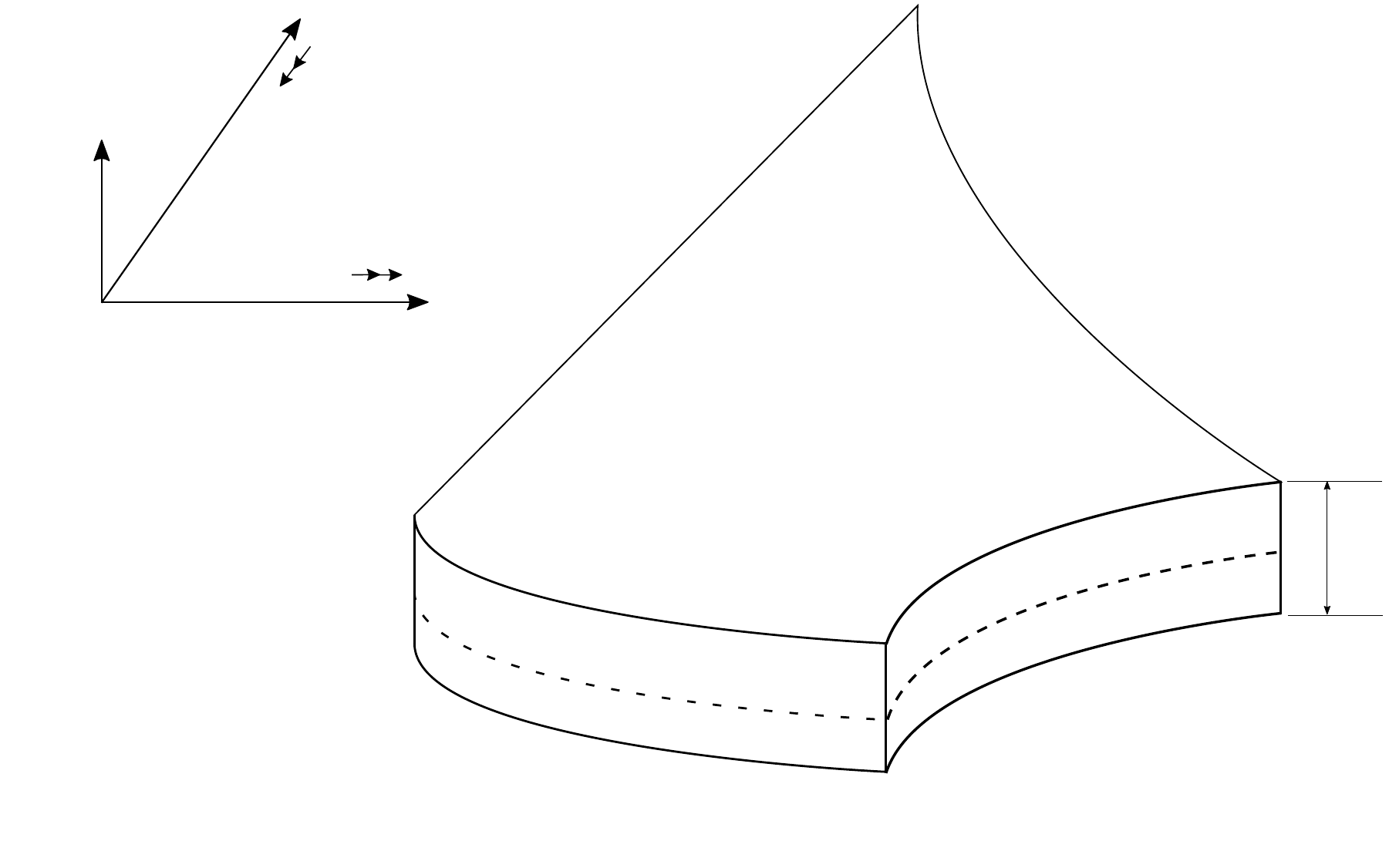
	\caption{Geometry and coordinate system of a classical Reissner-Mindlin plate.}
	\label{RMplate} 
\end{figure}

Let $\Omega\subset {\mathbb{R}}^2$ represent $x-y$ domain of the middle plane of a typical Mindlin plate, as shown in Fig.\ref{RMplate}. The boundary $\partial\Omega$ involves $\partial\Omega_{u}$, $\partial\Omega_{s}$ and $\partial\Omega_{m}$ such that: $\partial\Omega = \overline{\partial\Omega_{u}\cup\partial\Omega_{s}\cup\partial\Omega_{m}},~\partial\Omega_{u}\bigcap\partial\Omega_{s} = \emptyset, ~\partial\Omega_{u}\bigcap\partial\Omega_{m} = \emptyset$. The formal statement of governing equation can be expressed as
\begin{align}
	\begin{rcases}\label{govOmega}	
		-\dfrac{\rho h^{3}}{12}\ddot{\vect{\theta}} + \mathbb{L}^{T}\vect{M} + \vect{S} = 0\\
		-\rho\ddot{\vect{w}}h + \vect{\nabla}^{T}\vect{S} + \vect{q} = 0	
	\end{rcases}& \ \text{in} \ \Omega,\\
	\begin{rcases}
		w = \bar{w}\\
		\vect{\theta} = \vect{\bar{\theta}}&
	\end{rcases}& \ \text{on} \ \partial\Omega_{u},\\
	\vect{S} = \vect{\bar{S}} \ \ \text{on} \ \partial\Omega_{s},\\
	\vect{M} = \vect{\bar{M}} \ \ \text{on} \ \partial\Omega_{m}.
\end{align}
The $\rho$ is density and $h$ is thickness. The directions of deflection $w$ and rotation $\vect{\theta} = (\theta_{x}, \theta_{y})^{T}$ are presented in Fig.\ref{RMplate}. The $\vect{q}$ is transverse loading, and the operator $\mathbb{L}$ is defined as
\begin{align*}
\mathbb{L} =
\begin{bmatrix}
\dfrac{\partial}{\partial x} & 0\\\\
0 & \dfrac{\partial}{\partial y}\\\\
\dfrac{\partial}{\partial y} & \dfrac{\partial}{\partial x}
\end{bmatrix}.
\end{align*}
According to \cite{zienkiewicz2000finite}, the resultant shear force $\vect{S}$ and the moment $\vect{M}$ can be given as follows
\begin{align}\label{SM}
\begin{aligned}
\vect{M} = \matr{D}\mathbb{L}\vect{\theta}, \\\vect{S} = \mathcal{\kappa}Gh(\vect{\nabla}w - \vect{\theta}).
\end{aligned}
\end{align}
The elastic matrix $\matr{D}$ is defined by the assumption of plane stress such that
\begin{align*}
\matr{D} = D
\begin{bmatrix}
1 & \nu & 0\\
\nu & 1 & 0\\
0 & 0 & \dfrac{(1-\nu)}{2}
\end{bmatrix},
\end{align*}
where the parameter $D = \dfrac{Eh^{3}}{12(1 - \nu^{2})}$ denotes the bending stiffness of the plate, and $E,\nu$ express the Young's modulus, the Poisson's ratio, respectively. The shear elastic modulus is $G$, and the constant coefficient $\kappa$ set to be 5/6 in this study. Substituting Eq.\eqref{SM} into Eq.\eqref{govOmega}, we can rewrite Eq.\eqref{govOmega} as
\begin{align}\label{gov2}
	\begin{aligned}
		-\dfrac{\rho h^{3}}{12}\ddot{\vect{\theta}} + \mathbb{L}^{T}\matr{D}\mathbb{L}\vect{\theta} + \mathcal{\kappa}Gh(\vect{\nabla} w - \vect{\theta}) = 0,\\
		-\rho\ddot{w}h + \vect{\nabla}^{T}[\mathcal{\kappa}Gh(\vect{\nabla}w - \vect{\theta})] + \vect{q} = 0.
	\end{aligned}
\end{align}		
We now introduce the trial solution space $\mathscr{U}$ and test function space $\mathscr{V}$
\begin{align}
\mathscr{U} = \{u \in H^{2}(\Omega) : u = \bar{u} \mathrm{~on~} \partial\Omega_{u} \},\\
\mathscr{V} = \{v \in H^{2}(\Omega) : v = 0 \mathrm{~on~} \partial\Omega_{u} \},
\end{align}
and for all $[\delta w, \delta\vect{\theta}]^{T} \in \mathscr{V}$, $[w, \vect{\theta}]^{T} \in \mathscr{U}$, we can have the variational form of Eq.\eqref{gov2}
\begin{align}
	\begin{aligned}\label{weakform1}
	-\int_{\Omega}\delta\vect{\theta}^{T}\dfrac{\rho h^{3}}{12}\ddot{\vect{\theta}}d\Omega + \int_{\Omega}\delta\vect{\theta}^{T}\mathbb{L}^{T}\matr{D}\mathbb{L}\vect{\theta}d\Omega + \int_{\Omega}\delta\vect{\theta}^{T}\mathcal{\kappa}Gh(\vect{\nabla} w - \vect{\theta})d\Omega = 0,\\
	-\int_{\Omega}\delta w^{T}\rho\ddot{w}hd\Omega + \int_{\Omega}\delta w^{T}\vect{\nabla }^{T}[\mathcal{\kappa}Gh(\vect{\nabla}w - \vect{\theta})]d\Omega + \int_{\Omega}\delta w^{T}\vect{q}d\Omega = 0.		
	\end{aligned}
\end{align} 	
We integrate by parts two terms in Eq.\eqref{weakform1} such that
\begin{align}\label{integralbp}
	\begin{aligned}
	\int_{\Omega}\delta\vect{\theta}^{T}\mathbb{L}^{T}\matr{D}\mathbb{L}\vect{\theta}d\Omega &= -\int_{\Omega}(\mathbb{L}\delta\vect{\theta})^{T}\matr{D}\mathbb{L}\vect{\theta}d\Omega + \int_{\partial\Omega_{_{m}}}\delta\vect{\theta}^{T}\vect{\bar{M}}d\Omega,\\
	\int_{\Omega}\delta w^{T}\vect{\nabla}^{T}[\mathcal{\kappa}Gh(\vect{\nabla}w - \vect{\theta})]d\Omega &= -\int_{\Omega}(\vect{\nabla}\delta w)^{T}\mathcal{\kappa}Gh\vect{\nabla}wd\Omega+\int_{\Omega}(\vect{\nabla}\delta w)^{T}\mathcal{\kappa}Gh\vect{\theta}d\Omega+\int_{\partial\Omega_{s}}\delta w^{T}\vect{\bar{S}}d\Omega,
	\end{aligned}
\end{align}
where $ \vect{\bar{M}} $ and $ \vect{\bar{S}} $ denote the prescribed moment and shear force respectively. Since only the free vibration analysis is considered in this paper, $\vect{\bar{M}}, \vect{\bar{S}}, \vect{q}$ are all set to be zero. Substitution of Eq.\eqref{integralbp} into Eq.\eqref{weakform1} gives the weak form of the governing equation as
\begin{align}
	\begin{aligned}\label{weakformend}
		\int_{\Omega}\delta\vect{\theta}^{T}\dfrac{\rho h^{3}}{12}\ddot{\vect{\theta}}d\Omega + \int_{\Omega}(\mathbb{L}\delta\vect{\theta})^{T}\matr{D}\mathbb{L}\vect{\theta}d\Omega + \int_{\Omega}\delta\vect{\theta}^{T}\mathcal{\kappa}Gh\vect{\theta}d\Omega - \int_{\Omega}\delta\vect{\theta}^{T}\mathcal{\kappa}Gh\vect{\nabla}wd\Omega = 0, \\
		\int_{\Omega}\delta w^{T}\rho\ddot{w}hd\Omega + \int_{\Omega}(\vect{\nabla}\delta w)^{T}\mathcal{\kappa}Gh\vect{\nabla}wd\Omega - \int_{\Omega}(\vect{\nabla}\delta w)^{T}\mathcal{\kappa}Gh\vect{\theta}d\Omega = 0.
	\end{aligned}
\end{align}
The general time-dependent solution of the free vibration equation can be constructed by assuming
\begin{align}\label{uexp}
\vect{u}(x,t) = \vect{\phi}(x)\text{exp}(i\lambda t) = 
\begin{bmatrix}
w\\
\vect{\theta}
\end{bmatrix} = 
\begin{bmatrix}
\vect{\phi}_{w}\\
\vect{\phi}_{\vect{\theta}}
\end{bmatrix}
\text{exp}(i\lambda t),
\end{align}
where $\lambda$ is the frequency and $\vect{\phi}$ is the eigenvector. Substituting Eq.\eqref{uexp} in Eq.\eqref{weakformend}, the weak form becomes an eigenvalue problem as follows
\begin{align}\label{weakformfredomain}
	\begin{aligned}
	-\lambda^{2}\int_{\Omega}\delta\vect{\phi}_{\vect{\theta}}^{T}\dfrac{\rho h^{3}}{12}\vect{\phi}_{\vect{\theta}}d\Omega + \int_{\Omega}(\mathbb{L}\delta\vect{\phi}_{\vect{\theta}})^{T}\matr{D}\mathbb{L}\delta\vect{\phi}_{\vect{\theta}}d\Omega + \int_{\Omega}\delta\vect{\phi}_{\vect{\theta}}^{T}\mathcal{\kappa}Gh\vect{\phi}_{\vect{\theta}}d\Omega - \int_{\Omega}\delta\vect{\phi}_{\vect{\theta}}^{T}\mathcal{\kappa}Gh\vect{\nabla}\vect{\phi}_{w}d\Omega = 0, \\
	-\lambda^{2}\int_{\Omega}\delta\vect{\phi}_{w}^{T}\rho\vect{\phi}_{w}hd\Omega + \int_{\Omega}(\vect{\nabla}\delta \vect{\phi}_{w})^{T}\mathcal{\kappa}Gh\vect{\nabla}\vect{\phi}_{w}d\Omega - \int_{\Omega}(\vect{\nabla}\delta \vect{\phi}_{w})^{T}\mathcal{\kappa}Gh\vect{\phi}_{\vect{\theta}}d\Omega = 0.
	\end{aligned}
\end{align}

\section{Discretization of the eigenvalue problem using IGA and GIFT} \label{sec:discrete}
Let $\cP$ be the parametric domain, and the physical domain $\Omega$ is parametrized on $\cP$ by a geometrical mapping $\vF$
\begin{align}
\vF : \cP \rightarrow \Omega, \quad \vx = \vF(\vxi)
\label{eq:domain_param_1}.
\end{align}
We assume that the domain $\Omega$ consists of sub-domains $\Omega_{i}$, such that $ \Omega = \bigcup_{i=1}^{N}\Omega_{i}$. The geometrical map $\vF$ is given by a set of basis functions $N_{\vk} (\vxi)$ and a set of control points $\vP_{\vk}$ as
\begin{equation}
\vx (\vxi)= \sum \limits_{\vk \in \vI}N_{\vk} (\vxi)\vP_{\vk},
\label{eq:domain_param_2}
\end{equation}
where $N_{\vk} (\vxi)$ is a bivariate tensor product of the univariate basis functions $\vect{N}_{i}^{p}(\xi),\vect{N}_{j}^{q}(\eta)$ with the orders $p,q$, such that $N_{\vk} (\vxi) = \vect{N}_{i}^{p}(\xi)\vect{N}_{j}^{q}(\eta), (\xi,\eta) \in \vect{\xi}$. Moreover, $\vI$ is introduced as 2-dimensional multi-index $(i, j)$, and $\vk$ is interchangeably regarded as the collapsed notation for $\vI$. In what follows, we will refer to the set $\{N_{\vk}(\vxi)\}_{\vk \in \vI}$ as the \emph{geometry basis}, and introduce the discretization of the eigenvalue problem via the scheme of IGA and GIFT respectively.
\subsection{The framework of IGA}
In IGA, the solution field $\vect{\phi} = [\vect{\phi}_{w},\vect{\phi}_{\vect{\theta}}]^{T} $ is represented through the same spline functions which are used for the geometry, i.e.,
\begin{align}
\vect{\phi} = \sum \limits_{\vk \in \vI} N_{\vk}(\vxi) \bar{\vect{\phi}}_{\vk},
\end{align}
where $ \bar{\vect{\phi}}_{\vk}$ are unknown control variables. Then the deflection and rotations which serve as components of the solution variable can be denoted in the matrix form
\begin{align}\label{discretwtheta}
	\begin{bmatrix}
	\vect{\phi}_{w}\\
	\vect{\phi}_{\vect{\theta}}
	\end{bmatrix} = 
	\begin{bmatrix}
	\vect{N}_{w} & 0 \\
	0 & \vect{N}_{\vect{\theta}}
	\end{bmatrix}
	\begin{bmatrix}
	\bar{\vect{\phi}}_{w}\\
	\bar{\vect{\phi}}_{\vect{\theta}}
	\end{bmatrix} = 
	\begin{bmatrix}
	 \vect{N}_{w} & 0 & 0\\
	 0 & \vect{N}_{\theta_{x}} & 0\\
	 0 & 0 & \vect{N}_{\theta_{y}}
	 \end{bmatrix}
	 \begin{bmatrix}
	 \bar{\vect{\phi}}_{w}\\
	 \bar{\vect{\phi}}_{\theta_{x}}\\
	 \bar{\vect{\phi}}_{\theta_{y}}
	 \end{bmatrix}.
\end{align}
Supposed that the test functions $\delta\vect{\phi}_{w}$ and $\delta\vect{\phi}_{\vect{\theta}}$ are discretized using Eq.\eqref{discretwtheta}, the discrete form of Eq.\eqref{weakformfredomain} will become
\begin{align}
	\begin{aligned}\label{dicreteform1}
	\delta\bar{\vect{\phi}}_{\vect{\theta}}^{T}\left(-\lambda^{2}\int_{\Omega}\vect{N}_{\vect{\theta}}^{T}\dfrac{\rho h^{3}}{12}\vect{N}_{\vect{\theta}}d\Omega\right)\bar{\vect{\phi}}_{\vect{\theta}} + \delta\bar{\vect{\phi}}_{\vect{\theta}}^{T}\left(\int_{\Omega}(\mathbb{L}\vect{N}_{\vect{\theta}})^{T}\matr{D}\mathbb{L}\vect{N}_{\vect{\theta}}d\Omega + \int_{\Omega}\vect{N}_{\vect{\theta}}^{T}\mathcal{\kappa}Gh\vect{N}_{\vect{\theta}}d\Omega\right)\bar{\vect{\phi}}_{\vect{\theta}}\\ -\delta\bar{\vect{\phi}}_{\vect{\theta}}^{T}\left( \int_{\Omega}\vect{N}_{\vect{\theta}}^{T}\mathcal{\kappa}Gh\vect{\nabla}\vect{N}_{w}d\Omega\right)\bar{\vect{\phi}}_{w} = 0, \\
	\delta{\bar{\vect{\phi}}_{w}}^{T}\left( -\lambda^{2}\int_{\Omega}\vect{N}_{w}^{T}\rho h\vect{N}_{w}d\Omega\right)\bar{\vect{\phi}}_{w} - \delta\bar{\vect{\phi}}_{w}^{T}\left( \int_{\Omega}(\vect{\nabla}\vect{N}_{w})^{T}\mathcal{\kappa}Gh\vect{N}_{\vect{\theta}}d\Omega\right)\bar{\vect{\phi}}_{\vect{\theta}}\\  + \delta\bar{\vect{\phi}}_{w}^{T}\left( \int_{\Omega}(\vect{\nabla}\vect{N}_{w})^{T}\mathcal{\kappa}Gh\vect{\nabla}\vect{N}_{w}d\Omega\right)\bar{\vect{\phi}}_{w}  = 0.
	\end{aligned}	
\end{align}
The Jacobian matrix $J(\vxi)$ of the mapping $\vF$ is introduced as
\begin{equation}
\matr{J}(\vxi) = \dfrac{\partial\vect{x}}{\partial{\vect{\xi}}}
= \sum \limits_{\vk \in \vI} \vP_{\vk}
\dfrac{\partial N_{\vk}(\vxi)}{\partial\vect{\xi}}.
\label{jacobian}
\end{equation}
We can rewrite Eq.\eqref{dicreteform1} by defining the stiffness $\matr{K}$ and mass matrix $\matr{M}$ integrated in the parametric space $\mathcal{P}$ as follows:
\begin{align}\label{assembleIGA}
	\begin{aligned}
		&\matr{K} = \matr{K}_{b} + \matr{K}_{s},\\
		&\matr{K}_{b} = \int_{\mathcal{P}}(\mathbb{B}_{b}\vect{N})^{T}\matr{D}\mathbb{B}_{b}\vect{N}\left| \matr{J}(\vect{\xi})\right| d\mathcal{P}, \hspace{0.5cm}\text{bending stiffness}\\
		&\matr{K}_{s} = \int_{\mathcal{P}}(\mathbb{B}_{s}\vect{N})^{T}\matr{D}\mathbb{B}_{s}\vect{N}\left| \matr{J}(\vect{\xi})\right| d\mathcal{P}, \hspace{0.5cm}\text{shear stiffness}\\
		&\matr{M} = \int_{\mathcal{P}}\rho\vect{N}^{T}\matr{m}\vect{N}\left| \matr{J}(\vect{\xi})\right| d\mathcal{P},
	\end{aligned}
\end{align}
with 
\begin{align}\label{eq:BBm}
	\begin{aligned}
		\mathbb{B}_{b} = 
		\begin{bmatrix}
		0 & \dfrac{\partial}{\partial x} & 0\\\\
		0 & 0 & \dfrac{\partial}{\partial y}\\\\
		0 & \dfrac{\partial}{\partial y} & \dfrac{\partial}{\partial x}
		\end{bmatrix},
		\mathbb{B}_{s} =
		\begin{bmatrix}
		\dfrac{\partial}{\partial x} & -1 & 0\\\\
		\dfrac{\partial}{\partial y} & 0 & -1		
		\end{bmatrix},
		\matr{m} = 
		\begin{bmatrix}
		h & 0 & 0\\\\
		0 & \dfrac{h^{3}}{12} & 0\\\\
		0 & 0 & \dfrac{h^{3}}{12}
		\end{bmatrix}.
	\end{aligned}
\end{align}
Then Eq.\eqref{dicreteform1} can be compactly written into the final matrix form of eigenvalue problem by IGA scheme
\begin{align}\label{matrixformendIGA}
\delta\bar{\vect{\phi}}^{T}( \matr{K} - \lambda^{2}\matr{M})\bar{\vect{\phi}} = 0.
\end{align}

\subsection{The framework of GIFT}
A detailed exposition of GIFT is presented in \cite{atroshchenko2017weakening}. GIFT allows to choose a \emph{solution basis} $\{\varPsi_{\vk}(\vxi)\}_{\vk \in \vJ}$ that can be different from the geometry basis $\{N_{\vk}(\vxi)\}_{\vk \in \vI}$, but is defined on the physical domain with the help of the same mapping $\vect{F}^{-1}$ as in Eq.\eqref{eq:domain_param_1}, i.e.,
\begin{align}
	\varPsi_{\vk}(\vxi) = \varPsi_{\vk} \circ \vect{F}^{-1}(\vect{x}).
\end{align}
Hence, the solution variables $\vect{\phi} = [\vect{\phi}_{w},\vect{\phi}_{\vect{\theta}}]^{T} $ are described accordingly by
\begin{align}
	\vect{\phi} = \sum \limits_{\vk \in \vJ} \varPsi_{\vk}(\vxi) \bar{\vect{\phi}}_{\vk}, \label{discreteGIFT}	
\end{align}
Following the steps of derivation to obtain the statement of matrix form Eq.\eqref{matrixformendIGA} from the weak form Eq.\eqref{weakformfredomain} using IGA method, the discrete eigenvalue equation can be similarly acquired through GIFT as in Eq.\eqref{matrixformendIGA} with notations
\begin{align}\label{eq:assembleGIFT}
	\begin{aligned}
	&\matr{K} = \matr{K}_{b} + \matr{K}_{s},\\
	&\matr{K}_{b} = \int_{\mathcal{P}}(\mathbb{B}_{b}\vect{\varPsi})^{T}\matr{D}\mathbb{B}_{b}\vect{\varPsi}\left| \matr{J}(\vect{\xi})\right| d\mathcal{P},\\
	&\matr{K}_{s} = \int_{\mathcal{P}}(\mathbb{B}_{s}\vect{\varPsi})^{T}\matr{D}\mathbb{B}_{s}\vect{\varPsi}\left| \matr{J}(\vect{\xi})\right| d\mathcal{P},\\
	&\matr{M} = \int_{\mathcal{P}}\rho\vect{\varPsi}^{T}\matr{m}\vect{\varPsi}\left| \matr{J}(\vect{\xi})\right| d\mathcal{P}.
	\end{aligned}
\end{align}
\begin{remark}\label{rmk:fixF}
Despite the use of a different spline basis for assembling within GIFT and IGA (compare Eq.\eqref{assembleIGA} with Eq.\eqref{eq:assembleGIFT}), the mapping $\vect{F}$ is kept the same for both methods. Therefore, supposed that the geometry is exactly represented by spline basis $N_{\vk}(\vxi)$ on the initial (often very coarse) mesh with geometric mapping $\vect{F}_{0}$, the integrals on physical domain, e.g., Eq.\eqref{eq:assembleGIFT} can be computed by assuming that the mapping is fixed ($\vect{F} = \vect{F}_{0}$). It means the generation of new controls points and new geometric basis are not essential any more during the refinement, which is computationally advantageous. It is worth noting that this hypothesis of fixed mapping is satisfied, due to the fact that both $h-$refinement and $p-$refinement in IGA are with the geometry preserving.
\end{remark}

In this work, three kinds of spline basis functions, NURBS, PHT and RHT are applied for IGA and GIFT methods. They are introduced in detail in Appendix \ref{NURBS app}, \ref{PHT app} and \ref{app:RHT}, respectively. 
\subsection{Boundary conditions and multiple patch coupling} \label{secBC}
\begin{figure}
	\centering
	\includegraphics[width=0.5\columnwidth]{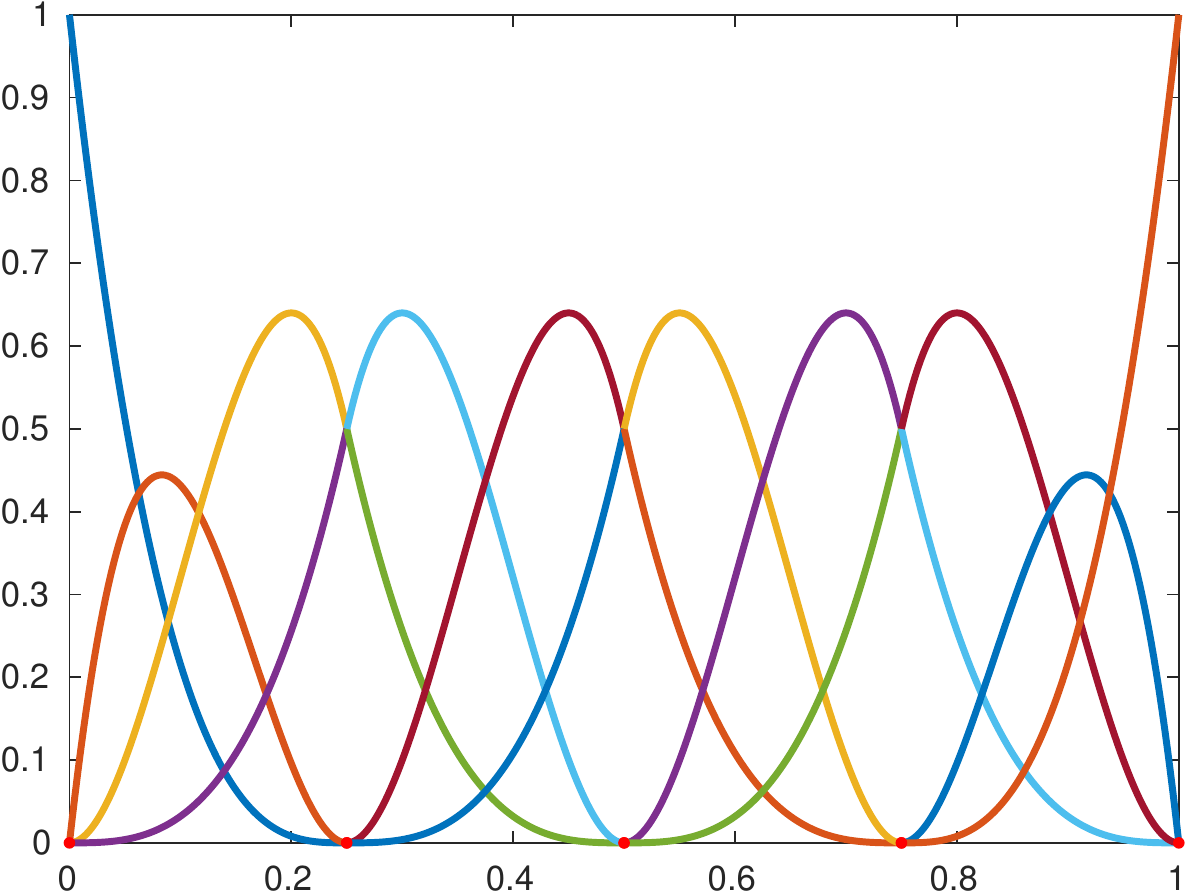}
	\caption{1D cubic B-spline basis functions defined in knot vector $\Xi = [0,0.25,0.5,0.75,1]$.}
	\label{fig:1DBspline}
\end{figure}

Two types of Dirichlet boundary conditions for free vibration are applied in this study
\begin{align}\label{BC}
	\begin{aligned}
	\bar{w} = 0, \ \bar{\theta}_{s} = 0, \ \bar{\theta}_{n} = 0 \hspace{0.5cm} &\text{clampled},\\
	\bar{w} = 0 \hspace{0.5cm} &\text{simply supported},
	\end{aligned}
\end{align}
where the subscripts $n$ and $s$ denote tangent and normal direction, as shown in Fig.\ref{RMplate}. Spline basis functions are generally not with interpolatory nature. This does not allow the imposition of Dirichlet boundary condition as straightforward as that in finite element method. Some strategies proposed for the mesh free method, e.g., \cite{liu2006meshfree, fernandez2004imposing} can be extended to the isogeometric framework, but we desire to seek a more simple method. As presented in Fig.\ref{fig:1DBspline}, at boundary points of the knot vector $\xi=0$ and $\xi=1$, only one B-spline basis function is equal to 1, while others are zero. So it is just required to find the control variables related to the non-zeros basis functions, and remove the relevant degrees of freedom. Thus, the boundary conditions in Eq.\eqref{BC} can be imposed. As the NURBS, PHT-splines and RHT-splines are all based on the B-splines, and refinements do not change the situations on the boundaries mentioned above, imposing the boundary conditions will be direct.
Regarding the coupling of multiple patches, we utilize a convenient and robust approach that patches are conforming at interfaces through the imposition which allows the $C^{0}$ continuity by identifying the corresponding degrees of freedom. Owing to the property of local refinement possessed by PHT, it is simple to realize that this process is with minimal additional refinements and computing resources. Besides, other weak coupling methods \cite{nguyen2017isogeometric,nguyen2014nitsche} could also be used without additional difficulties.

\section{Adaptivity for one mode} \label{sec:adpone}
As mentioned in Remark \ref{rmk:fixF}, in the GIFT scheme, the refinement is only required for solution space. Therefore, in this section, we are going to present the procedure of PHT mesh adaptivity when a particular mode is targeted.
\subsection{Error estimator}\label{sec:Errestm}
Supposed that $i$ is the mode of interest on \emph{current (coarse) mesh} $\mathbb{T}$ ($\mathbb{T}$ is a hierarchical T-mesh), we define a corresponding mode $\tilde{i}$ on \emph{refined mesh} $\tilde{\mathbb{T}}$, wherein the elements are created by dividing each element in $\mathbb{T}$ into $2^{d \cdot L_{e}}$ elements, where $d$ is the dimension of the problem, and $L_{e}$ is the level of refinement. Assuming that $\vect{\phi}_{i}^{h}$ and $\lambda_{i}^{h}$ denote the eigenvector and frequency obtained using $\mathbb{T}$ for mode $\tilde{i}$, and $\tilde{\vect{\phi}}_{\tilde{i}}$ and $\tilde{\lambda}_{\tilde{i}}$ indicate solutions acquired on $\tilde{\mathbb{T}}$ for mode $\tilde{i}$, then the error estimators for frequency and mode shape can be defined as
\begin{align}\label{eq:errestmor}
\left| e_{i}^{\lambda} \right| = \left| \text{log}\tilde{\lambda}_{\tilde{i}} - \text{log}\lambda_{i}\right|, \hspace{0.2cm} \delta_{i}^{\vect{\phi}} = \dfrac{\left\| e_{i}^{\vect{\phi}}\right\| _{E}}{\left\|\tilde{\vect{\phi}}_{\tilde{i}}\right\|_{E}}= \dfrac{\left\| \tilde{\vect{\phi}}_{\tilde{i}} - \vect{\phi}_{i}^{h}\right\|_{E} }{\left\|\tilde{\vect{\phi}}_{\tilde{i}}\right\|_{E}},
\end{align}	
where 100\textcelsius $\left\| \cdot \right\|_{E} := \left[ \int_{\Omega}\mathbb{B}_{b}^{T}(\cdot)\matr{D}\mathbb{B}_{b}(\cdot)d\Omega + \int_{\Omega}\mathbb{B}_{s}^{T}(\cdot)\matr{D}\mathbb{B}_{s}(\cdot)d\Omega\right]^{\frac{1}{2}}$ is the energy norm. Since $\vect{\phi}_{i}^{h}$ and $\tilde{\vect{\phi}}_{\tilde{i}}$ are discretized by
\begin{align*}
\vect{\phi}_{i}^{h} = \vect{T}\bar{\vect{\mathcal{\phi}}}_{i}^{h}, ~\tilde{\vect{\phi}}_{\tilde{i}} = \tilde{\vect{T}}\tilde{\bar{\vect{\mathcal{\phi}}}}_{\tilde{i}},
\end{align*}
where $\vect{T}$ and $\tilde{\vect{T}}$ are PHT-spline basis functions defined over spaces $\mathbb{T}$ and $\tilde{\mathbb{T}}$. In order to compute $\left\| \tilde{\vect{\phi}}_{\tilde{i}} - \vect{\phi}_{i}^{h}\right\|_{E}$, the control variables $ \bar{\vect{\mathcal{\phi}}}_{i}^{h} $ should be prolongated onto the $\tilde{\mathbb{T}}$
\begin{align}
\bar{\vect{\mathcal{\phi}}}_{i}^{h} \longrightarrow \mathbb{P}\bar{\vect{\mathcal{\phi}}}_{i}^{h}
\end{align}
where $\mathbb{P}$ is the prolongation operator. Two strategies are introduced to compute this prolongation in Appendix.\ref{app:prolong}. One is based on the insertion of control points in PHT refinement, and the other one is based upon the projection. Owing to the features of the isogeometric system that refinement does not change the field, the solution $\vect{\mathcal{\phi}}_{i}^{h}$ is preserved exactly after the prolongation, which reads 
\begin{align}
\vect{\mathcal{\phi}}_{i}^{h} = \vect{T}\bar{\vect{\mathcal{\phi}}}_{i}^{h} = \tilde{\vect{T}}\mathbb{P}\bar{\vect{\mathcal{\phi}}}_{i}^{h}.
\end{align}
Hence, it yields that
\begin{align}
\left\| \tilde{\vect{\phi}}_{\tilde{i}} - \vect{\phi}_{i}^{h}\right\|_{E} = \left[(\tilde{\bar{\vect{\mathcal{\phi}}}}_{\tilde{i}}-\mathbb{P}\bar{\vect{\mathcal{\phi}}}_{i}^{h})^{T}\tilde{\matr{K}}(\tilde{\bar{\vect{\mathcal{\phi}}}}_{\tilde{i}}-\mathbb{P}\bar{\vect{\mathcal{\phi}}}_{i}^{h})\right]^{\frac{1}{2}},
\end{align}
where the stiffness matrix $\tilde{\matr{K}}$ is obtained by the GIFT method
\begin{align}
\tilde{\matr{K}} = \int_{\mathcal{P}}\left[ (\mathbb{B}_{b}\tilde{\vect{T}})^{T}\matr{D}\mathbb{B}_{b}\tilde{\vect{T}} + (\mathbb{B}_{s}\tilde{\vect{T}})^{T}\matr{D}\mathbb{B}_{s}\tilde{\vect{T}} \right]\left|\matr{J}(\vect{\xi})\right|d\mathcal{P}.
\end{align} 
For the reason that the adaptive mesh requires a local criterion, by referring to $\Omega_{e}$ as an element-wise physical domain, the local error estimator of eigenvector is posed, i.e.
\begin{align}\label{eq:e_Kloc}
\left\| e_{i}^{\vect{\phi}}(\Omega_{e})\right\|_{E} = \left[ \int_{\Omega_{e}}(\mathbb{B}_{b}e_{i}^{\vect{\phi}})^{T}\matr{D}\mathbb{B}_{b}e_{i}^{\vect{\phi}}d\Omega_{e} + \int_{\Omega_{e}}(\mathbb{B}_{s}e_{i}^{\vect{\phi}})^{T}\matr{D}\mathbb{B}_{s}e_{i}^{\vect{\phi}}d\Omega_{e}\right]^{\frac{1}{2}}.
\end{align}
Letting $N$ denote the total number of elements, it yields
\begin{align}\label{eq:e_K}
\left\| e_{i}^{\vect{\phi}}\right\|_{E}^{2} = \sum \limits_{j = 1}^{N} \left\| e_{i}^{\vect{\phi}}(\Omega_{e}^{j})\right\|_{E}^{2}.
\end{align}
\subsection{Marking method}\label{sec:Markstrategy} 
In order to take the error contribution by each cell into consideration, the marking strategy is proposed based on the approach \cite{dorfler1996convergent}. To be specific, we sort the values of $\left\| e_{i}^{\vect{\phi}}(\Omega_{e}^{j})\right\|_{E}^{2} (j=1, 2, \ldots, N)$ from large to small. Then mark the set of $N^{*}$ elements to be refined, if the following criterion is satisfied
\begin{align}\label{eq:criterion}
\sum \limits_{j = 1}^{N^{*}} \left\| e_{i}^{\vect{\phi}}(\Omega_{e}^{j})\right\|_{E}^{2} \geqslant \tau\left\| e_{i}^{\vect{\phi}}\right\|_{E}^{2},
\end{align}
where $\tau \in (0,1]$ is the percentage, $N^{*}$ is the minimum number of elements to satisfy Eq.\eqref{eq:criterion}. Each marked element will be subdivided into 4 elements in case of $d =2, L_{e}=1$. The interpretation for the adaptive PHT refinement process is presented in Fig.\ref{fig:adppros}, and more details on PHT refinement can be found in Appendix \ref{app:PHTrefine}. The adaptivity for mode $i$ will proceed until both of the following criteria are fulfilled
\begin{align}\label{eq:tol}
\left| e_{i}^{\lambda} \right| \leqslant \tau_{\lambda}, ~\delta_{i}^{\vect{\phi}} \leqslant \tau_{\vect{\phi}}, 
\end{align}
where $\tau_{\lambda}$ and $\tau_{\vect{\phi}}$ are error tolerances for the frequency and the eigenvector respectively.
\begin{remark}
In terms of the option for $\tau$ by a given accuracy, there is always a compromise between refinement steps and number of elements to be refined at each step. When $\tau\ll1$, it may achieve an optimal final mesh, however, it would sacrifice the computational cost due to too many refinement steps. Whilst if $\tau$ is too large, the effect of adaptivity would be diminished, since $\tau=1$ leads to the uniform refinement. We have found through some experimentations that, $\tau = 0.3$ offers a good balance in the context of this work.
\end{remark}
\begin{figure} 
	\centering
	\def\svgwidth{0.8\columnwidth}
	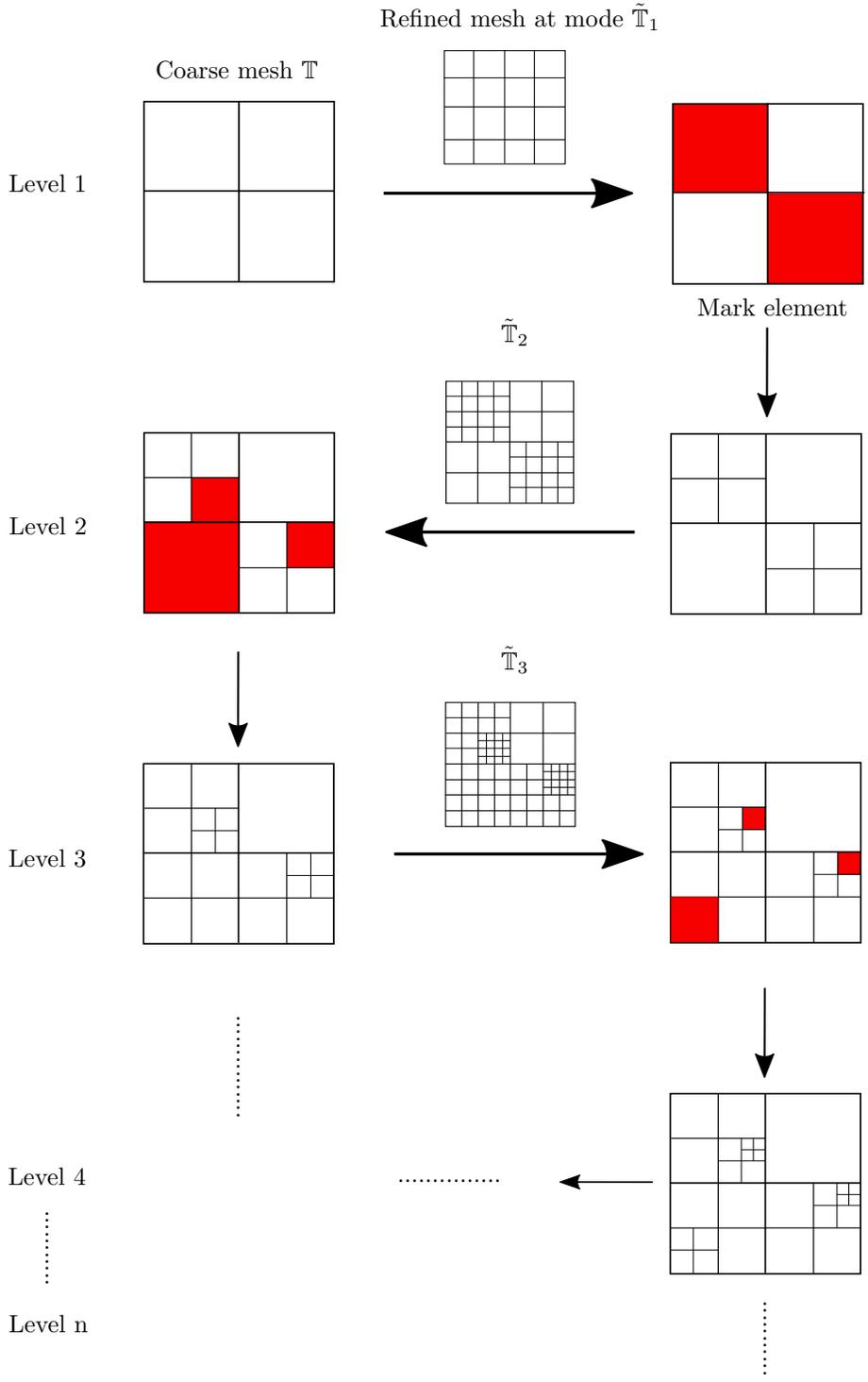
	\caption{The schematic illustration for adaptive PHT refinement procedure in parametric domain at mode $i$ in case of $d=2, L_{e}=1$.}
	\label{fig:adppros} 
\end{figure}

\begin{algorithm}
	\caption{Adaptivity process for the single mode $i$}
	\label{Algm:siglemode}
	\hspace*{-0.5cm}\textbf{Input:} $\left| e_{i}^{\lambda} \right|$ and $\delta_{i}^{\vect{\phi}}$ on coarse mesh $\mathbb{T}$.\\
	\hspace*{-0.5cm}\textbf{Output:} Updated $\mathbb{T}$ after refinement.\\
	\While{$\left| e_{i}^{\lambda} \right| \leqslant \tau_{\lambda}, ~\delta_{i}^{\vect{\phi}} \leqslant \tau_{\vect{\phi}}$}{
		\For{$j\leftarrow 1$ \KwTo $N$}{
			Compute $\left\| e_{i}^{\vect{\phi}}(\Omega_{e}^{j})\right\|_{E}^{2}$ by Eq.\eqref{eq:e_Kloc}.
		}
		
		Sort values of $\left\| e_{i}^{\vect{\phi}}(\Omega_{e}^{j})\right\|_{E}^{2}$ from large to small.	
				
		\For{$j\leftarrow 1$ \KwTo $N$}{
			\If{$\sum \limits_{j^{*} = 1}^{j} \left\| e_{i}^{\vect{\phi}}(\Omega_{e}^{j*})\right\|_{E}^{2} \geqslant \tau\left\| e_{i}^{\vect{\phi}}\right\|_{E}^{2}$}{
				Mark $N^{*} = j$\\ 
				\textbf{break}
				}		
		}
		
		\For{$j\leftarrow 1$ \KwTo $N^{*}$}{
			Refine element $j$ to update $\mathbb{T}$.
		}
		Renew $\left| e_{i}^{\lambda} \right|$ and $\delta_{i}^{\vect{\phi}}$ by Eq.\eqref{eq:errestmor}.
	}
\end{algorithm}

\section{Adaptivity for a range of frequencies} \label{sec:adpfreq}
Following the section above, we will discuss how to deal with the adaptivity when modes are inside a band of frequencies of interest.

\subsection{Algorithm of adaptivity by sweeping modes}
\begin{figure} 
	\centering
	\def\svgwidth{1\columnwidth}
	\hspace*{2cm}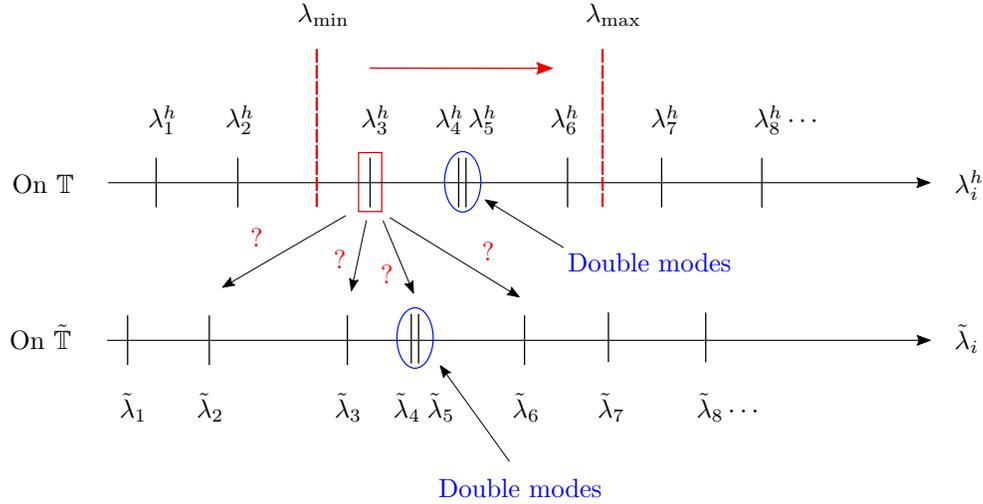
	\caption{The schematic of adaptivity algorithm for an interval of frequencies of interest by sweeping modes from low to high.}
	\label{fig:sweepmodes} 
\end{figure}

As it is shown in Fig.\ref{fig:sweepmodes}, suppose that frequencies of interest are inside in a band, that is, $\lambda_{i}^{h} \in [\lambda_{min}, \lambda_{max}]$ (marked with red dash line), and four modes are involved, with frequencies are  $\lambda_{3}^{h}, \lambda_{4}^{h}, \lambda_{5}^{h}$ and $\lambda_{6}^{h}$. We start with the lowest mode (mode 3). Since mode 3 is a single mode by checking the multiplicity, calling the Algorithm \ref{Algm:siglemode} directly, the adaptivity for mode 3 will be proceeded. Afterwards, we move to the mode 4. Through multiplicity identification, mode 4 and mode 5 are double modes. It is worth noting that there are no actual double (or multiple) modes. The so-called double (or multiple) modes in our study are numerical. For example, if $|\lambda_{4}^{h} - \lambda_{5}^{h}| \leqslant \tau_{\lambda}^{\text{mul}}$, where $\tau_{\lambda}^{\text{mul}}$ is a threshold, the mode 4 and mode 5 are considered as a double mode. Exploiting the Algorithm \ref{Algm:mulmode}, the adaptivity of double mode $(4,5)$ can be realized. Thus, by sweeping the modes until mode 6, the adaptivity for the frequencies in the window can be delivered. This algorithm is summarized in Algorithm \ref{Algm:adpFoI}.

\begin{algorithm}
	\caption{Adaptivity for a range of frequencies of interest $[\lambda_{\text{min}}, \lambda_{\text{max}}]$}
	\label{Algm:adpFoI}
	\While{$\lambda_{\text{min}} \leqslant \lambda_{i} \leqslant \lambda_{\text{max}}$ }{
		\textbf{Step 1.} \vspace*{0.1cm}Call Algorithm \ref{Algm:FEC} or Algorithm \ref{Algm:MAC} to find the multiplicity $n$ of mode $i$ on $\mathbb{T}$, and the related modes on $\tilde{\mathbb{T}}$.\\
		\textbf{Step 2.}\\
		\uIf{$n=1$}{
			Call Algorithm \ref{Algm:siglemode} to perform adaptivity for mode $i$.
		}
		\Else{
			Call Algorithm \ref{Algm:mulmode} to process adaptivity for $n$ multiple modes $\{i, \ldots, i+n-1\}$.
		}
		\textbf{Step 3.} 
		$\begin{cases}
		i = i+1, ~\text{for single mode}.\\
		i = i+n, ~\text{for $n$ multiple modes}.	
		\end{cases}$ 
	}	
\end{algorithm}

\subsection{Location of modal correspondence} \label{sec:Locmode}
As we focus on the error estimation and adaptivity in Section \ref{sec:adpone}, mode $i$ and mode $\tilde{i}$ are assumed to be related. But in fact, for example, as illustrated in Fig.\ref{fig:sweepmodes}, $\lambda_{3}^{h}$ could be related to any mode on $\tilde{\mathbb{T}}$, such as $\tilde{\lambda}_{2}$,  $\tilde{\lambda}_{3}$, $\tilde{\lambda}_{4}$, $\tilde{\lambda}_{5}$ or $\tilde{\lambda}_{6}$. Therefore, it is required to find the method to recognize this modal resemblance. Two approaches are introduced in the following sections.

\subsubsection{Frequency Error Criterion (FEC)}
The FEC strategy is to regard the modes, with the closest frequencies, as the related modes. The algorithm is summarized in Algorithm \ref{Algm:FEC}. Specifically, when using FEC strategy, first of all, we check the multiplicity of mode $i$ with $\lambda_{i}^{h} \in [\lambda_{\text{min}}, ~\lambda_{\text{max}}] $. If \vspace*{0.1cm} $|\lambda_{i+1}^{h} - \lambda_{i}^{h}| > \tau_{\lambda}^{\text{mul}}$, then mode $i$ is considered as a single mode. While, if the modes satisfy the following conditions \vspace*{0.2cm}\\
\fbox{
	\parbox{\textwidth}{
		\begin{align*}
		|\lambda_{i+1} - \lambda_{i}| \leqslant \tau_{\lambda}^{\text{mul}}, |\lambda_{i+2} - \lambda_{i+1}| \leqslant \tau_{\lambda}^{\text{mul}}, \cdots,|\lambda_{i+n-1} - \lambda_{i+n-2}| \leqslant \tau_{\lambda}^{\text{mul}}, 
		\end{align*}
		the multiplicity of mode $i$ is $n$.
	}
}\vspace*{0.2cm}
Then, we have the set of multiple modes $\{i,i+1,\ldots,i+n-1\}$.

Now, we need to find the corresponding modes on $\tilde{\mathbb{T}}$. If mode $i$ is a single mode, we compute the set of absolute values of errors
\begin{align}\label{eq:errFEC}
\{ | e_{i, \tilde{i}} | \} = \{| \lambda_{i} - \tilde{\lambda}_{\tilde{i}} |\}, ~\forall\tilde{\lambda}_{\tilde{i}} \in [\lambda_{\text{min}} - a, ~\lambda_{\text{max}} + a], 
\end{align}
where the constant $a$ is to make sure the interval $[\lambda_{\text{min}} - a, ~\lambda_{\text{max}} + a]$ \vspace*{0.1cm} is wide enough to include the corresponding mode inside. Select the minimum of $\{ | e_{i, \tilde{i}} | \}$, and then the relevant $\tilde{i}$ is the corresponding mode number. 

\begin{figure} 
	\centering
	\def\svgwidth{1\columnwidth}
	\hspace*{2.5cm}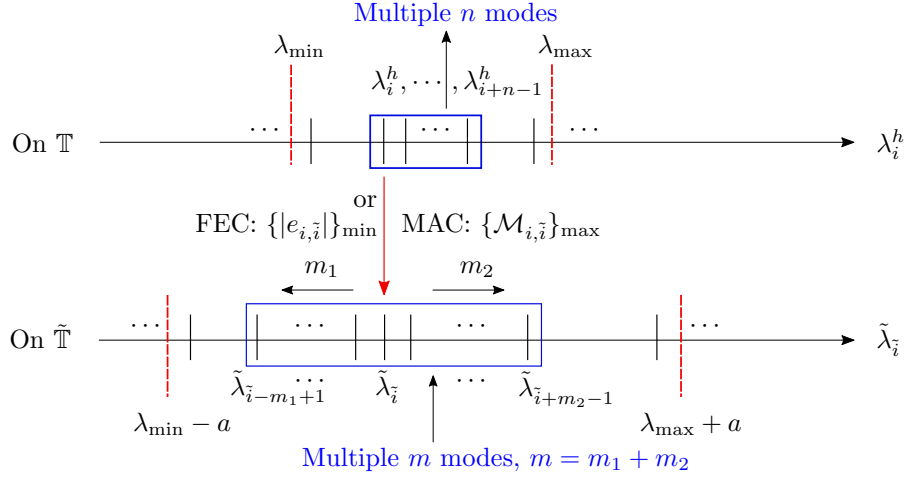
	\caption{The interpretation of the modal resemblance location for multiple modes by FEC and MAC scheme.}
	\label{fig:mulmodesLoc} 
\end{figure}

When modes $\{i, i+1,\ldots, i+n-1\}$ are $n$ multiple modes, the approach is presented in the Fig.\ref{fig:mulmodesLoc}. \vspace*{0.1cm} We firstly still find the mode $\tilde{i}$ by obtaining the minimum of $\{ | e_{i, \tilde{i}} | \}$ by Eq.\eqref{eq:errFEC}. Next, we check the multiplicity of mode $\tilde{i}$. It is required to check the modes lower than mode $\tilde{i}$, as well the modes higher than mode $\tilde{i}$ (look at Fig.\ref{fig:mulmodesLoc}). By this way, we can have that\vspace*{0.2cm}\\
\fbox{
	\parbox{\textwidth}{If
		\begin{align*}
		|\tilde{\lambda}_{\tilde{i}-m_{1}+2} - \tilde{\lambda}_{\tilde{i}-m_{1}+1}| \leqslant \tau_{\lambda}^{\text{mul}},\cdots,|\tilde{\lambda}_{\tilde{i}} - \tilde{\lambda}_{\tilde{i}-1}| \leqslant \tau_{\lambda}^{\text{mul}}, ~|\tilde{\lambda}_{\tilde{i}+1} - \tilde{\lambda}_{\tilde{i}}| \leqslant \tau_{\lambda}^{\text{mul}}, \cdots,|\tilde{\lambda}_{\tilde{i}+m_{2}-1} - \tilde{\lambda}_{\tilde{i}+m_{2}-2}| \leqslant \tau_{\lambda}^{\text{mul}}, 
		\end{align*}
		then, the multiplicity of mode $\tilde{i}$ is $m = m_{1} + m_{2}$.
	}
}\vspace*{0.2cm}
Reset $\tilde{i} = \tilde{i}-m_{1}+1$, and thus we obtain the related multiple modes $\{\tilde{i},\tilde{i}+1,\ldots,\tilde{i}+m-1\}$.

\begin{algorithm}
	\caption{Identification of mode multiplicity and location of mode correspondence by FEC}
	\label{Algm:FEC}
	\hspace*{-0.5cm}\textbf{Input:} $\lambda_{i}$ and mode $i$ on $\mathbb{T}$.\\
	\hspace*{-0.5cm}\textbf{Output:} Multiplicity $n$ of mode $i$ and the related set of modes $\left\lbrace \tilde{i}, \ldots, \tilde{i}+m-1\right\rbrace $  on $\tilde{\mathbb{T}}$.\\
	$n =1, ~j = i+1$.\\
	\While{$ \lambda_{j}^{h} \leqslant \lambda_{\text{max}}$ \tcc*{Find multiplicity $n$ of mode $i$}}{
		\uIf{$\left| \lambda_{j}^{h} - \lambda_{j-1}^{h} \right| \leqslant \tau_{\lambda}^{\text{mul}}$}{
			\vspace*{0.1cm}
			$ j = j+1$.\\
			$ n = n+1$.
			}
		\Else{ \textbf{break}
		}
	}
	The multiplicity of mode $i$ is $n$.\\
	\While{$\lambda_{\text{min}}-a \leqslant \tilde{\lambda}_{\tilde{i}} \leqslant \lambda_{\text{max}}+a$}{
		Compute $ \{ | e_{i, \tilde{i}}| \} = |\lambda_{i}^{h} - \tilde{\lambda}_{\tilde{i}}|$.\\
		$ \tilde{i} = \tilde{i}+1$.
	}
	Select the minimum of $ \{ | e_{i, \tilde{i}}| \}$ and mark the relevant $\tilde{i}$.\\
	\uIf{$n=1$}{
		 The mode $i$ is related to mode $\tilde{i}$.
		}
	\Else{$m_{1} = 1, ~j = \tilde{i}-1$. \tcc*{Find multiplicity $m$ of mode $\tilde{i}$}
		\While{$\lambda_{\text{min}}-a \leqslant \tilde{\lambda}_{j}$}{
			\uIf{$\left| \tilde{\lambda}_{j} - \tilde{\lambda}_{j+1} \right| \leqslant \tau_{\lambda}^{\text{mul}}$}{
				$ j = j-1$.\\
				$ m_{1} = m_{1}+1$.
			}
			\Else{ \textbf{break}
			}
		}
		$m_{2} = 1, ~j = \tilde{i}+1.$ \vspace*{0.1cm}\\
		\While{$\tilde{\lambda}_{j} \leqslant \lambda_{\text{max}}+a $ }{
			\uIf{$\left| \tilde{\lambda}_{j} - \tilde{\lambda}_{j-1} \right| \leqslant \tau_{\lambda}^{\text{mul}}$}{
				$ j = j+1$.\\
				$ m_{2} = m_{2}+1$.
			}
			\Else{ \textbf{break}
			}		
		}
		The multiplicity of mode $\tilde{i}$ is $m = m_{1} + m_{2}$.\\
		\vspace*{0.1cm}
		Reset $\tilde{i} = \tilde{i}-m_{1}+1$. Then, the multiple modes $\{ i, \ldots, i+n-1 \}$ are related to multiple modes $\{\tilde{i}, \ldots, \tilde{i}+m-1 \}$.
		}	
	
\end{algorithm}

\subsubsection{Modal Assurance Criterion (MAC)}\label{sec:MAC}
The MAC method has been widely used to build correlation between analytical and experimental modal vectors for validation of experiment \cite{allemang2003modal,piegl2012nurbs}. In this paper, we utilize it to locate the correspondence between two computational modal shapes. The values of MAC are computed by computational eigenvectors, and then they are assembled into MAC matrix $\matr{\mathcal{M}}$ by using the formula
\begin{align}\label{eq:MAC}
\matr{\mathcal{M}}_{i,j}(\vect{\phi}_{i},\vect{\phi}_{j}) = \dfrac{\int_{\Omega}\vect{\phi}_{i}^{T}\matr{m}\vect{\phi}_{j}d\Omega}{\left\| \vect{\phi}_{i}\right\|_{m}^{2} \| \vect{\phi}_{j} \|_{m}^{2}},
\end{align}
where $\left\|\cdot\right\|_{m}$ is the mass norm and defined by
\begin{align}\label{massnorm}
\left\|\cdot\right\|_{m} := \left[ \int_{\Omega}(\cdot)^{T}\matr{m}(\cdot)d\Omega\right]^{\frac{1}{2}},
\end{align}
and $\matr{m}$ is defined in Eq.\eqref{eq:BBm}. Note that, if the eigenvectors $\vect{\phi}_{i}$ and $\vect{\phi}_{j}$ are obtained by the same mesh, the term $\int_{\Omega}\vect{\phi}_{i}^{T}\matr{m}\vect{\phi}_{j}d\Omega$ can be computed straightforward. If not, we should use projection to ensure integral is processed in the same domain. The details of projection can be found in Appendix \ref{app:prolong}. The values of the MAC are located in the interval $[0,1]$, where 0 means no consistent resemblance whereas 1 means a consistent correspondence. Generally, it is accepted that large values denote relatively consistent correlation whilst small value represents poor association.

In the MAC method, for mode $i$ with $\lambda_{i}^{h} \in [\lambda_{\text{min}}, ~\lambda_{\text{max}}]$, if the MAC value $\matr{\mathcal{M}}_{i,i+1} < \tau_{\text{MAC}}$ , where $\tau_{\text{MAC}}$ is the tolerance, mode $i$ is interpreted as a single mode. If the modes are multiple, such that \vspace*{0.2cm}\\
\fbox{
	\parbox{\textwidth}{
		\begin{align*}
		\matr{\mathcal{M}}_{i,i+1} \geqslant \tau_{\text{MAC}}, ~\matr{\mathcal{M}}_{i+1,i+2} \geqslant \tau_{\text{MAC}}, \cdots, ~\matr{\mathcal{M}}_{i+n-2,i+n-1} \geqslant \tau_{\text{MAC}}, 
		\end{align*}
		the multiplicity of mode $i$ is $n$, and the set of multiple modes are $\{i, i+1, \ldots i-n+1\}$.
	}
}\vspace*{0.2cm}
The strategy to deal with the multiple modes is similar to the FEC scheme, as shown in Fig.\ref{fig:mulmodesLoc}. The slight difference is that, in MAC method, we find the mode $\tilde{i}$ by the maximum of $\{ \matr{\mathcal{M}}_{i,\tilde{i}} \}, ~\forall\tilde{\lambda}_{\tilde{i}} \in [\lambda_{\text{min}} - a, ~\lambda_{\text{max}} + a]$. Also, we use MAC values to recognize the multiplicity of mode $\tilde{i}$ by following: \vspace*{0.2cm}\\
\fbox{
	\parbox{\textwidth}{If
		\begin{align*}
		\matr{\mathcal{M}}_{\tilde{i}-m_{1}+2, ~\tilde{i}-m_{1}+1} \geqslant \tau_{\text{MAC}},\cdots,\matr{\mathcal{M}}_{\tilde{i}, \tilde{i}-1} \geqslant \tau_{\text{MAC}}, ~\matr{\mathcal{M}}_{\tilde{i}+1, \tilde{i}} \geqslant \tau_{\text{MAC}}, \cdots,\matr{\mathcal{M}}_{\tilde{i}+m_{2}-1, ~\tilde{i}+m_{2}-2} \geqslant \tau_{\text{MAC}},
		\end{align*}
		the multiplicity of mode $\tilde{i}$ is $m = m_{1} + m_{2}$.
	}
}\vspace*{0.2cm}
Reset $\tilde{i} = \tilde{i}-m_{1}+1$, and then we obtain the related multiple modes $\{\tilde{i},\tilde{i}+1,\ldots,\tilde{i}+m-1\}$.

\begin{figure} 
	\centering
	\hspace*{4.5cm}\includegraphics[width=1\columnwidth]{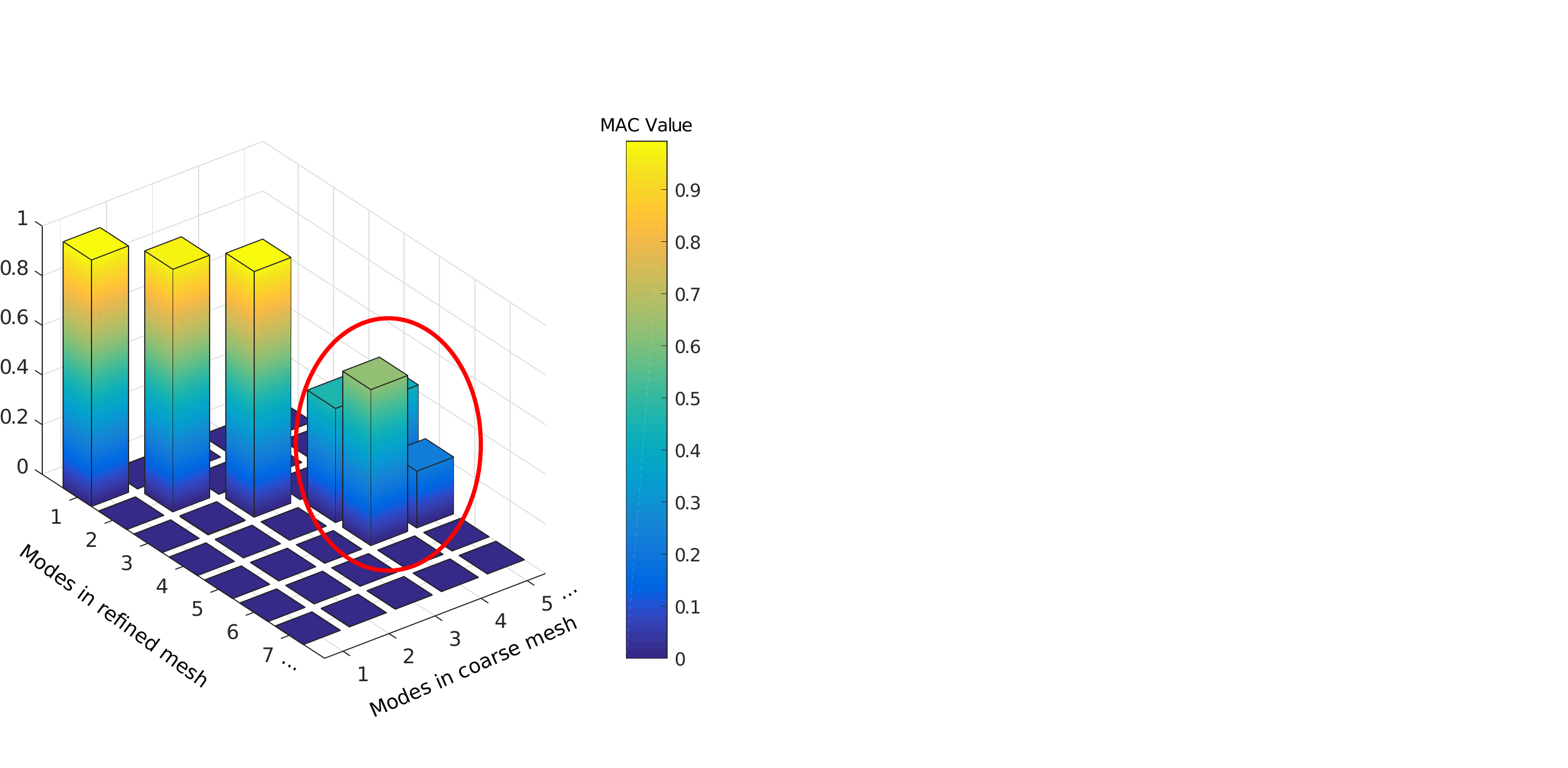}
	\caption{An example of 3D view for a MAC matrix obtained between coarse and refined meshes. Specially, the blocks of MAC values in red circle mean that double modes arise.}
	\label{fig:MAC}
\end{figure}

For instance, an example of MAC matrix $\matr{\mathcal{M}}$ is illustrated with 3D view in Fig.\ref{fig:MAC}. It is obvious that, the single modes 1,2 and 3 on coarse mesh are correlated to the single modes 1,2 and 3 on refined mesh, respectively. In addition, double modes 4 and 5 on $\mathbb{T}$ are associated to the double modes 4 and 5 on $\tilde{\mathbb{T}}$.

\begin{algorithm}
	\caption{Identification of mode multiplicity and location of mode correspondence by MAC}
	\label{Algm:MAC}
	\hspace*{-0.5cm}\textbf{Input:} $\lambda_{i}$ and mode $i$ on $\mathbb{T}$.\\
	\hspace*{-0.5cm}\textbf{Output:} Multiplicity $n$ of mode $i$ and the related set of modes $\left\lbrace \tilde{i}, \ldots, \tilde{i}+m-1\right\rbrace $  on $\tilde{\mathbb{T}}$.\\
	$n =1, ~j =i+1$\\
	\While{$\lambda_{\text{min}} \leqslant \lambda_{j} \leqslant \lambda_{\text{max}}$ \tcc*{Find multiplicity $n$ of mode $i$}}{
		Compute $\matr{\mathcal{M}}(\vect{\phi}_{j}^{h},\vect{\phi}_{j-1}^{h})$ by Eq.\eqref{eq:MAC}.\\
		\vspace*{0.1cm}
		\uIf{$\matr{\mathcal{M}}(\vect{\phi}_{j}^{h},\vect{\phi}_{j-1}^{h}) \geqslant \tau_{\text{MAC}}$}{
			\vspace*{0.1cm}
			$ j = j+1$.\\
			$ n = n+1$.
		}
		\Else{ \textbf{break}
		}
	}
	The multiplicity of mode $i$ is $n$.\\
	\While{$\lambda_{\text{min}}-a \leqslant \tilde{\lambda}_{\tilde{i}} \leqslant \lambda_{\text{max}}+a$}{
		\vspace*{0.1cm}
		Compute $ \left\lbrace  \matr{\mathcal{M}}_{i, \tilde{i}}(\vect{\phi}_{i}^{h},\tilde{\vect{\phi}}_{\tilde{i}}) \right\rbrace$ by Eq.\eqref{eq:MAC}.\\
		$ \tilde{i} = \tilde{i}+1$.
	}
	Select the maximum of $ \left\lbrace  \matr{\mathcal{M}}_{i, \tilde{i}} \right\rbrace $ and mark the relevant $\tilde{i}$.\\
	\uIf{$n=1$}{
		The mode $i$ is related to mode $\tilde{i}$.
	}
	\Else{
		$m_{1} = 1, ~j = \tilde{i}-1.$\\
		\While{$\lambda_{\text{min}}-a \leqslant \tilde{\lambda}_{j}$ \tcc*{Find multiplicity $m$ of mode $\tilde{i}$}}{
			Compute $\matr{\mathcal{M}}(\tilde{\vect{\phi}}_{j},\tilde{\vect{\phi}}_{j+1})$ by Eq.\eqref{eq:MAC}.\\
			\uIf{$\matr{\mathcal{M}}(\tilde{\vect{\phi}}_{j},\tilde{\vect{\phi}}_{j+1}) \geqslant \tau_{\text{MAC}}$}{
				\vspace*{0.1cm}
				$ j = j-1$.\\
				$ m_{1} = m_{1}+1$.
			}
			\Else{ \textbf{break}
			}
		}
		$m_{2} = 1, ~j = \tilde{i}+1.$ \vspace*{0.1cm}\\
		\While{$\tilde{\lambda}_{j} \leqslant \lambda_{\text{max}}+a $ }{
			\uIf{$\matr{\mathcal{M}}(\tilde{\vect{\phi}}_{j},\tilde{\vect{\phi}}_{j-1}) \geqslant \tau_{\lambda}^{\text{mul}}$}{
				\vspace*{0.1cm}
				$ j = j+1$.\\
				$ m_{2} = m_{2}+1$.
			}
			\Else{ \textbf{break}
			}		
		}
		The multiplicity of mode $\tilde{i}$ is $m = m_{1} + m_{2}$.\\
		Reset $\tilde{i} = \tilde{i}-m_{1}+1$. Then, The multiple modes $\{ i, \ldots, i+n-1 \}$ are related to multiple modes $\{\tilde{i}, \ldots, \tilde{i}+m-1 \}$.
	}		
\end{algorithm}

These two schemes, that is, FEC and MAC, will be compared in the numerical examples Section \ref{sec:eye}.

\subsection{Adaptivity for multiple modes}\label{sec:adpmul}
In section \ref{sec:adpone}, the adaptivity strategy for a single mode is established. In this section, we intend to propose a methodology to deal with the adaptivity for multiple modes.

Suppose that the $n$ multiple modes $\{i,i+1,\ldots,i+n-1\}$ on coarse mesh $\mathbb{T}$ are related to $m$ multiple modes  $\{\tilde{i},\tilde{i}+1,\ldots,\tilde{i}+m-1\}$ on refined mesh $\tilde{\mathbb{T}}$. Then all the vectors for modes $\{i,\ldots,i+n-1\}$ are actually defined by the linear combinations of basis eigenvectors $\vect{\Phi}$ in a eigenspace $\mathscr{P}$
\begin{align}\label{eq:eigenspace}
\mathscr{P} = \{\vect{\varphi}: \vect{\varphi} =  \vect{\Phi}\vect{\alpha},~~ \vect{\alpha}\in\mathbb{R}^n \},
\end{align}
where $\vect{\alpha} = (\alpha_{i} ~\cdots~ \alpha_{i+n-1})^T$ are arbitrary coefficients. The basis eigenvectors $\vect{\Phi}$ are defined by 
\begin{align*}
\vect{\Phi} = (\vect{\phi}_{i}^{h} ~~\cdots~~ \vect{\phi}_{i+n-1}^{h}) \in \mathbb{R}^{3\times n},
\end{align*}
where $\vect{\phi}_{i}^{h} = [\vect{\phi}_{i,w}^{h} ~~\vect{\phi}_{i, \theta_{x}}^{h} ~~\vect{\phi}_{i, \theta_{y}}^{h}]^{T}$ is the solution field for mode $i$. Likewise, the eigenspace for multiple modes $\{\tilde{i},\ldots,\tilde{i}+m-1\}$ is defined by
\begin{align}
\tilde{\mathscr{P}} = \{\tilde{\vect{\varphi}}: \tilde{\vect{\varphi}} =  \tilde{\vect{\Phi}}\tilde{\vect{\alpha}},~~ \tilde{\vect{\alpha}}\in\mathbb{R}^m \},
\end{align}
where
\begin{align*}
\tilde{\vect{\Phi}} = (\tilde{\vect{\phi}}_{\tilde{i}} ~~\cdots~~ \tilde{\vect{\phi}}_{\tilde{i}+m-1}) \in \mathbb{R}^{3\times m}, ~\tilde{\vect{\phi}}_{\tilde{i}} = [\tilde{\vect{\phi}}_{\tilde{i},w} ~~\tilde{\vect{\phi}}_{\tilde{i}, \theta_{x}} ~~\tilde{\vect{\phi}}_{\tilde{i}, \theta_{y}}]^{T}.
\end{align*}

\begin{figure} 
	\centering
	\def\svgwidth{2\columnwidth}
	\hspace*{0.5cm}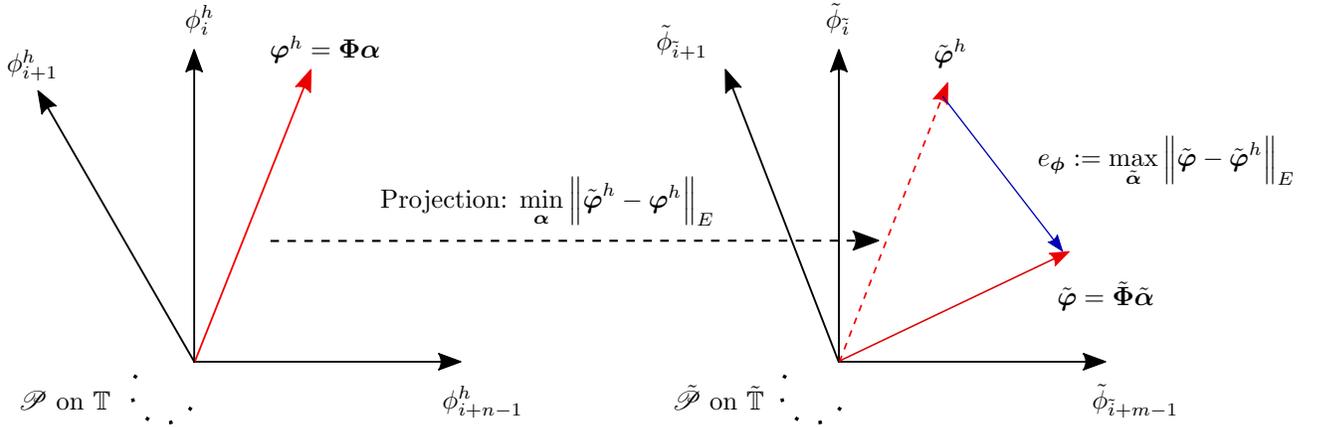
	\caption{The schematic of the measurement of the error estimator of eigenvectors for multiple modes between $\mathscr{P}$ and $\tilde{\mathscr{P}}$.}
	\label{fig:mulproject} 
\end{figure}

Now we aim to define error estimator of eigenvectors $e_{\vect{\phi}}$, between $\mathscr{P}$ and $\tilde{\mathscr{P}}$. The strategy, briefly summarized, is to project all the vectors in $\mathscr{P}$ onto $\tilde{\mathscr{P}}$, and then compute all the errors between the projected vectors and the vectors in $\tilde{\mathscr{P}}$. Among these errors, the maximum will be regarded as the $e_{\vect{\phi}}$. This process is presented in Fig.\ref{fig:mulproject}. As shown in Fig.\ref{fig:mulproject}, The minimization of the projection and the maximization of the errors can be defined as the following problem \vspace*{0.2cm}\\
\fbox{
	\parbox{\textwidth}{
		Find $\vect{\alpha}$ and $\tilde{\vect{\alpha}}$ such that
		\begin{align}\label{eq:maxmin}
		e_{\vect{\phi}} :=  \underset{\tilde{\vect{\alpha}}}{\text{max}},\underset{\vect{\alpha}}{\text{min}}\left\| \tilde{\vect{\Phi}}\tilde{\vect{\alpha}} -   \vect{\Phi}\vect{\alpha}\right\|_{E},
		\end{align}
		with $\left\| \tilde{\vect{\alpha}}\right\|_{2} = \sqrt{\tilde{\vect{\alpha}}^{T}\tilde{\vect{\alpha}}}= 1$.
	}	
}\vspace*{0.2cm}

Aiming to solve the problem in Eq.\eqref{eq:maxmin}, we build a Lagrange function $\mathscr{L}(\cdot)$ such that 
\begin{align}\label{eq:Lagrag}
\mathscr{L}(\vect{\alpha}, \tilde{\vect{\alpha}}, \mu) = \left\| \tilde{\vect{\Phi}}\vect{\tilde{\alpha}} - \vect{\Phi}\vect{\alpha} \right\|_{E}^{2} + \mu\left( \left\| \vect{\tilde{\alpha}}\right\|_{2}^{2}-1\right),
\end{align}
where $\left\|\vect{\tilde{\alpha}}\right\|_{2} = 1$ specifies the range of vectors $ \left\| \tilde{\vect{\varphi}}\right\| $ and $ \left\| \vect{\varphi}^{h} \right\|$. Then the problem defined in Eq.\eqref{eq:maxmin} can also be mathematically understood as follows \vspace*{0.2cm}\\
\fbox{
	\parbox{\textwidth}{
	Find $(\vect{\alpha}, \tilde{\vect{\alpha}}, \mu)$, such that,
	\begin{align}
	     &\delta_{\vect{\alpha}}\mathscr{L} = 0,\ \Rightarrow \delta_{\vect{\alpha}}\left\| \tilde{\vect{\Phi}}\vect{\tilde{\alpha}} - \vect{\Phi}\vect{\alpha} \right\|_{E}^{2}= 0,\label{eq:Lder1}\\
	     &\delta_{\tilde{\vect{\alpha}}}\mathscr{L} = 0,\ \Rightarrow \ \delta_{\tilde{\vect{\alpha}}}\left\| \tilde{\vect{\Phi}}\vect{\tilde{\alpha}} - \vect{\Phi}\vect{\alpha} \right\|_{E}^{2} + \mu\cdot\delta_{\tilde{\vect{\alpha}}}\left(\left\| \vect{\tilde{\alpha}}\right\| _{2}^{2}-1\right) = 0,\label{eq:Lder2}\\
	    &\delta_{\mu}\mathscr{L} = 0, \Rightarrow \left\| \vect{\tilde{\alpha}}\right\|_{2}  = 1. \label{eq:Lder3}
	\end{align}
	}
}\vspace*{0.2cm}
The basis eigenvectors $\vect{\Phi}$ and $\tilde{\vect{\Phi}}$ are discretized by PHT-spline basis functions $\vect{T}$ and $\tilde{\vect{T}}$, that is,
\begin{align}\label{eq:discret}
\vect{\Phi} = \vect{T}\bar{\vect{\Phi}}, ~
\tilde{\vect{\Phi}} = \tilde{\vect{T}}\tilde{\bar{\vect{\Phi}}},
\end{align}
where $\bar{\vect{\Phi}} \in \mathbb{R}^{3k \times n}$ and $\tilde{\bar{\vect{\Phi}}} \in \mathbb{R}^{3\tilde{k} \times m}$ are control variables of eigenvectors, and $k$ and $\tilde{k}$ are the numbers of control variables over $\mathbb{T}$ and $\tilde{\mathbb{T}}$, respectively. As discussed in Section \ref{sec:Errestm}, since refined mesh $\tilde{\mathbb{T}}$ is obtained by the hierarchical refinement of coarse mesh $\mathbb{T}$, we can have the prolongation of $\bar{\vect{\Phi}}$ in the $\tilde{\mathbb{T}}$, i.e.,
\begin{align}\label{eq:prolongbasiseigvec}
\vect{\Phi} = \vect{T}\bar{\vect{\Phi}} = \tilde{\vect{T}}\mathbb{P}\bar{\vect{\Phi}}.
\end{align}
Substitution of Eq.\eqref{eq:discret} and Eq.\eqref{eq:prolongbasiseigvec} into Eq.\eqref{eq:Lder1} and solving the equation, yields
\begin{align}\label{eq:alpharelation}
\bar{\vect{\Phi}}\matr{K}\bar{\vect{\Phi}}\vect{\alpha} = (\mathbb{P}\bar{\vect{\Phi}})^{T}\tilde{\matr{K}}\tilde{\bar{\vect{\Phi}}}\vect{\tilde{\alpha}},
\end{align}
where $\matr{K}$ and $\tilde{\matr{K}}$ are the stiffness matrices on $\mathbb{T}$ and $\tilde{\mathbb{T}}$,  defined in Eq.\eqref{eq:assembleGIFT}. Defining a $n\times m$ matrix $\vect{\mathcal{A}}$, that is,
\begin{align*}
\vect{\mathcal{A}} := \left((\mathbb{P}\bar{\vect{\Phi}})^{T}\tilde{\matr{K}}\tilde{\bar{\vect{\Phi}}}\right)^{-1}\bar{\vect{\Phi}}\matr{K}\bar{\vect{\Phi}},
\end{align*}
we can rewrite Eq.\eqref{eq:alpharelation}, which reads
\begin{align}\label{eq:alpharelationA}
 \vect{\alpha} = \vect{\mathcal{A}}\vect{\tilde{\alpha}}.
\end{align}
Similarly, substituting Eq.\eqref{eq:alpharelationA} into Eq.\eqref{eq:Lder2} and solving the equation, we have an eigenvalue problem for $\vect{\tilde{\alpha}}$
\begin{align}\label{eq:eigenalpha}
	\matr{K}^{*}\vect{\tilde{\alpha}} = \mu\vect{\tilde{\alpha}},
\end{align}
where
\begin{align*}
\matr{K}^{*} = \tilde{\bar{\vect{\Phi}}}^{T}\tilde{\matr{K}}\tilde{\bar{\vect{\Phi}}} + \vect{\mathcal{A}}^{T}\bar{\vect{\Phi}}\matr{K}\bar{\vect{\Phi}}\vect{\mathcal{A}} - \vect{\mathcal{A}}(\mathbb{P}\bar{\vect{\Phi}})^{T}\tilde{\matr{K}}\tilde{\bar{\vect{\Phi}}}.
\end{align*}
Solving the eigenvalue problem in Eq.\eqref{eq:eigenalpha}, $\vect{\tilde{\alpha}}$ can be obtained. Then, substituting $\vect{\tilde{\alpha}}$ into Eq.\eqref{eq:alpharelationA}, $\vect{\alpha}$ will be obtained.
\vspace*{0.1cm}Furthermore, the error estimator $ e_{\vect{\phi}}$ in Eq.\eqref{eq:maxmin} can be determined. Additionally, the local error estimator $ e_{\vect{\phi}}(\Omega_{e})$ can be given by
\begin{align}\label{eq:e_Kloc_mul}
e_{\vect{\phi}}(\Omega_{e}) = \left[ \int_{\Omega_{e}}(\mathbb{B}_{b}\vect{e}_{\vect{\varphi}})^{T}\matr{D}\mathbb{B}_{b}\vect{e}_{\vect{\varphi}}d\Omega_{e} + \int_{\Omega_{e}}(\mathbb{B}_{s}\vect{e}_{\vect{\varphi}})^{T}\matr{D}\mathbb{B}_{s}\vect{e}_{\vect{\varphi}}d\Omega_{e}\right]^{\frac{1}{2}},
\end{align}
where $\vect{e}_{\vect{\varphi}} = \vect{\tilde{\varphi}}-\vect{\varphi}$. The marking strategy has been discussed in Section \ref{sec:Markstrategy}. Define the error estimator of frequency and relative error estimator of eigenvector
\begin{align}\label{eq:errestmmul}
|e_{\lambda}|  = \left| \dfrac{1}{n}\sum_{j=i}^{i+n-1}\lambda_{j} -\dfrac{1}{m}\sum_{j=\tilde{i}}^{\tilde{i}+m-1}\tilde{\lambda}_{j} \right|, ~\delta_{\vect{\phi}} = \dfrac{e_{\phi}}{\left\|\tilde{\vect{\varphi}}\right\|_{E}},
\end{align}
and subsequently, the adaptivity can be conducted until they satisfy the given thresholds
\begin{align}
|e_{\lambda}| \leqslant \tau_{\lambda}, ~\delta_{\vect{\phi}} \leqslant \tau_{\phi}.
\end{align}
The methodology proposed above is summarized in Algorithm \ref{Algm:mulmode}. 

\begin{algorithm}
	\caption{Adaptivity process for the $n$ multiple modes $\{i, \ldots, i+n-1\}$}
	\label{Algm:mulmode}
	\hspace*{-0.5cm}\textbf{Input:} Multiple modes $\{i,\ldots\,i-n+1\}$ on $\mathbb{T}$ and $\{\tilde{i},\ldots\,\tilde{i}-m+1\}$ on $\tilde{\mathbb{T}}$.\\
	\hspace*{-0.5cm}\textbf{Output:} Updated $\mathbb{T}$ after refinement.\\
	\textbf{Step 1.} Define eigenspaces $\mathscr{P}$ and $\tilde{\mathscr{P}}$, and vectors $\vect{\varphi} = \vect{\Phi}\vect{\alpha} \in \mathscr{P}$ and $\tilde{\vect{\varphi}} = \tilde{\vect{\Phi}}\tilde{\alpha} \in \tilde{\mathscr{P}}$.\\
	\textbf{Step 2.} Define the error estimator of eigenvector $e_{\vect{\phi}}$ in Eq.\eqref{eq:maxmin}.\\
	\textbf{Step 3.} Build a Lagrange function $\mathscr{L}(\cdot)$ in Eq.\eqref{eq:Lagrag}, and solve it by Eq.\eqref{eq:Lder1}-\eqref{eq:Lder3} to obtain vectors $\vect{\varphi}$, $\tilde{\vect{\varphi}}$ and $e_{\vect{\phi}}$.\\
	\textbf{Step 4.} Compute $\left|e_{\lambda}\right|$ and $\delta_{\vect{\phi}}$ by Eq.\eqref{eq:errestmmul}.\\
	\While{$\left|e_{\lambda}\right| \leqslant \tau_{\lambda}, ~\delta_{\vect{\phi}} \leqslant \tau_{\vect{\phi}}$}{
		\For{$j\leftarrow 1$ \KwTo $N$}{
			Compute $\left\| e_{\vect{\phi}}(\Omega_{e}^{j})\right\|_{E}^{2}$ by Eq.\eqref{eq:e_Kloc_mul}.
		}
		
		Sort values of $\left\| e_{\vect{\phi}}(\Omega_{e}^{j})\right\|_{E}^{2}$ from large to small.	
		
		\For{$j\leftarrow 1$ \KwTo $N$}{
			\If{$\sum \limits_{j^{*} = 1}^{j} \left\| e_{\vect{\phi}}(\Omega_{e}^{j*})\right\|_{E}^{2} \geqslant \tau e_{\vect{\phi}}^{2}$}{
				Mark $N^{*} = j$\\ 
				\textbf{break}
			}		
		}
		
		\For{$j\leftarrow 1$ \KwTo $N^{*}$}{
			\vspace*{0.1cm}Refine element $j$ to update $\mathbb{T}$.
		}
		Repeat \textbf{Step 1 -- Step 4}.
	}
\end{algorithm}

\section{Numerical examples} \label{sec:Numtests}
In this section, three numerical examples are carried out for the following purposes. The example in Section\ref{sec:cirplate} aims to show that the GIFT method (NURBS for design and PHT splines for analysis) delivers the results with good accuracy as those obtained in IGA framework. Afterwards, we use GIFT method for examples in both Section \ref{sec:eye} and Section \ref{sec:plateholes}. The example in Section \ref{sec:eye} presents a comparison between two strategies that is, MAC and FEC proposed in Section \ref{sec:Locmode}, for modal resemblance determinations. The results illustrate that the shorter the width of frequencies of interest is, the better the MAC method is to be used. Finally, in Section \ref{sec:plateholes}, we apply MAC method within GIFT framework to study the local adaptivity for structural vibration by sweeping the modes from low to high frequencies.
\subsection{Homogeneous circular plate}\label{sec:cirplate}
\begin{figure}
	\centering
	\begin{subfigure}[t]{0.4\textwidth}
		\centering
		\includegraphics[width=1\columnwidth]{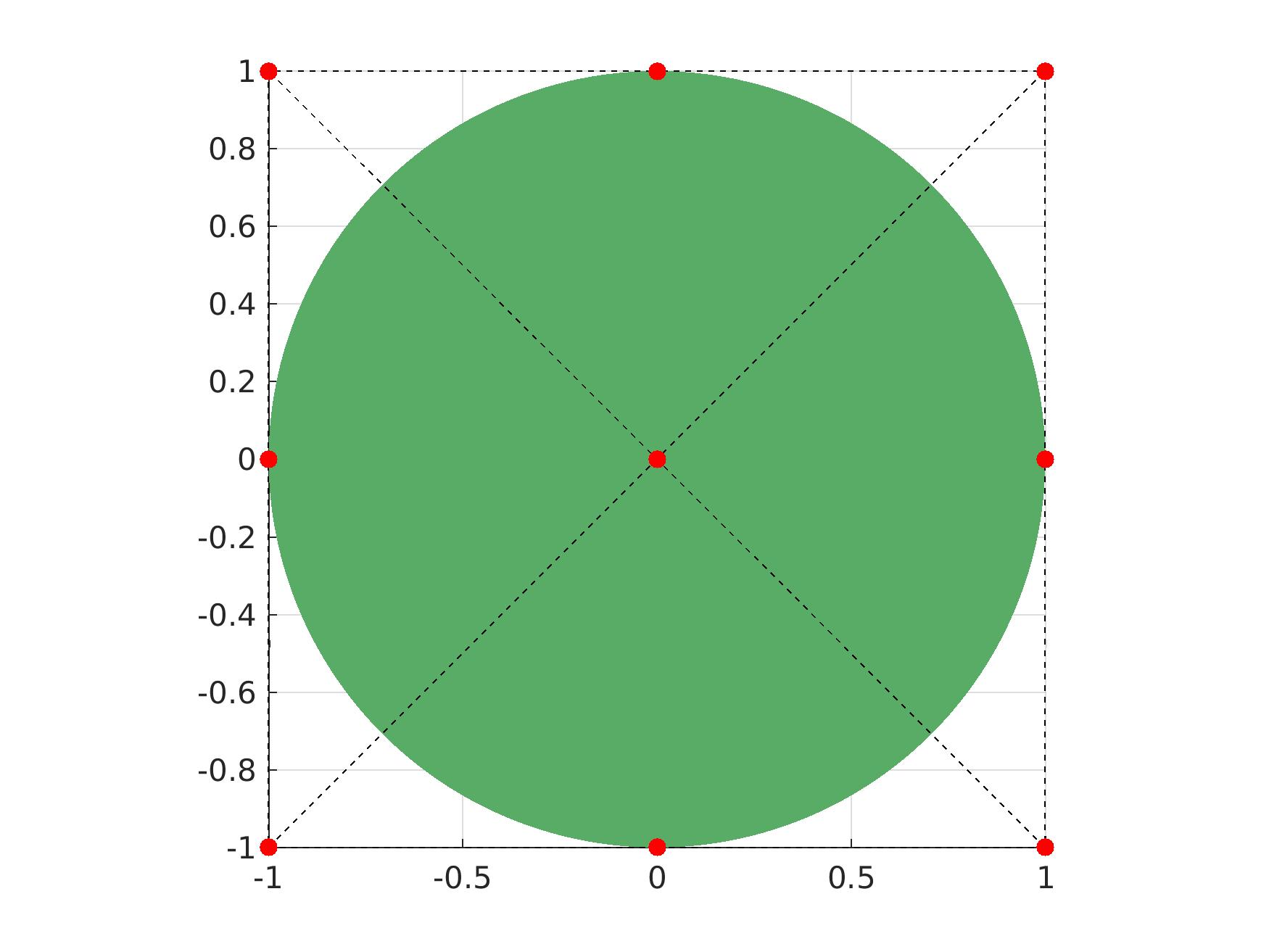}
		\caption{Quadratic NURBS}
	\end{subfigure}
	\begin{subfigure}[t]{0.4\textwidth}
		\centering
		\includegraphics[width=1\columnwidth]{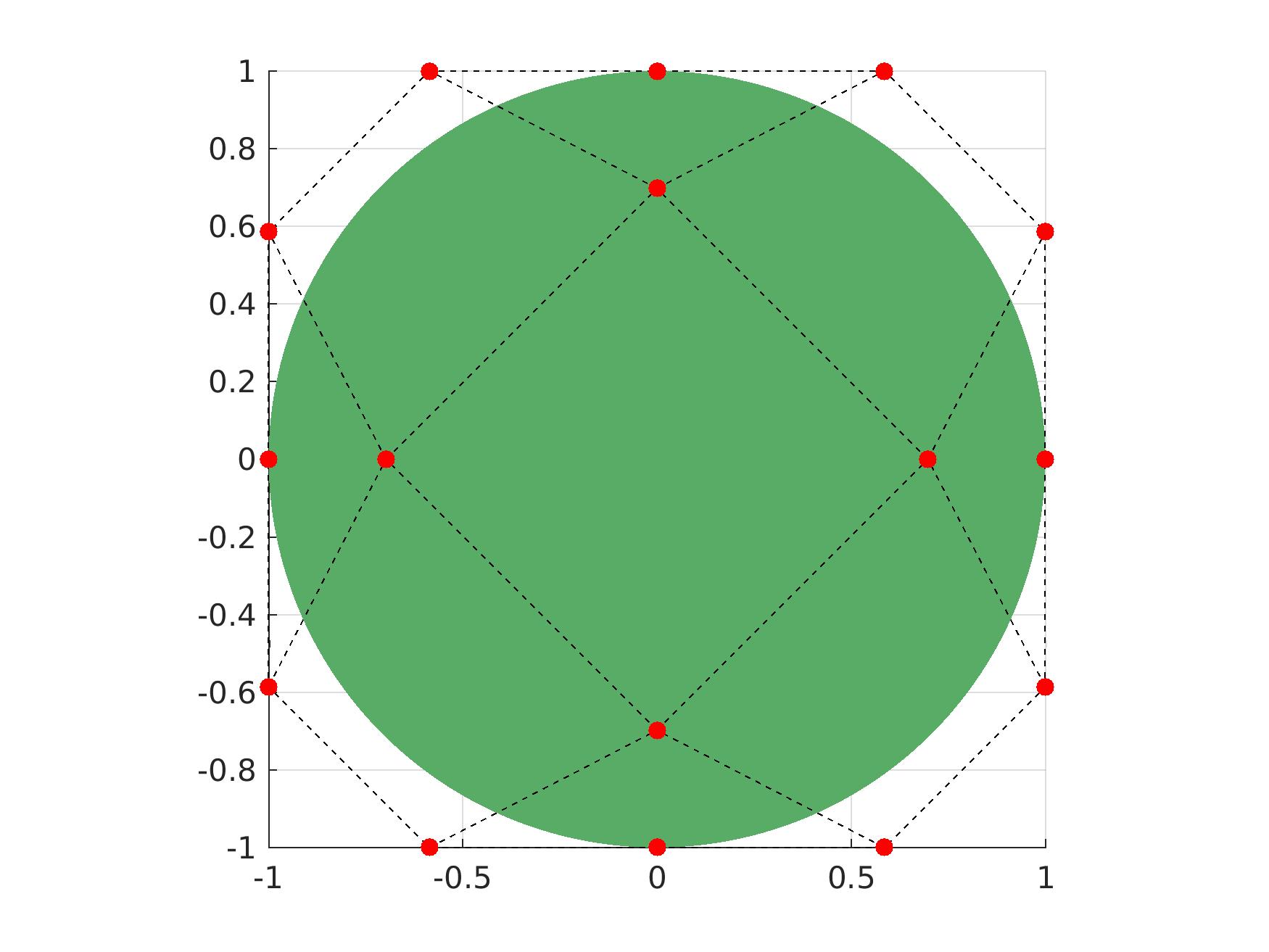}
		\caption{Cubic NURBS and RHT}
	\end{subfigure}
	\caption{Geometry and intial control points for disk generated by GIFT(a) and IGA(b).}
	\label{fig:controlpts}
\end{figure}

In this example, we compare normalized natural frequency $\lambda_{N}$ obtained by IGA(NURBS), IGA(RHT) and GIFT(NURBS+PHT) in case of the vibration of the disk. The $\lambda_{N}$ can be expressed as $\lambda_{N} = \lambda^{h}/\lambda_{\text{ext}}$, where $\lambda_{\text{ext}}$ is the exact solution obtained from \cite{liew1993transverse}. The material parameters are as follows: Young's modulus $E=1$, density $\rho = 1$, Poisson's ratio $\nu=0.3$, thickness-span ratios $h/r$ ($h$ is the thickness and $r$ is the radius). As mentioned in Appendix \ref{app:RHT}, without any refinement, the initial RHT splines are NURBS so that bi-cubic IGA(RHT) and IGA(NURBS) share the same control points, as shown in Fig.\ref{fig:controlpts}(b). Whilst in GIFT method, the quadratic NURBS is adopted for geometry as it is precise enough to generate the circular shape, as seen in Fig.\ref{fig:controlpts}(a), and the cubic PHT splines are exploited to represent solution fields. In terms of the uniform refinement, as it can be seen in Fig\ref{fig:PHTRHTmesh}, the PHT and RHT mesh are exactly the same. From the results illustrated in Tab.\ref{Sdisk}, Tab.\ref{Cdisk} and Fig.\ref{fig:Neigval3method}, we can see the results obtained by GIFT method has an excellent agreement with those computed by IGA method, for the first 6 modes with both simply supported and clamped boundary conditions, in case of $h/r=0.1$ and $h/r=0.2$. 

\begin{table}
	\caption{Comparison of normalized frequency $\lambda_{N}$ for simply supported circular plates.}
	\label{Sdisk}
	\centering
	\begin{tabular}{ccccccccc}\hline
		\multirow{2}*{$h/r$} & \multirow{2}*{Method} & \multirow{2}*{Dof} & \multicolumn{6}{l}{Mode number}\\\cline{4-9}
		& & & 1 & 2 & 3 & 4 & 5 & 6\\\hline
		0.1 & NURBS & & 1.0216 & 1.0801 & 1.0801 & 1.1147 & 1.7991 & 1.7236\\ 
		& RHT & 108 & 1.0020 & 1.1569 & 1.1569 & 1.2127 & 1.4319 & 1.4026\\
		& GIFT &  & 1.0051 & 1.2026 & 1.2026 & 1.2418 & 1.5794 & 1.5266 \\\cline{2-9}
		& NURBS &  & 1.0002 & 1.0024 & 1.0024 & 1.0052 & 1.0147 & 1.0155\\
		& RHT &  300 & 1.0006 & 1.0057 & 1.0057 & 1.0036 & 1.0116 & 1.0030\\
		& GIFT & & 1.0012 & 1.0077 & 1.0077 & 1.0041 & 1.0143 & 1.0037 \\\cline{2-9}
		& NURBS &  & 1.0000 & 1.0003 & 1.0003 & 1.0005 & 1.0005 & 1.0005\\
		& RHT & 972 & 1.0001 & 1.0004 & 1.0004 & 1.0008 & 1.0010 & 1.0009\\
		& GIFT &  & 1.0001 & 1.0005 & 1.0005 & 1.0008 & 1.0010 & 1.0010\\\cline{2-9}
		\\
		0.2  & NURBS & & 1.0139	& 1.0405 & 1.0405 & 1.0633 & 1.4913 & 1.4551\\
		& RHT & 108 & 1.0016 & 1.0572 & 1.0572 & 1.0823 & 1.2744 & 1.2616\\
		& GIFT &  & 1.0036 & 1.0737 & 1.0737 & 1.0948 & 1.3626 & 1.3372\\\cline{2-9}
		& NURBS &  & 1.0004 & 1.0014 & 1.0014 & 1.0029 & 1.0063 & 1.0070\\
		& RHT & 300 & 1.0006 & 1.0026 & 1.0026 & 1.0036 & 1.0051 & 1.0035\\
		& GIFT &  & 1.0008 & 1.0033 & 1.0033 & 1.0039 & 1.0060 & 1.0040\\\cline{2-9}
		& NURBS &  & 1.0003 & 1.0008 & 1.0008 & 1.0014 & 1.0014 & 1.0015\\
		& RHT & 972 & 1.0004 & 1.0009 & 1.0009 & 1.0015 & 1.0016 & 1.0017\\
		& GIFT & & 1.0004 & 1.0009 & 1.0009 & 1.0015 & 1.0016 & 1.0017\\
		\hline
	\end{tabular}
\end{table}

\begin{table}
	\caption{Comparisons of normalized frequency $\lambda_{N}$ for fully clamped circular plates.}
	\label{Cdisk}
	\centering
	\begin{tabular}{ccccccccc}\hline
		\multirow{2}*{$h/r$} & \multirow{2}*{Method} & \multirow{2}*{Dof} & \multicolumn{6}{l}{Mode number}\\\cline{4-9}
		& & & 1 & 2 & 3 & 4 & 5 & 6\\\hline
		0.1 & NURBS & & 1.0984 & 1.1702 & 1.1702 & 1.1958 & 2.6834 & 2.4558\\ 
		& RHT & 108 & 1.0253 & 1.2405 & 1.2405 & 1.2387 & 2.3174 & 2.1222\\
		& GIFT &  & 1.0489 & 1.2664 & 1.2664 & 1.2635 & 2.4701 & 2.2603\\\cline{2-9}
		& NURBS &  & 1.0013 & 1.0036 & 1.0036 & 1.0048 & 1.0254 & 1.0302\\
		& RHT &  300 & 1.0020 & 1.0101 & 1.0101 & 1.0027 & 1.0137 & 1.0094\\
		& GIFT & & 1.0026 & 1.0134 & 1.0134 & 1.0059 & 1.0173 & 1.0125\\\cline{2-9}
		& NURBS &  & 1.0003 & 0.9978 & 0.9978 & 0.9946 & 0.9947 & 1.0010\\
		& RHT & 972 & 1.0004 & 0.9981 & 0.9981 & 0.9949 & 0.9954 & 1.0016\\
		& GIFT & & 1.0004 & 0.9981 & 0.9981 & 0.9950 & 0.9955 & 1.0017\\\cline{2-9}
		\\
		0.2  & NURBS & & 1.0556 & 1.0835 & 1.0835 & 1.0953 & 1.6798 & 1.6112\\
		& RHT & 108 & 1.0153 & 1.0810 & 1.0810 & 1.0872 & 1.4816 & 1.4279\\
		& GIFT & & 1.0285 & 1.0973 & 1.0973 & 1.1043 & 1.5594 & 1.4996\\\cline{2-9}
		& NURBS &  & 1.0013	& 0.9991 & 0.9991 & 0.9971 & 1.0041 & 1.0125\\
		& RHT & 300 & 1.0017 & 1.0012 & 1.0012 & 0.9983 & 1.0000 & 1.0065\\
		& GIFT &  & 1.0019 & 1.0022 & 1.0022 & 0.9993 & 1.0011 & 1.0076\\\cline{2-9}
		& NURBS &  & 1.0011 & 0.9975 & 0.9975 & 0.9941 & 0.9942 & 1.0024\\
		& RHT & 972 & 1.0011 & 0.9976 & 0.9976 & 0.9942 & 0.9944 & 1.0026\\
		& GIFT & &1.0011 & 0.9976 & 0.9976 & 0.9943 & 0.9945 & 1.0027\\
		\hline
	\end{tabular}
\end{table}

\begin{figure}
	\centering
	\begin{subfigure}[t]{0.2\textwidth}
		\centering
		\includegraphics[width=1\columnwidth]{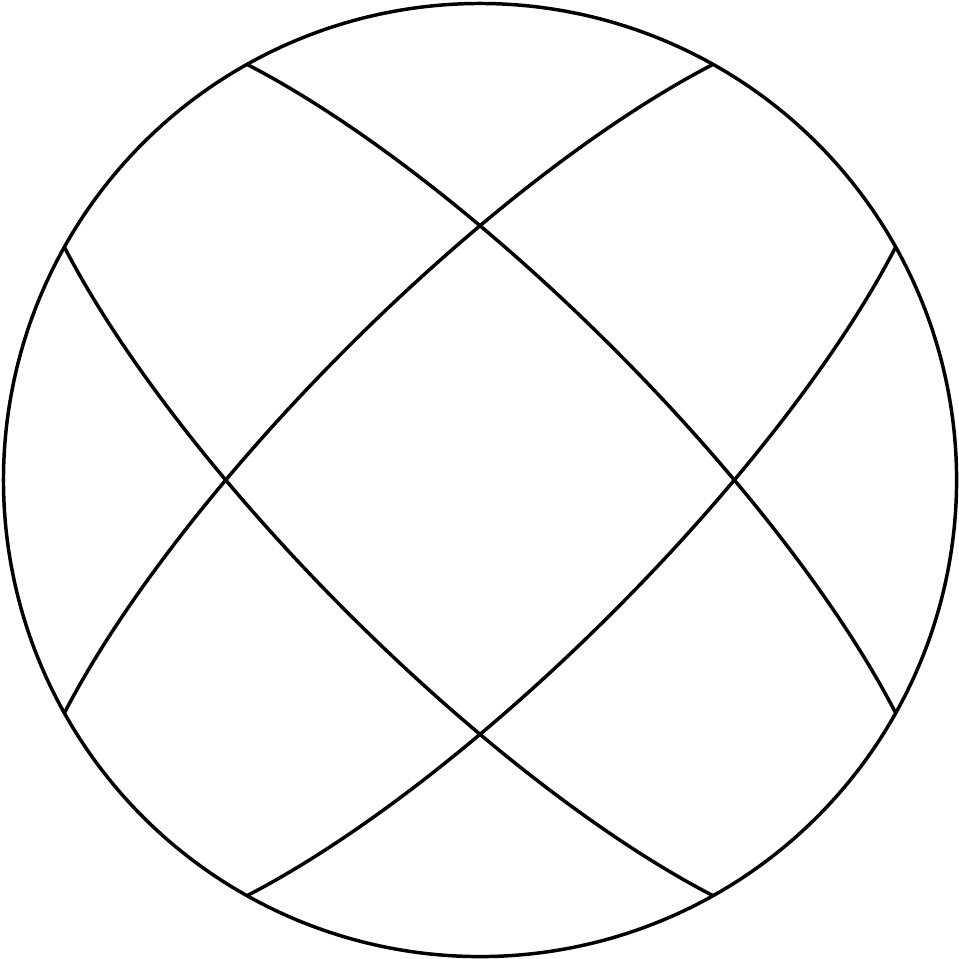}
		\caption{108 dofs}
	\end{subfigure}
	\hspace{1cm}
	\begin{subfigure}[t]{0.2\textwidth}
		\centering
		\includegraphics[width=1\columnwidth]{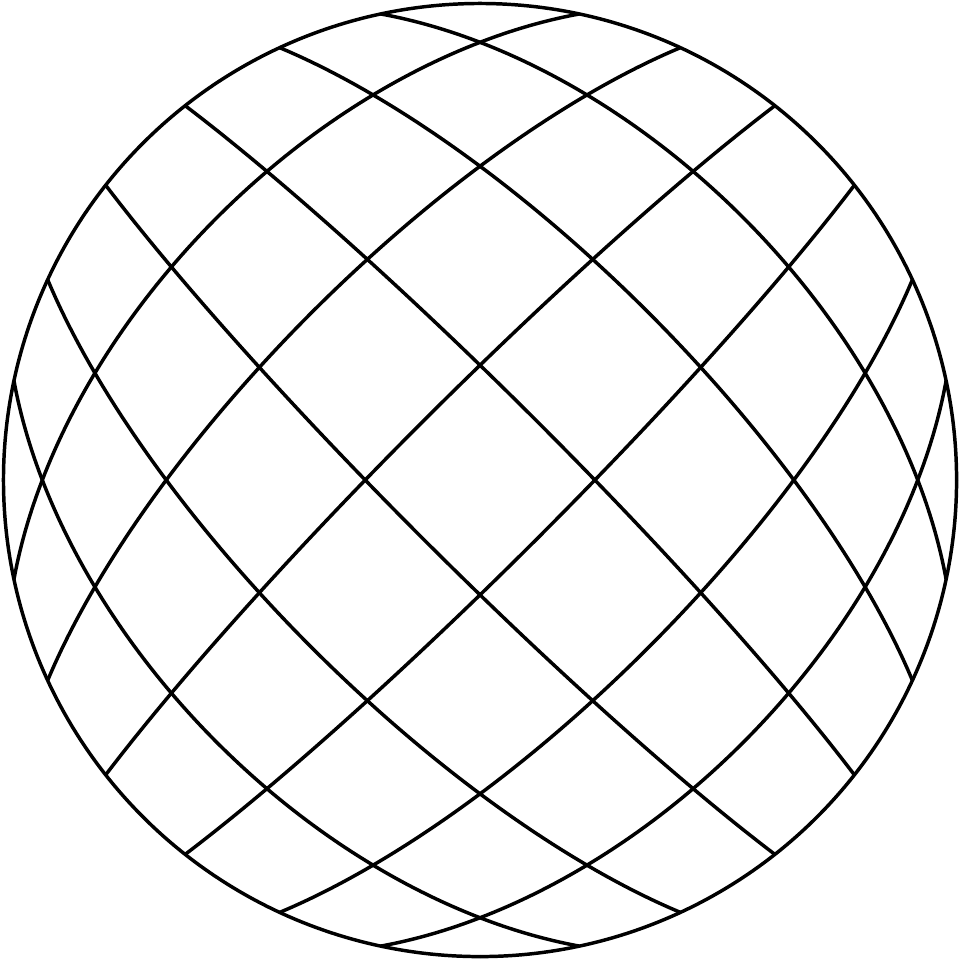}
		\caption{300 dofs}
	\end{subfigure}
	\hspace{1cm}
	\begin{subfigure}[t]{0.2\textwidth}
		\centering
		\includegraphics[width=1\columnwidth]{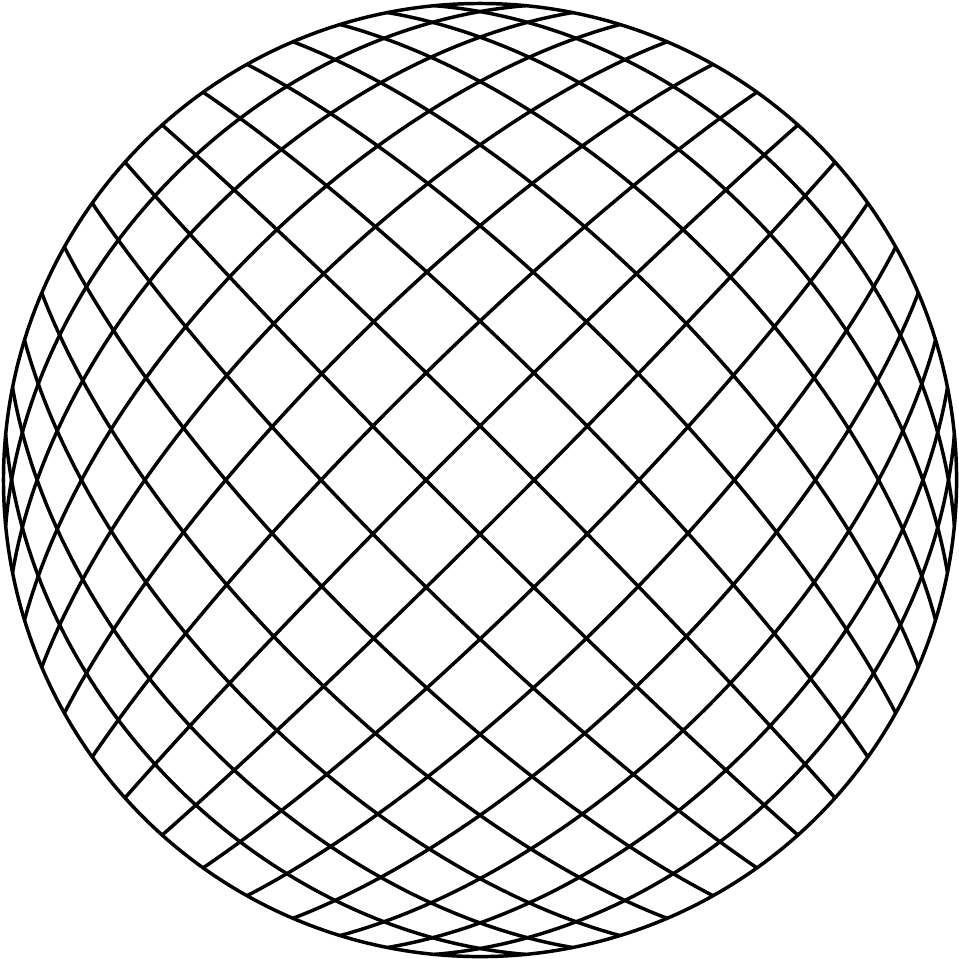}
		\caption{972 dofs}
	\end{subfigure}
	\caption{Mesh of solution field with NURBS $(p=3,q=3)$.}
	\label{control points}	
\end{figure}

\begin{figure}
	\centering
	\begin{subfigure}[t]{0.2\textwidth}
		\centering
		\includegraphics[width=1\columnwidth]{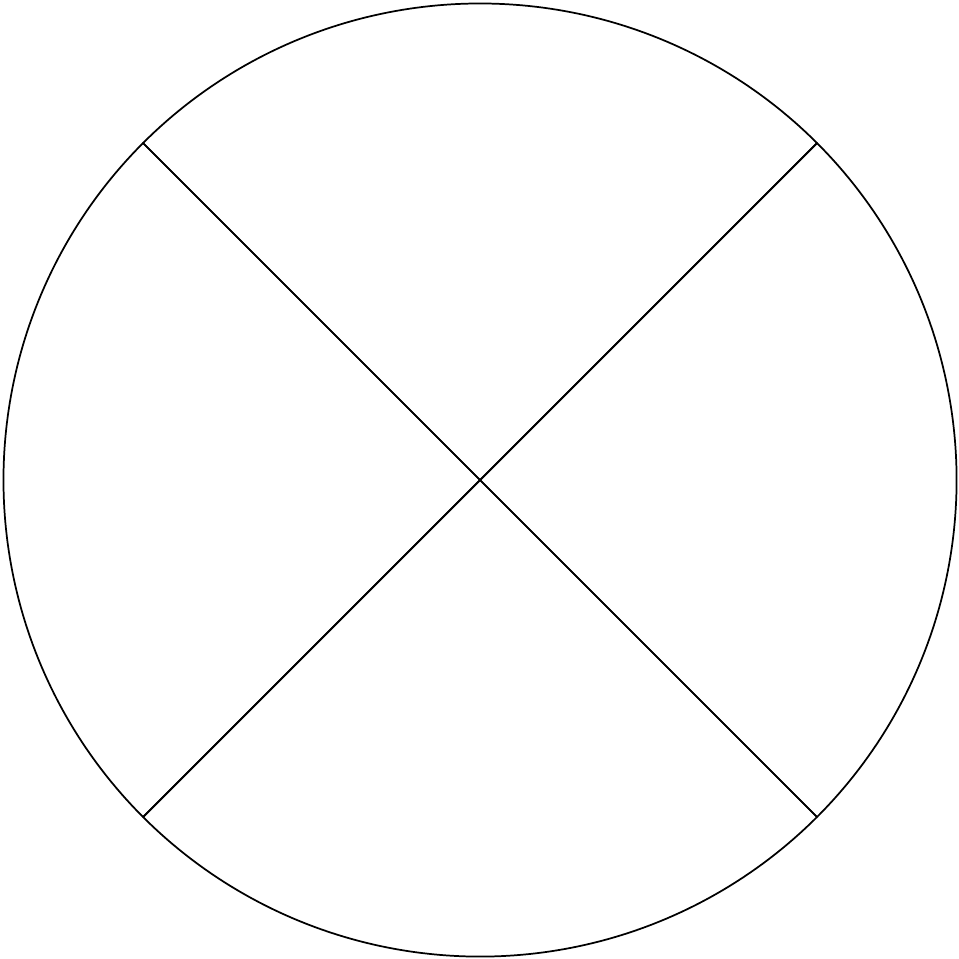}
		\caption{108 dofs}
	\end{subfigure}
	\hspace{1cm}
	\begin{subfigure}[t]{0.2\textwidth}
		\centering
		\includegraphics[width=1\columnwidth]{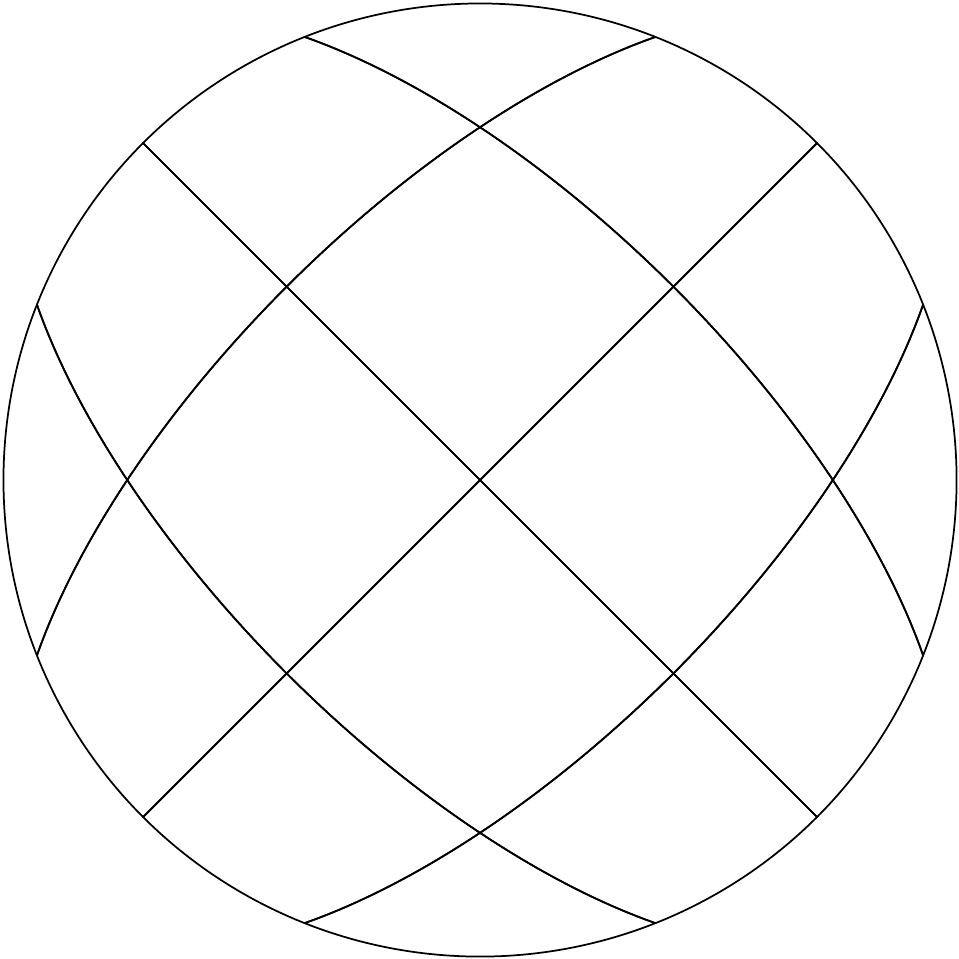}
		\caption{300 dofs}
	\end{subfigure}
	\hspace{1cm}
	\begin{subfigure}[t]{0.2\textwidth}
		\centering
		\includegraphics[width=1\columnwidth]{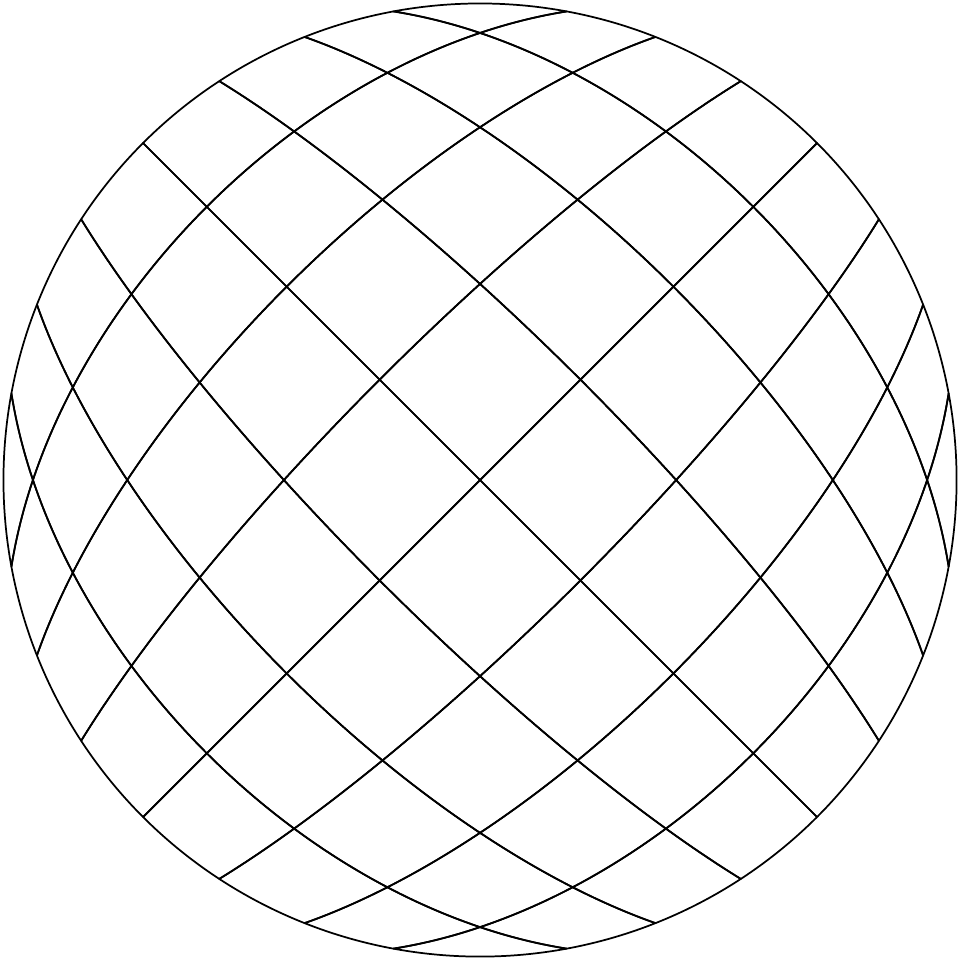}
		\caption{972 dofs}
	\end{subfigure}
	\caption{Mesh of solution field with PHT $(p=3,q=3)$ and RHT $(p=3,q=3)$.}
	\label{fig:PHTRHTmesh}	
\end{figure}

\begin{figure}
	\centering
	\begin{subfigure}[t]{0.25\textwidth}
		\centering
		\includegraphics[width=1\columnwidth]{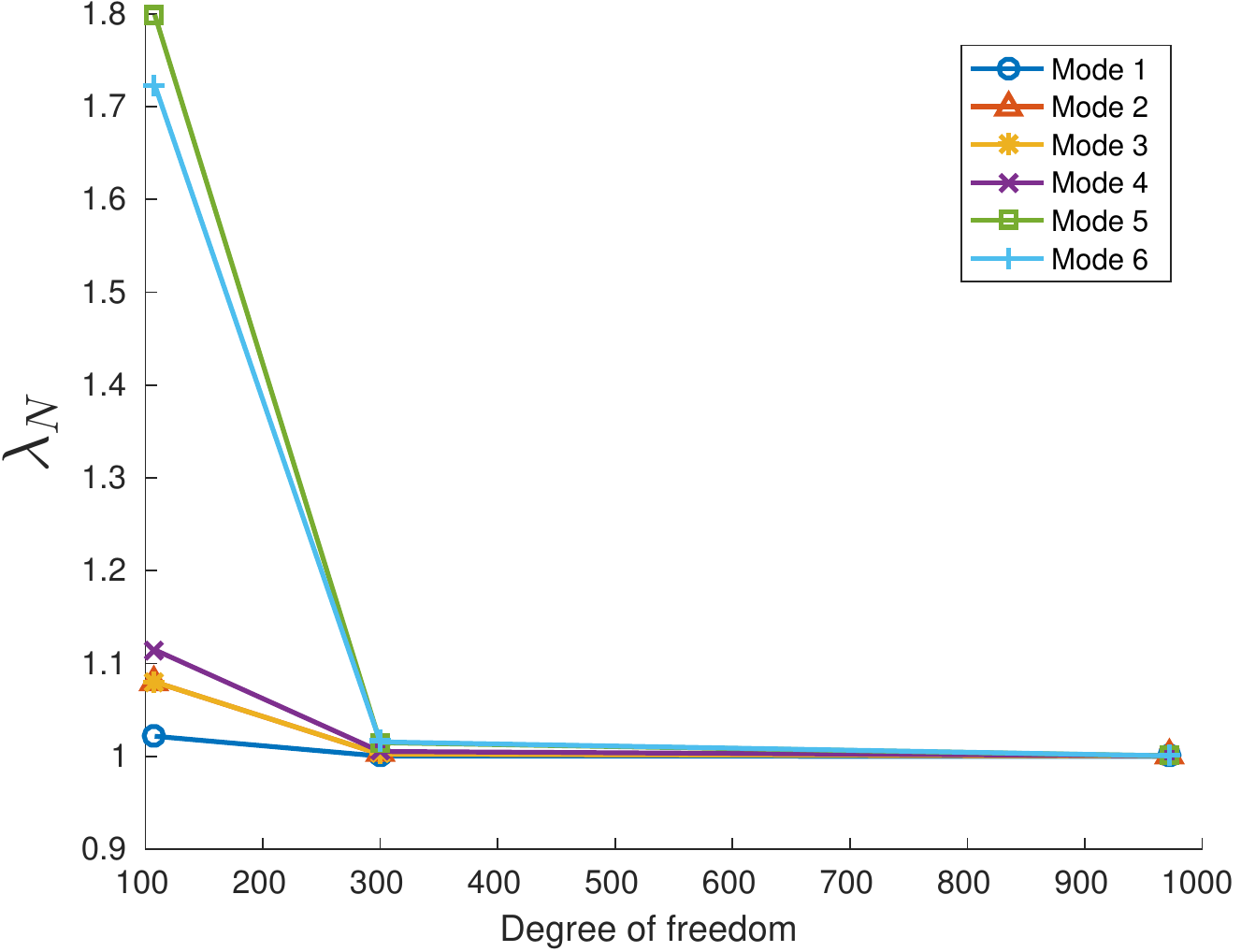}
		\caption{NURBS}
	\end{subfigure}
	\hspace{1cm}
	\begin{subfigure}[t]{0.25\textwidth}
		\centering
		\includegraphics[width=1\columnwidth]{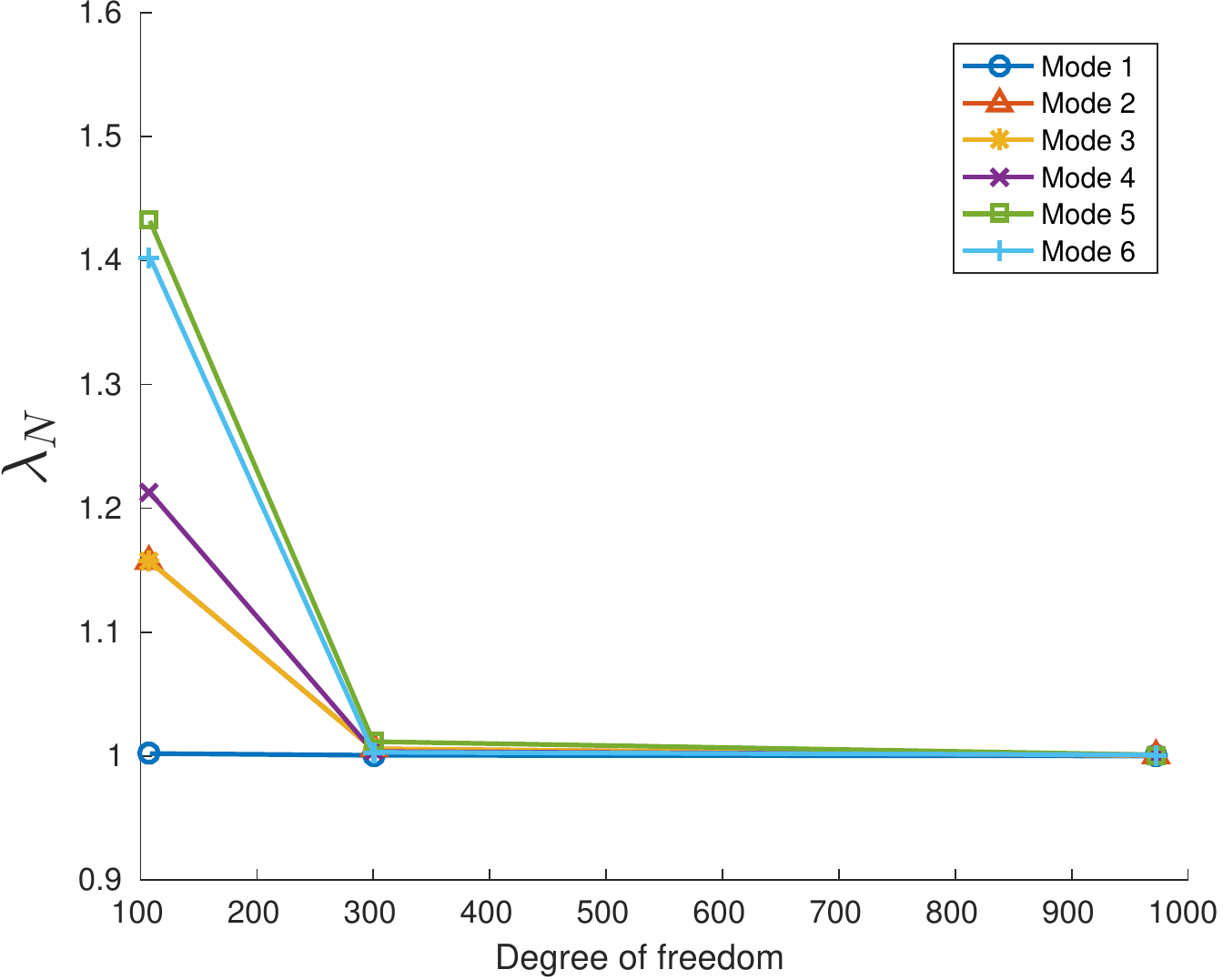}
		\caption{RHT}
	\end{subfigure}
	\hspace{1cm}
	\begin{subfigure}[t]{0.25\textwidth}
		\centering
		\includegraphics[width=1\columnwidth]{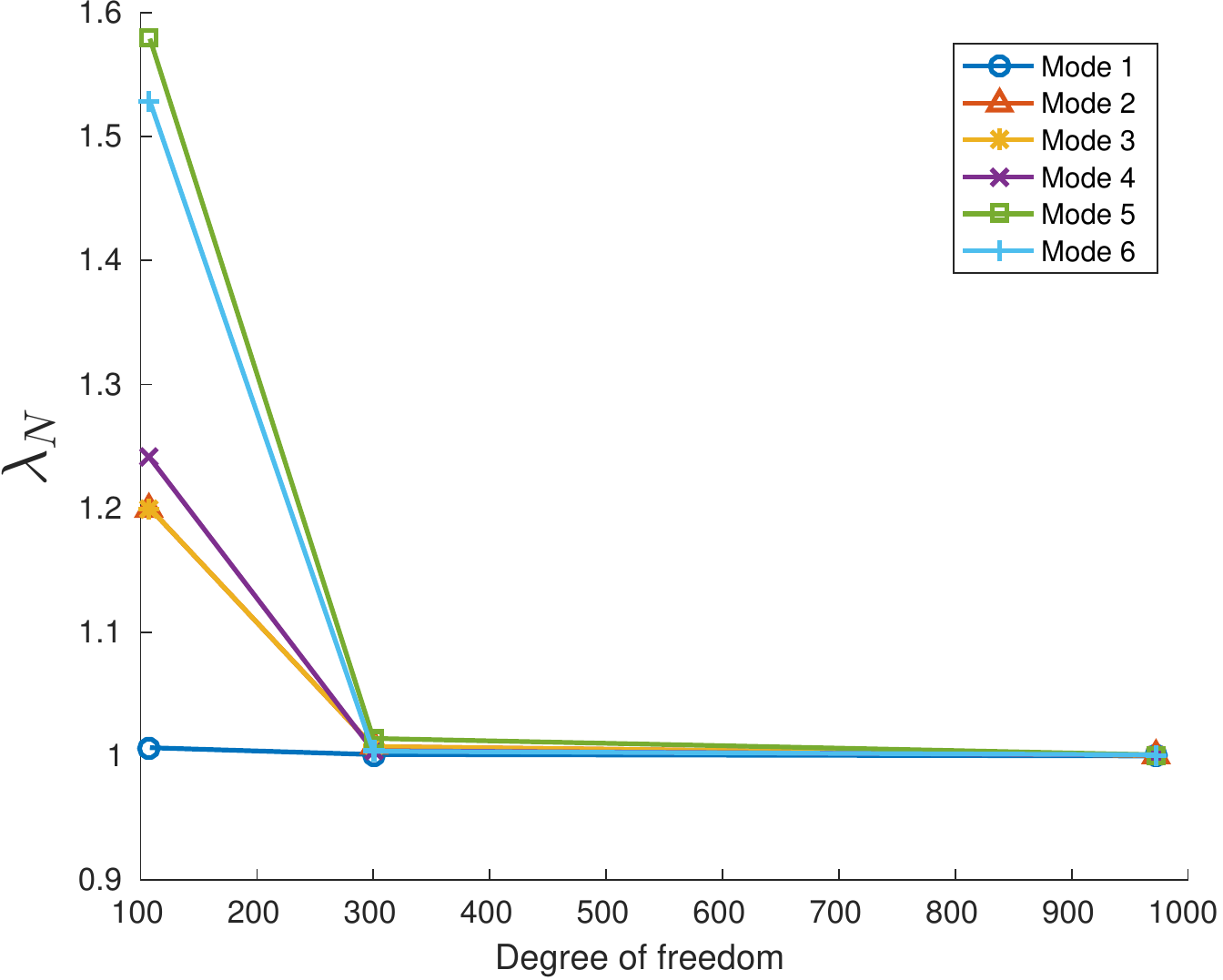}
		\caption{GIFT}
	\end{subfigure}
	\caption{Convergence of normalized eigenvalue $\lambda_{N}$ computed by NURBS, RHT and GIFT for vibration of the simply supported circular plate with $h/r = 0.1$.}
	\label{fig:Neigval3method}	
\end{figure}

\subsection{Heterogeneous eye shape with a hole}\label{sec:eye}
\begin{figure}
	\centering
	\begin{subfigure}[t]{0.3\textwidth}
		\centering
		\includegraphics[width=1\columnwidth]{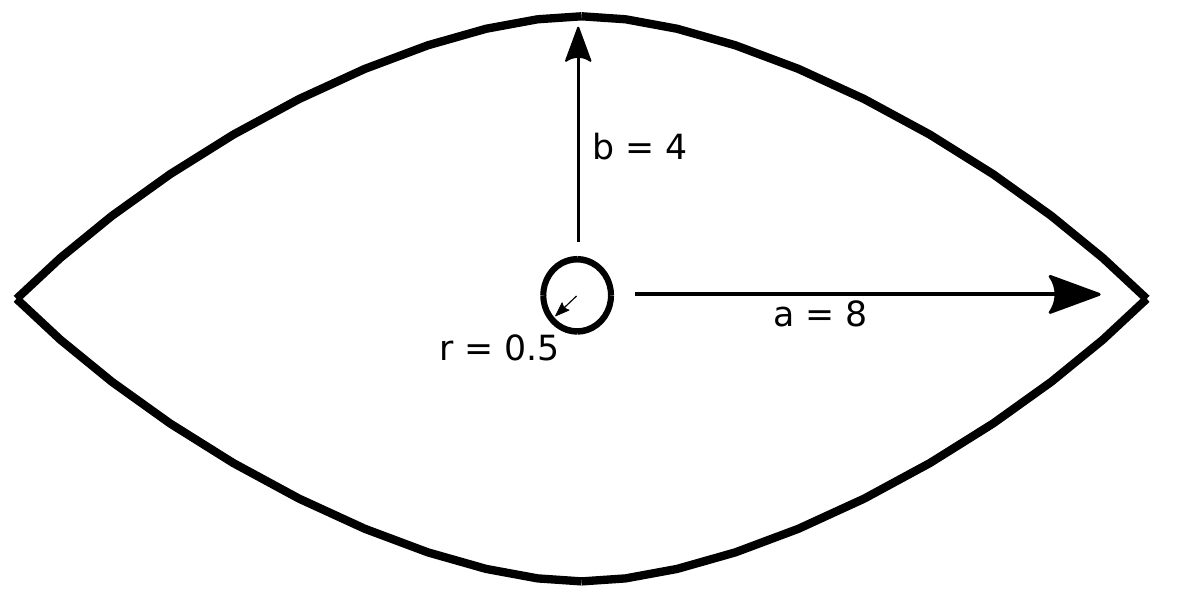}
		\caption{Geometry of structure}
	\end{subfigure}
	\hspace{1cm}
	\begin{subfigure}[t]{0.3\textwidth}
		\centering
		\includegraphics[width=1\columnwidth]{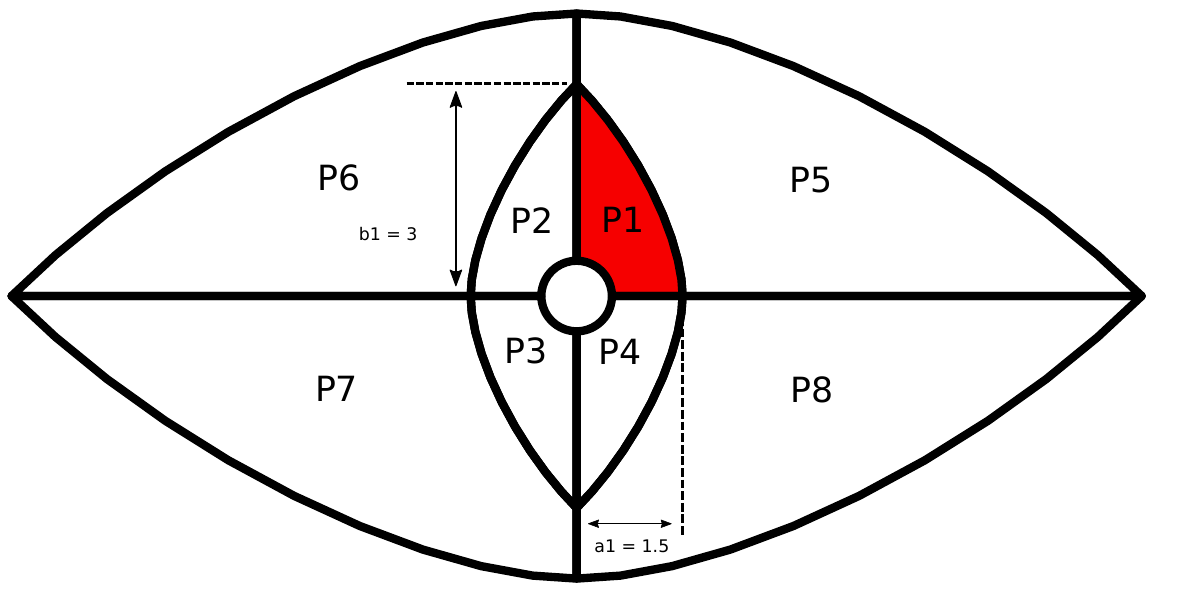}
		\caption{Initial discretization of patches}
	\end{subfigure}
	\caption{The simply supported eye shape with a hole is discretized by 8 patches, with the material parameters that $E_{1} = 0.03, ~\rho_{1} = 0.7.~E_{i} = 1, \rho_{i} = 1, (i=2\ldots8)$. The geometry is represented by NURBS ($ p = 2, ~q = 2$), and the solution field is approximated by PHT splines ($ p = 3, ~q = 3 $).}
	\label{eyeshapeGeo}
\end{figure}

\begin{figure}
	\centering
	\includegraphics[width=0.4\columnwidth]{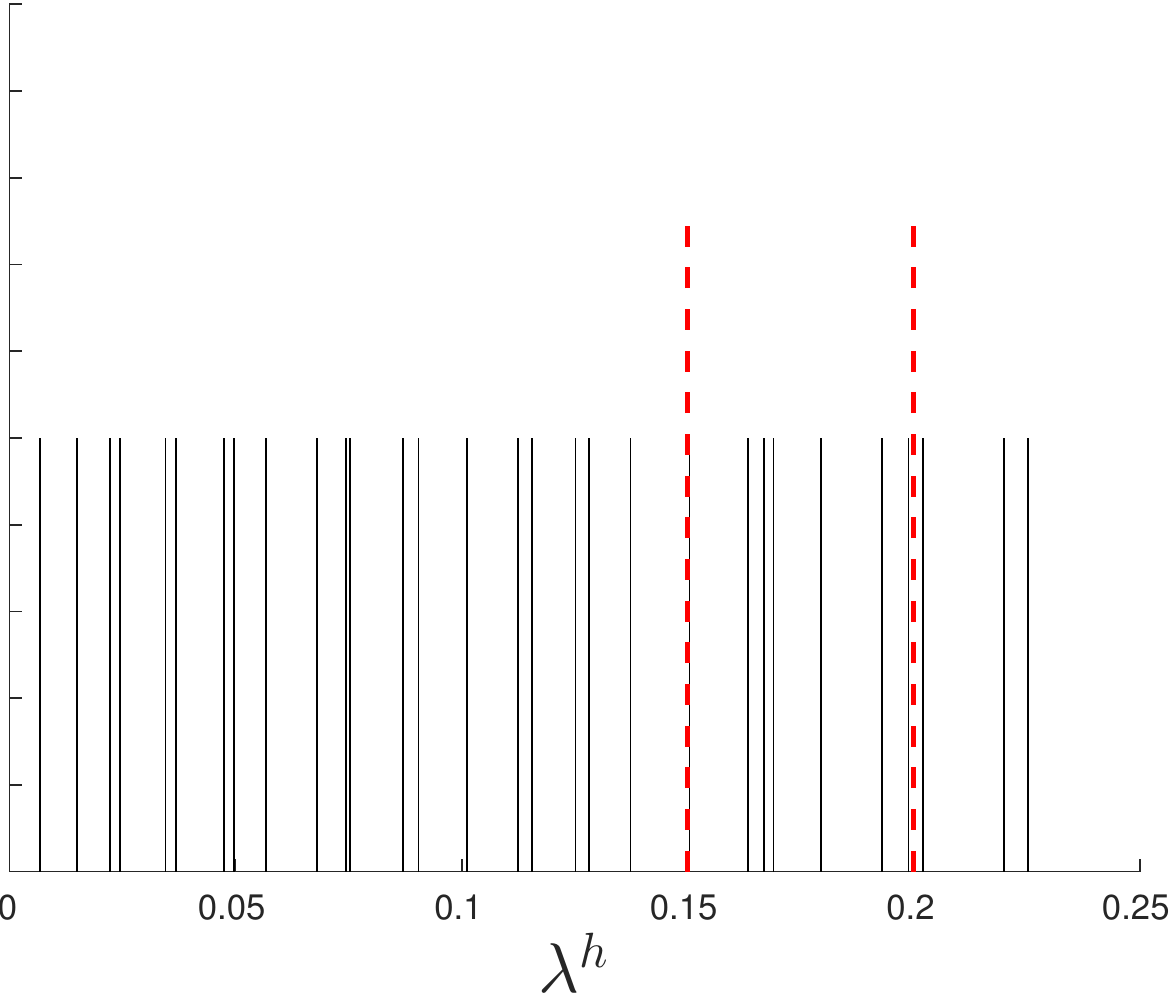}	
	\caption{Frequencies of interest in the window with interval $[0.15,  0.2]$.}
	\label{fig:winsmall}
\end{figure}

The geometry of eye shape with a hole and its material property are presented in Fig.\ref{eyeshapeGeo}. The boundary condition of the outer edge is simply supported, and the edge of the hole is free. This structure is built by 8 patches, where patch 1 with $ E_{1} = 0.03, \rho_{1} = 0.7$ is softer and lighter than other patches with $ E_{i} = 1, \rho_{i} = 1 ~(i = 2\ldots8)$. For all patches, the Poisson's rate is $\nu = 0.3$ and thickness is $h = 0.1$. The frequencies of interest are in the range of $\lambda_{i}^{h}\in[\lambda_{min}, \lambda_{max}]$ (see Fig.\ref{fig:winsmall}, $[\lambda_{min}, \lambda_{max}] = [0.15, 0.20]$, where the range is set by a window marked by red dash lines). Define the $\vect{\mathcal{I}}$, i.e.,
\begin{align}\label{eq:mode_mul}
\vect{\mathcal{I}}(n) = 
\begin{cases}
i, ~~~n=1\\
\{i,\ldots, i+n-1\}, ~~~n>1,
\end{cases}
\end{align}
where $n$ is the multiplicity of mode $i$. The adaptive process is based on the Algorithm \ref{Algm:adpFoI}, but there is a slight difference. In this example, the adaptivity is not conducted by sweeping modes from low to high. Instead, at each step of adaptivity, we pick up the mode(s) $\vect{\mathcal{I}}$ with the maximum of $\delta_{\vect{\mathcal{I}}}^{\vect{\phi}}$ to deliver the adaptivity, where $\delta_{\vect{\mathcal{I}}}$ reads
\begin{align}\label{eq:errvect_mul}
\delta_{\vect{\mathcal{I}}(n)}^{\vect{\phi}} = 
\begin{cases}
\delta_{i}^{\vect{\phi}} ~~~\text{in Eq.\eqref{eq:errestmor}}, ~~~n=1\\
\delta_{\vect{\phi}} ~~~\text{in Eq.\eqref{eq:errestmmul}}, ~~~n>1.
\end{cases}
\end{align}
Accordingly, the error estimator of frequencies is referred as
\begin{align}\label{eq:errfreq_mul}
|e_{\vect{\mathcal{I}}(n)}^{\lambda}| = 
\begin{cases}
|e_{i}^{\lambda}| ~~~\text{in Eq.\eqref{eq:errestmor}}, ~~~n=1\\
|e_{\lambda}| ~~~\text{in Eq.\eqref{eq:errestmmul}}, ~~~n>1.
\end{cases}
\end{align}
 
\begin{figure}
	\centering		
	\begin{subfigure}[t]{0.3\textwidth}
		\centering
		\includegraphics[width=1\columnwidth]{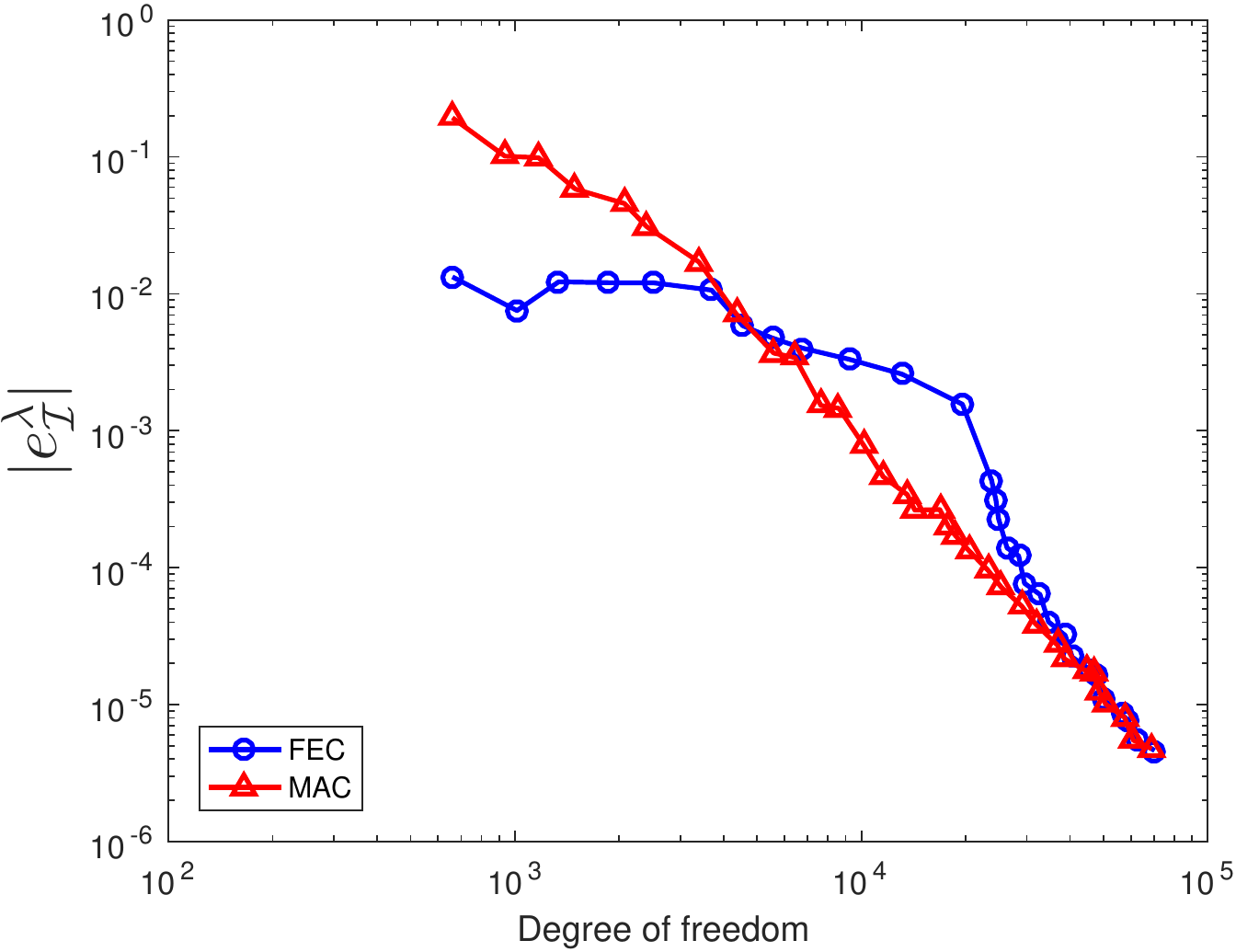}
		\caption{}
	\end{subfigure}
	\hspace{0.2cm}
	\begin{subfigure}[t]{0.3\textwidth}
		\centering
		\includegraphics[width=1\columnwidth]{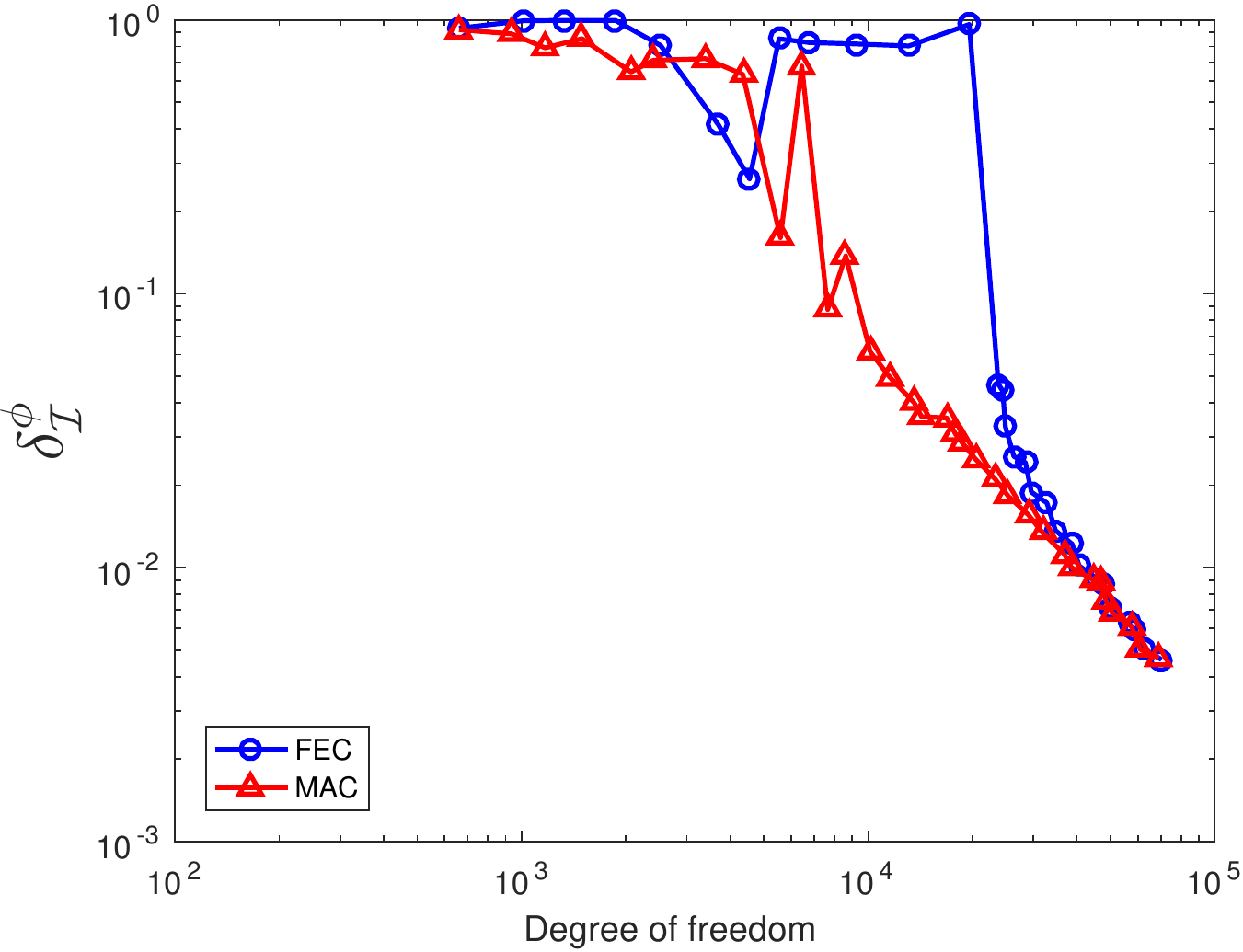}
		\caption{}
	\end{subfigure}
	\hspace{0.2cm}
	\begin{subfigure}[t]{0.31\textwidth}
		\centering
		\includegraphics[width=1\columnwidth]{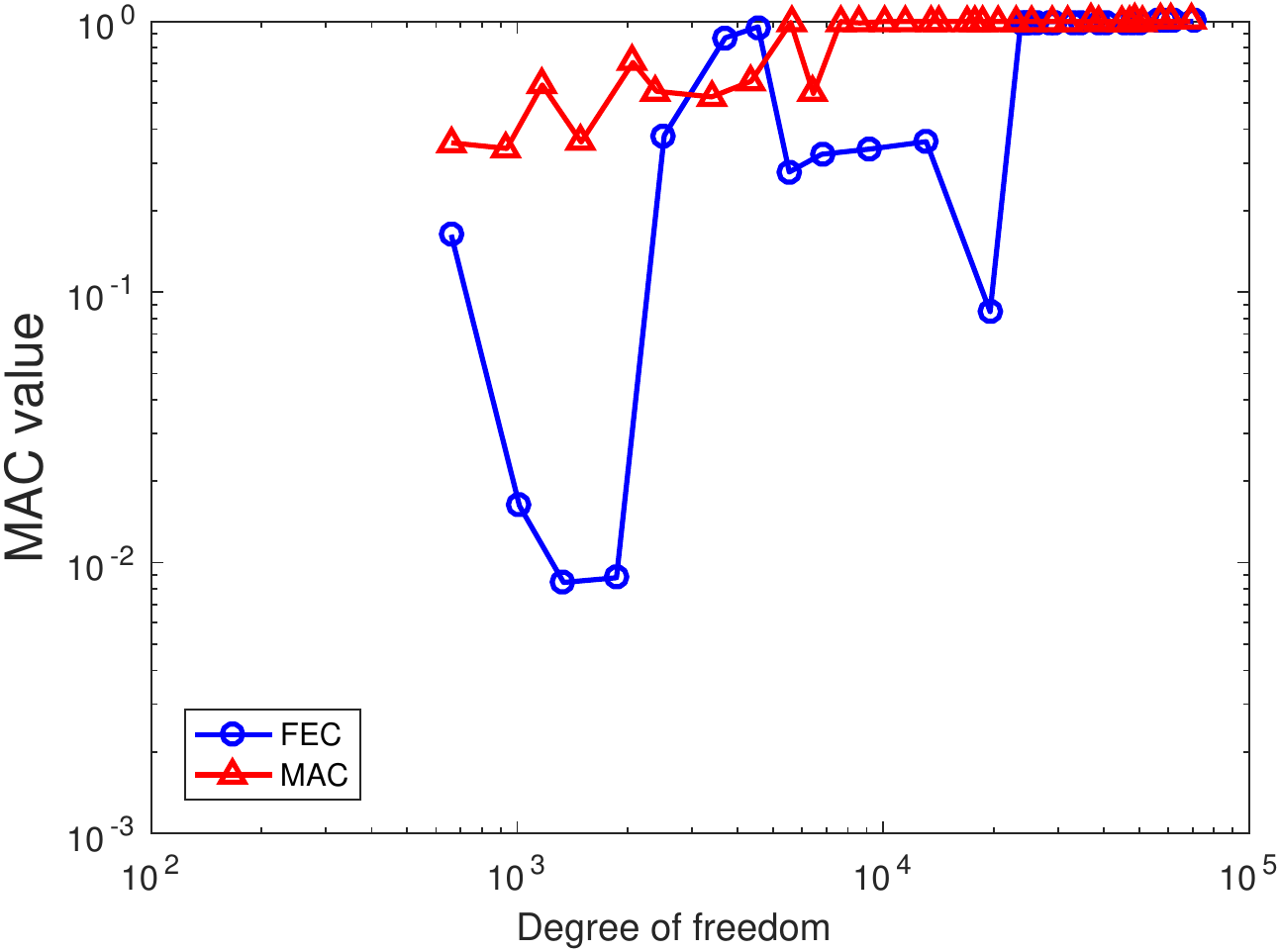}
		\caption{}
	\end{subfigure}
	\caption{Comparisons between MAC and FEC methods in a frequency range [0.15,0.2] in respect of (a) $|e_{\vect{\mathcal{I}}}^{\lambda}|$, (b) $ \delta_{\vect{\mathcal{I}}}^{\phi}$, (c) MAC value. Refinement level $L_{e}$ is chosen to be 1.}
	\label{fig:MACvsFECsmall}		
\end{figure}

In this case, two schemes, FEC and MAC, proposed in Section \ref{sec:Locmode} are compared in order to investigate theirs effects on the adaptivity. As it can be seen in Fig.\ref{fig:MACvsFECsmall}(a),(b), it is clear that MAC method has a better convergence than FEC method. The reason can be obtained by tracing the results in Tab.\ref{tab:FECvmode} and Tab.\ref{tab:MACvmode}. To be specific, the local adaptivity is driven by the error estimation of mode shapes, FEC can not always guarantee to locate the right corresponding modal vector on refined mesh (see Tab.\ref{tab:FECvmode} that $\vect{\phi}_{\vect{\mathcal{I}}}^{h}$ and $\tilde{\vect{\phi}}_{\tilde{\vect{\mathcal{I}}}}$ are not consistent until adaptive step 13). Here, $\vect{\phi}_{\vect{\mathcal{I}}}$ is expressed as
\begin{align}
\vect{\phi}_{\vect{\mathcal{I}}(n)} = 
\begin{cases}
\vect{\phi}_{i} ~~~\text{eigenvector}, ~~~n=1\\
\vect{\varphi} ~~~\text{vector in Eq.\eqref{eq:eigenspace}}, ~~~n>1.
\end{cases}
\end{align}
Therefore, It leads to the inefficient adaptive mesh, and furthermore causes that the error estimator $\delta_{\vect{\mathcal{I}}}^{\vect{\phi}}$ is divergent, as shown in Fig.\ref{fig:MACvsFECsmall}(b). Only when FEC scheme is able to identify the related mode correctly (from the step 13 and forwards in Tab.\ref{tab:FECvmode}), the error estimators are just convergent accordingly. In contrast, MAC method can find the associated mode accurately at very early stage of adaptivity (at around 6th step displayed in Tab.\ref{tab:MACvmode} and Fig.\ref{fig:MACvsFECsmall}(c)). Therefore, for the given accuracy such as $\left| e_{\vect{\mathcal{I}}}^{\lambda}\right| \leqslant 10^{-4}$ and $\delta_{\vect{\mathcal{I}}}^{\vect{\phi}} \leqslant 10^{-2}$, MAC method is more efficient than FEC method.

\begin{figure}
	\centering
	\includegraphics[width=0.4\columnwidth]{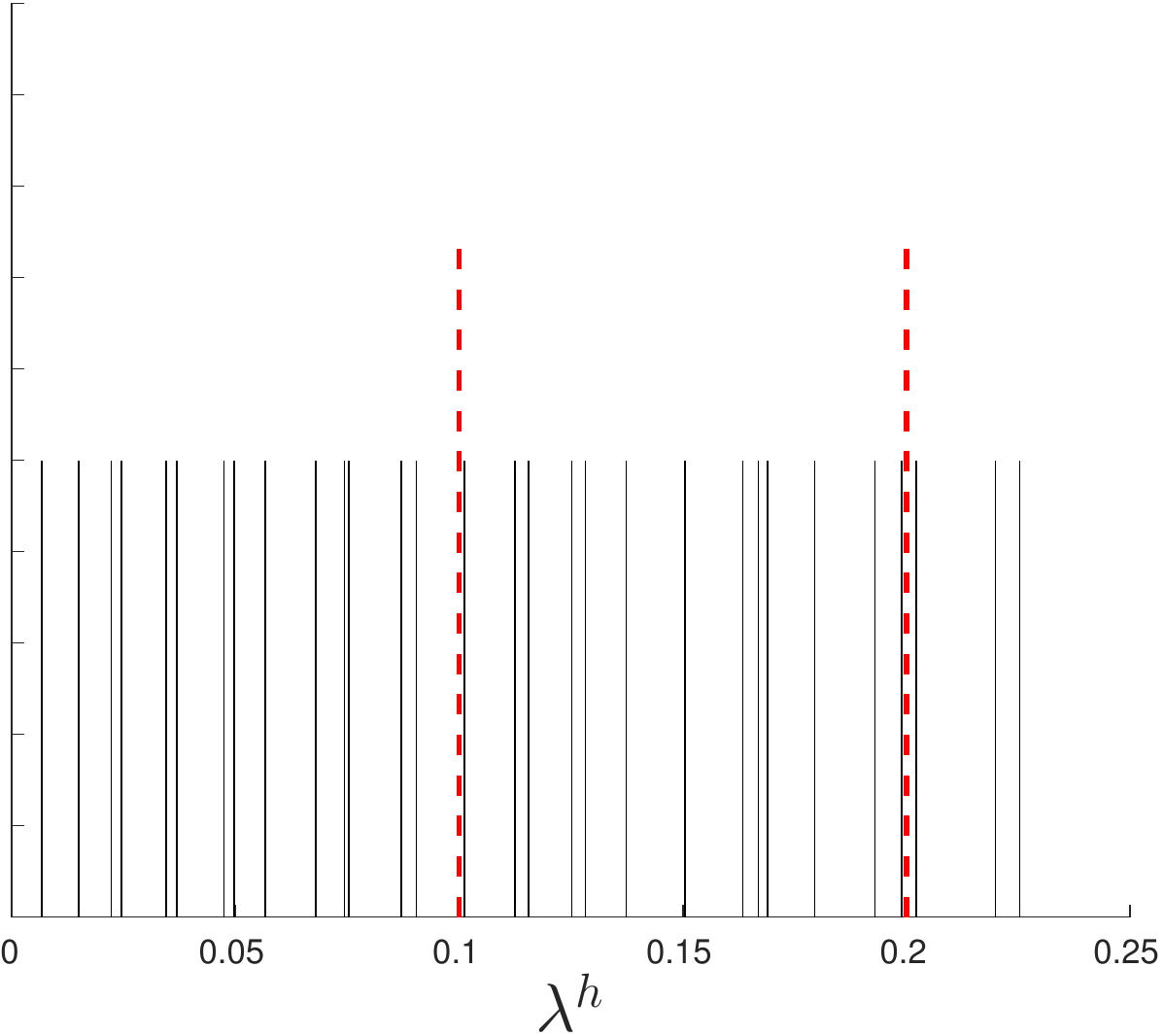}
	\caption{Frequencies of interest in the window with interval $[0.1, 0.2]$.}
	\label{fig:winlarge}
\end{figure}

\begin{figure}
	\centering			
	\begin{subfigure}[t]{0.32\textwidth}
		\centering
		\includegraphics[width=1\columnwidth]{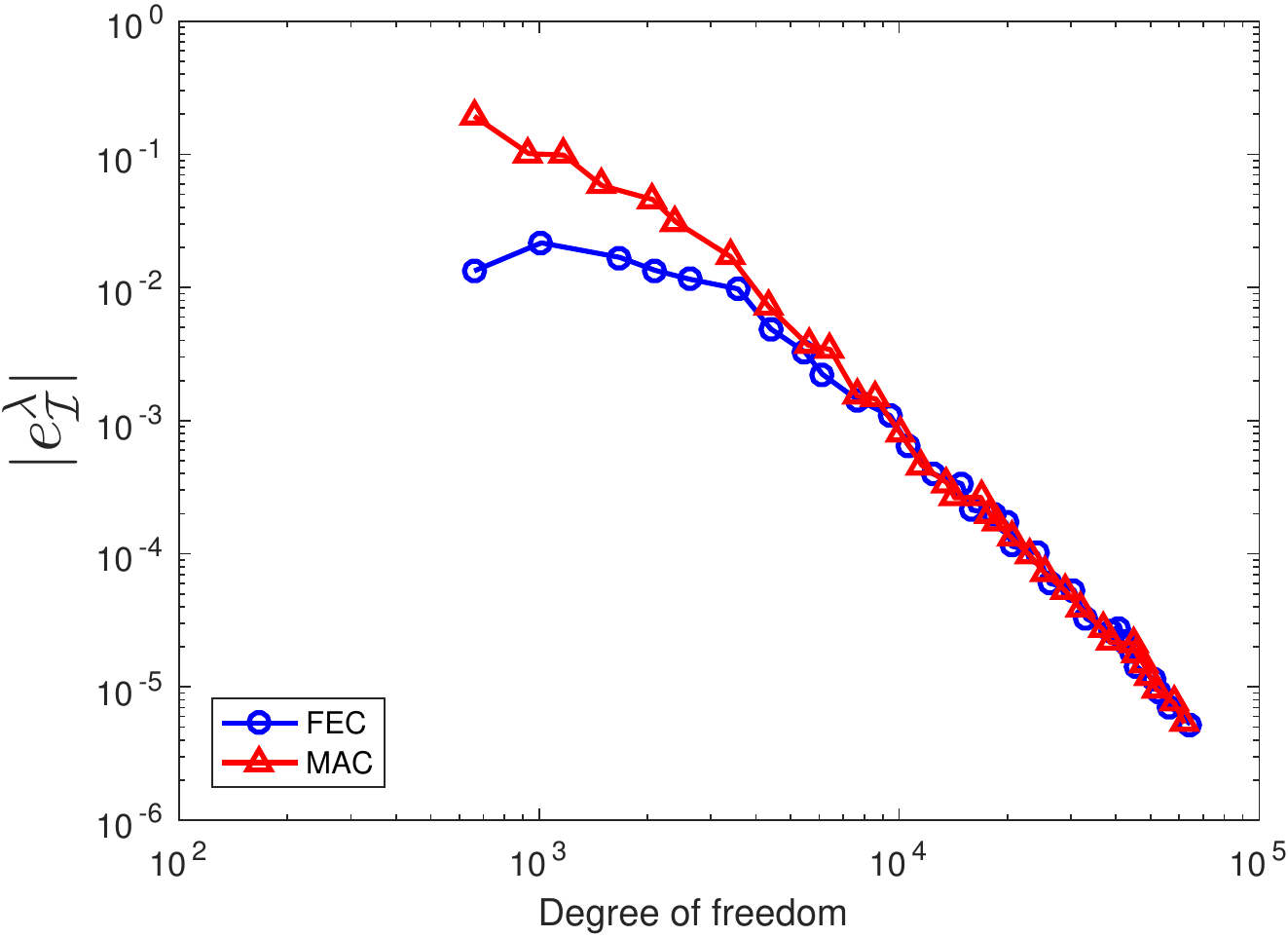}
		\caption{}
	\end{subfigure}
	\hspace{0.2cm}
	\begin{subfigure}[t]{0.3\textwidth}
		\centering
		\includegraphics[width=1\columnwidth]{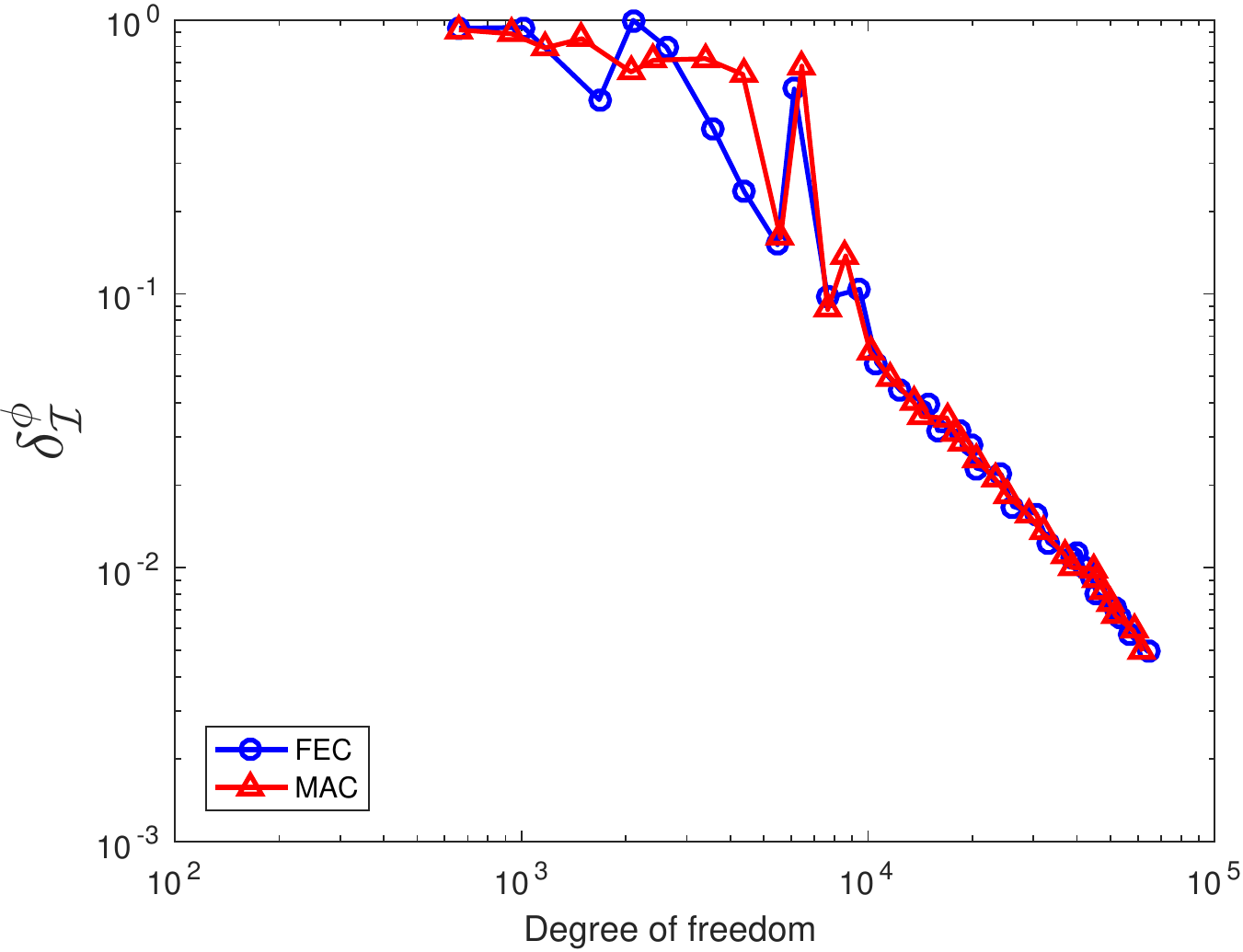}
		\caption{}
	\end{subfigure}
	\hspace{0.2cm}
	\begin{subfigure}[t]{0.31\textwidth}
		\centering
		\includegraphics[width=1\columnwidth]{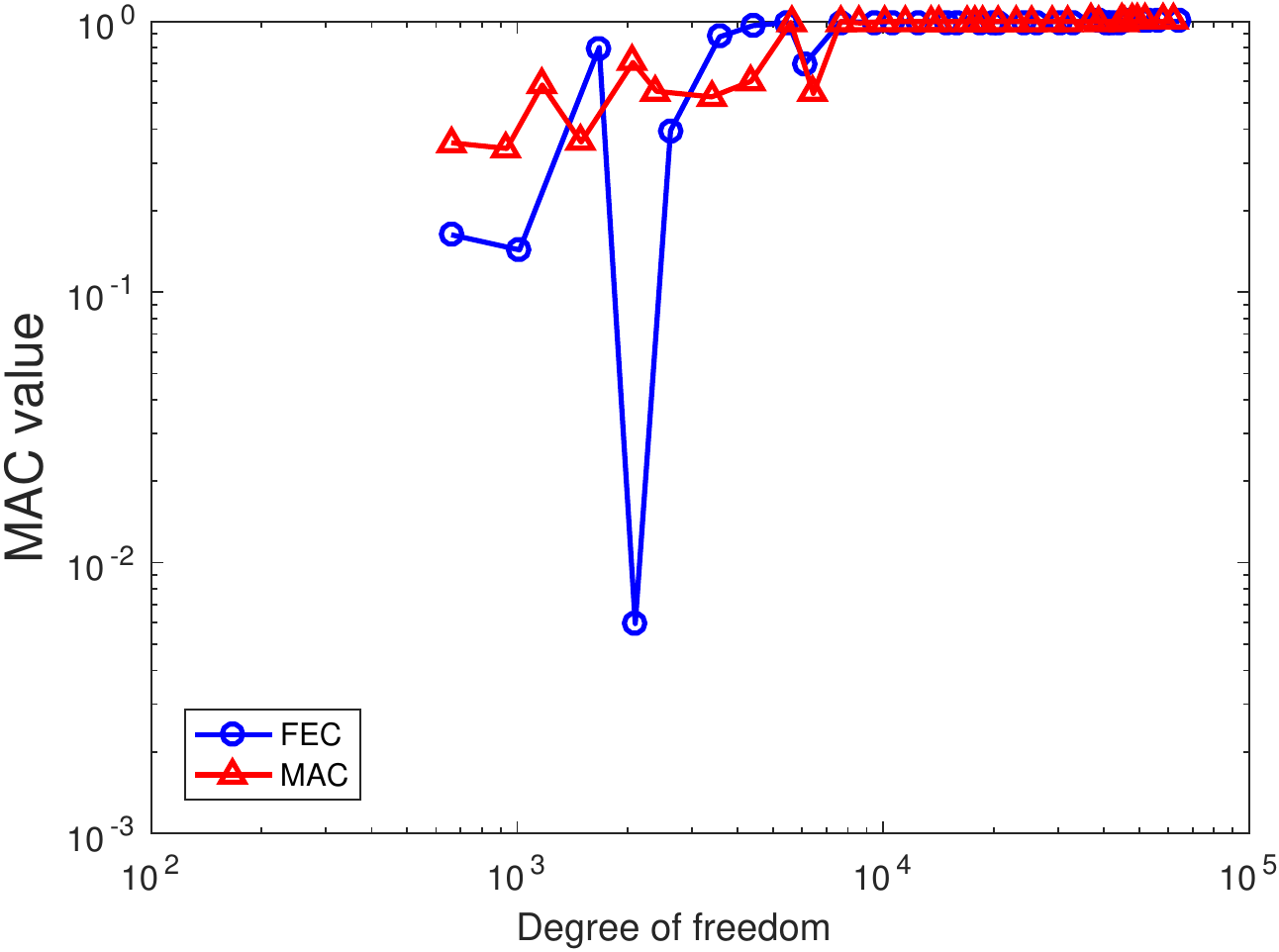}
		\caption{}
	\end{subfigure}
	\caption{Comparisons between MAC and FEC methods in a frequency band [0.1,0.2] in terms of (a) $|e_{\vect{\mathcal{I}}}^{\lambda}|$, (b) $ \delta_{\vect{\mathcal{I}}}^{\phi}$, (c) MAC value.}
	\label{fig:MACvsFEClarge}		
\end{figure}

As the window expands to cover $\lambda_{i}^{h} \in [0.1,0.2]$ in Fig.\ref{fig:winlarge}, the advantage of using MAC is not so noticeable that good convergence is achieved by both MAC and FEC, as presented in Fig.\ref{fig:MACvsFEClarge}(a),(b). This is because that the modal shapes at low frequency modes are often distinct. If the adaptivity starts from low frequency modes, it is easy to locate the corresponding mode even by FEM scheme. Thus, the precisely adaptive refinement in low modes will help to accurately locate modal correspondence for high frequency modes. Regardless of that, the MAC will be our preference to the remaining computations as it will not make mistakes on modal resemblance recognition in any case. Simultaneously, it is cheap to implement and execute.

\begin{figure}
	\centering
	\begin{subfigure}[t]{0.3\textwidth}
		\centering
		\includegraphics[width=1\columnwidth]{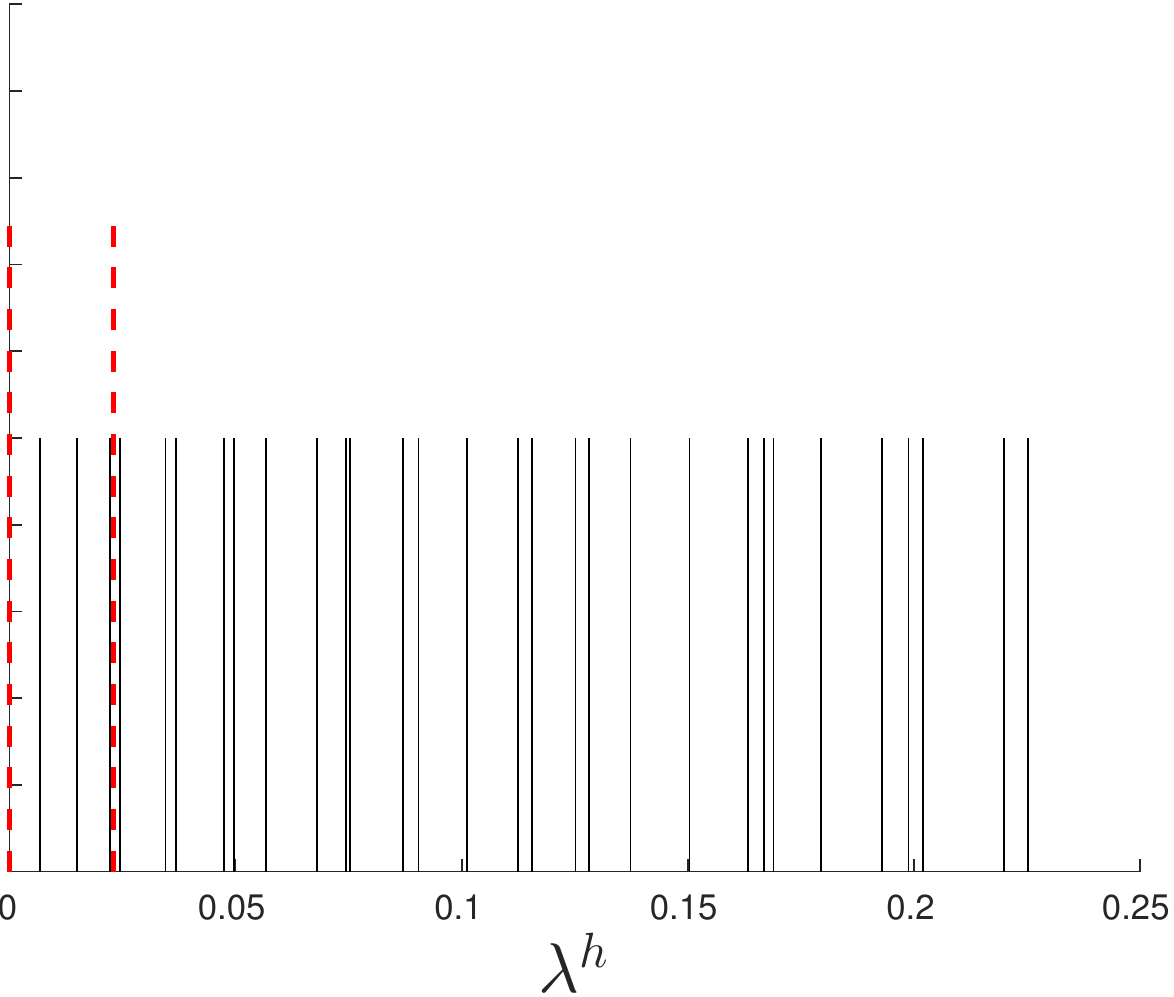}
		\caption{[0, 0.023]}
	\end{subfigure}
	\hspace{0.2cm}
	\begin{subfigure}[t]{0.3\textwidth}
		\centering
		\includegraphics[width=1\columnwidth]{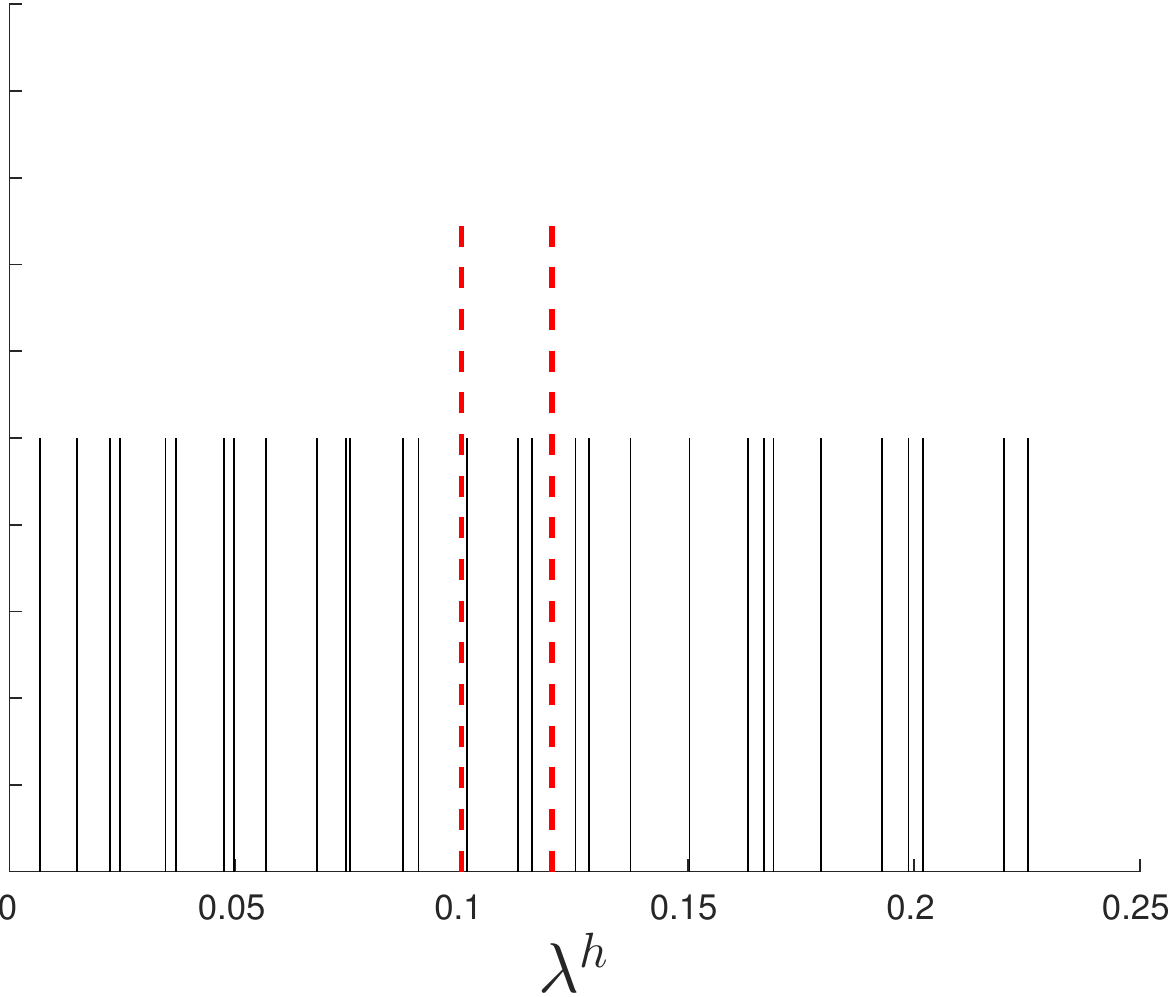}
		\caption{[0.1, 0.12]}
	\end{subfigure}
	\hspace{0.2cm}
	\begin{subfigure}[t]{0.3\textwidth}
		\centering
		\includegraphics[width=1\columnwidth]{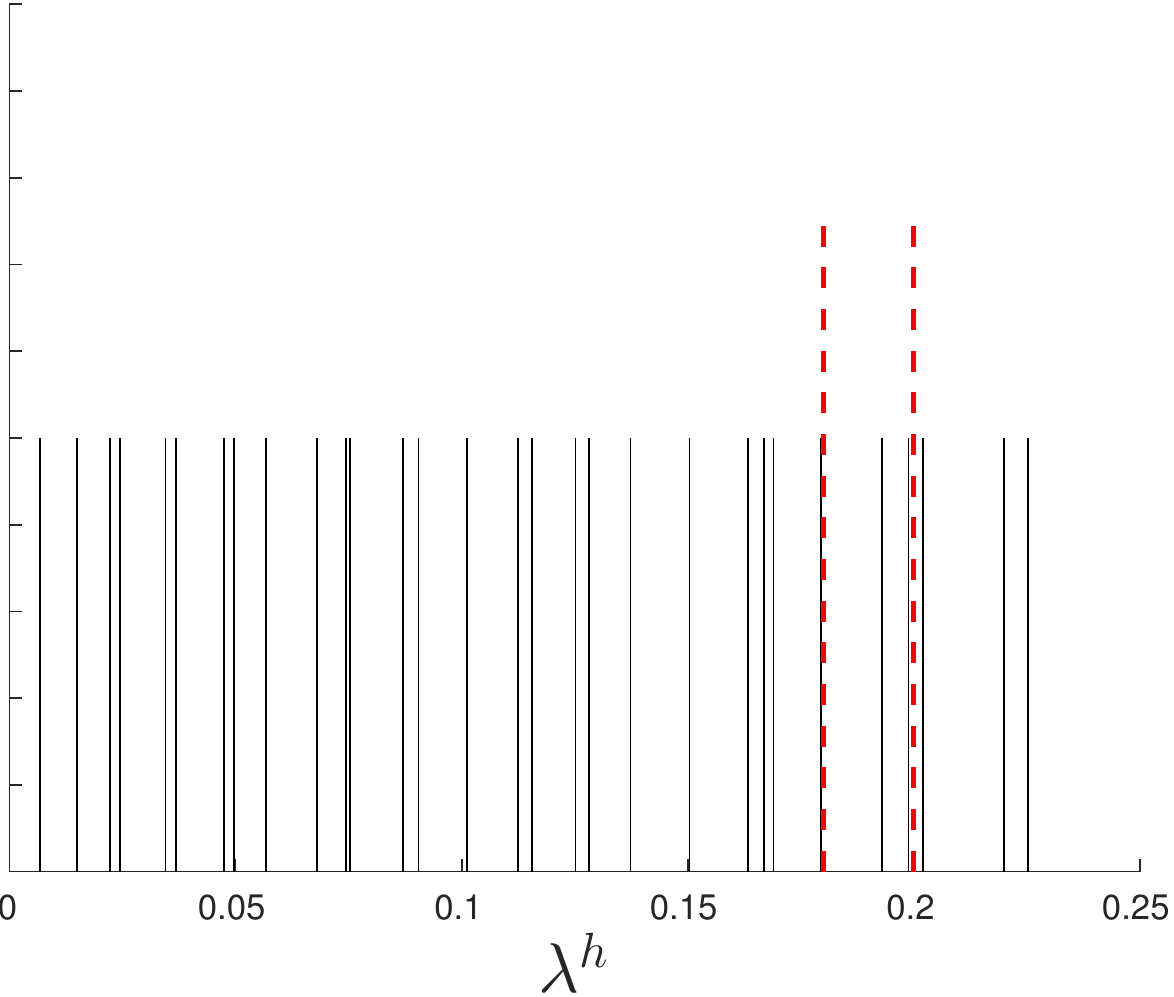}
		\caption{[0.18 0.2]}
	\end{subfigure}\\
	\begin{subfigure}[t]{0.3\textwidth}
		\centering
		\includegraphics[width=1\columnwidth]{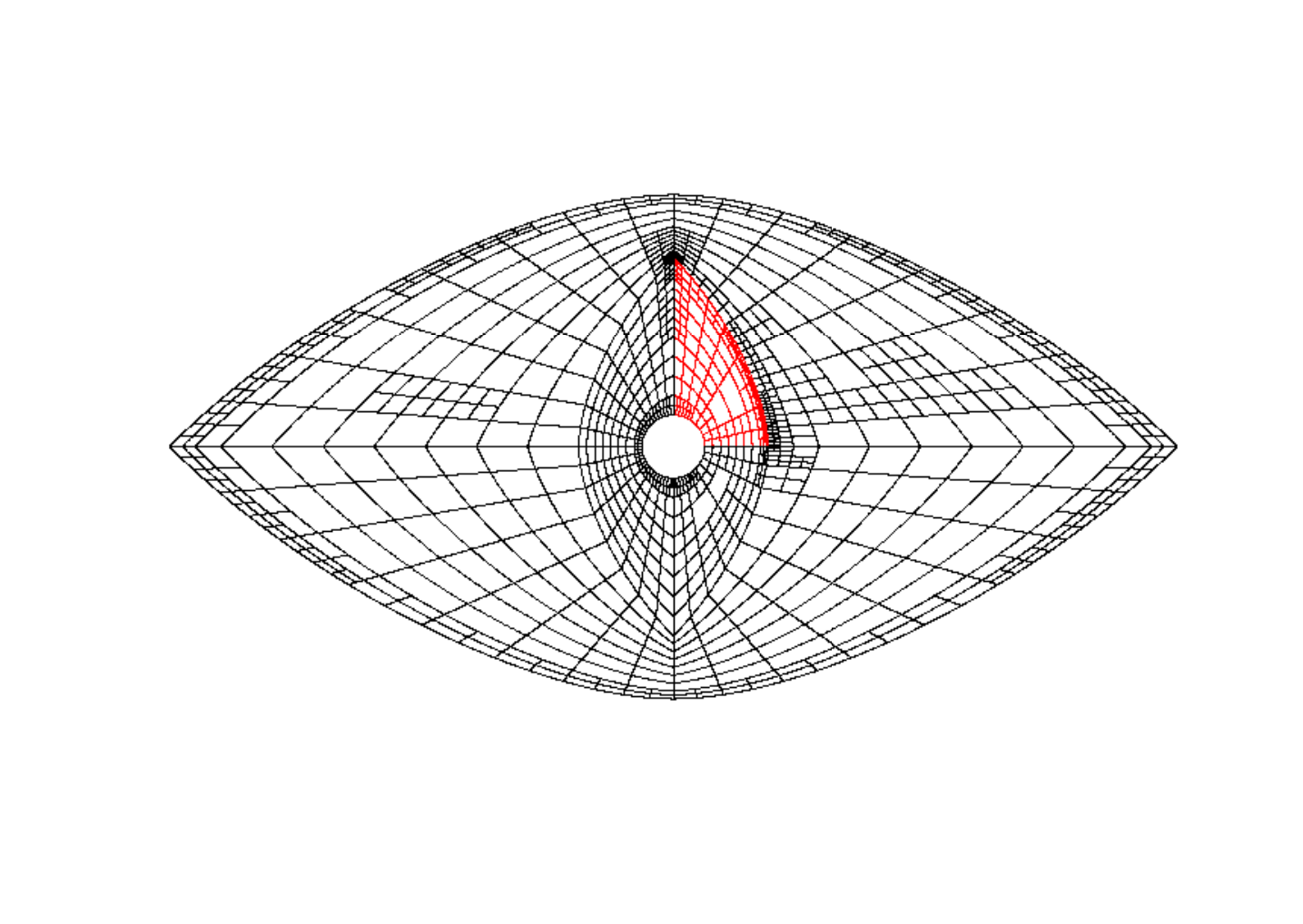}
		\caption{Refinement at $[0, 0.02]$, with $\left| e_{\vect{\mathcal{I}}}^{\lambda} \right| = 2 \times 10^{-5}, ~\delta_{\vect{\mathcal{I}}}^{\vect{\phi}} = 0.0096$}
	\end{subfigure}
	\hspace{0.2cm}
	\begin{subfigure}[t]{0.3\textwidth}
		\centering
		\includegraphics[width=1\columnwidth]{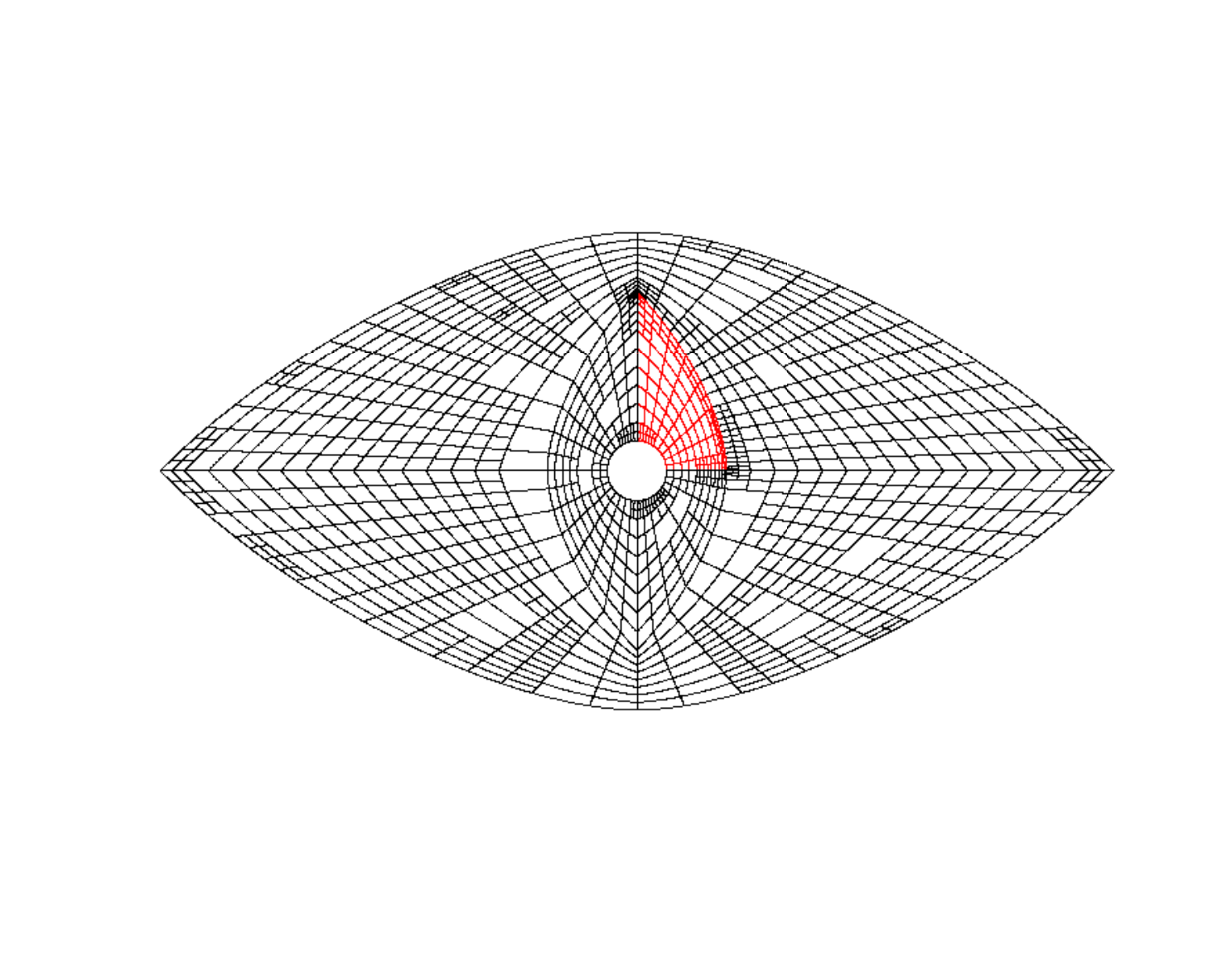}
		\caption{Refinement at [0.1, 0.12], with $\left| e_{\vect{\mathcal{I}}}^{\lambda} \right| = 1.2 \times 10^{-4}, ~\delta_{\vect{\mathcal{I}}}^{\vect{\phi}} = 0.024$}
	\end{subfigure}
	\hspace{0.2cm}
	\begin{subfigure}[t]{0.3\textwidth}
		\centering
		\includegraphics[width=1\columnwidth]{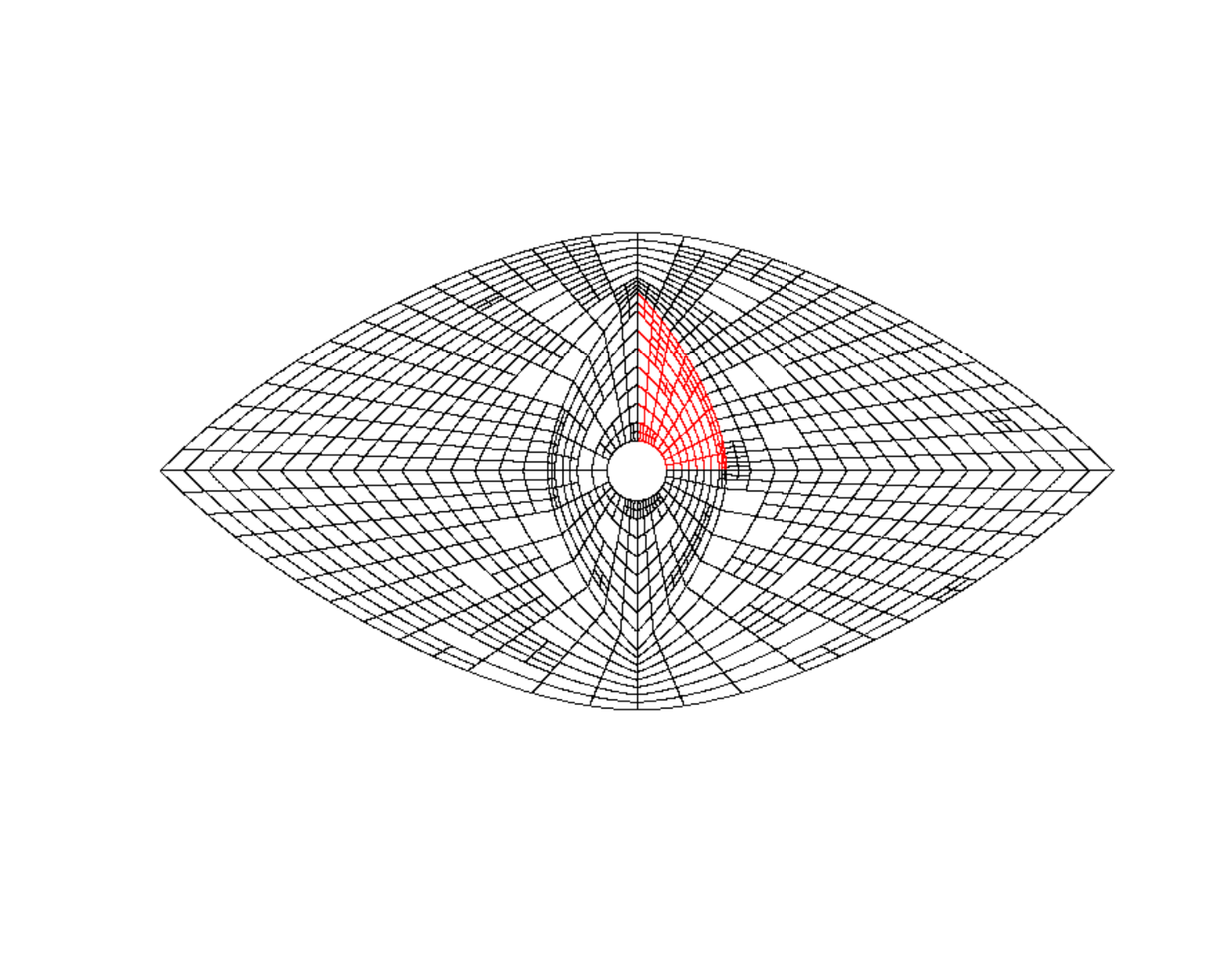}
		\caption{Refinement at [0.18 0.20], with $\left| e_{\vect{\mathcal{I}}}^{\lambda} \right| = 3.2 \times 10^{-4}, ~\delta_{\vect{\mathcal{I}}}^{\vect{\phi}} = 0.042$}
	\end{subfigure}		
	\caption{The Comparison of adaptive refinements among different intervals of frequencies of interest (a)-(c) based on the MAC method. The numbers of elements for the meshes in (d),(e) and (f) are all 1150.}
	\label{fig:compmesh3}		
\end{figure}

Note that the final meshes in both Tab.\ref{tab:FECvmode} (at step 28) and Tab.\ref{tab:MACvmode} (at step 24) are close to uniform meshes. This is because the modes of interest are with high frequencies (see Fig.\ref{fig:MACvsFECsmall}), the structural vibrations are normally global. If modes of interest are low (as in Fig.\ref{fig:compmesh3}(a)), the refinement will localize around the patch with soft material and the hole (see Fig.\ref{fig:compmesh3}(d)). While as the modes of interest become higher (as shown in Fig.\ref{fig:compmesh3}(b),(c)), with the same number of elements, the adaptive refinements will get closer to uniform refinements, as in Fig.\ref{fig:compmesh3}(e),(f), and the error estimators $\left| e_{\vect{\mathcal{I}}}^{\lambda}\right| $ and $\delta_{\vect{\mathcal{I}}}^{\vect{\phi}} $ are larger.

\begin{table}
	\caption{Targeted mode shapes $\vect{\phi}_{\vect{\mathcal{I}}}^{h}$ at coarse mesh and related mode shapes $\tilde{\vect{\phi}}_{\tilde{\vect{\mathcal{I}}}}$ over refined mesh, and the adaptive refinement at different steps obtained through FEC.}
	\label{tab:FECvmode}
	\centering
	\begin{tabular}{cccc}\hline
		\textbf{Step} & $\vect{\phi}_{\vect{\mathcal{I}}}^{h}$ & $\tilde{\vect{\phi}}_{\tilde{\vect{\mathcal{I}}}}$ & \textbf{Adaptive mesh}\\
		\raisebox{1cm}{2} & \includegraphics[width=0.2\columnwidth]{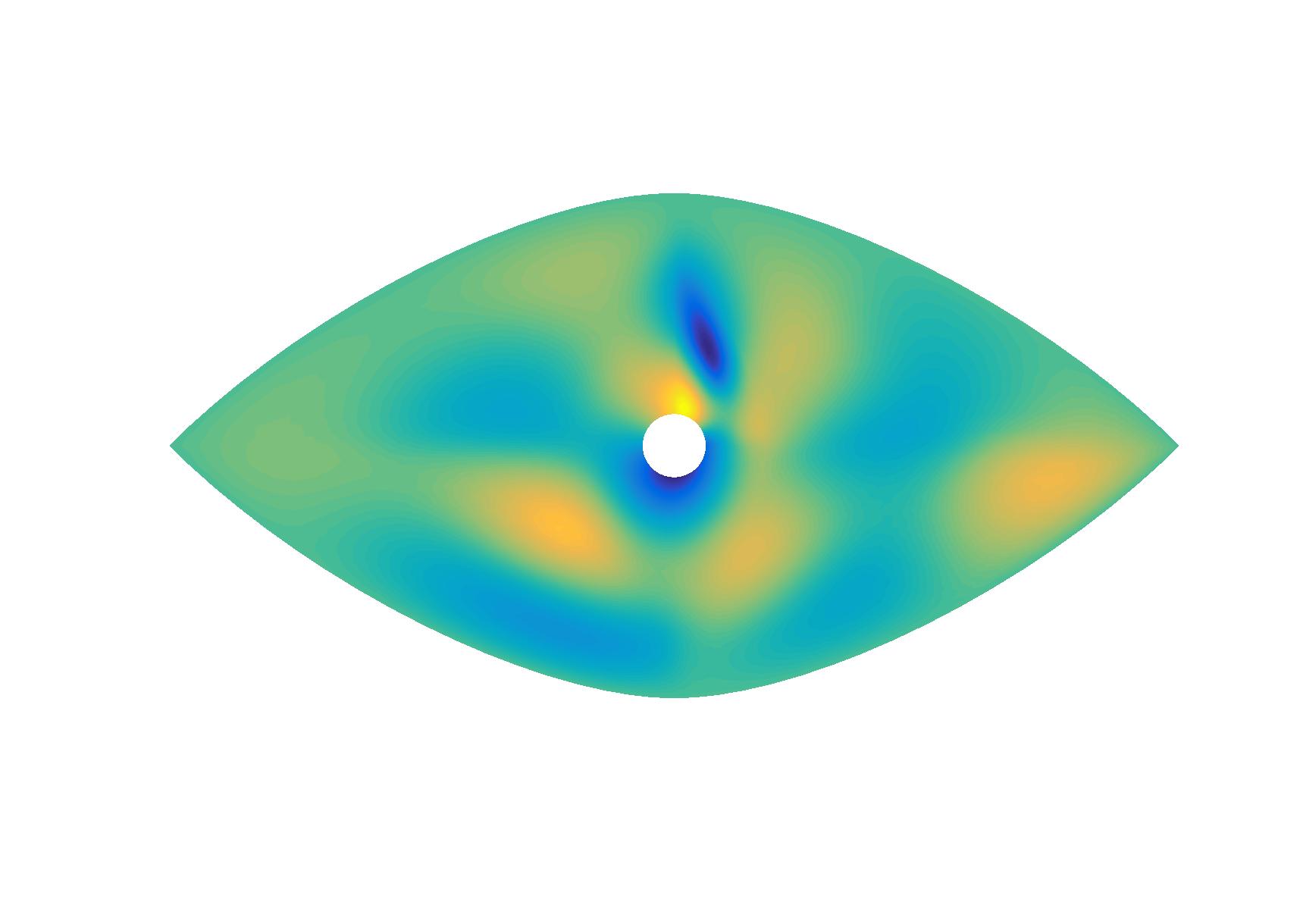}&  \includegraphics[width=0.2\columnwidth]{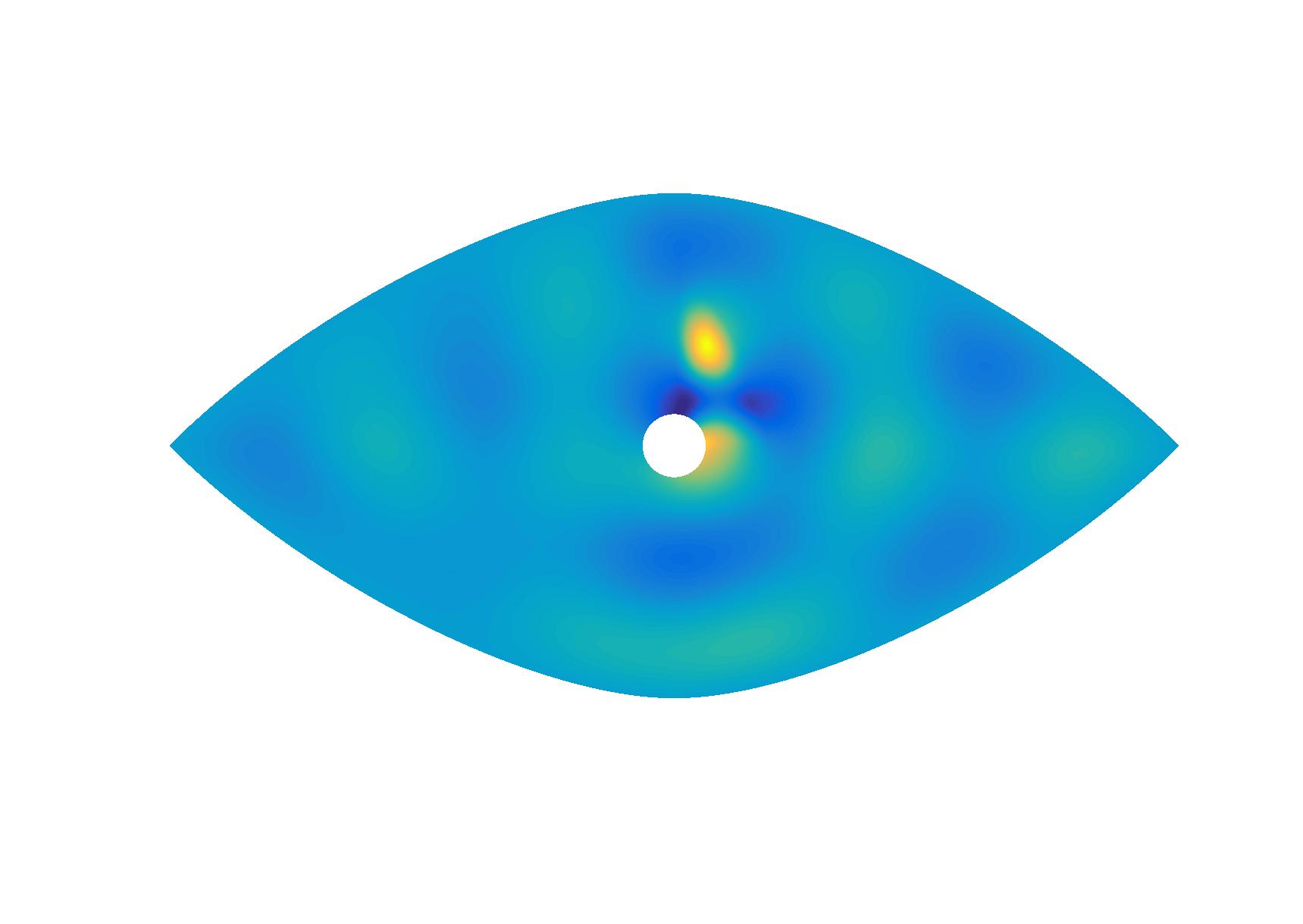} & \raisebox{0.55cm}{\includegraphics[width=0.17\columnwidth]{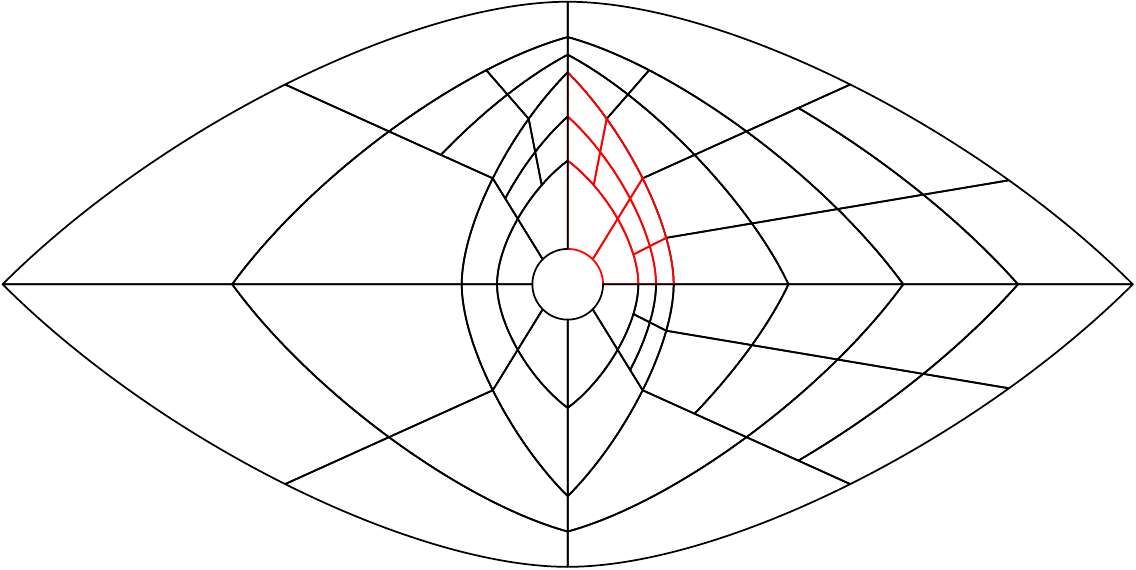}}\\
		\raisebox{1cm}{4} & \includegraphics[width=0.2\columnwidth]{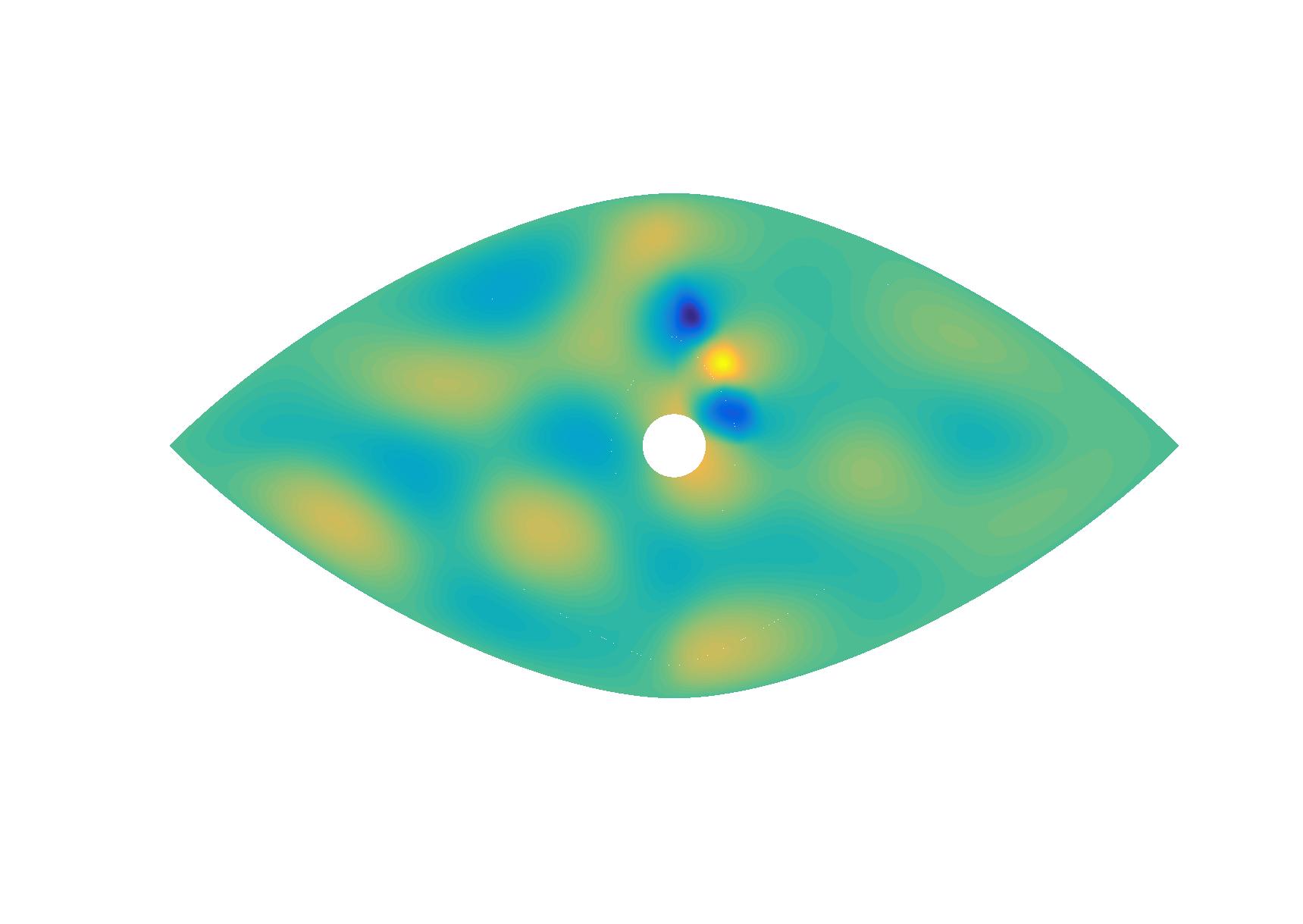} & \includegraphics[width=0.2\columnwidth]{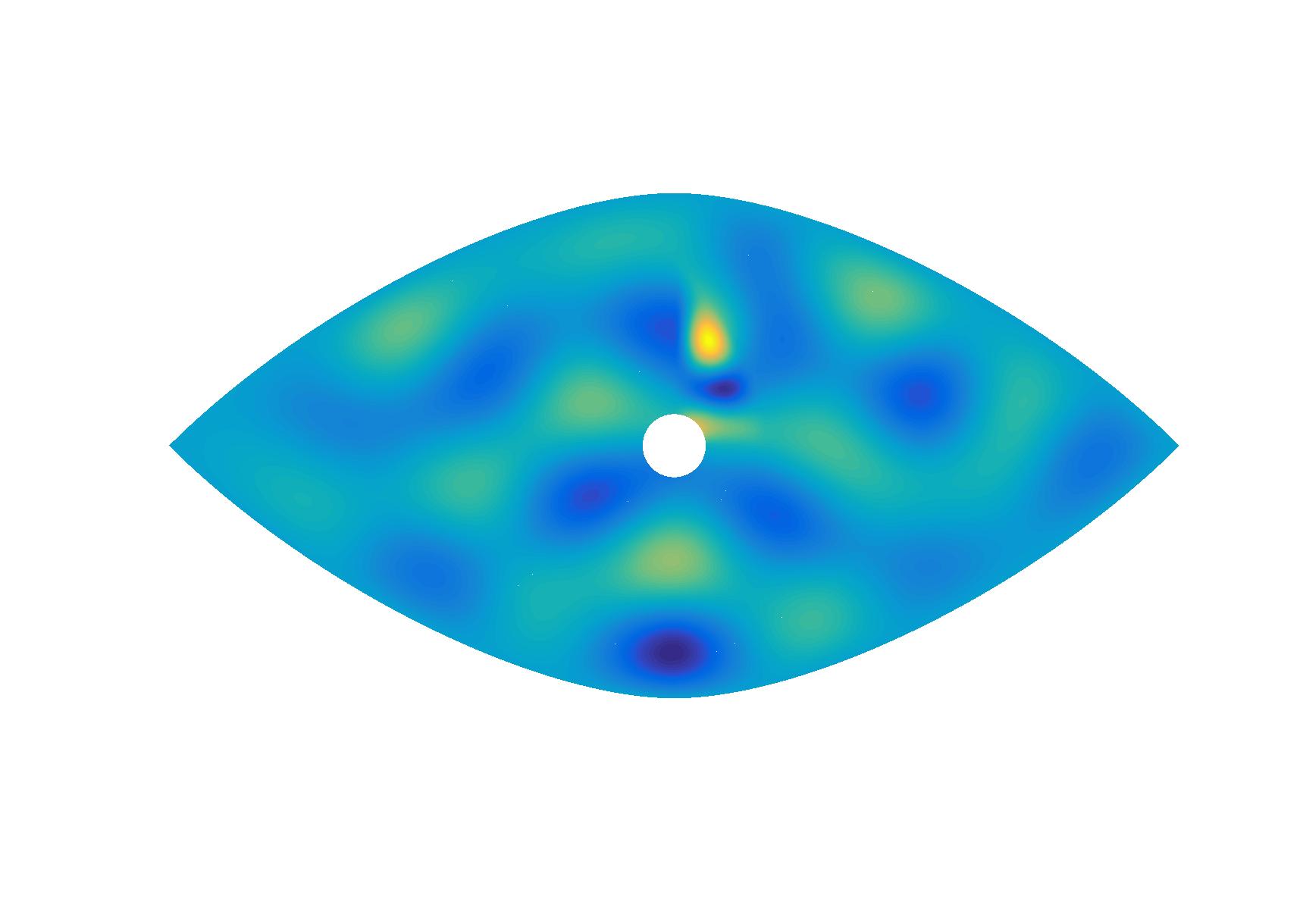} & \raisebox{0.55cm}{\includegraphics[width=0.17\columnwidth]{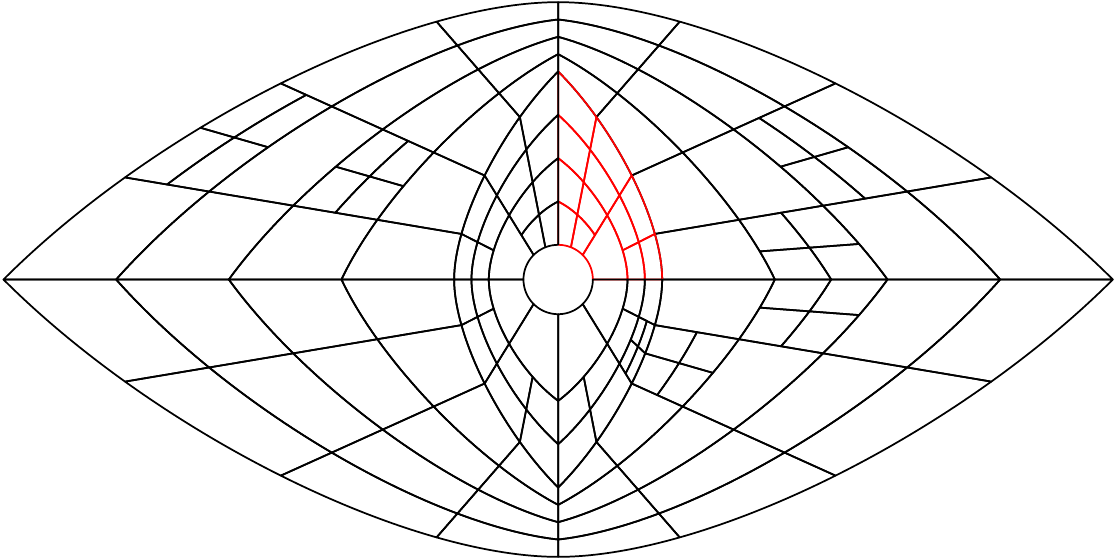}} \\
		\raisebox{1cm}{6} & \includegraphics[width=0.2\columnwidth]{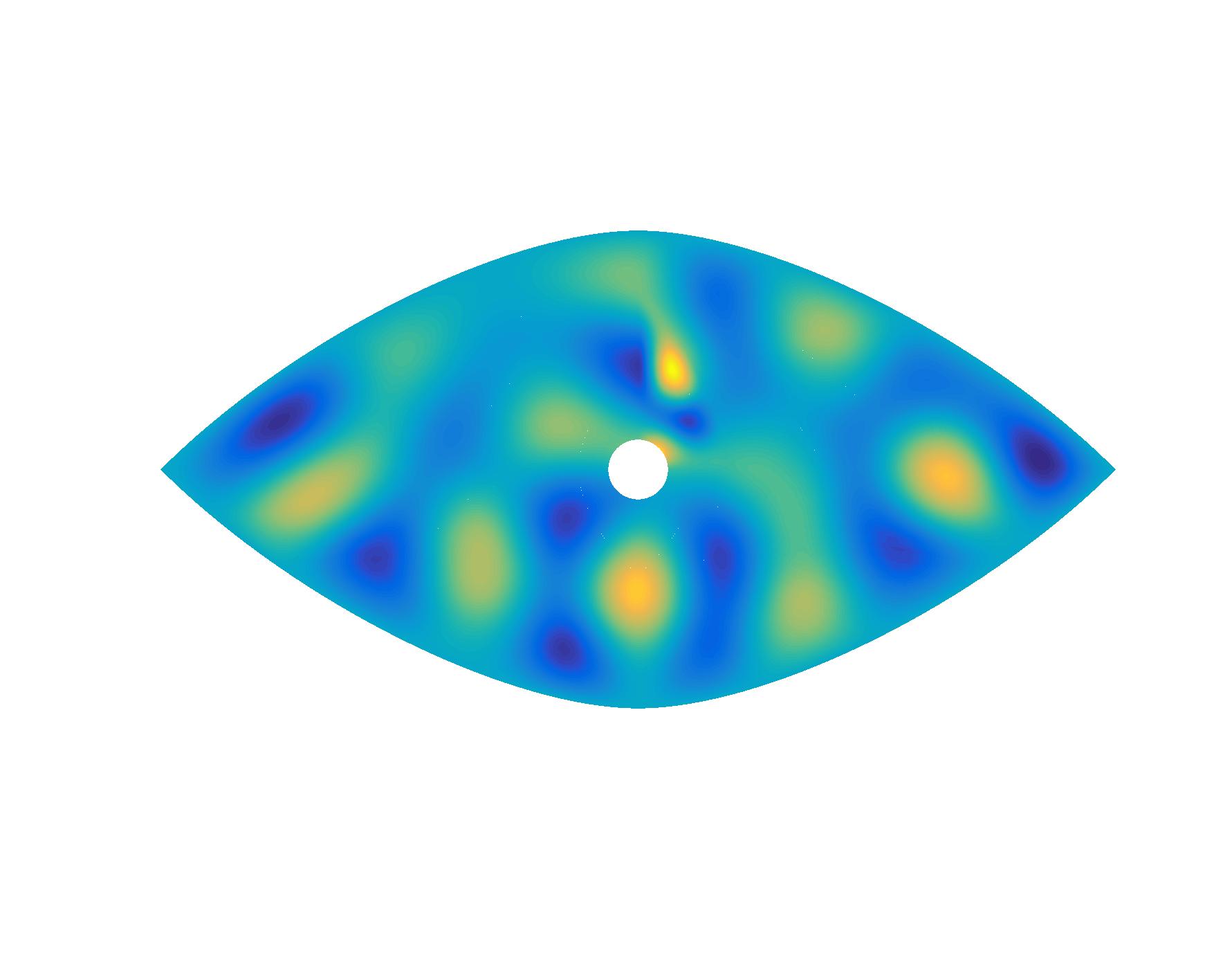} & \includegraphics[width=0.2\columnwidth]{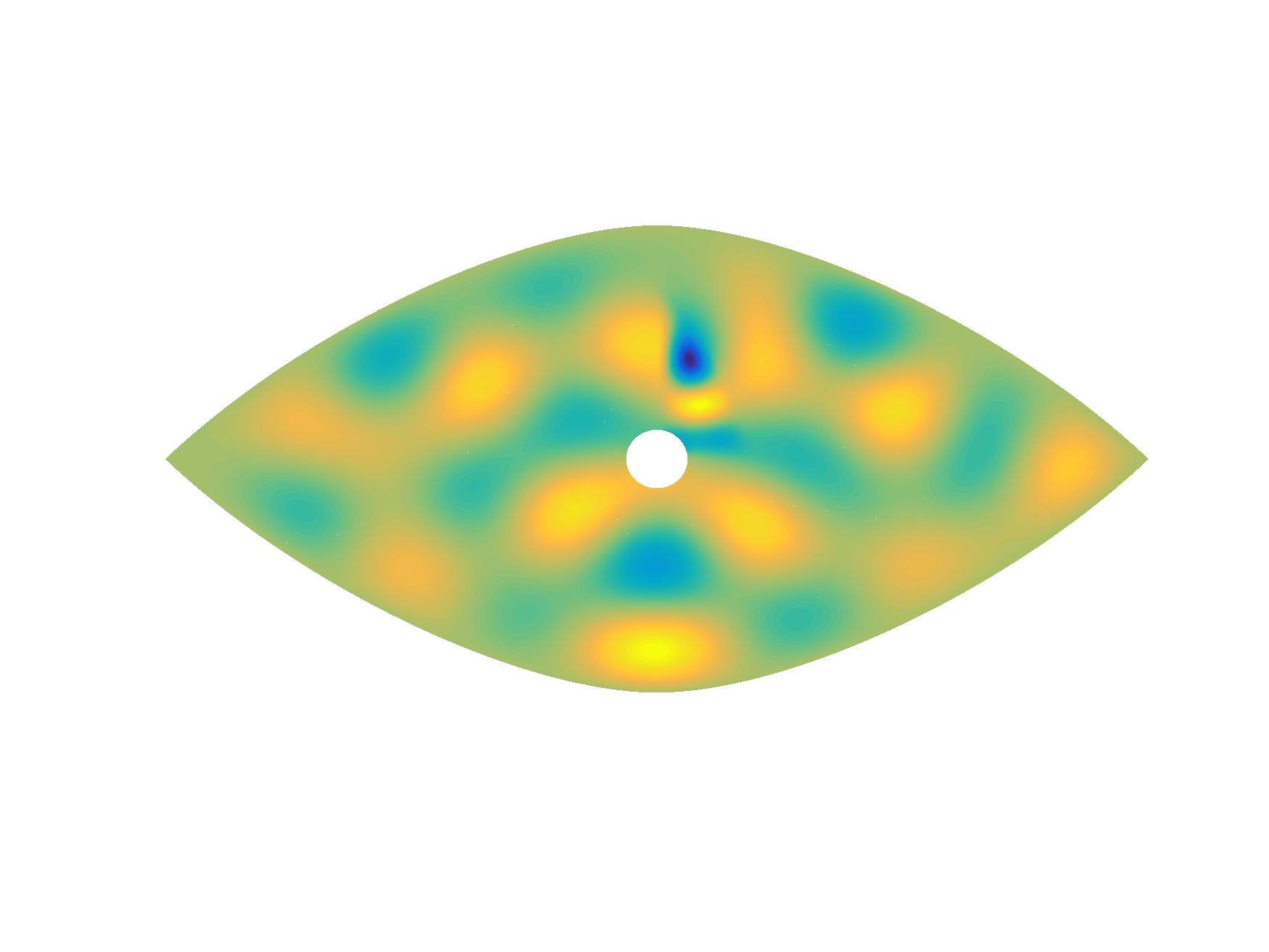} & \raisebox{0.55cm}{\includegraphics[width=0.17\columnwidth]{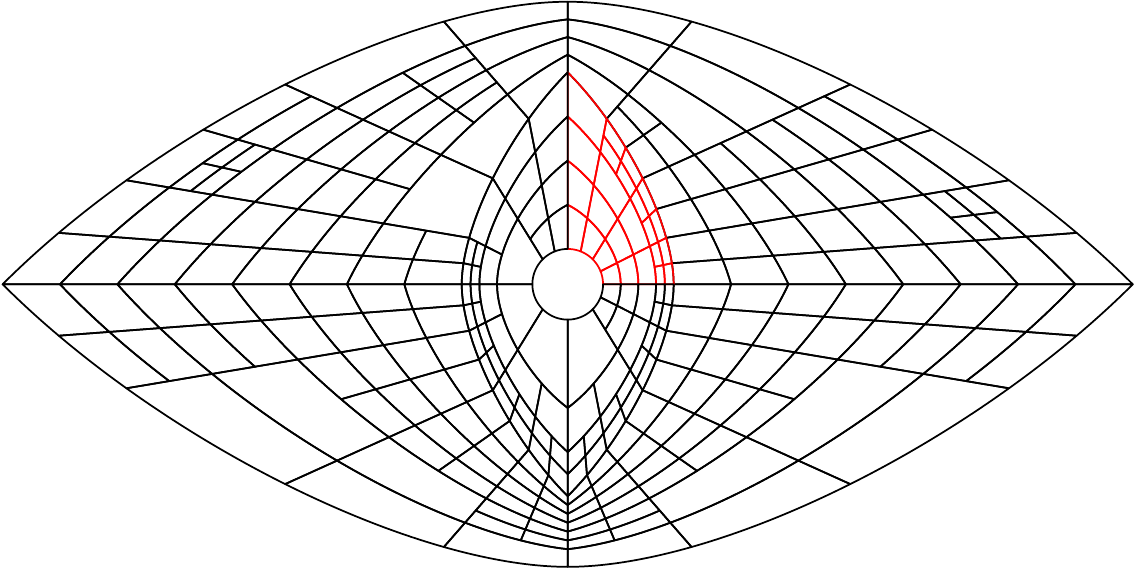}}\\
		\raisebox{1cm}{12} & \includegraphics[width=0.2\columnwidth]{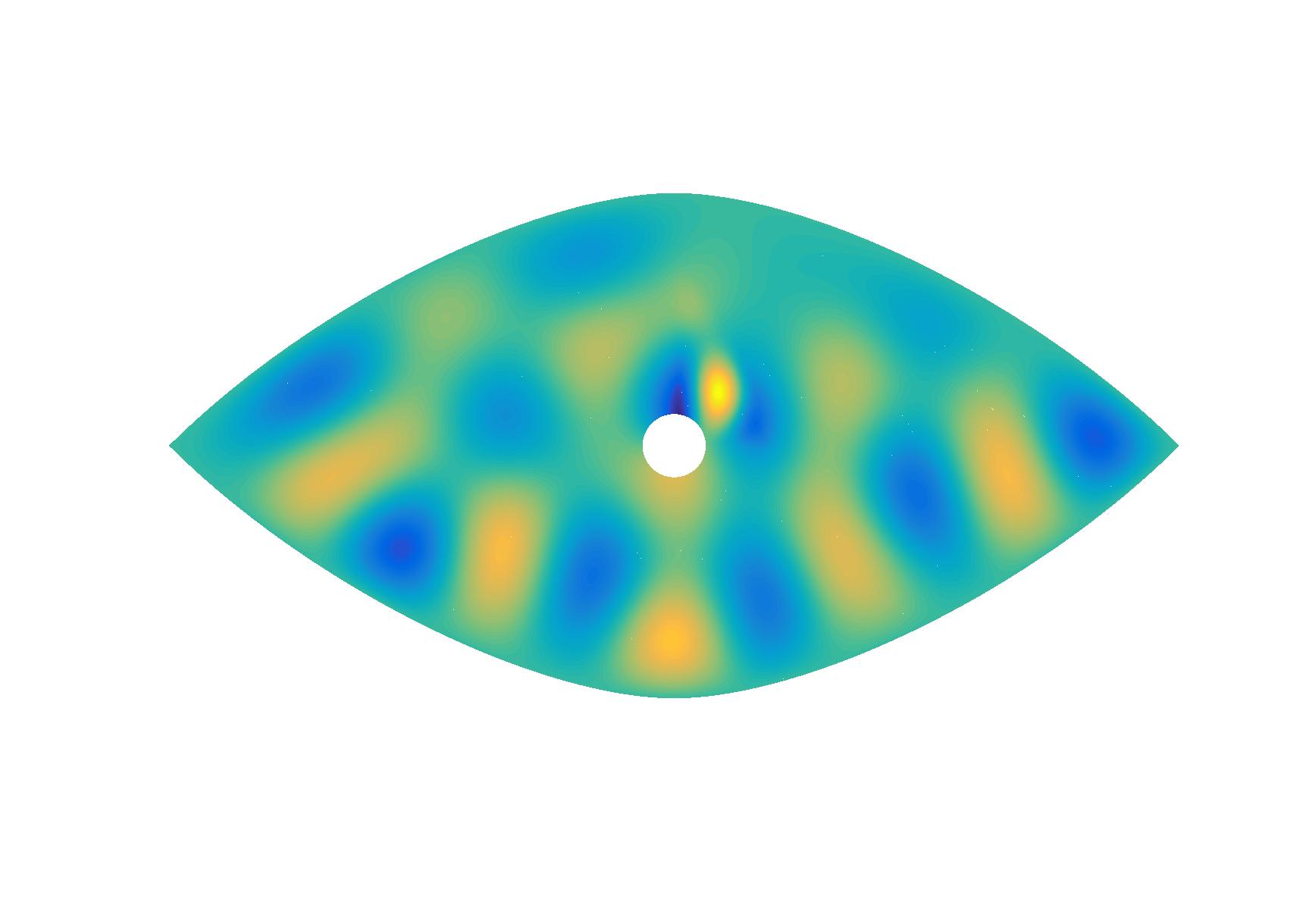} & \includegraphics[width=0.2\columnwidth]{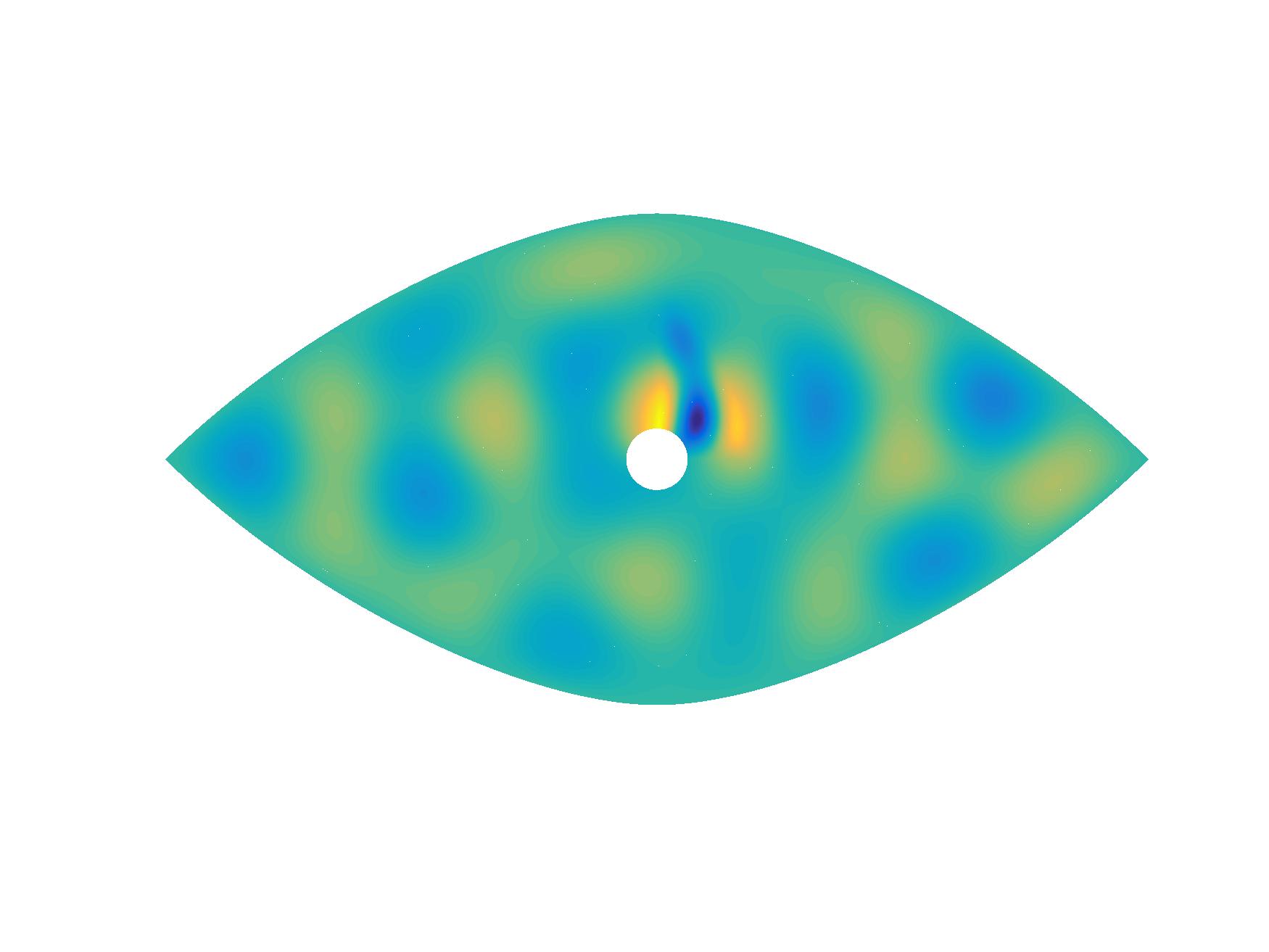} & \raisebox{-0.12cm}{\includegraphics[width=0.22\columnwidth]{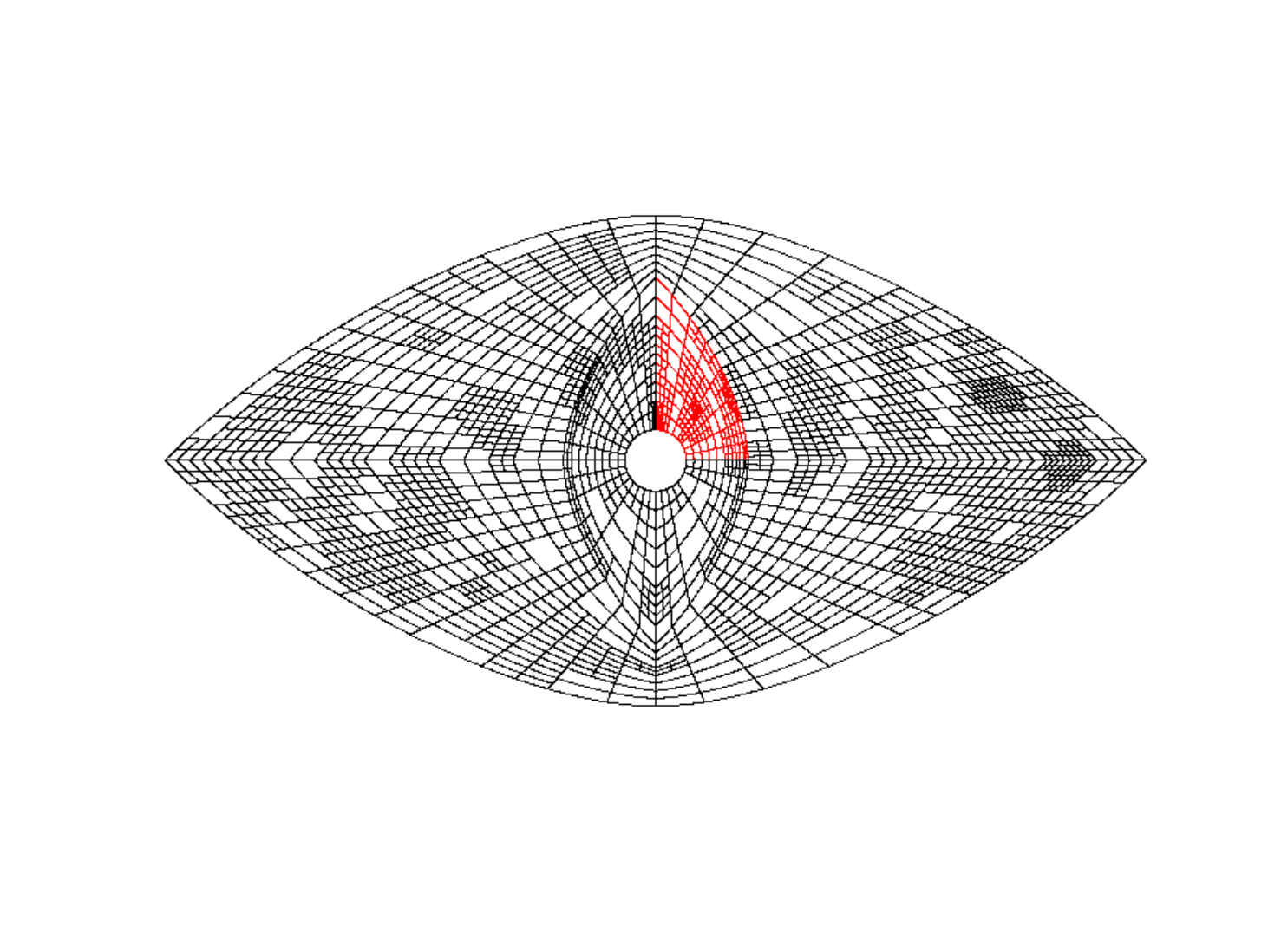}} \\
		\raisebox{1cm}{13} & \includegraphics[width=0.2\columnwidth]{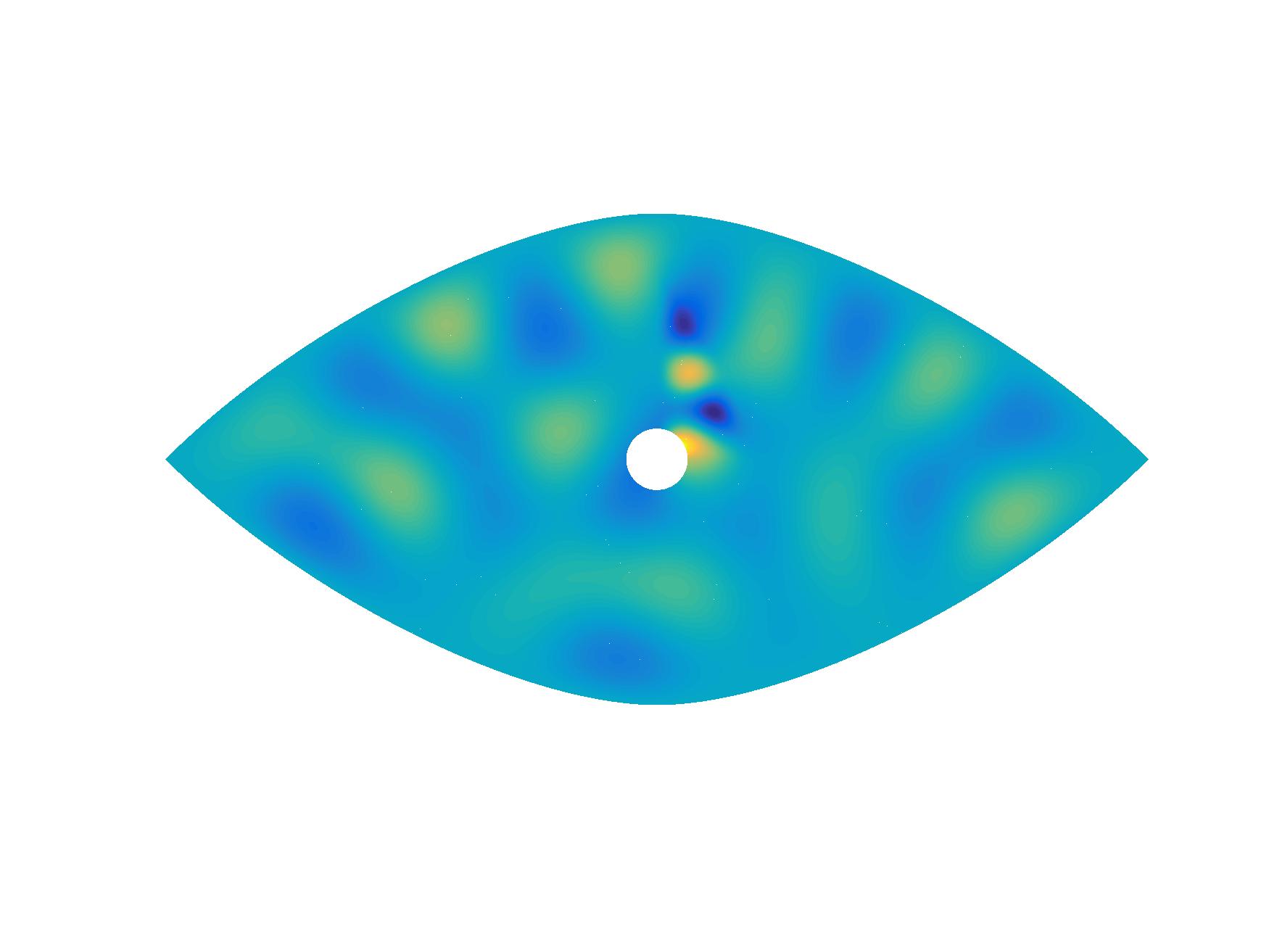} & \includegraphics[width=0.2\columnwidth]{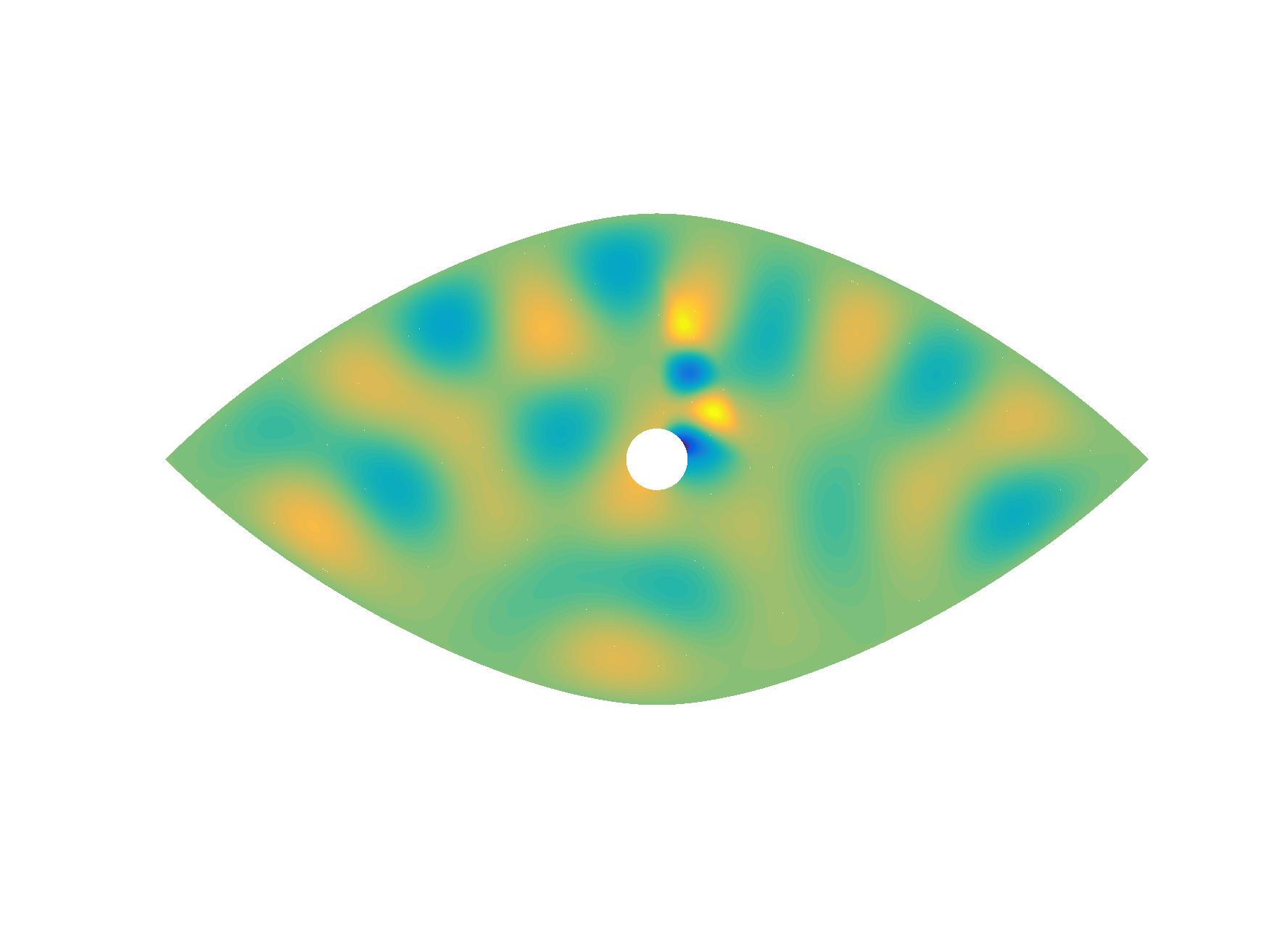} & \raisebox{-0.12cm}{\includegraphics[width=0.22\columnwidth]{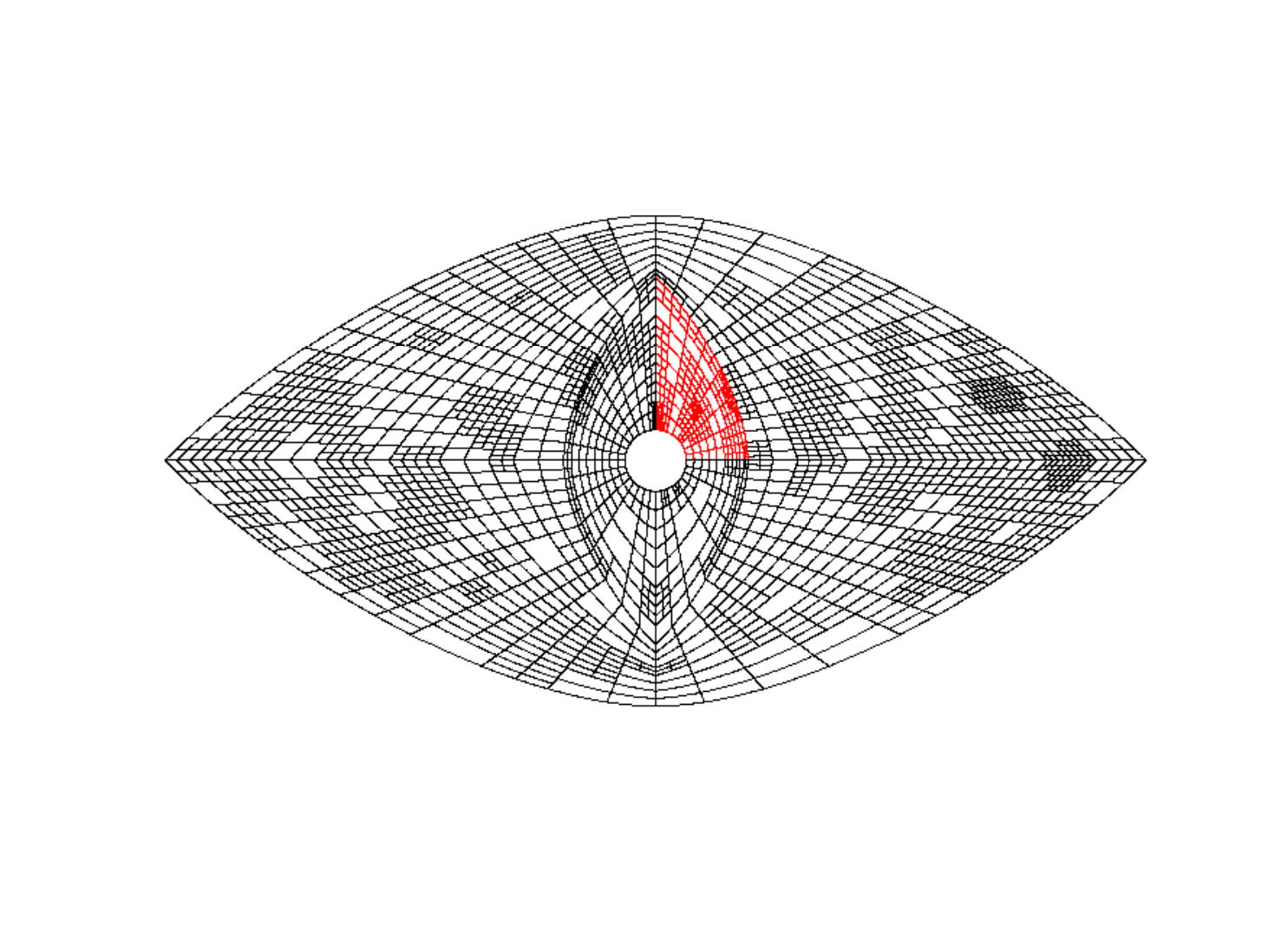}} \\
		\raisebox{1cm}{28} & \includegraphics[width=0.2\columnwidth]{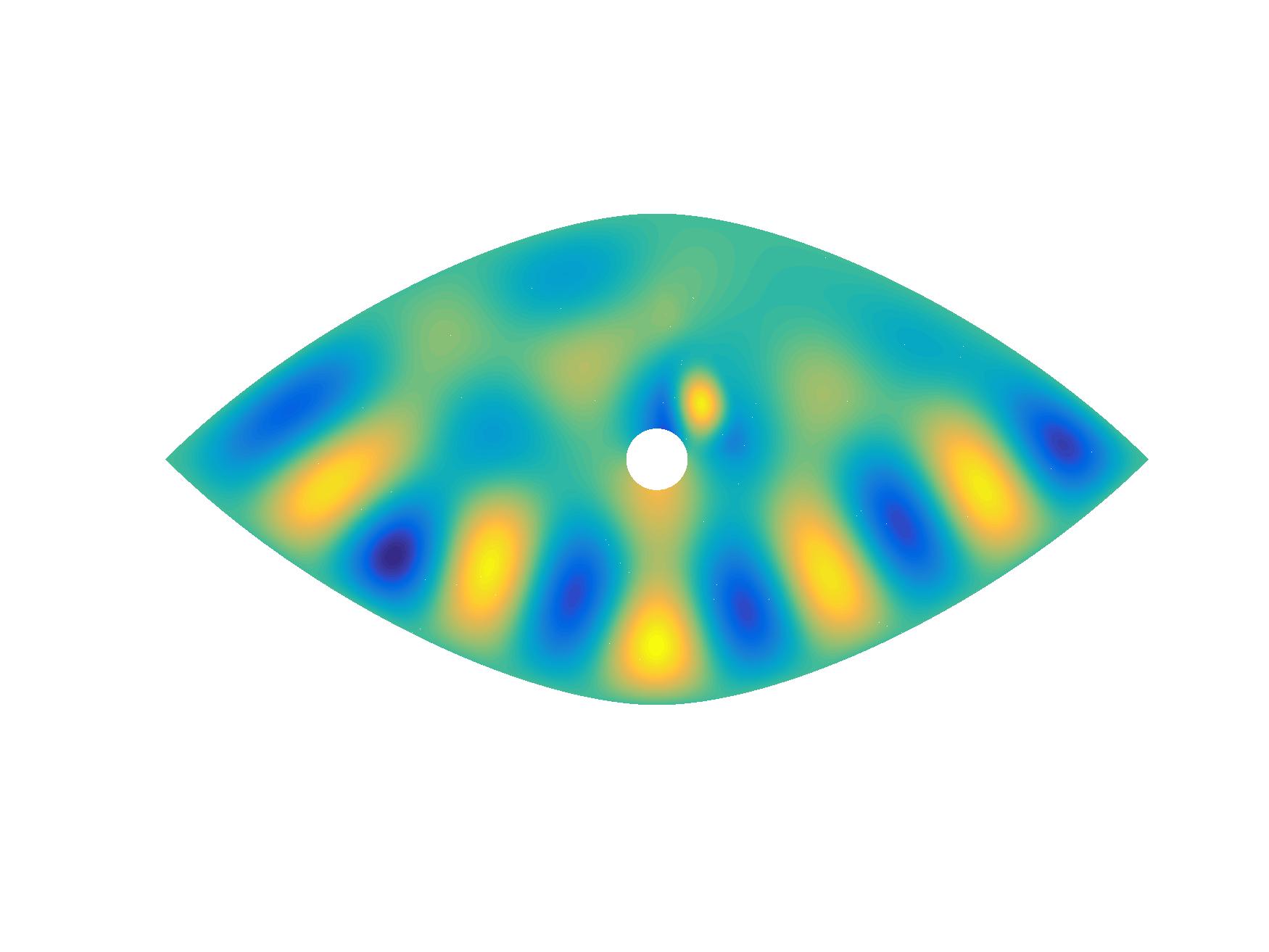} &  \includegraphics[width=0.2\columnwidth]{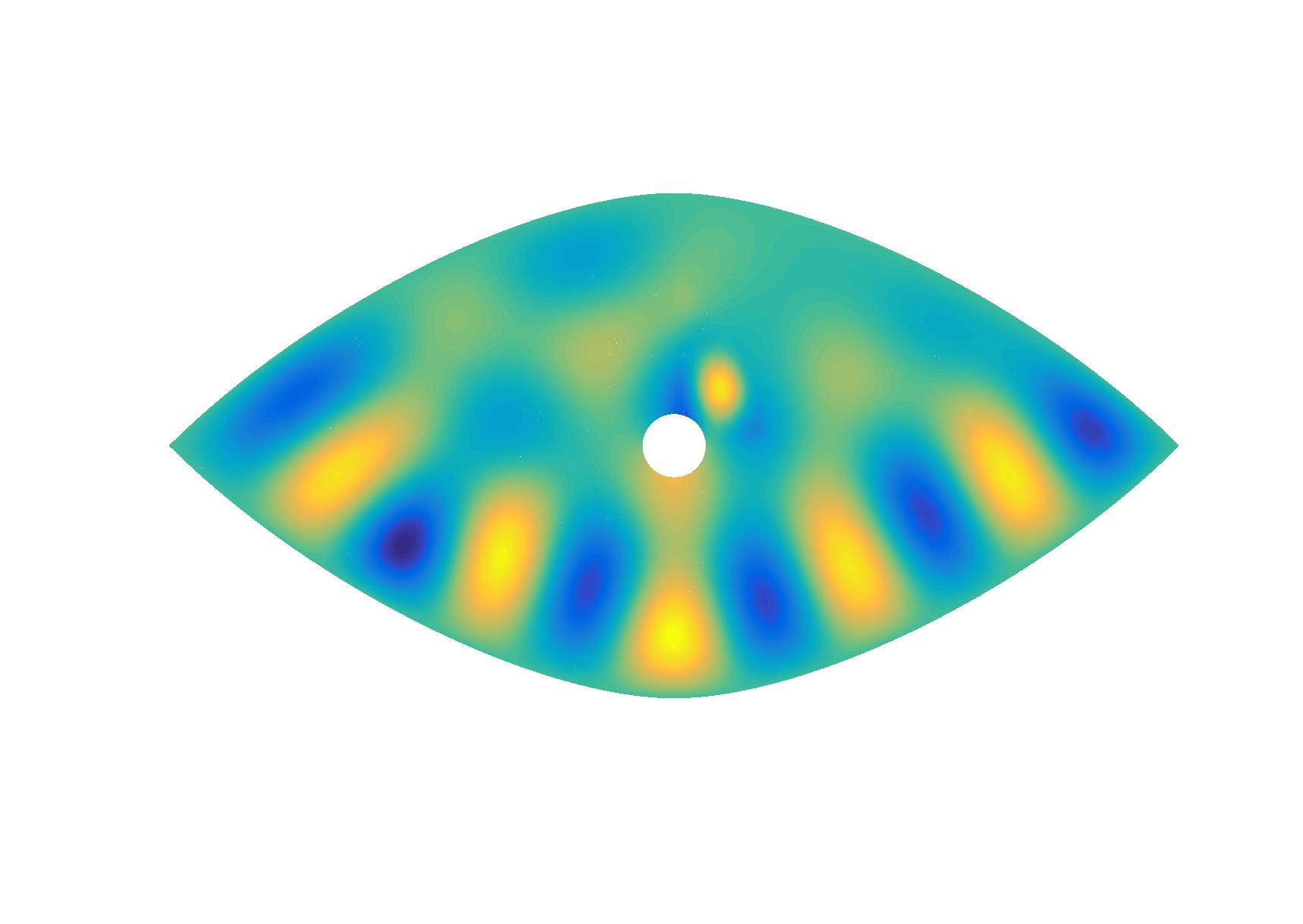} & \raisebox{-0.12cm}{\includegraphics[width=0.21\columnwidth]{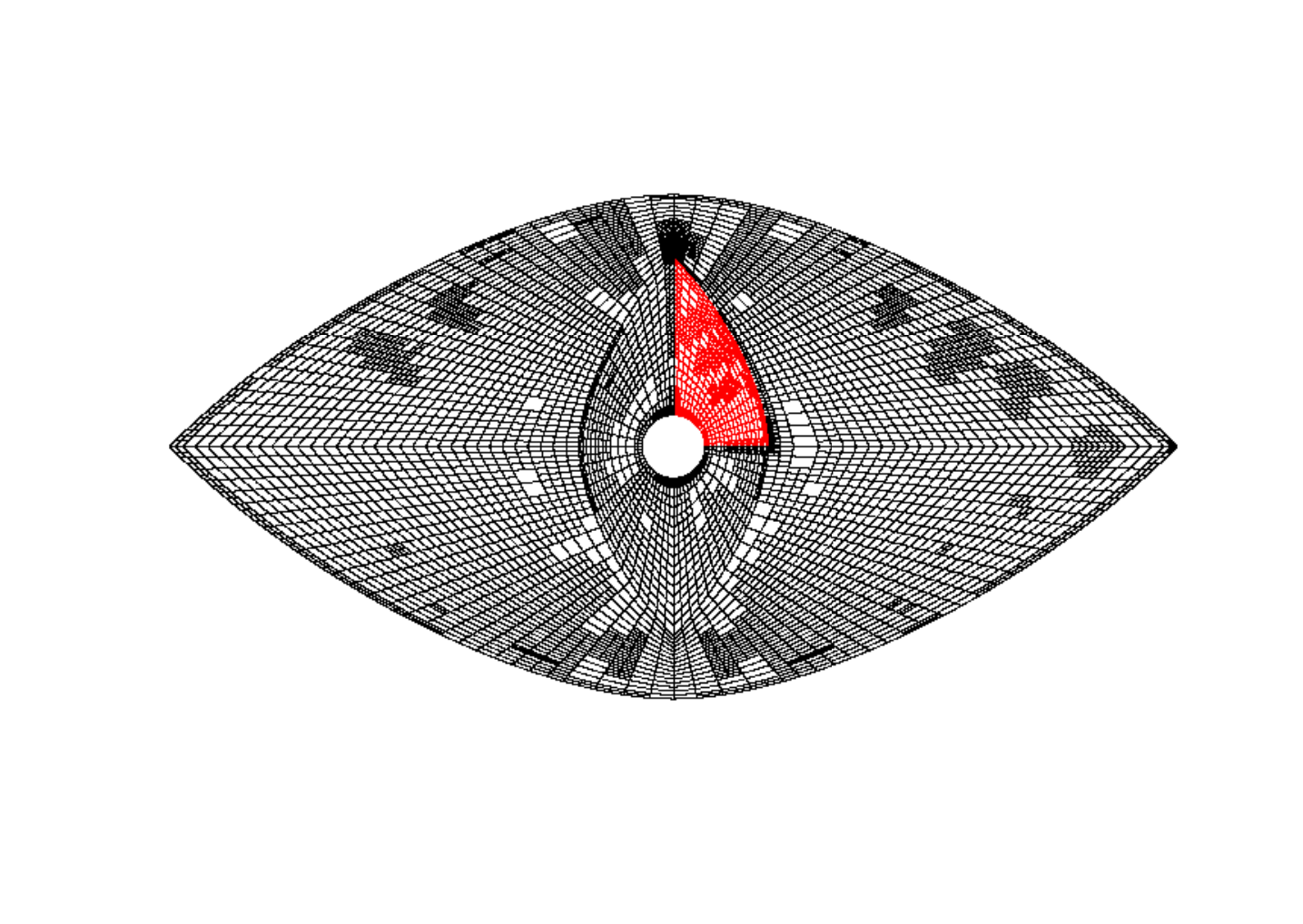}} \\
		\hline
	\end{tabular}
\end{table}

\begin{table}
	\caption{Targeted mode shapes $\vect{\phi}_{\vect{\mathcal{I}}}^{h}$ at coarse mesh and related mode shapes $\tilde{\vect{\phi}}_{\tilde{\vect{\mathcal{I}}}}$ over refined mesh, and the adaptive refinement at different steps acquired through MAC.}
	\label{tab:MACvmode}
	\centering
	\begin{tabular}{cccc}\hline
		\textbf{Step} & $\vect{\phi}_{\vect{\mathcal{I}}}^{h}$ & $\tilde{\vect{\phi}}_{\tilde{\vect{\mathcal{I}}}}$ & \textbf{Adaptive mesh}\\
		\raisebox{1cm}{2} & \includegraphics[width=0.2\columnwidth]{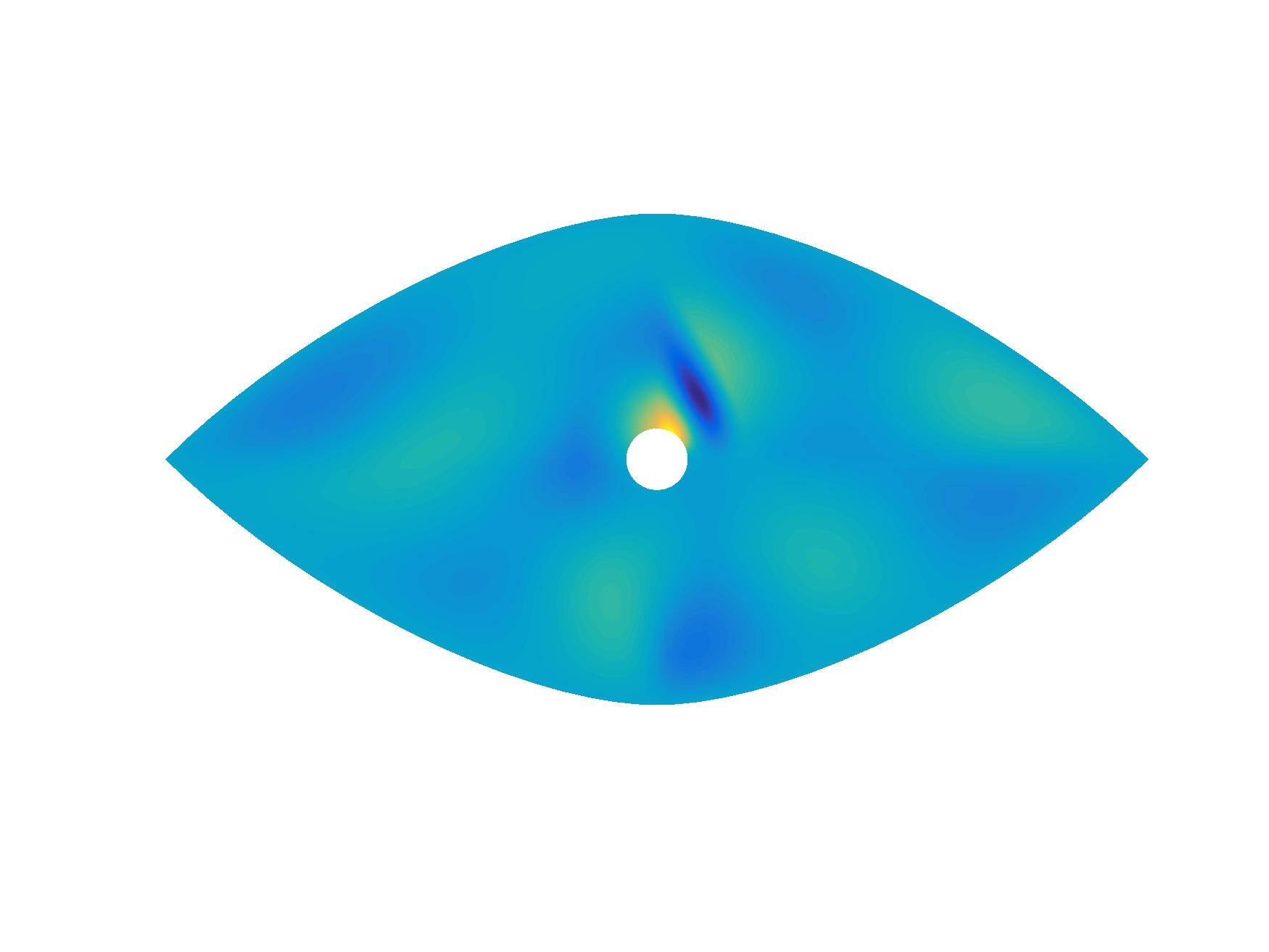}  &  \includegraphics[width=0.2\columnwidth]{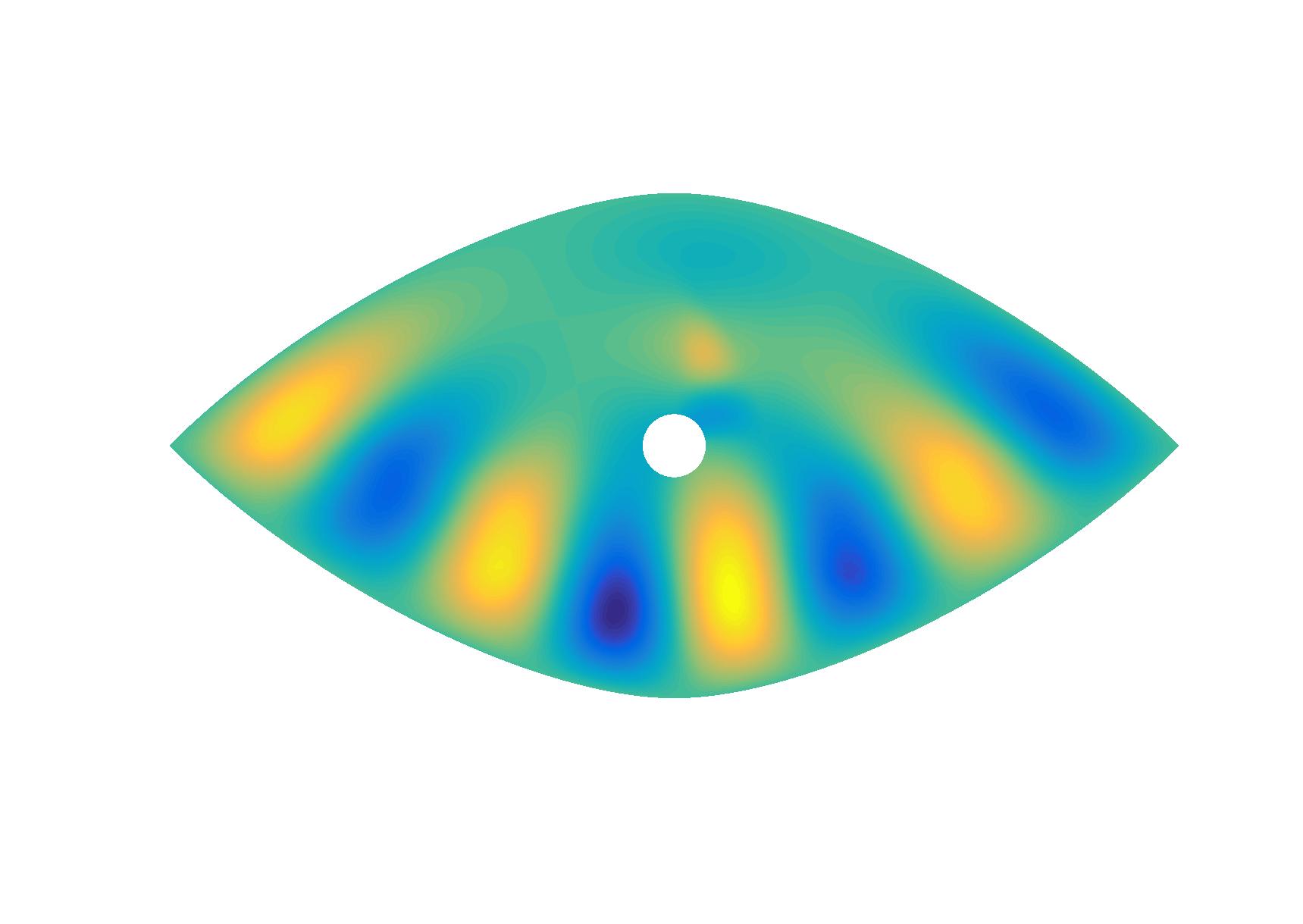} & \raisebox{0.55cm}{\includegraphics[width=0.17\columnwidth]{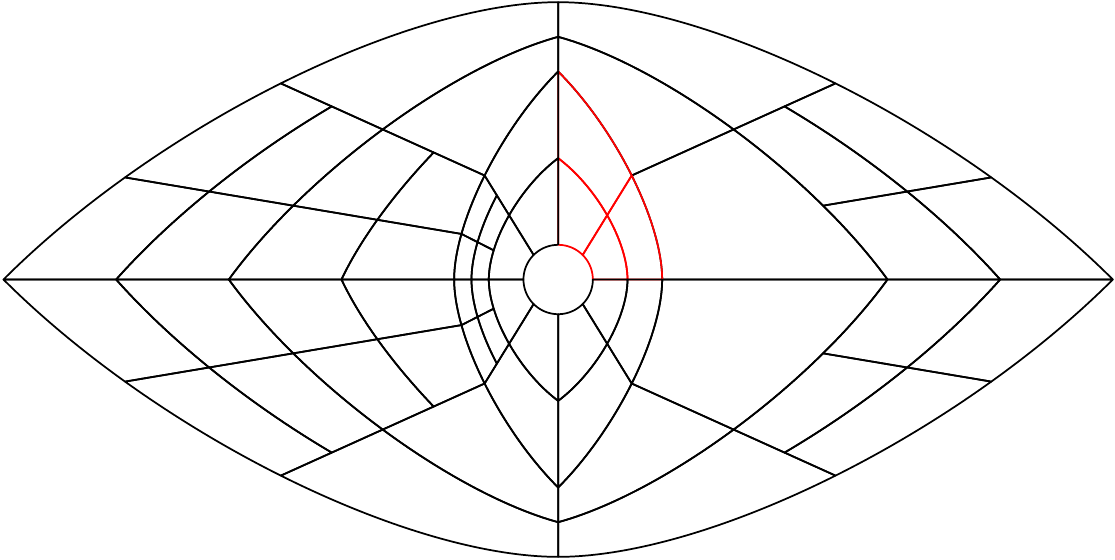}}\\
		\raisebox{1cm}{3} & \includegraphics[width=0.2\columnwidth]{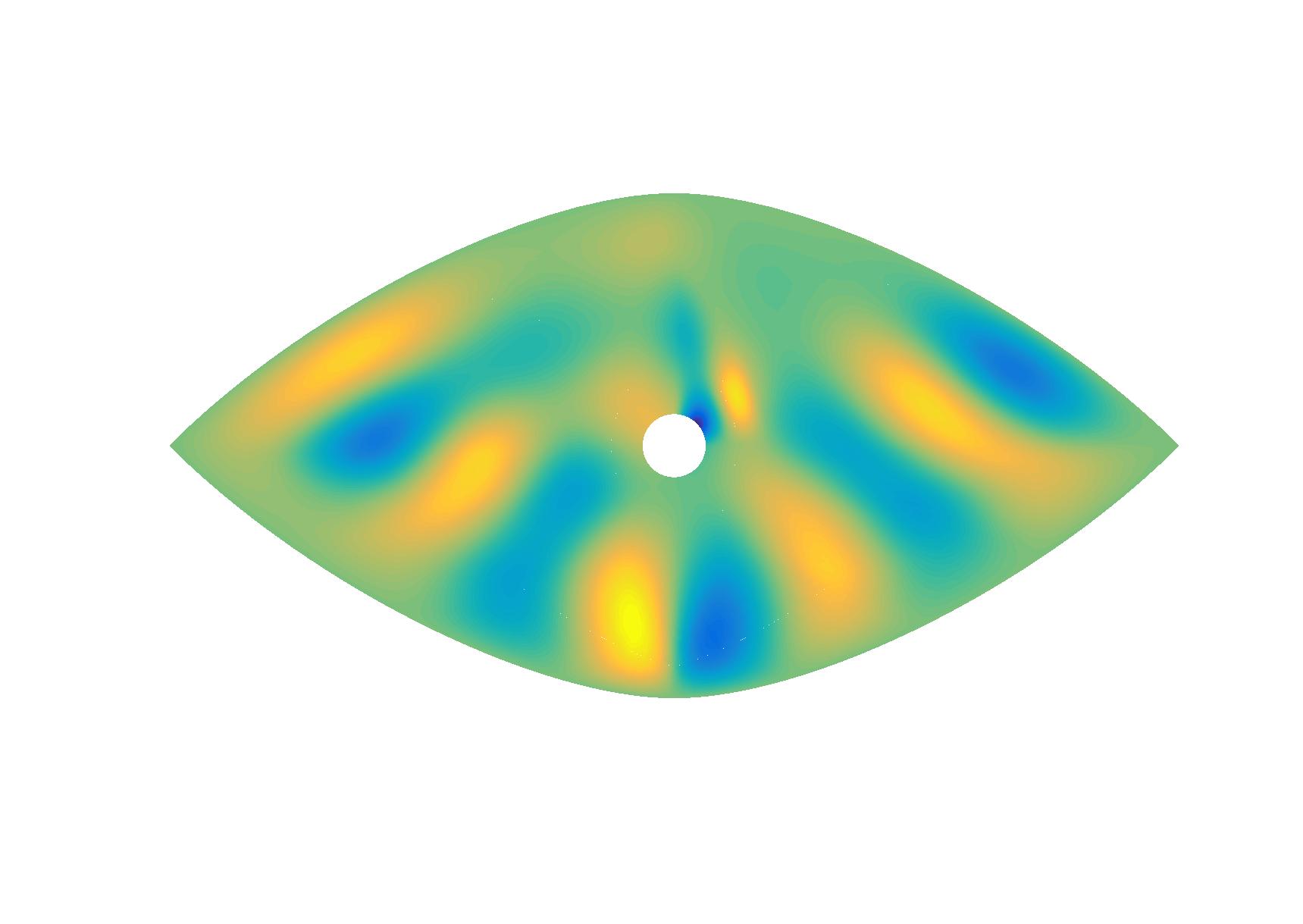} &  \includegraphics[width=0.2\columnwidth]{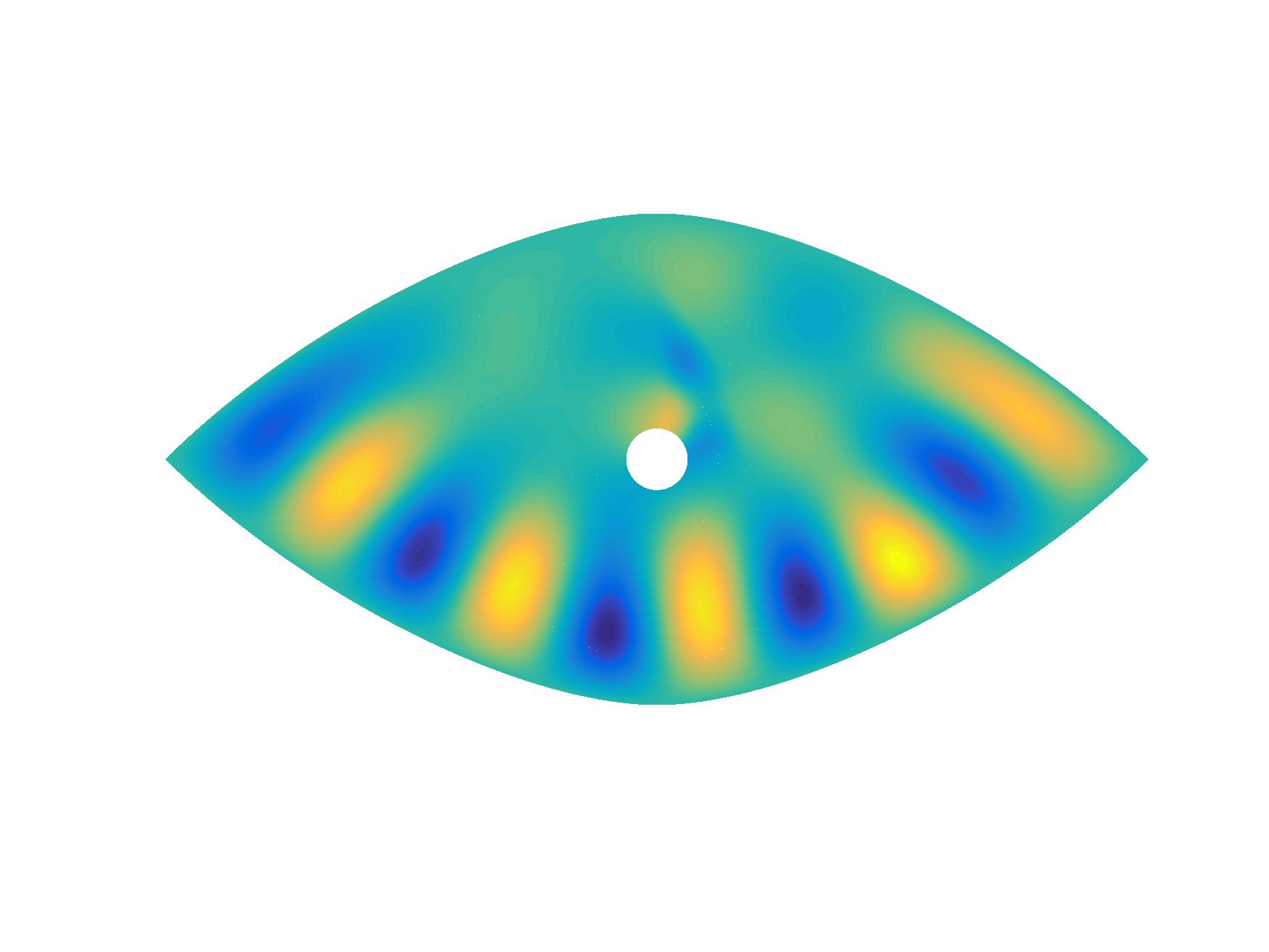}  & \raisebox{0.55cm}{\includegraphics[width=0.17\columnwidth]{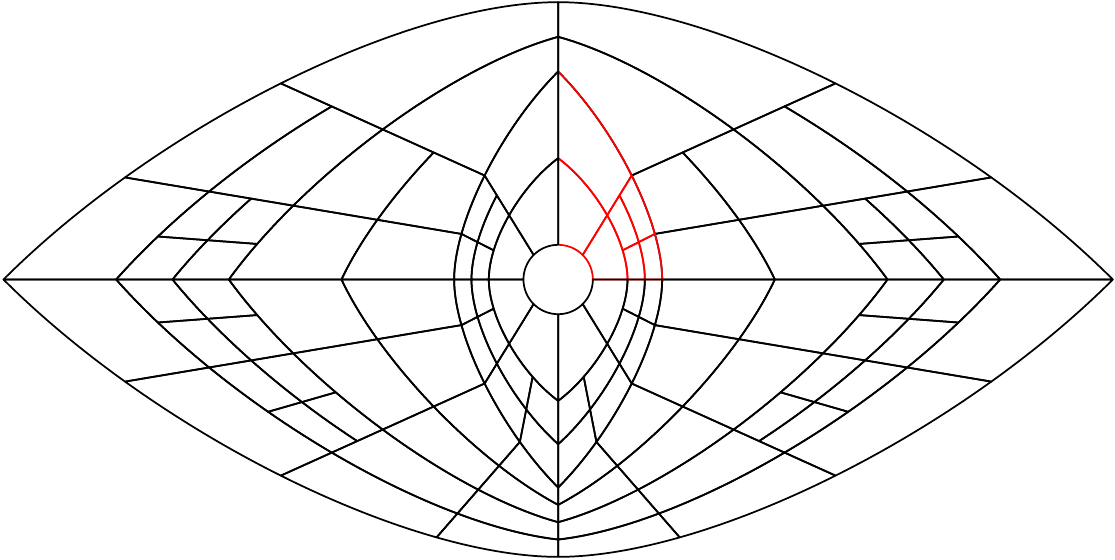}}\\
		\raisebox{1cm}{6} & \includegraphics[width=0.2\columnwidth]{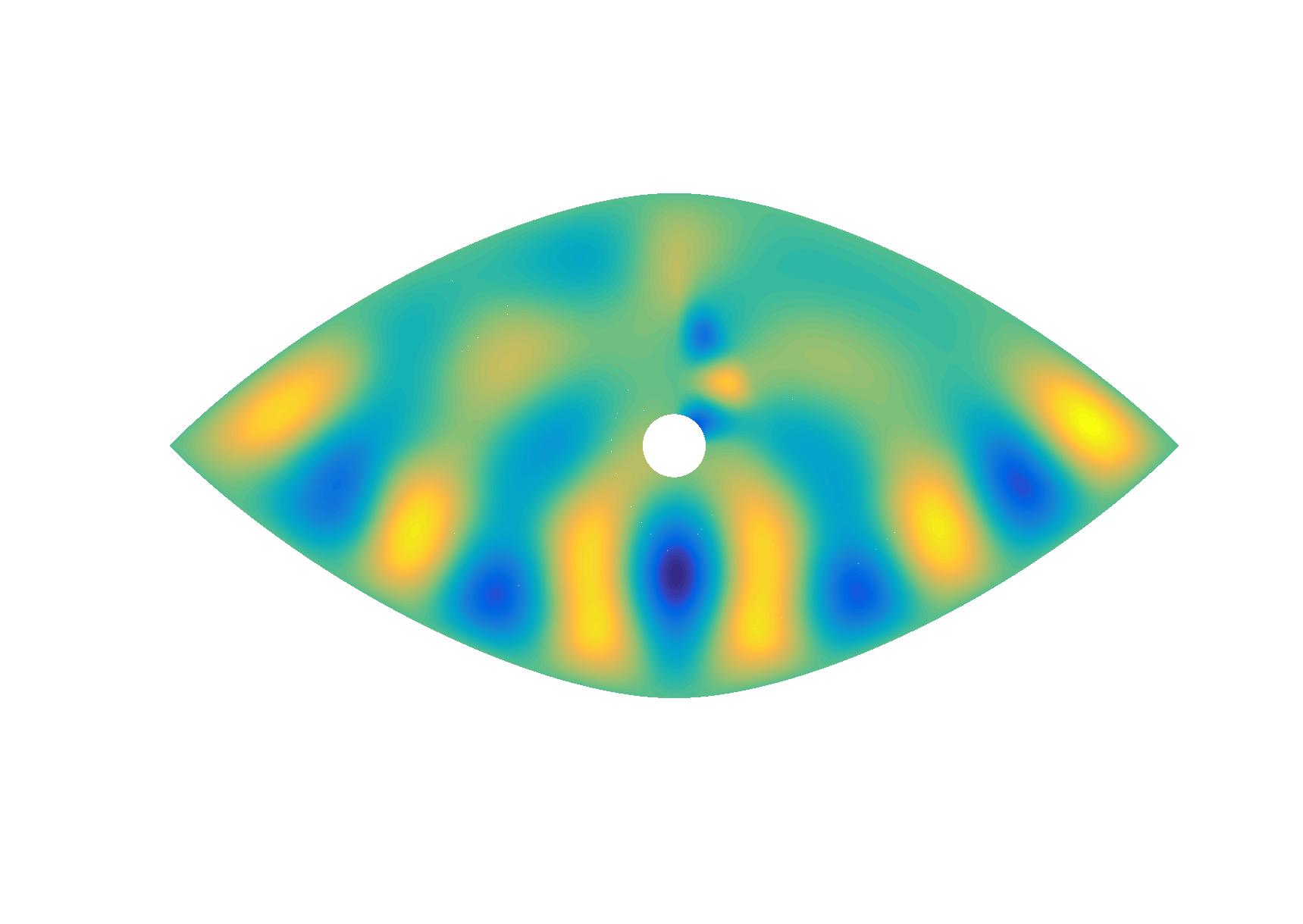} & \includegraphics[width=0.2\columnwidth]{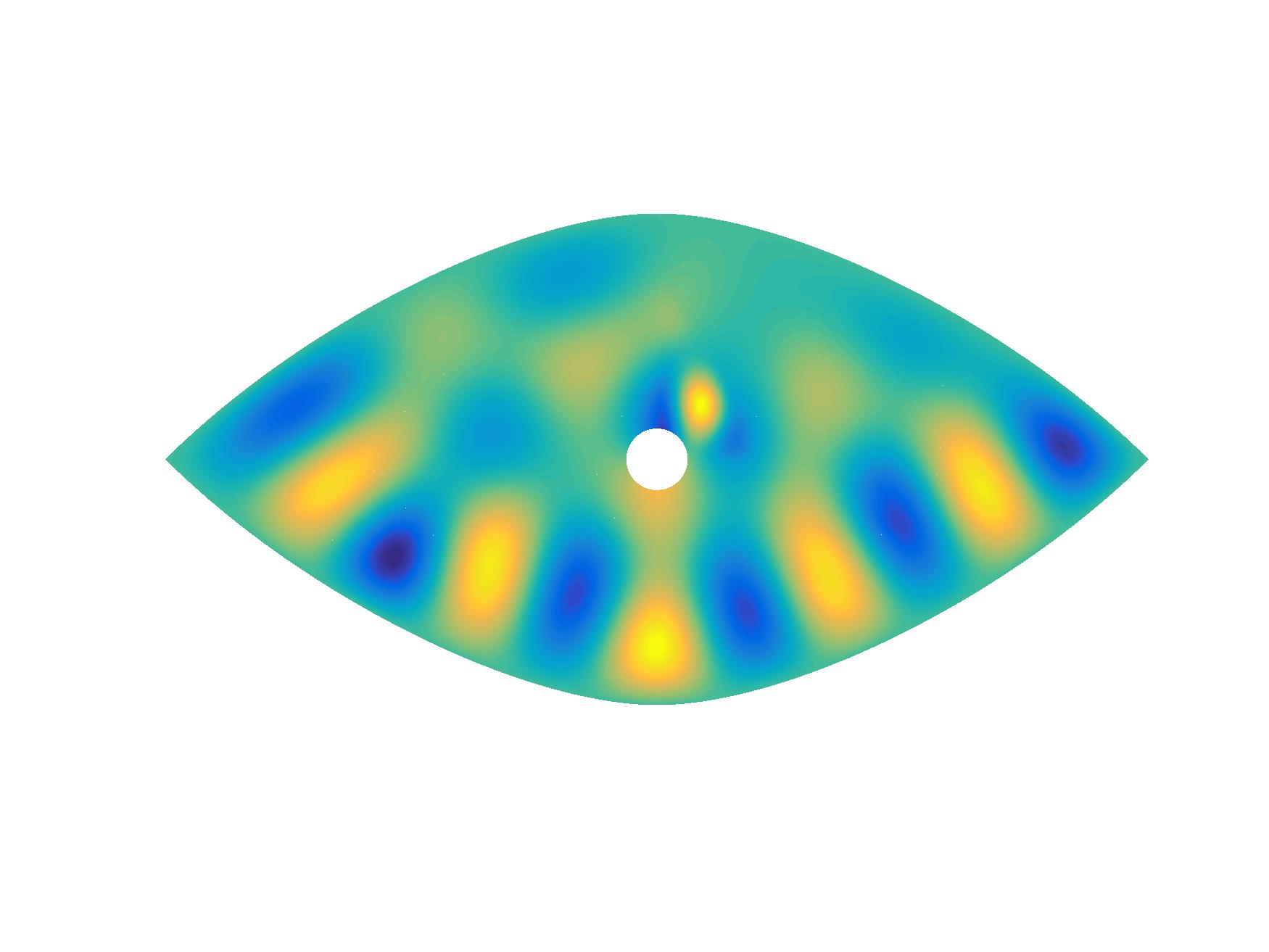} & \raisebox{0.55cm}{\includegraphics[width=0.17\columnwidth]{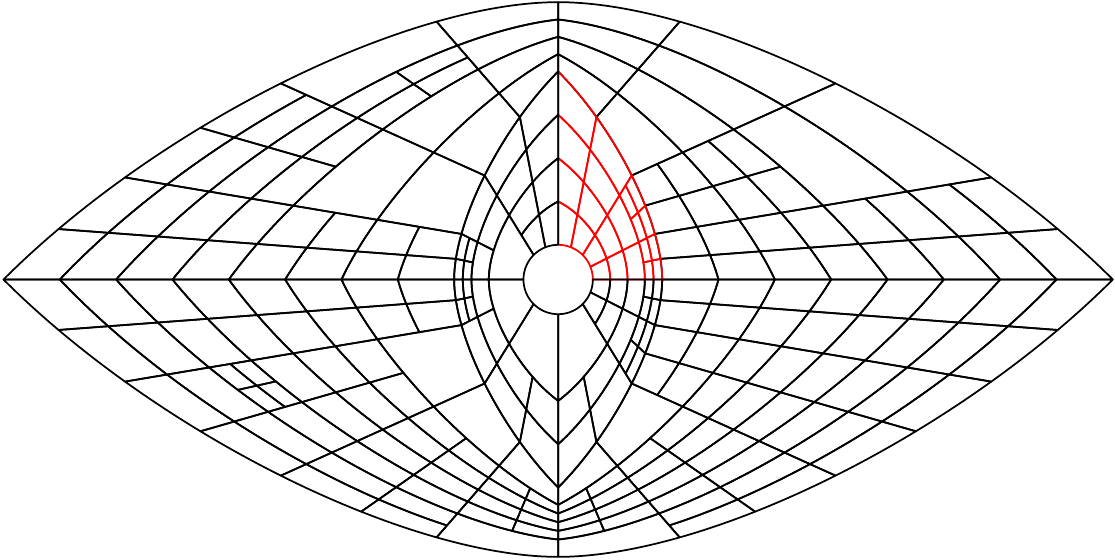}} \\
		\raisebox{1cm}{13} & \includegraphics[width=0.2\columnwidth]{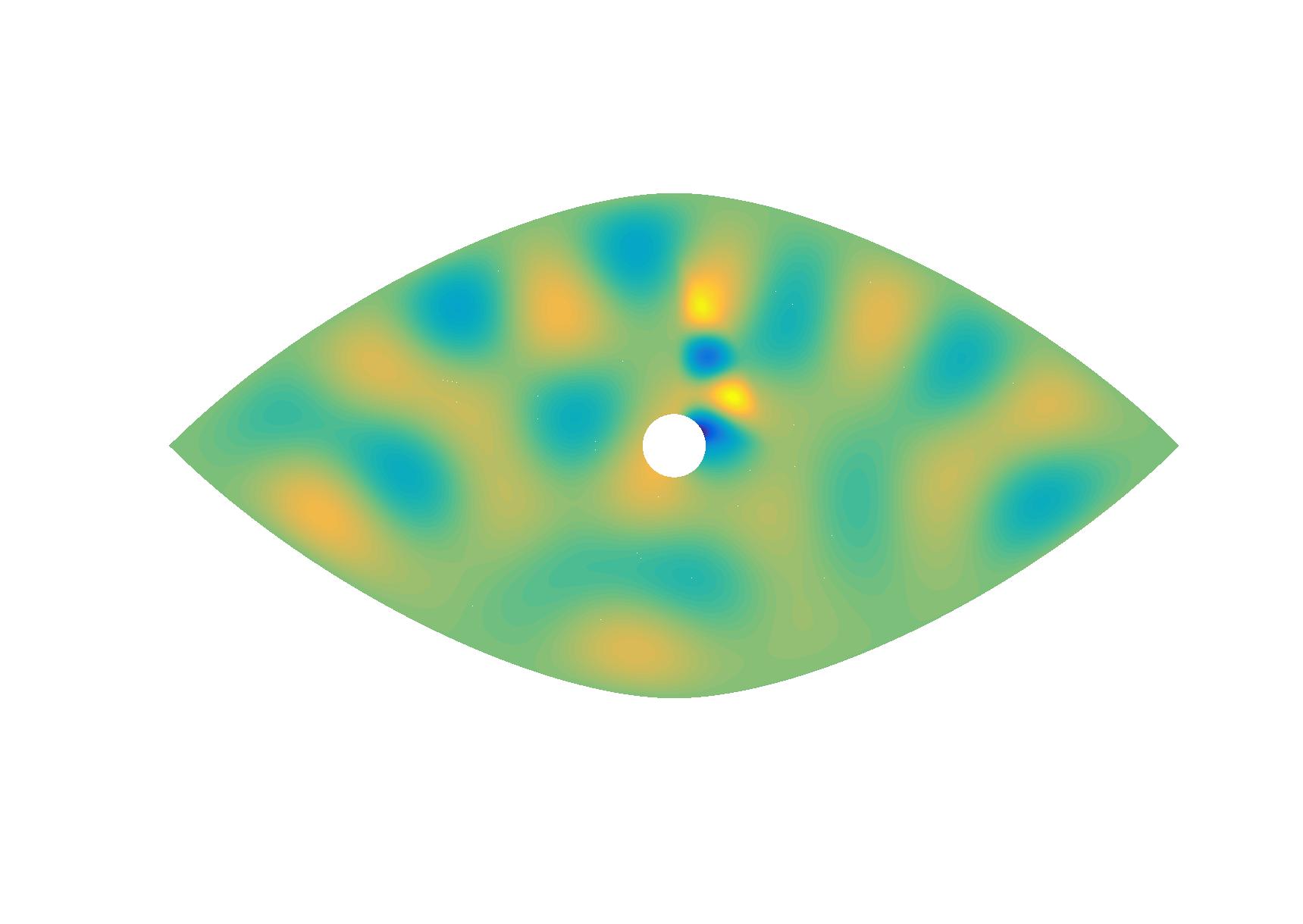} & \includegraphics[width=0.2\columnwidth]{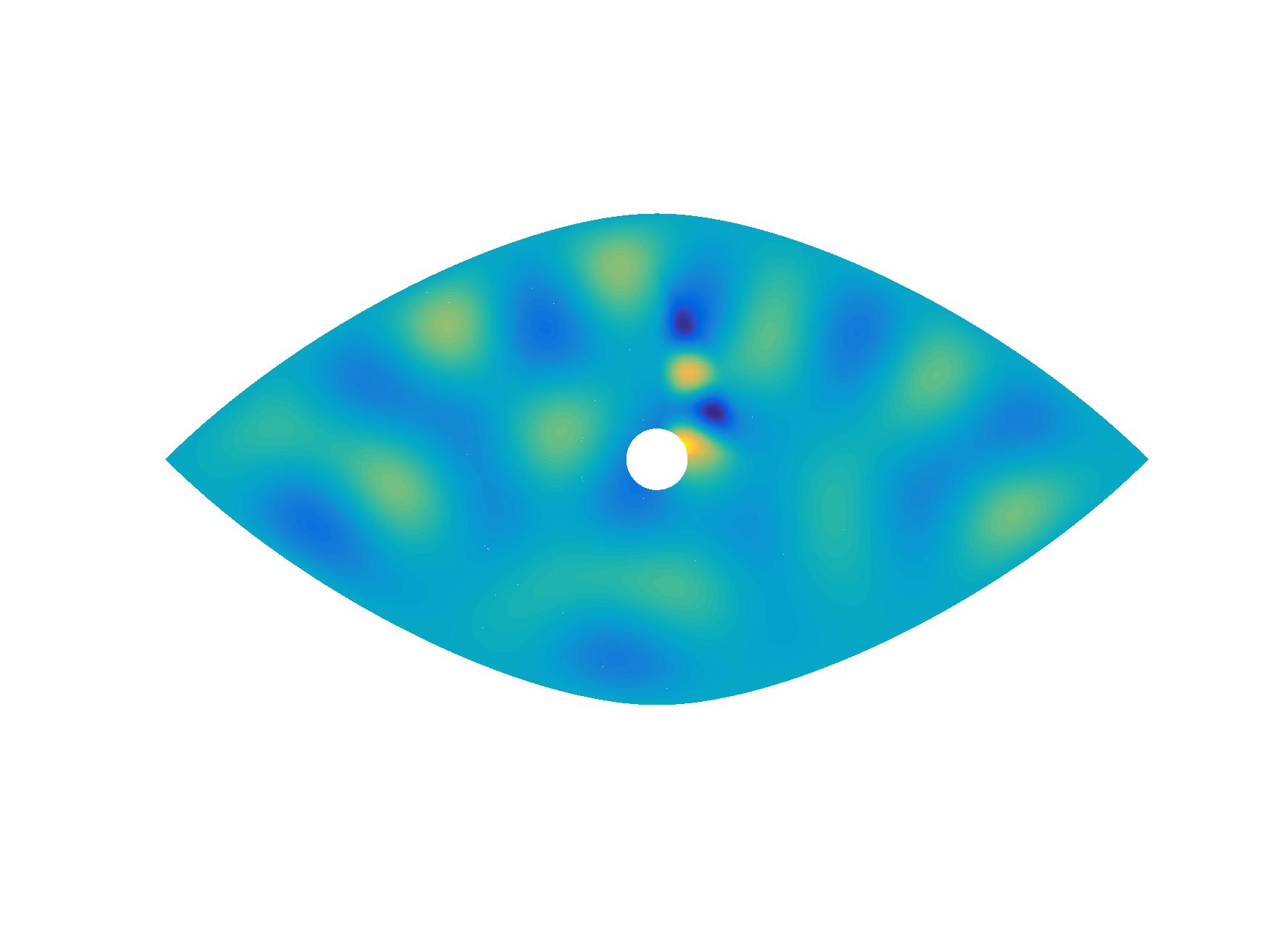} & \raisebox{-0.12cm}{\includegraphics[width=0.22\columnwidth]{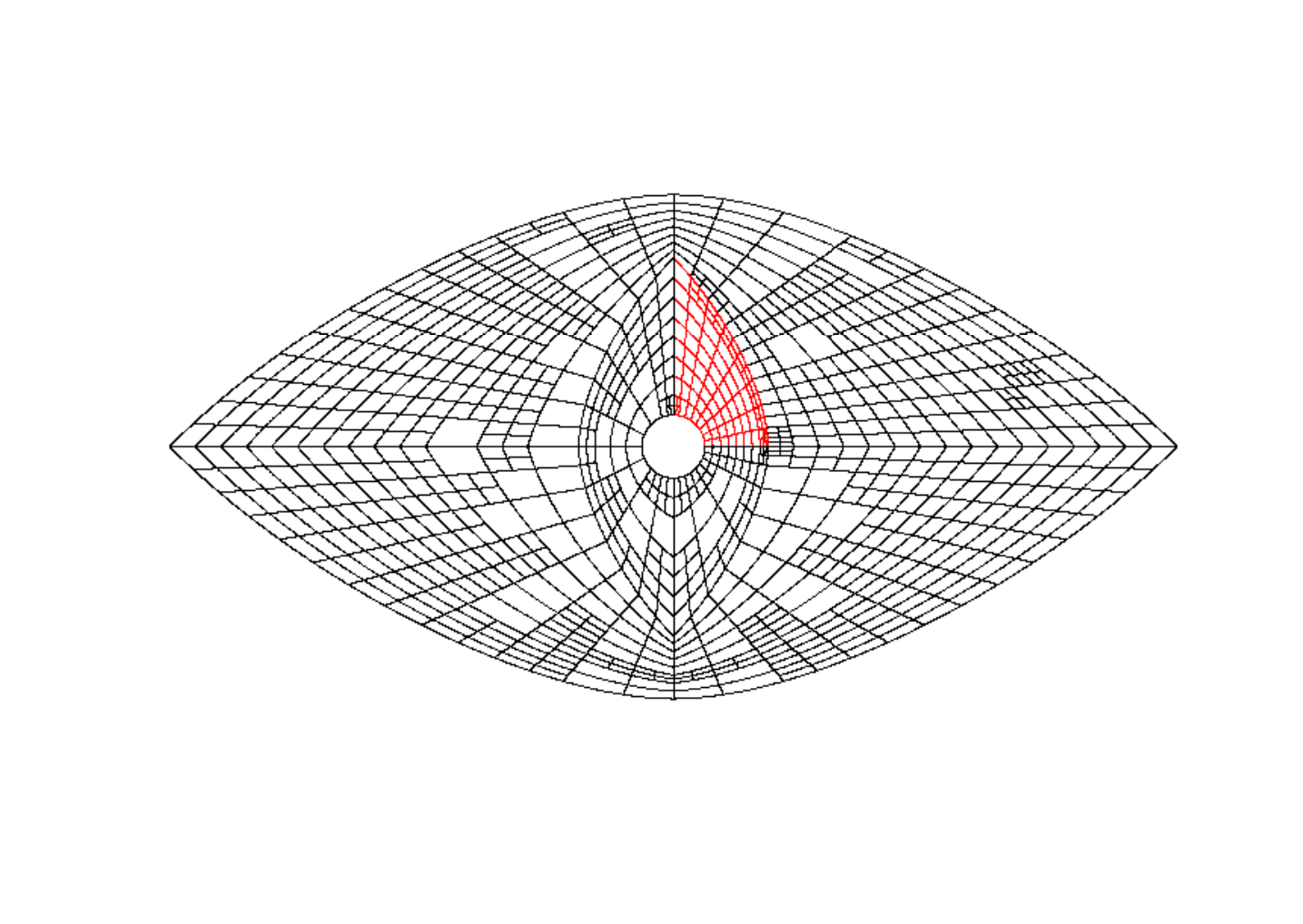}} \\
		\raisebox{1cm}{24} & \includegraphics[width=0.2\columnwidth]{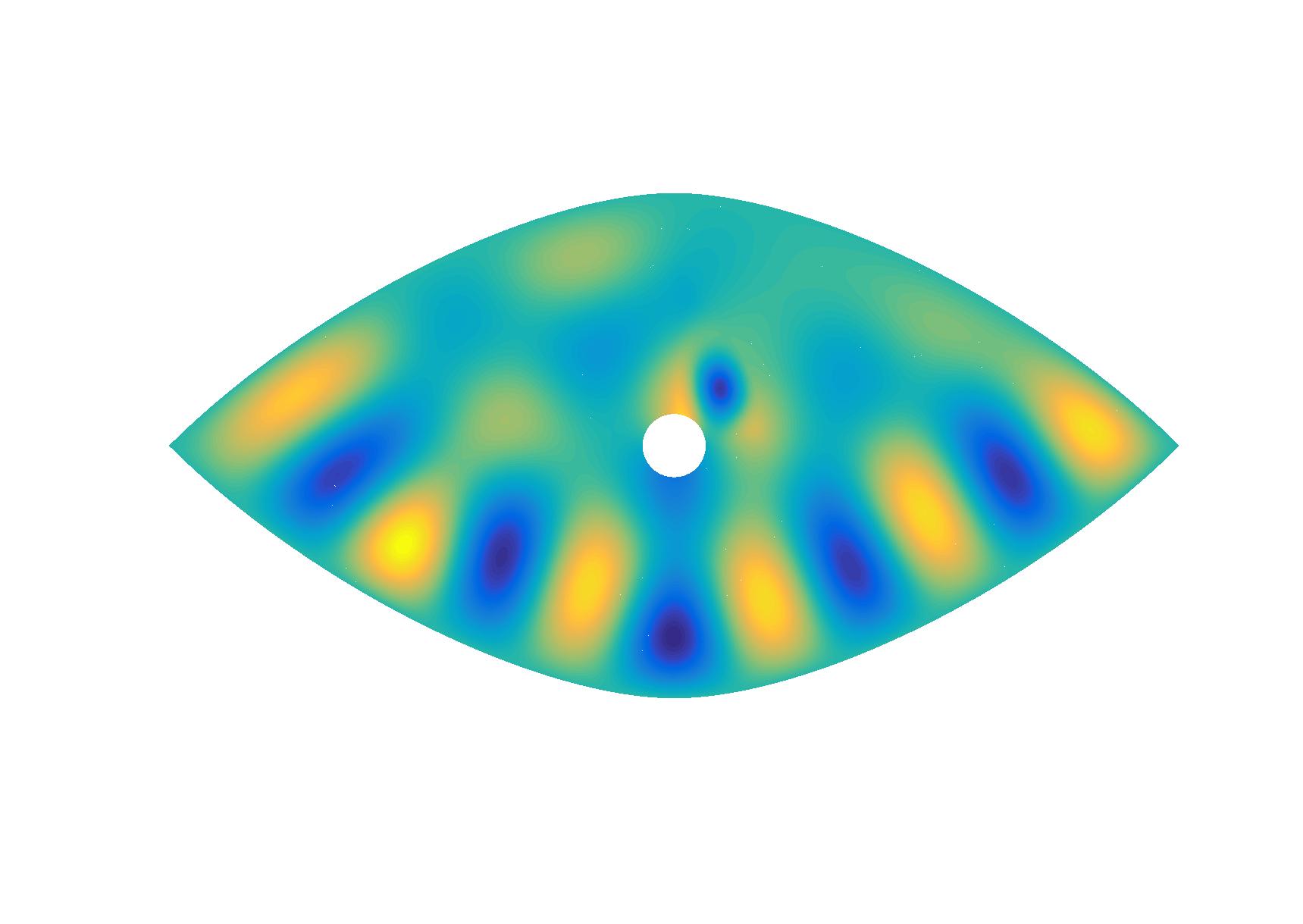} & \includegraphics[width=0.2\columnwidth]{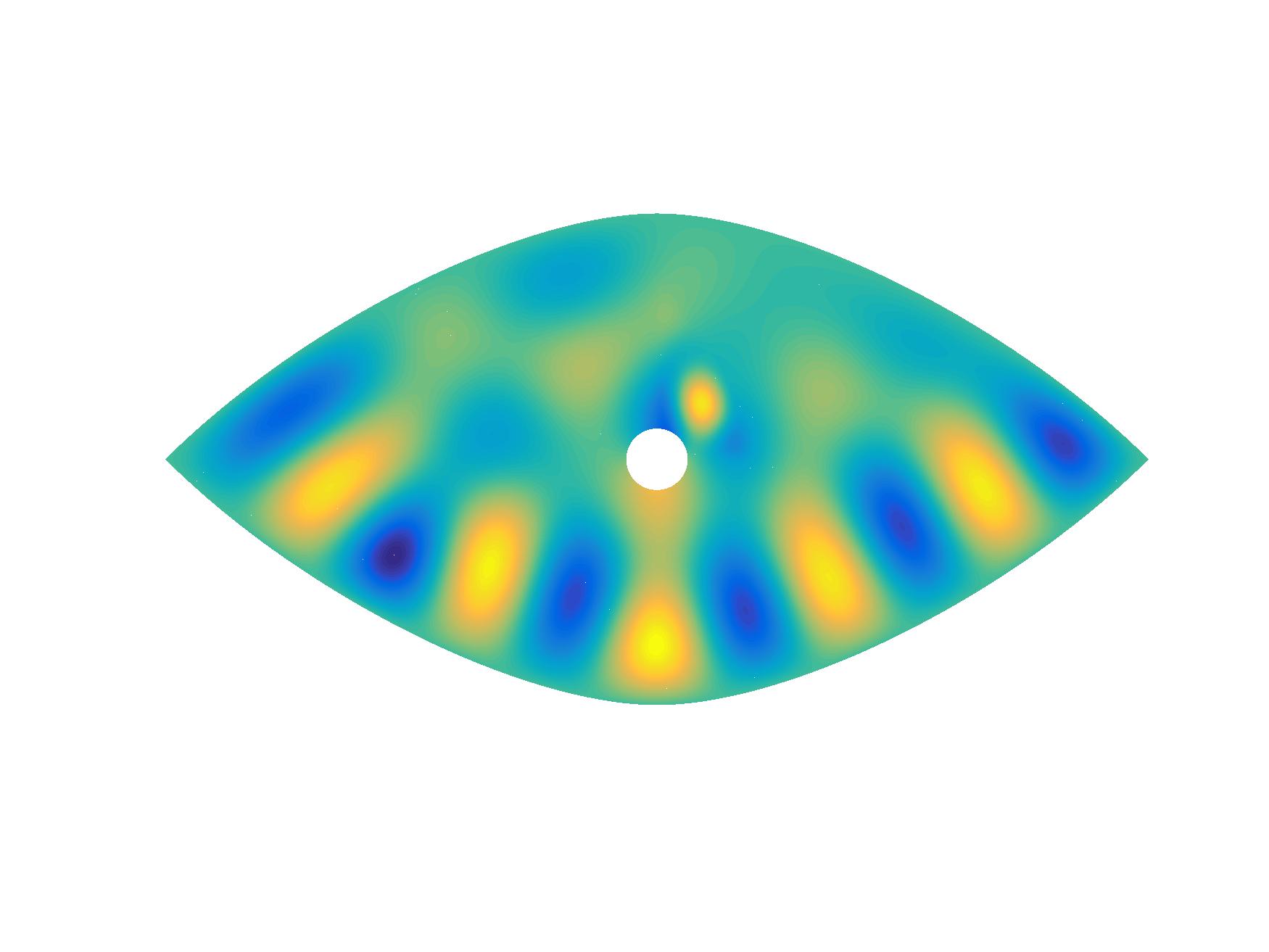} & \raisebox{-0.1cm}{\includegraphics[width=0.22\columnwidth]{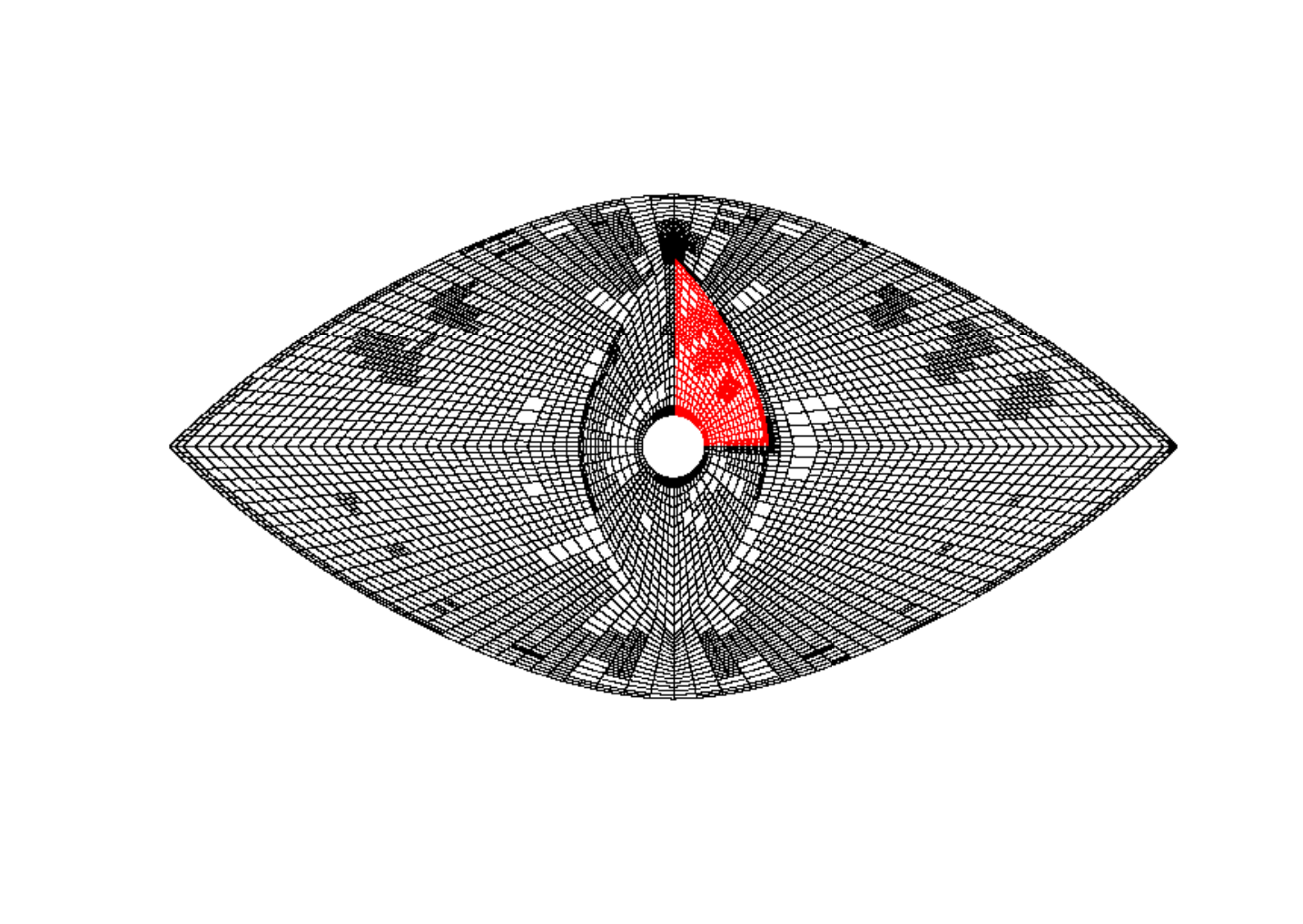}} \\
		\hline
	\end{tabular}
\end{table}

\subsection{Heterogeneous square plate with holes}\label{sec:plateholes}
\begin{figure}
	\centering
	\begin{subfigure}[t]{0.4\textwidth}
		\centering
		\includegraphics[width=1\columnwidth]{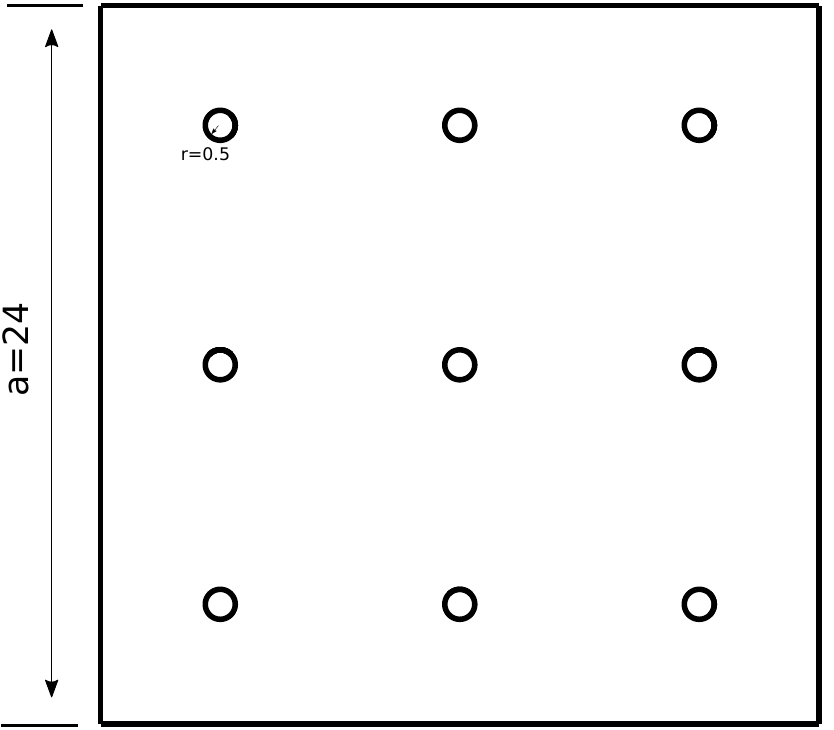}
		\caption{Geometry of structure}
	\end{subfigure}
	\hspace{1cm}
	\begin{subfigure}[t]{0.4\textwidth}
		\centering
		\includegraphics[width=1\columnwidth]{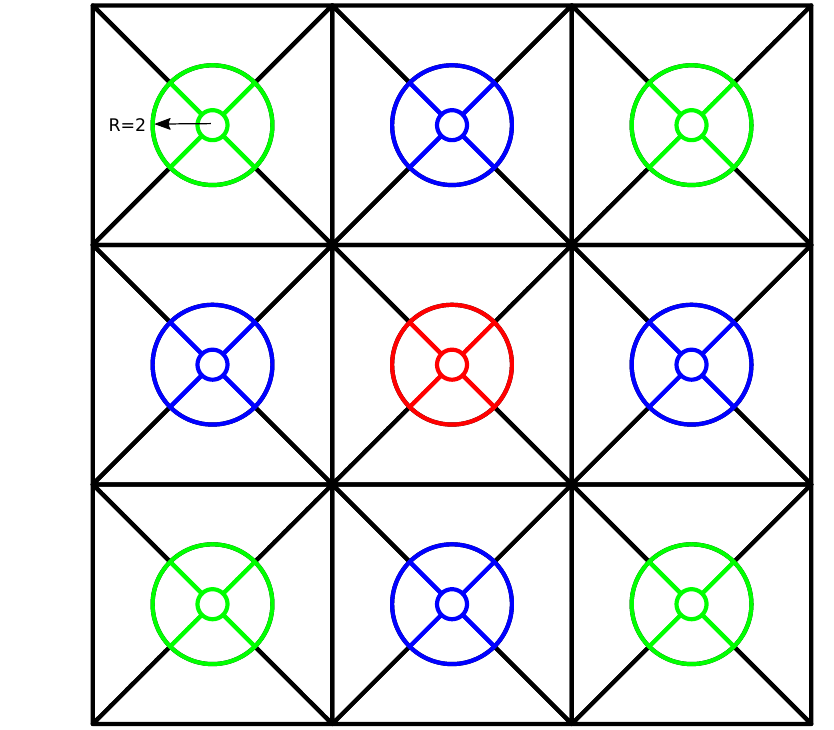}
		\caption{Discretization of patches}
	\end{subfigure}
	\caption{The simply supported square plate with holes is discretized by 72 patches with the Young's modulus $E_{\text{red}} = 0.05, ~E_{\text{blue}} = 0.08, ~E_{\text{green}} = 0.12, ~E_{\text{black}} = 1$. The geometry is represented by NURBS $(p = 2, q = 2)$, and the solution field is approximated by PHT splines $(p = 3, q = 3)$.}
	\label{plateholesGeo}
\end{figure}

\begin{figure}
	\centering
	\includegraphics[width=0.4\columnwidth]{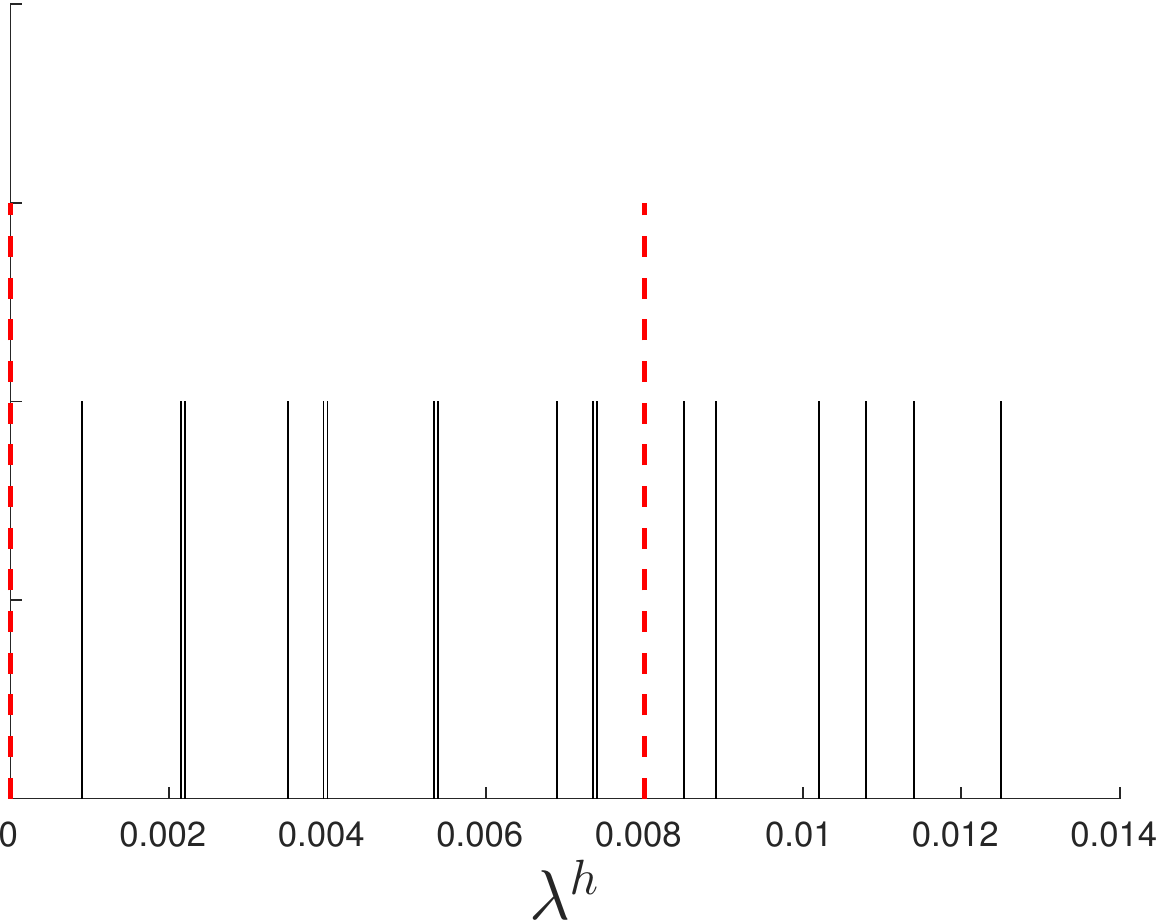}
	\caption{Frequencies of interest in the window with interval $[0, 0.008]$.}
	\label{fig:winplatehole}
\end{figure}

\begin{figure}
	\centering
	\begin{subfigure}[t]{0.3\textwidth}
		\centering
		\includegraphics[width=1\columnwidth]{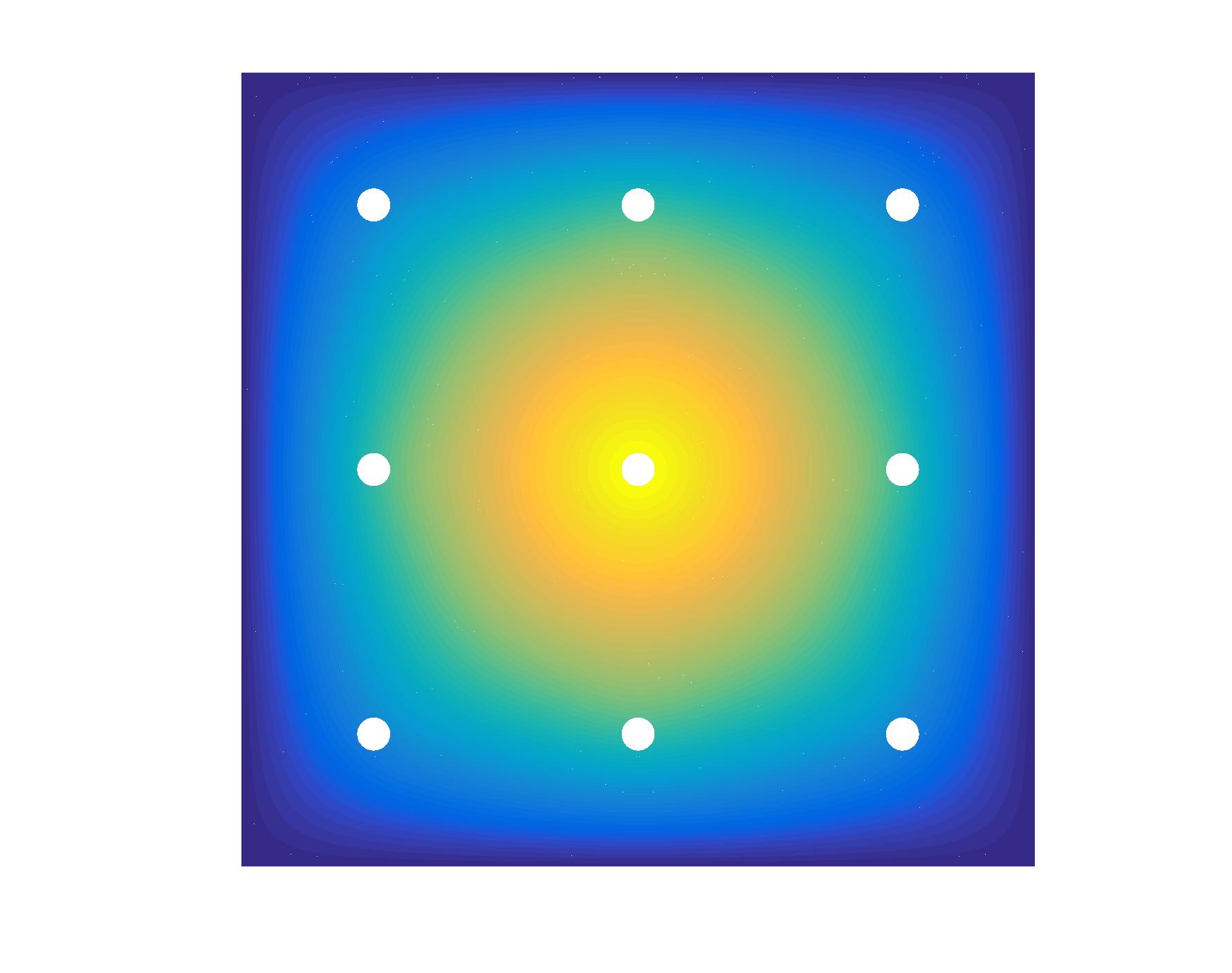}
		\caption{1st mode, $ \lambda_{N}^{1} = 1 $}
	\end{subfigure}
	\begin{subfigure}[t]{0.3\textwidth}
		\centering
		\includegraphics[width=1\columnwidth]{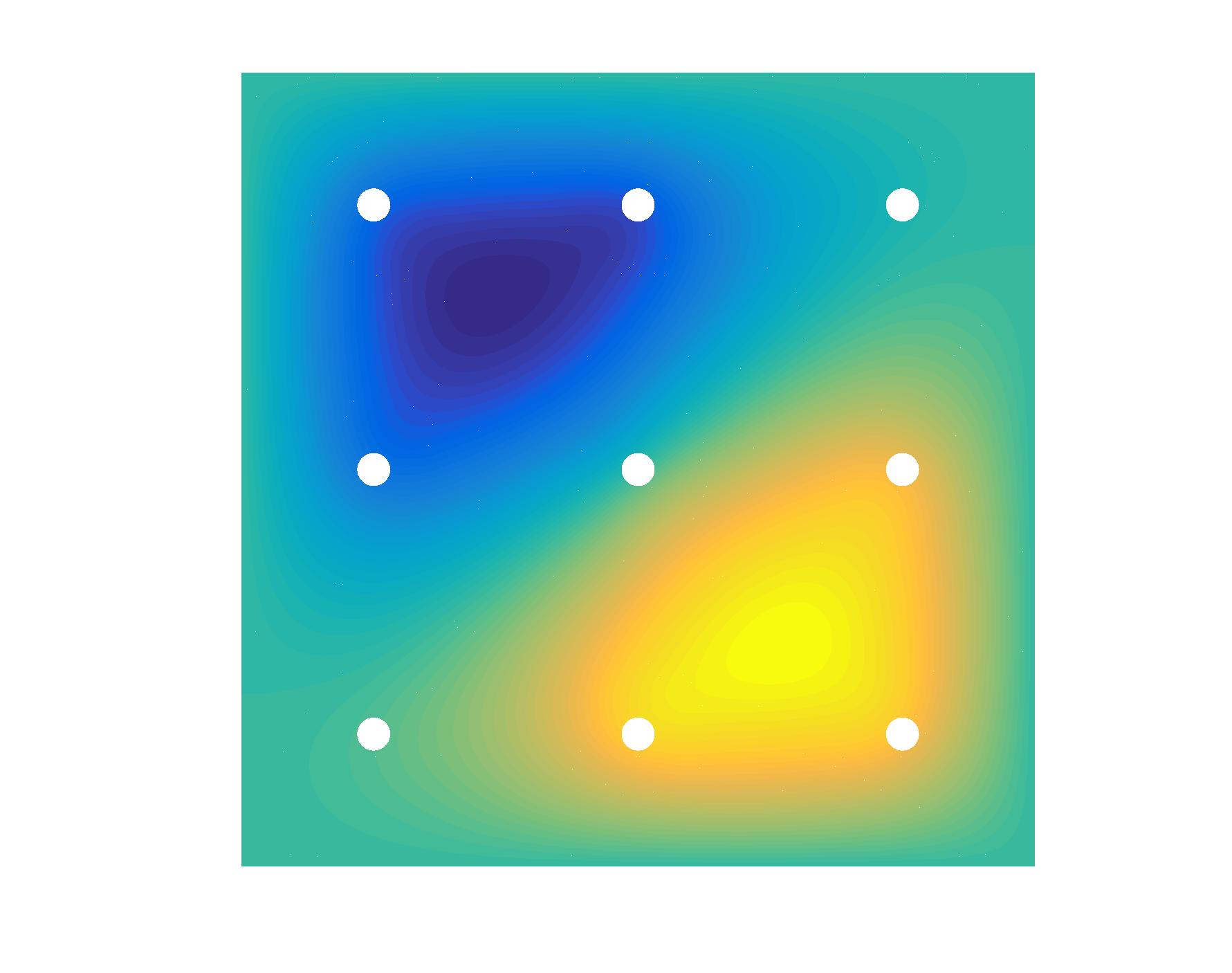}
		\caption{2nd mode, $ \lambda_{N}^{2} = 2.498 $}
	\end{subfigure}
	\begin{subfigure}[t]{0.3\textwidth}
		\centering
		\includegraphics[width=1\columnwidth]{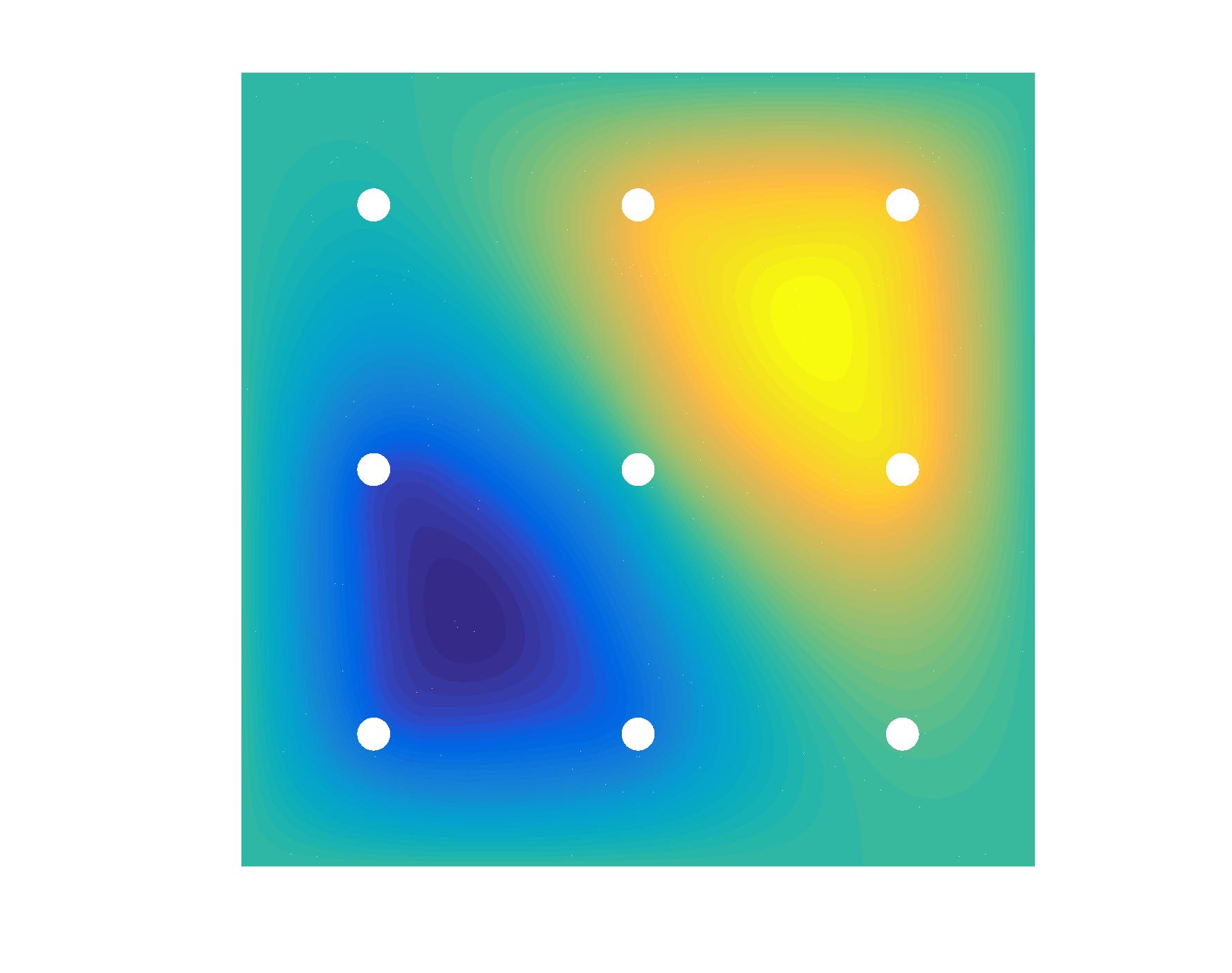}
		\caption{3rd mode, $ \lambda_{N}^{3} = 2.498 $}
	\end{subfigure}
	\begin{subfigure}[t]{0.3\textwidth}
		\centering
		\includegraphics[width=1\columnwidth]{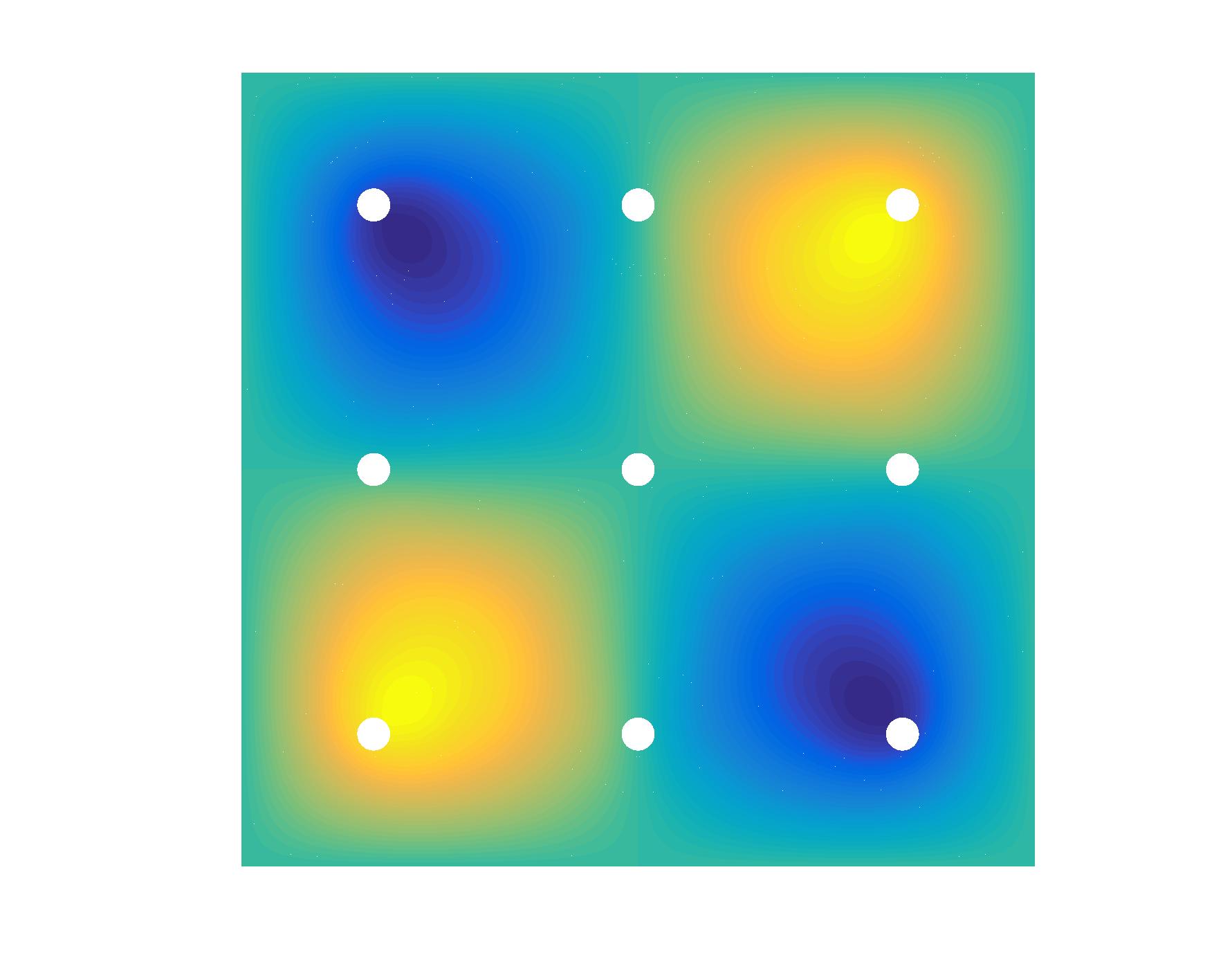}
		\caption{4th mode, $ \lambda_{N}^{4} = 4.013 $}
	\end{subfigure}
	\begin{subfigure}[t]{0.3\textwidth}
		\centering
		\includegraphics[width=1\columnwidth]{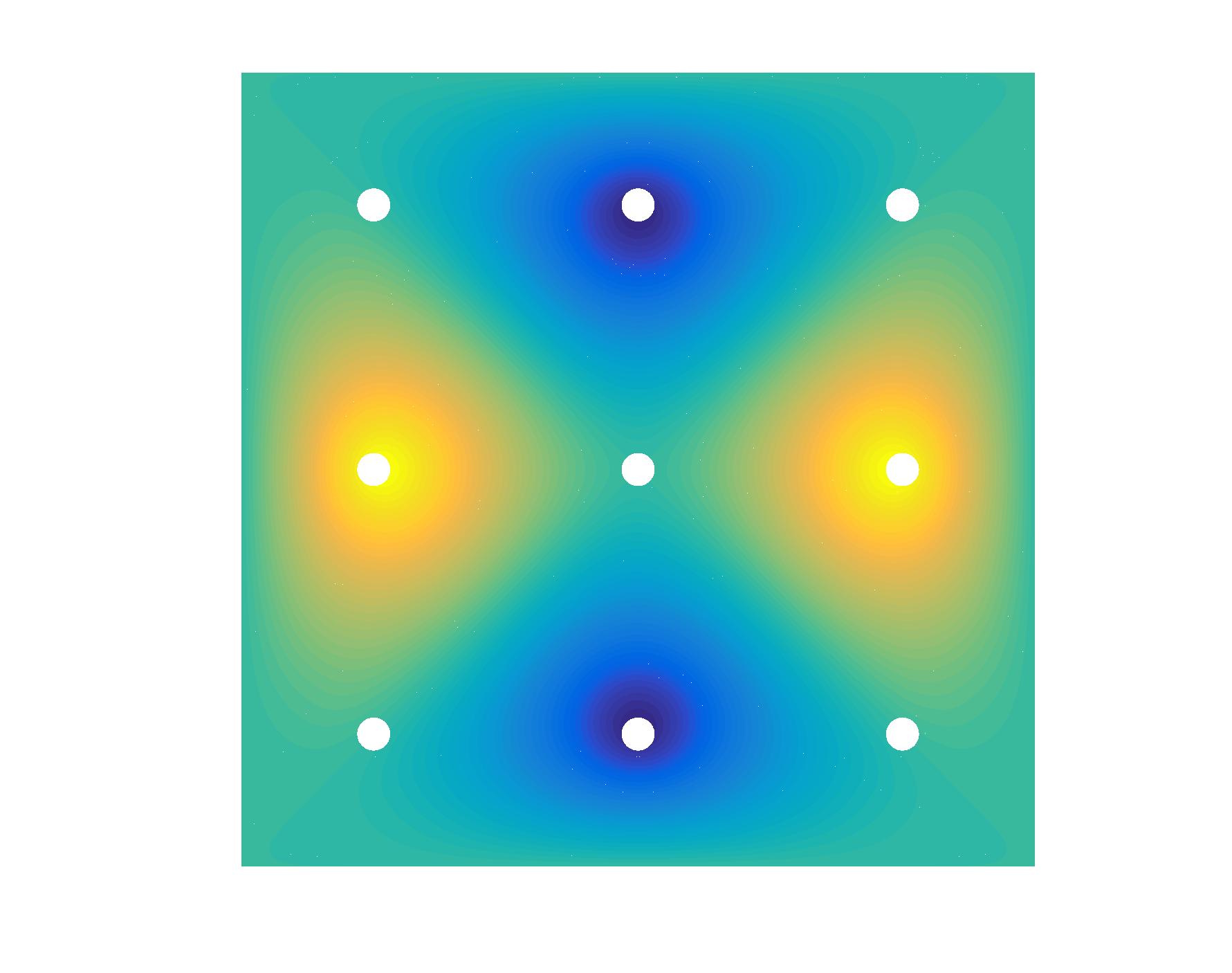}
		\caption{5th mode, $ \lambda_{N}^{5} = 4.616 $}
	\end{subfigure}
	\begin{subfigure}[t]{0.3\textwidth}
		\centering
		\includegraphics[width=1\columnwidth]{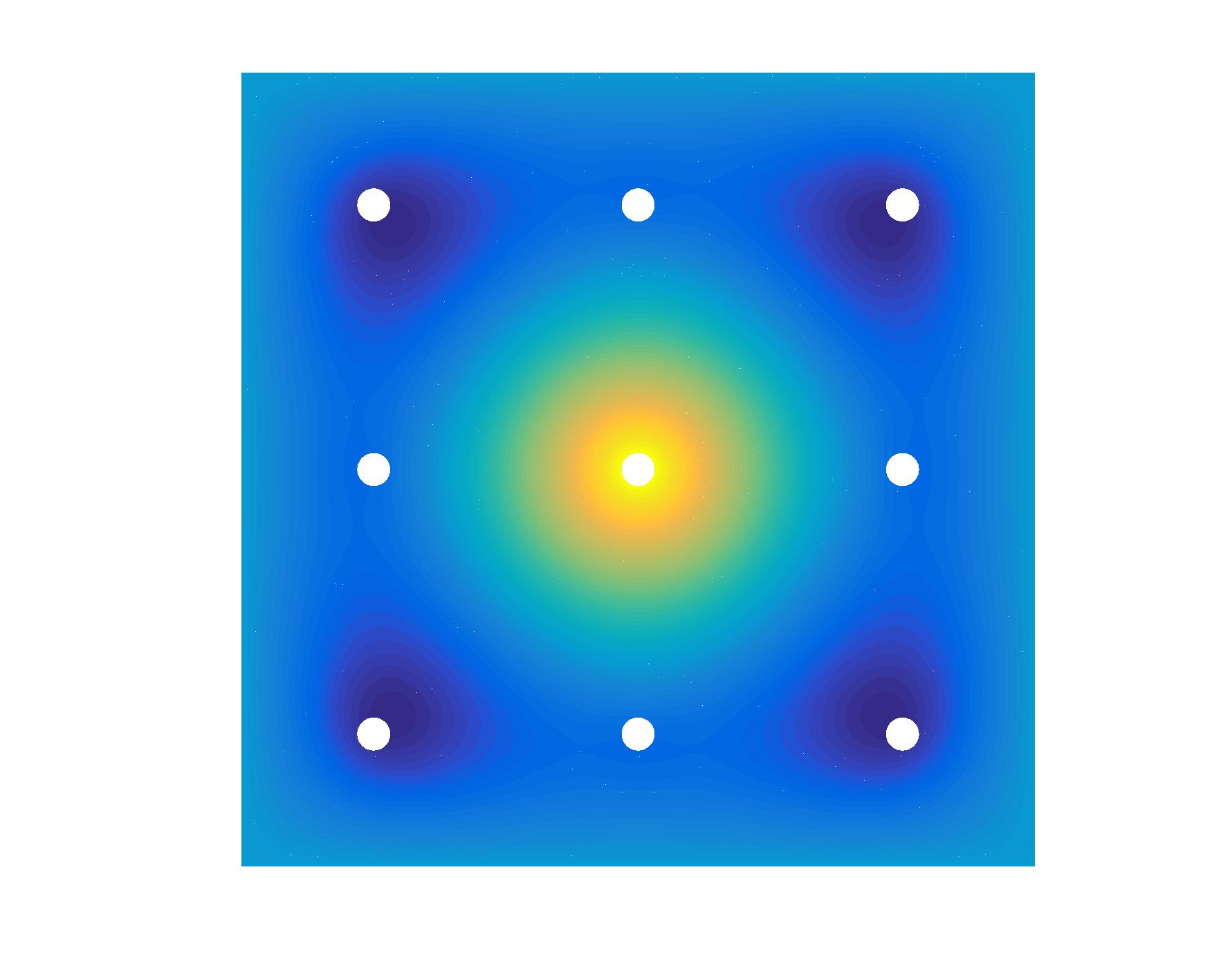}
		\caption{6th mode, $ \lambda_{N}^{6} = 4.618 $}
	\end{subfigure}
	\begin{subfigure}[t]{0.3\textwidth}
		\centering
		\includegraphics[width=1\columnwidth]{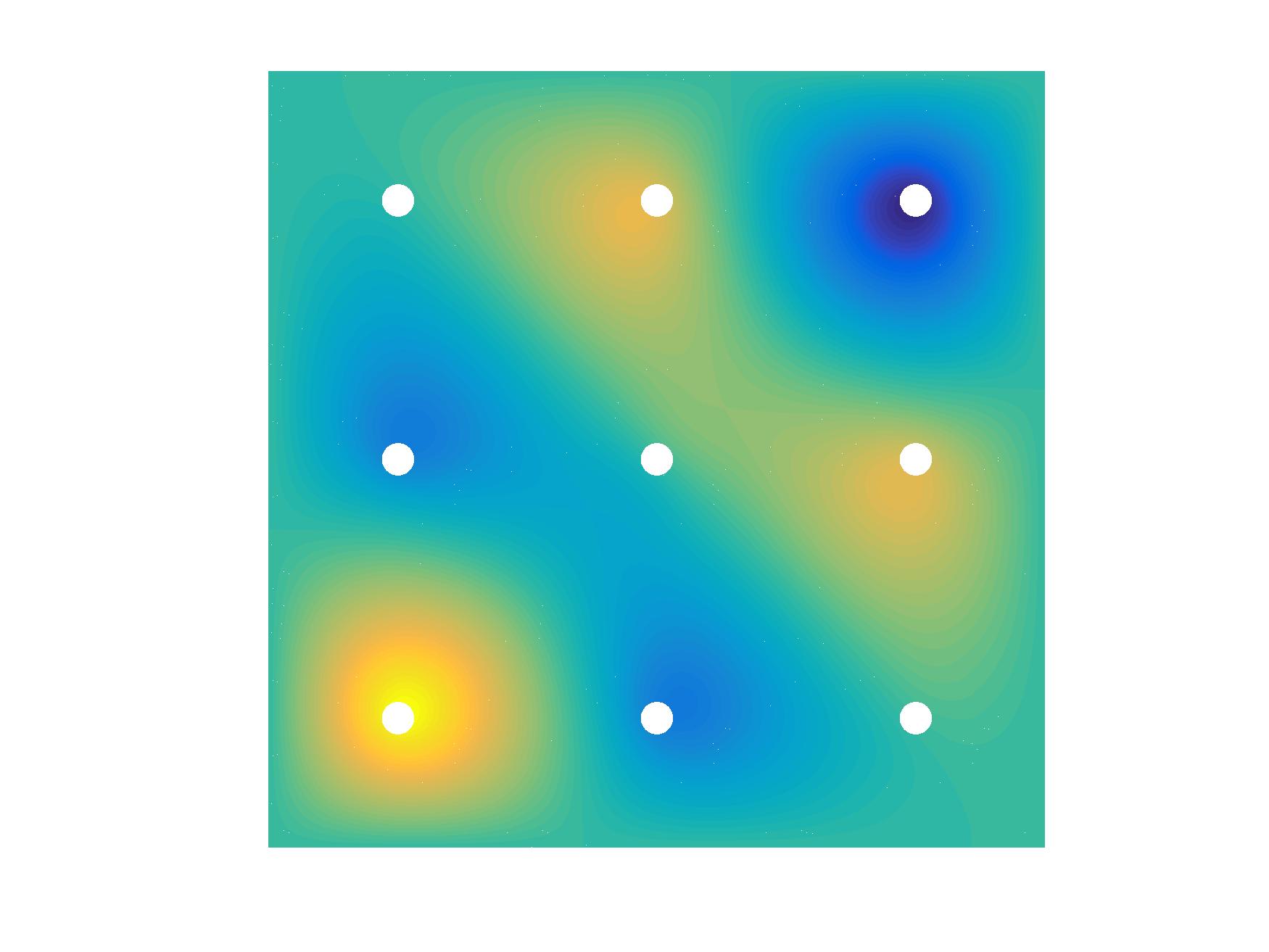}
		\caption{7th mode, $ \lambda_{N}^{7} = 6.166 $}
	\end{subfigure}
	\begin{subfigure}[t]{0.3\textwidth}
		\centering
		\includegraphics[width=1\columnwidth]{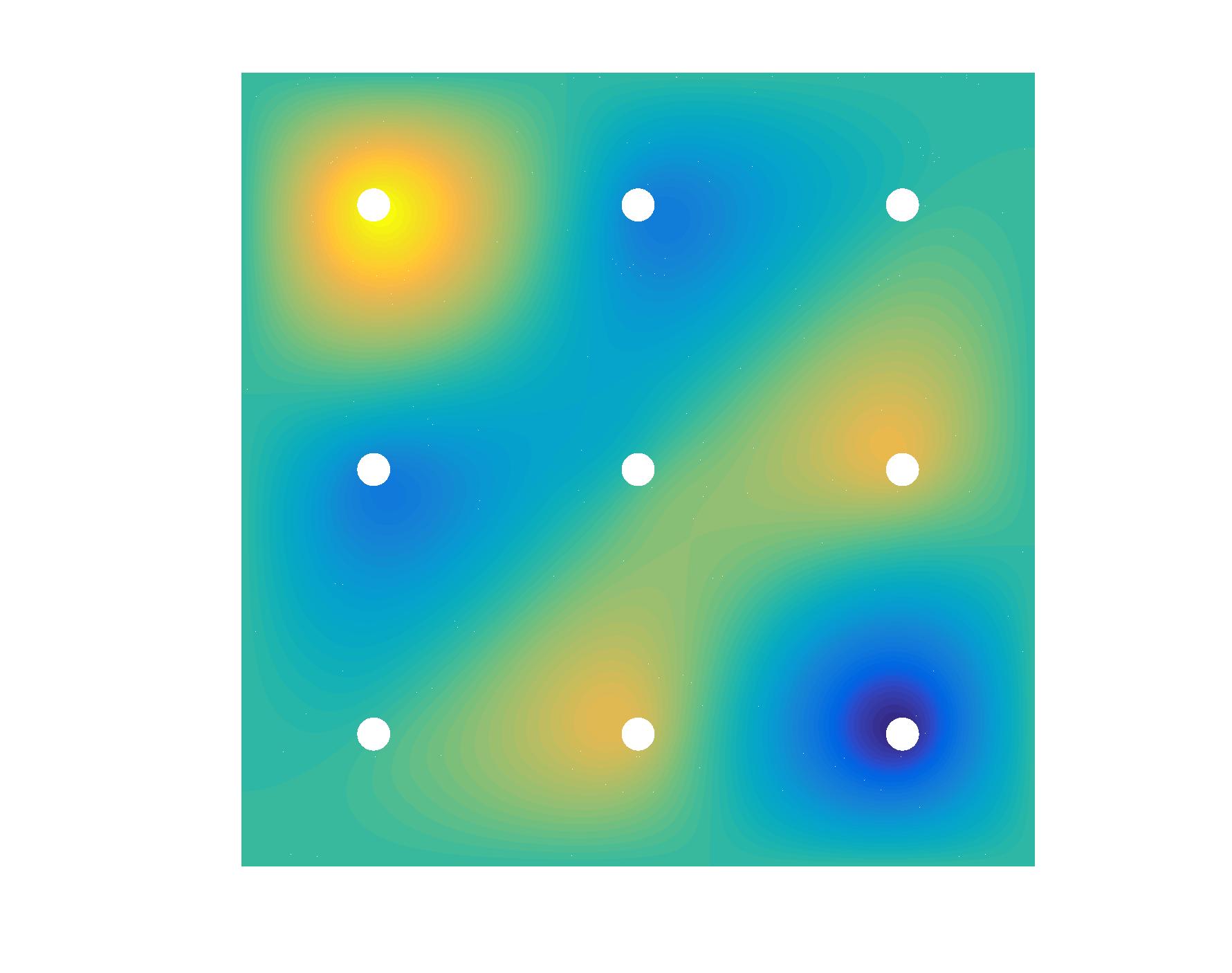}
		\caption{8th mode, $ \lambda_{N}^{8} = 6.166 $}
	\end{subfigure}
	\begin{subfigure}[t]{0.3\textwidth}
		\centering
		\includegraphics[width=1\columnwidth]{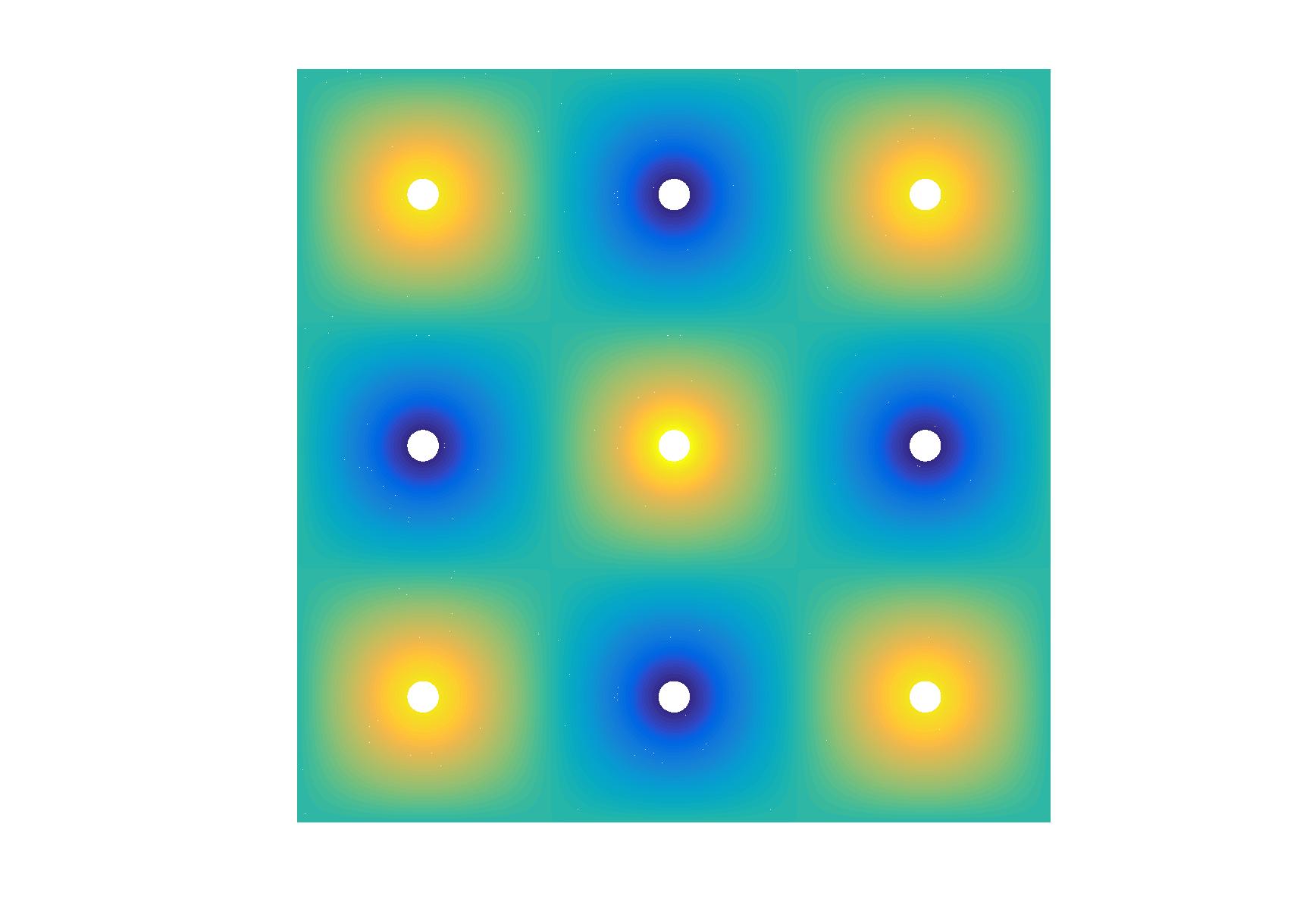}
		\caption{9th mode, $ \lambda_{N}^{9} = 7.874 $}
	\end{subfigure}	
	\begin{subfigure}[t]{0.3\textwidth}
		\centering
		\includegraphics[width=1\columnwidth]{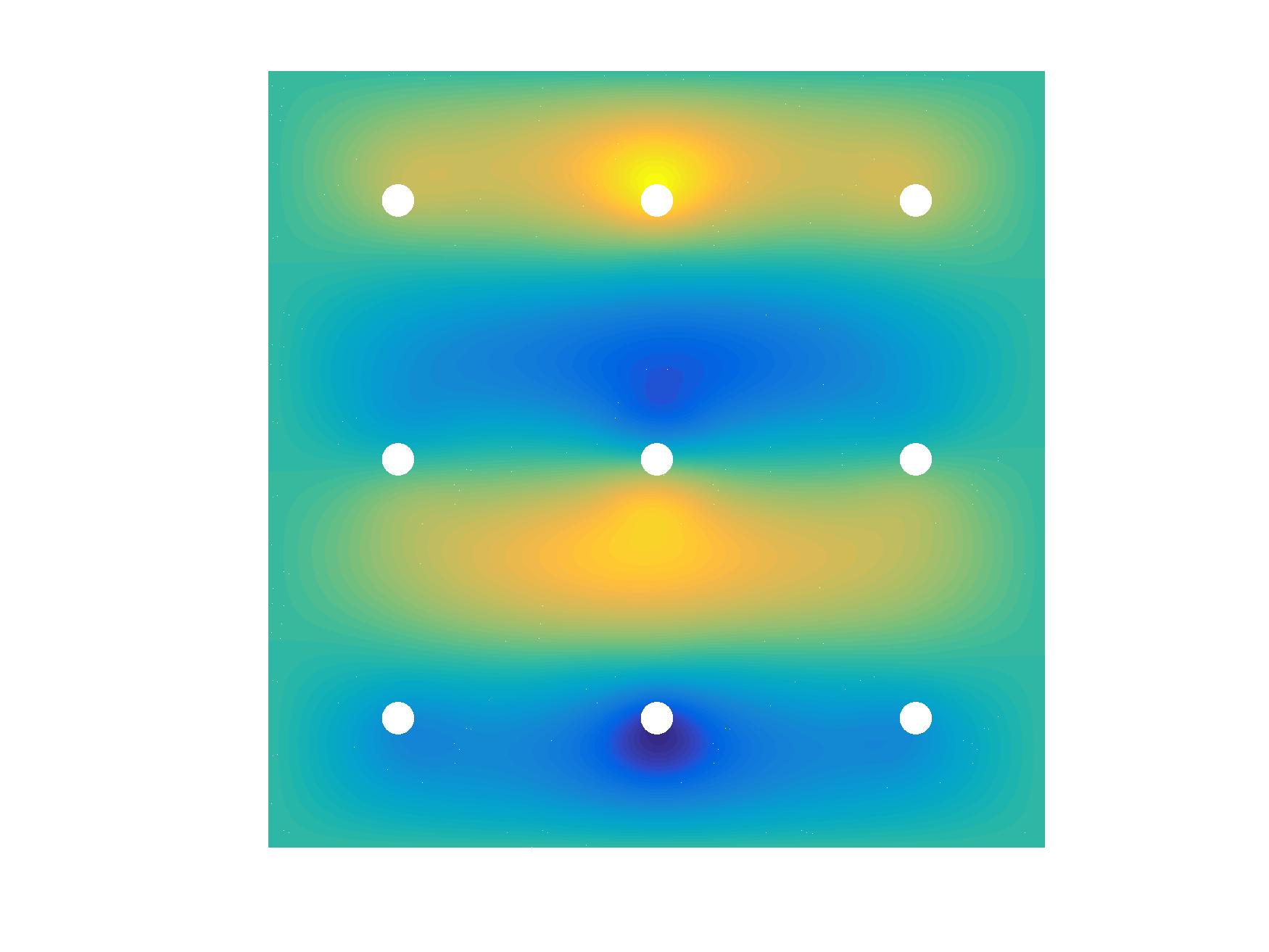}
		\caption{10th mode, $ \lambda_{N}^{10} = 8.383 $}
	\end{subfigure}
	\begin{subfigure}[t]{0.3\textwidth}
		\centering
		\includegraphics[width=1\columnwidth]{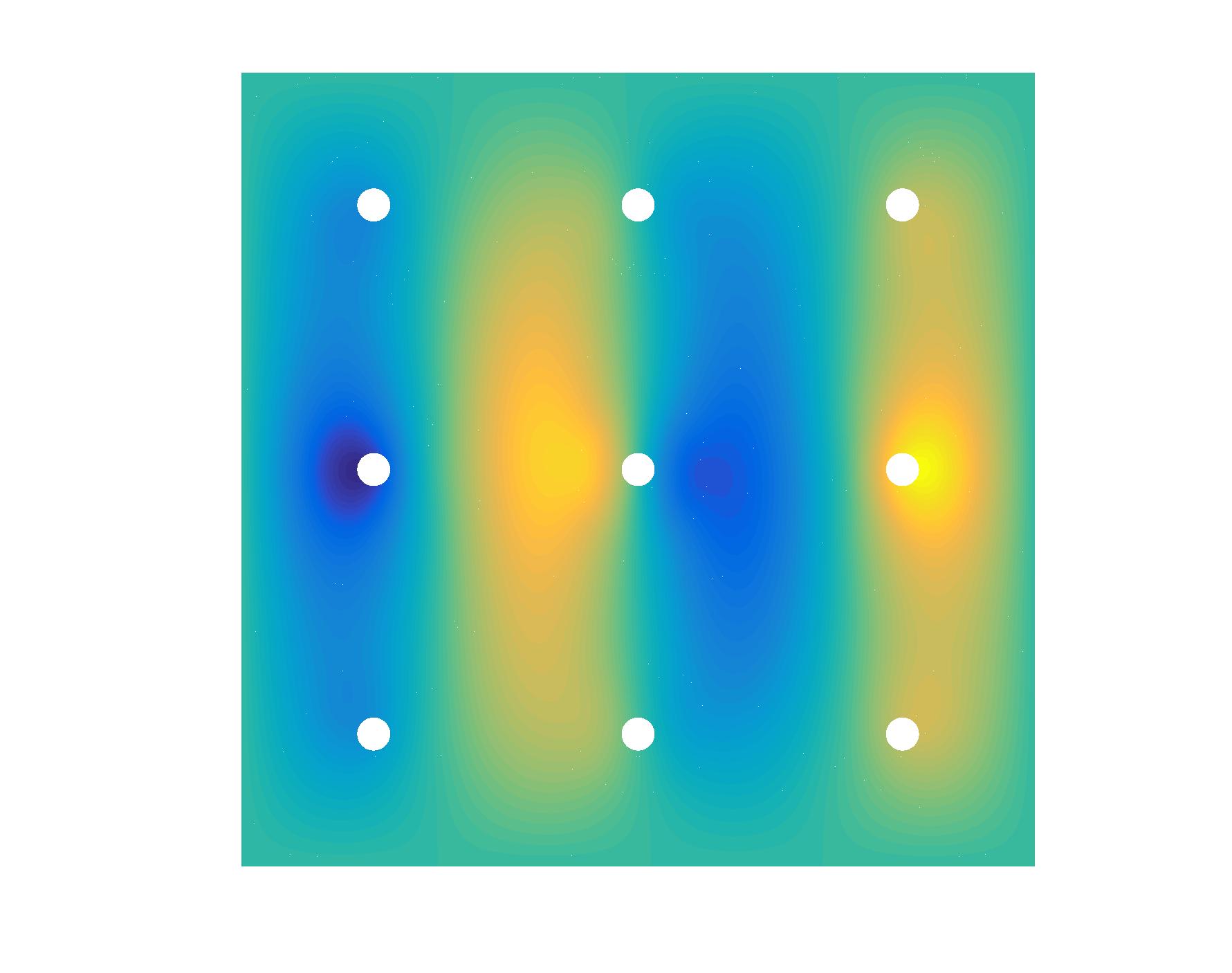}
		\caption{11th mode, $ \lambda_{N}^{11} = 8.383 $}
	\end{subfigure}
	\caption{Vibration of mode shapes for structure of plate with 9 holes.}
	\label{fig:plateholesVmode}
\end{figure}

\begin{figure}
	\centering
	\begin{subfigure}[t]{0.3\textwidth}
		\centering
		\raisebox{.4cm}{\includegraphics[width=0.65\columnwidth]{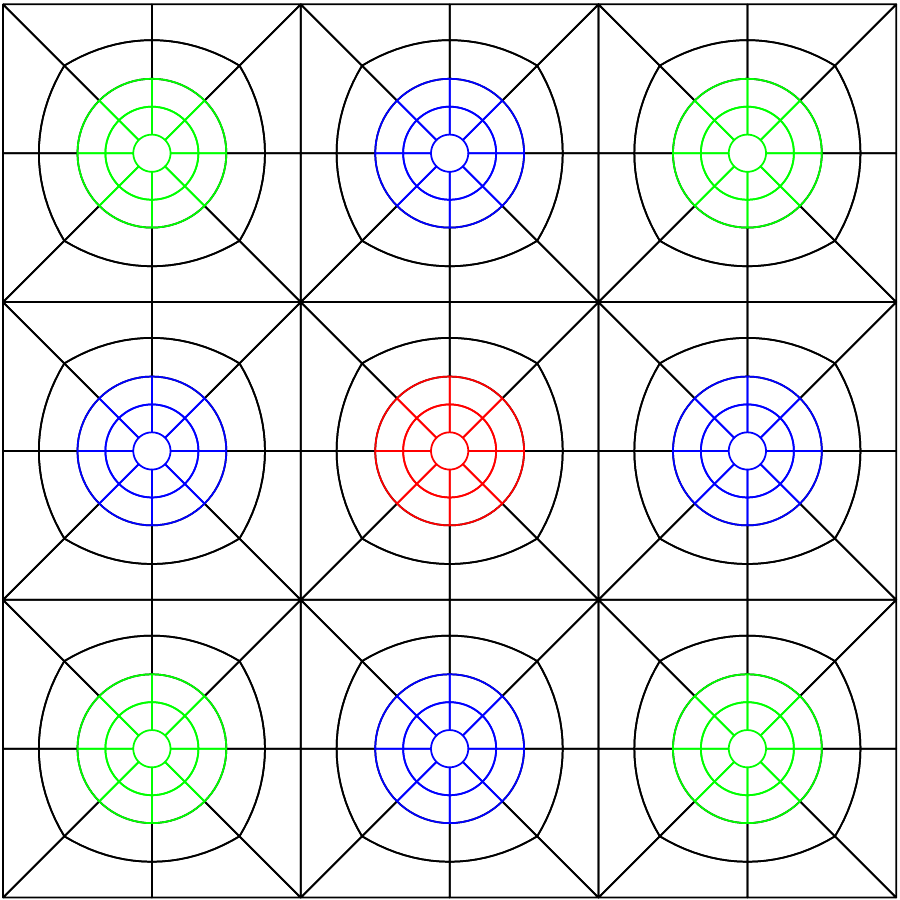}}
		\caption{Initial}
	\end{subfigure}
	\begin{subfigure}[t]{0.3\textwidth}
		\centering
		\raisebox{0cm}{\includegraphics[width=1.05\columnwidth]{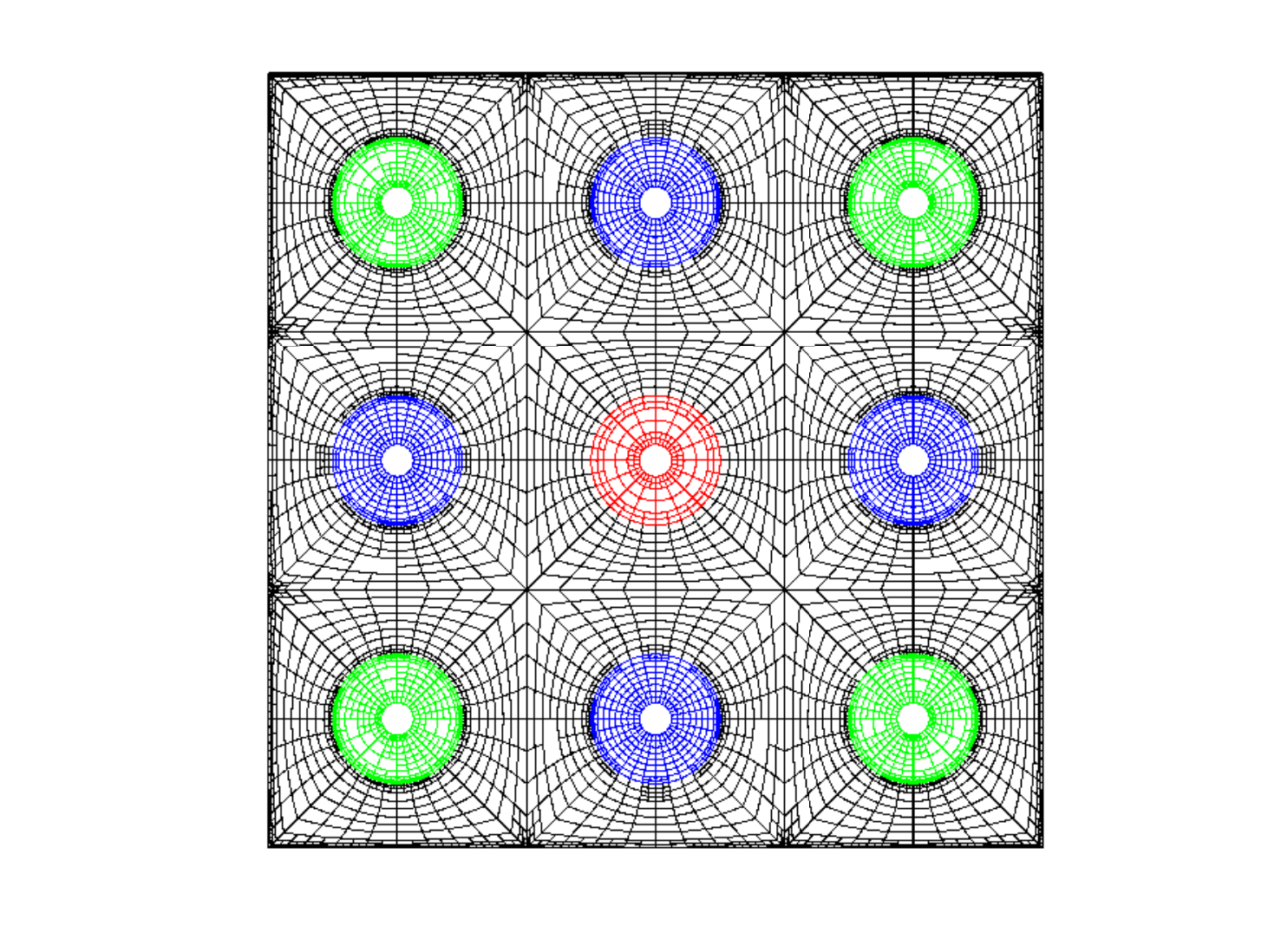}}
		\caption{1th mode}
	\end{subfigure}
	\begin{subfigure}[t]{0.3\textwidth}
		\centering
		\includegraphics[width=1\columnwidth]{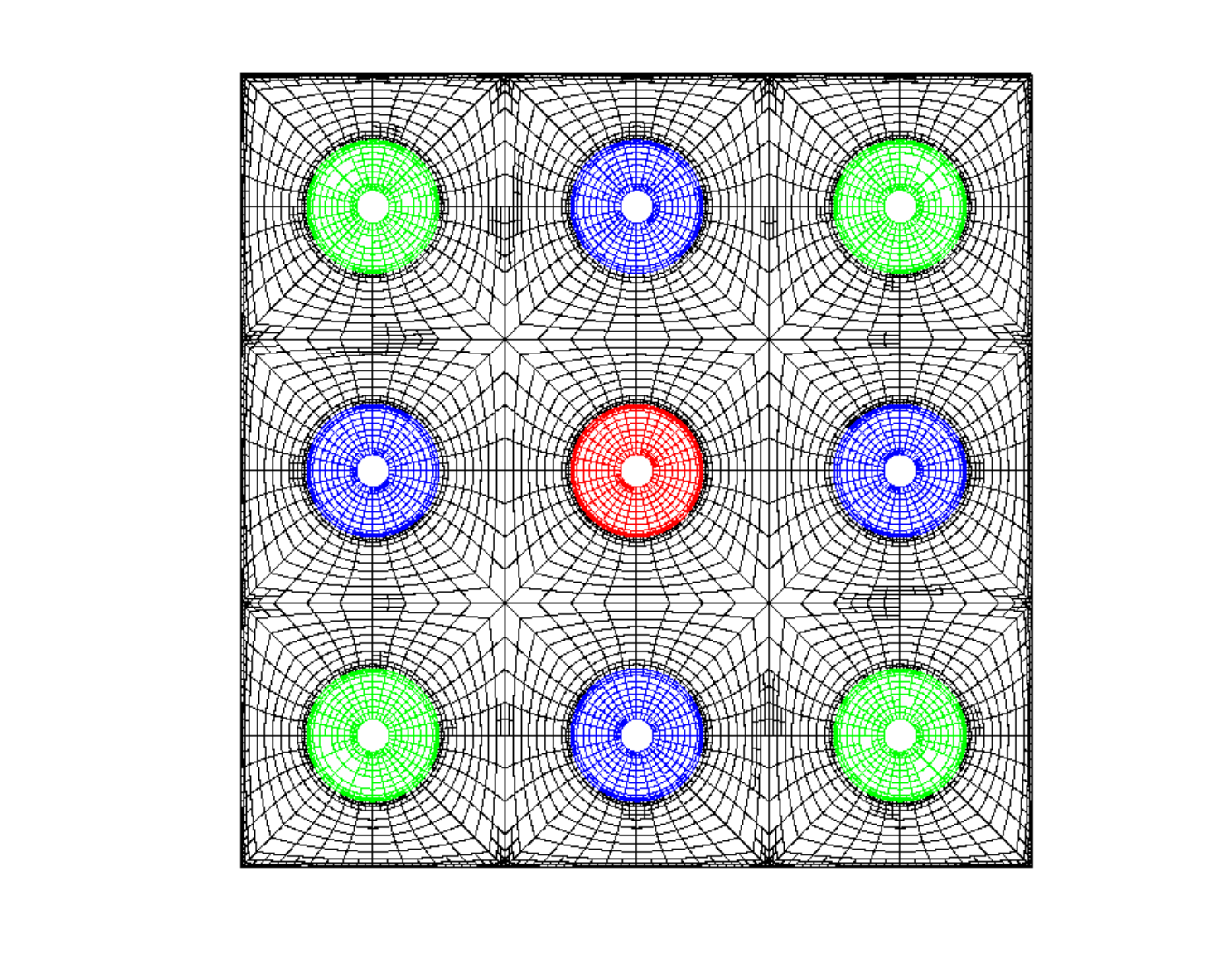}
		\caption{2-3, 4th modes}
	\end{subfigure}
	\begin{subfigure}[t]{0.3\textwidth}
		\centering
		\includegraphics[width=1\columnwidth]{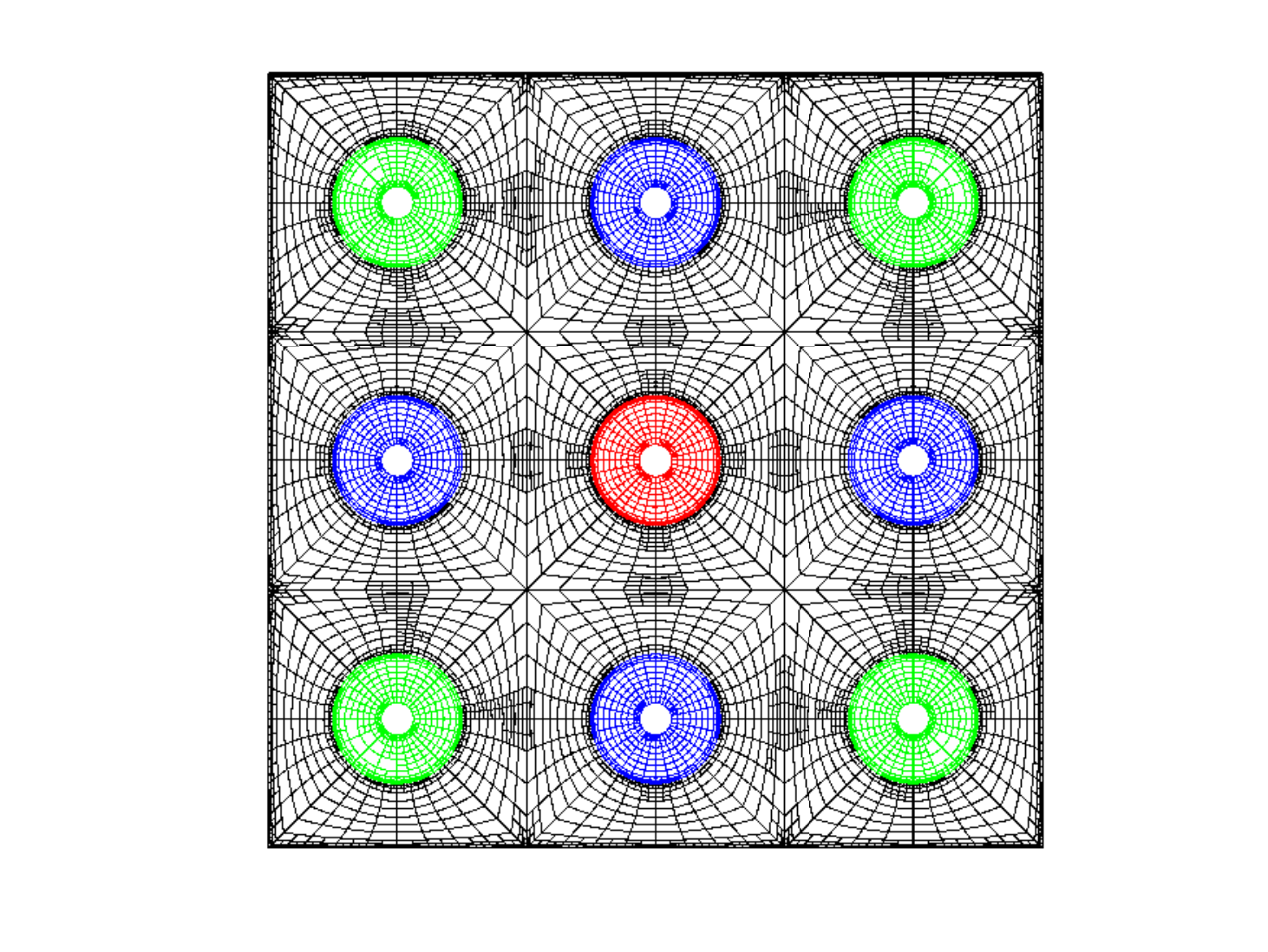}
		\caption{5, 6th modes}
	\end{subfigure}
	\begin{subfigure}[t]{0.3\textwidth}
		\centering
		\raisebox{-0.05cm}{\includegraphics[width=1.08\columnwidth]{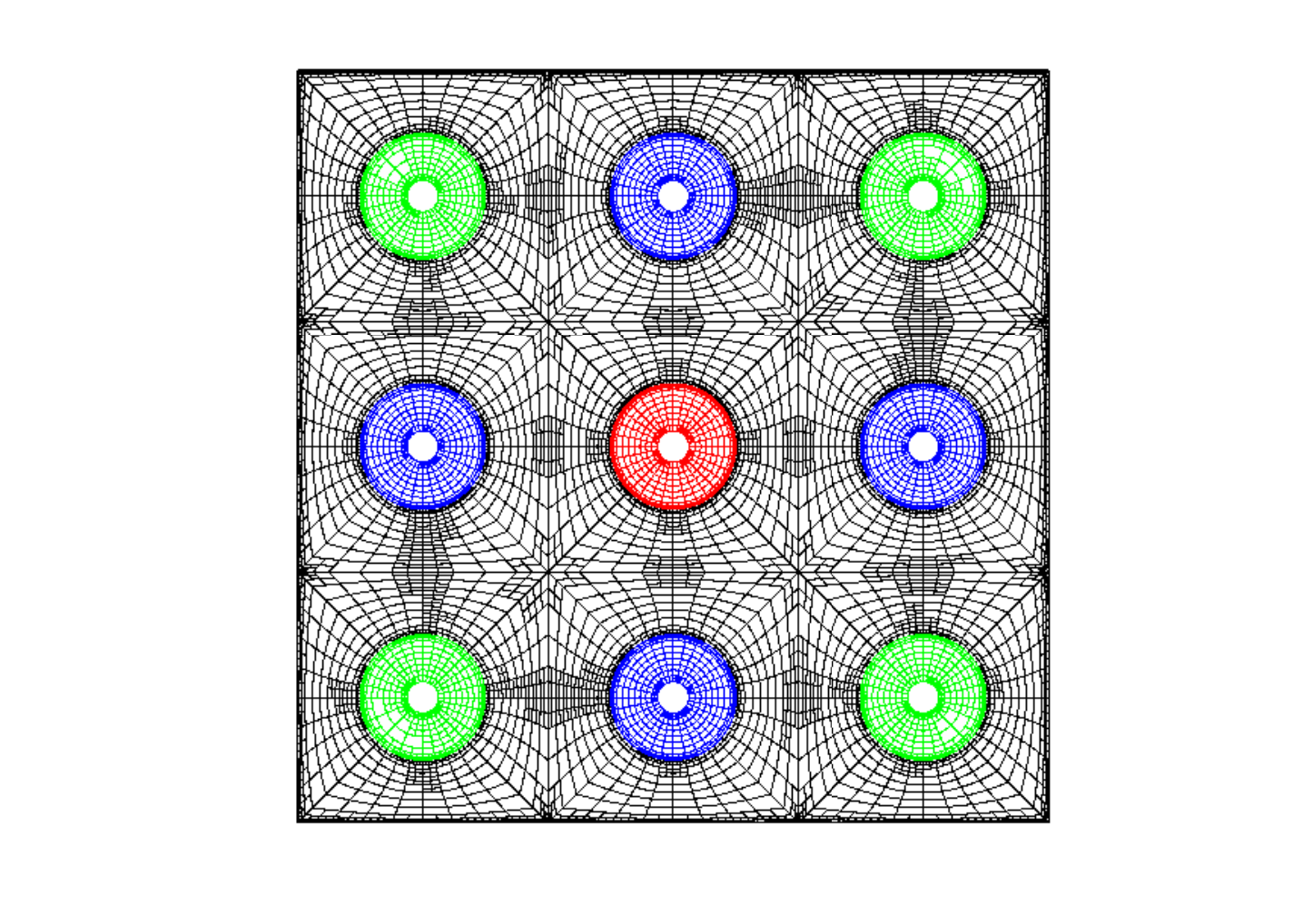}}
		\caption{7-8, 9th modes}
	\end{subfigure}
	\begin{subfigure}[t]{0.3\textwidth}
		\centering
		\raisebox{-0.1cm}{\includegraphics[width=1\columnwidth]{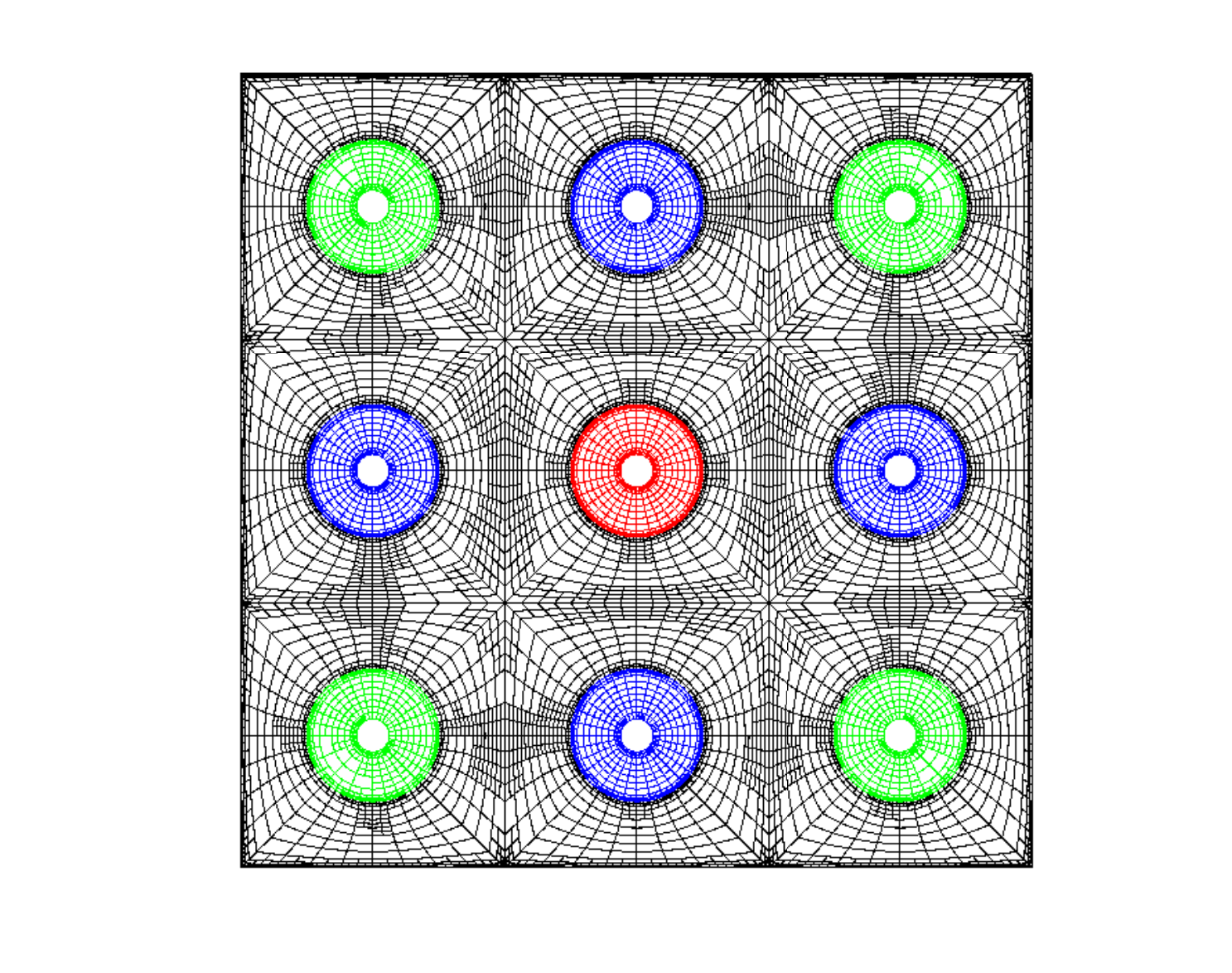}}
		\caption{10-11th modes}
	\end{subfigure}					
	\caption{Adaptive refinement process of 1-11th modes for vibration of the plate with holes, with refinement level $L_{e} = 1$.}
	\label{fig:adpplateholes}
\end{figure}

\begin{figure}
	\centering
	\begin{subfigure}[t]{0.4\textwidth}
		\centering
		\includegraphics[width=1\columnwidth]{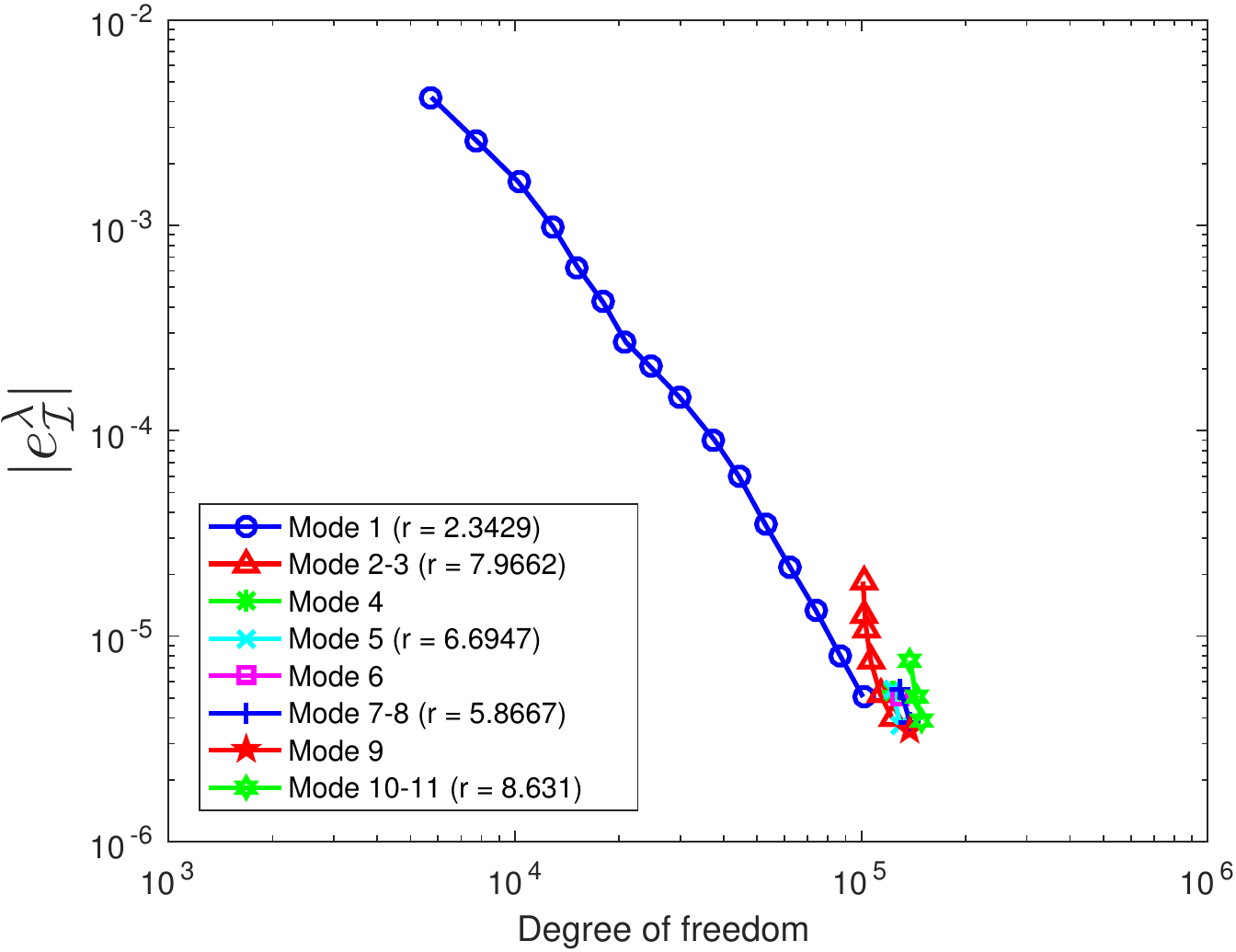}
		\caption{}
	\end{subfigure}
	\hspace{1cm}
	\begin{subfigure}[t]{0.43\textwidth}
		\centering
		\includegraphics[width=1\columnwidth]{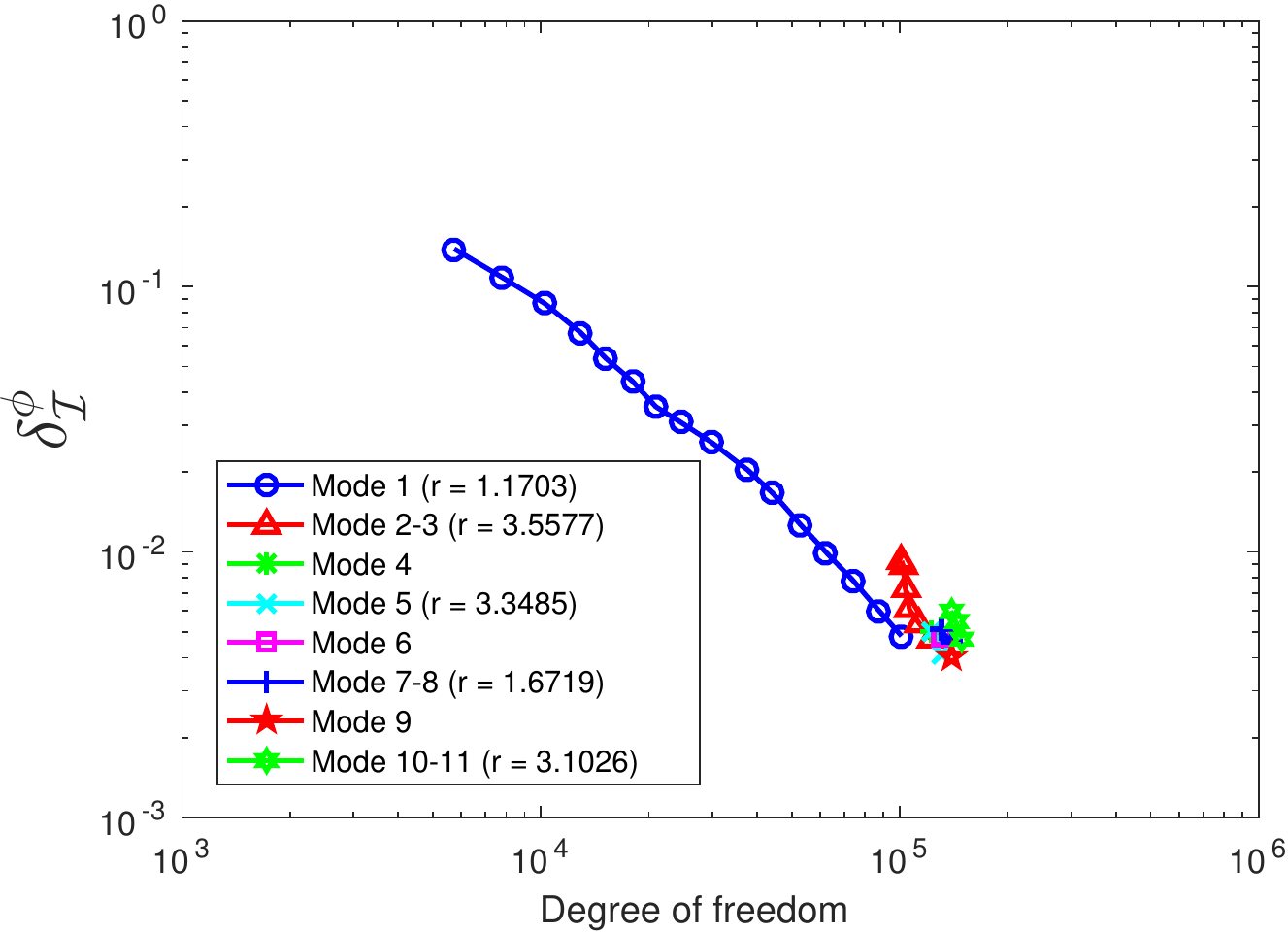}
		\caption{}
	\end{subfigure}
	\caption{The study of convergence of error estimators (a)\ $|e_{\vect{\mathcal{I}}}^{\lambda}|$ and (b)\ $\delta_{\vect{\mathcal{I}}}^{\phi}$ from the 1st to the 11th mode of the plate with holes.}
	\label{fig:conmodesltoh}
\end{figure}

\begin{figure}	
	\centering
	\begin{subfigure}[t]{0.4\textwidth}
		\centering
		\includegraphics[width=1\columnwidth]{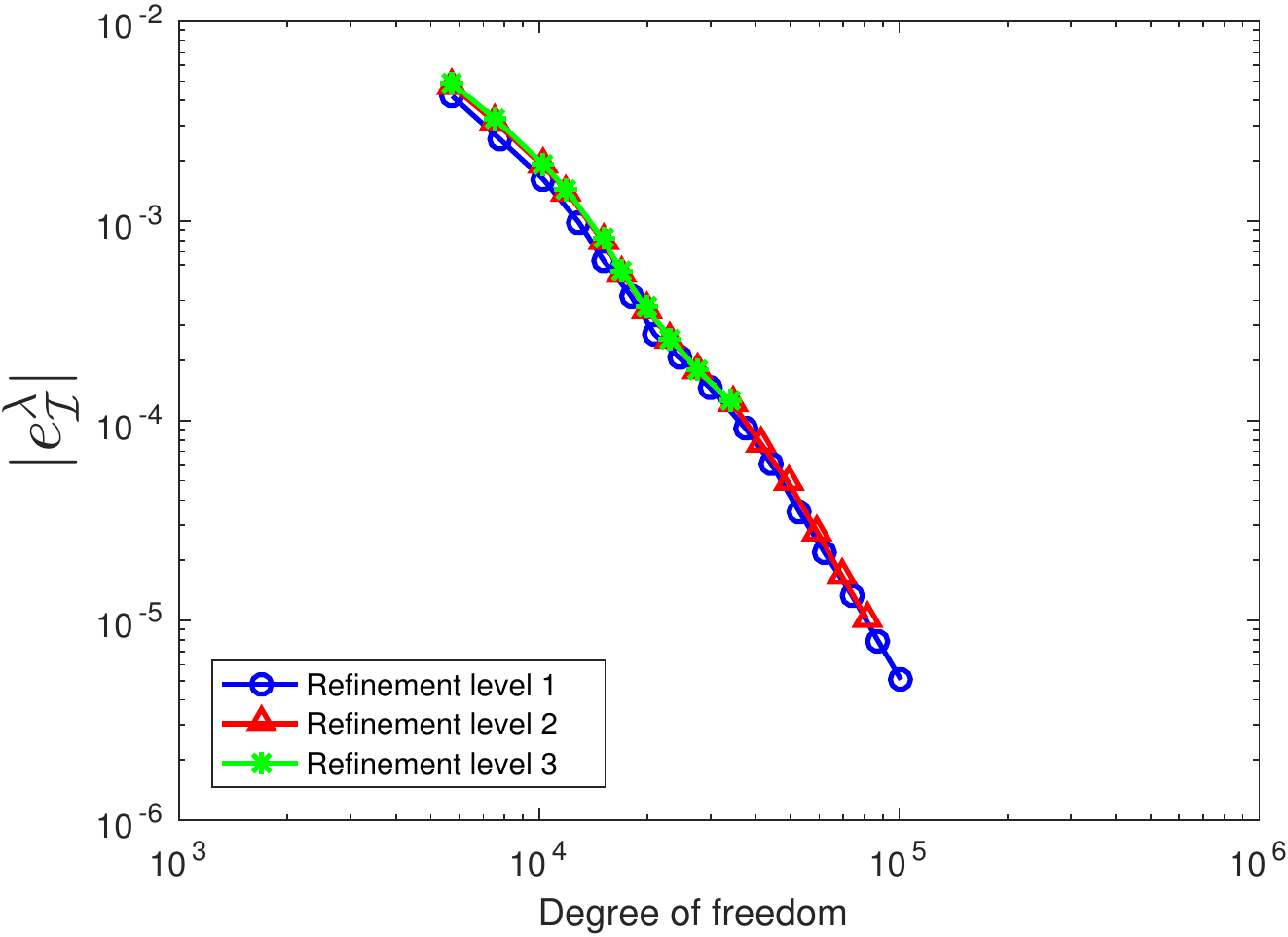}
		\caption{eigenvalue}
	\end{subfigure}
	\hspace{1cm}
	\begin{subfigure}[t]{0.4\textwidth}
		\centering
		\includegraphics[width=1\columnwidth]{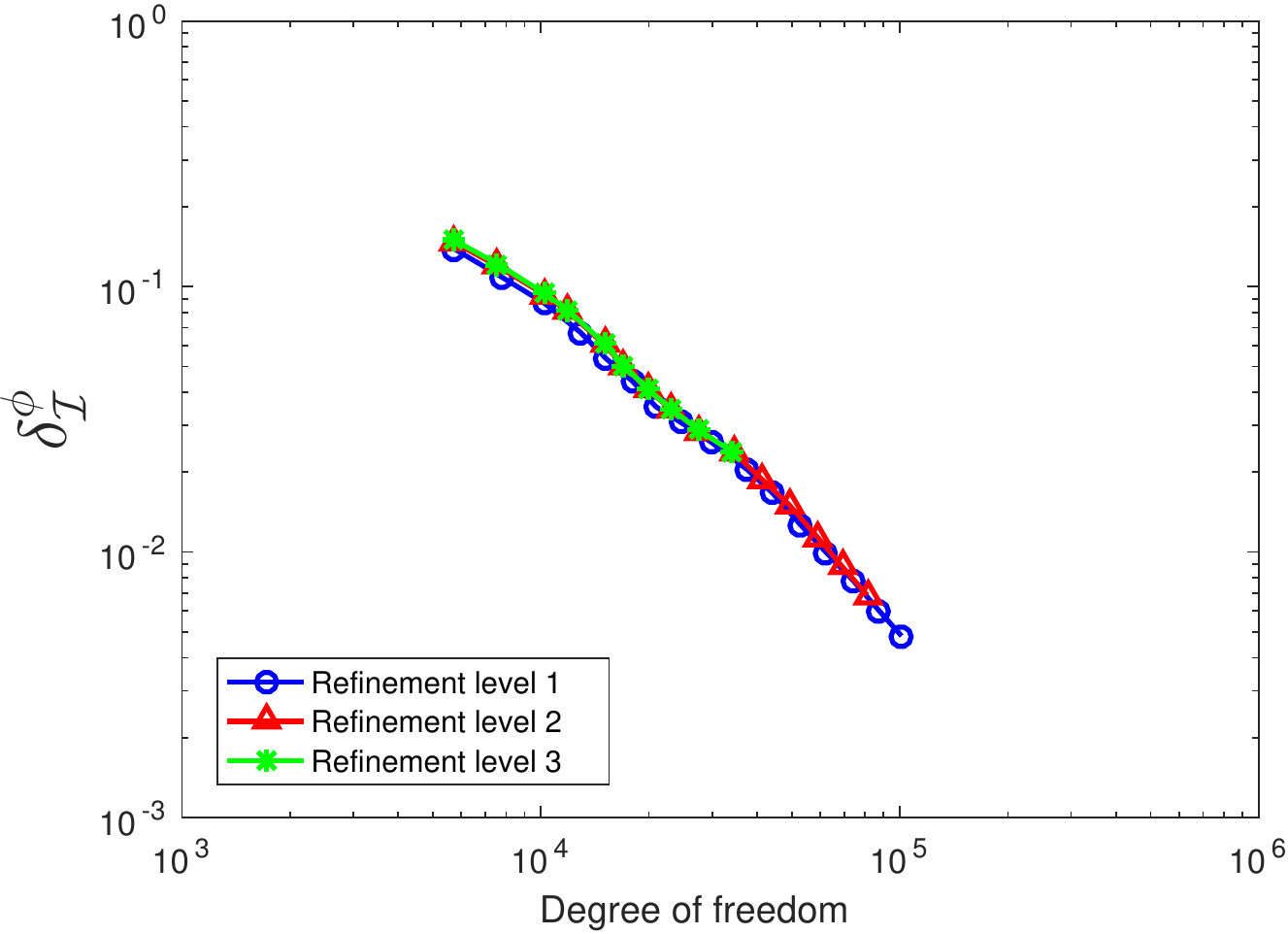}
		\caption{eigenvector}
	\end{subfigure}
	\caption{The Comparison of error estimators (a)\ $|e_{\vect{\mathcal{I}}}^{\lambda}|$ and (b)\ $\delta_{\vect{\mathcal{I}}}^{\phi}$ at 1st mode obtained by different refinement levels: $L_{e} = 1, L_{e} = 2, L_{e} = 3$.}
	\label{fig:conlevel}
\end{figure}

\begin{figure}	
	\centering
	\begin{subfigure}[t]{0.4\textwidth}
		\centering
		\includegraphics[width=1\columnwidth]{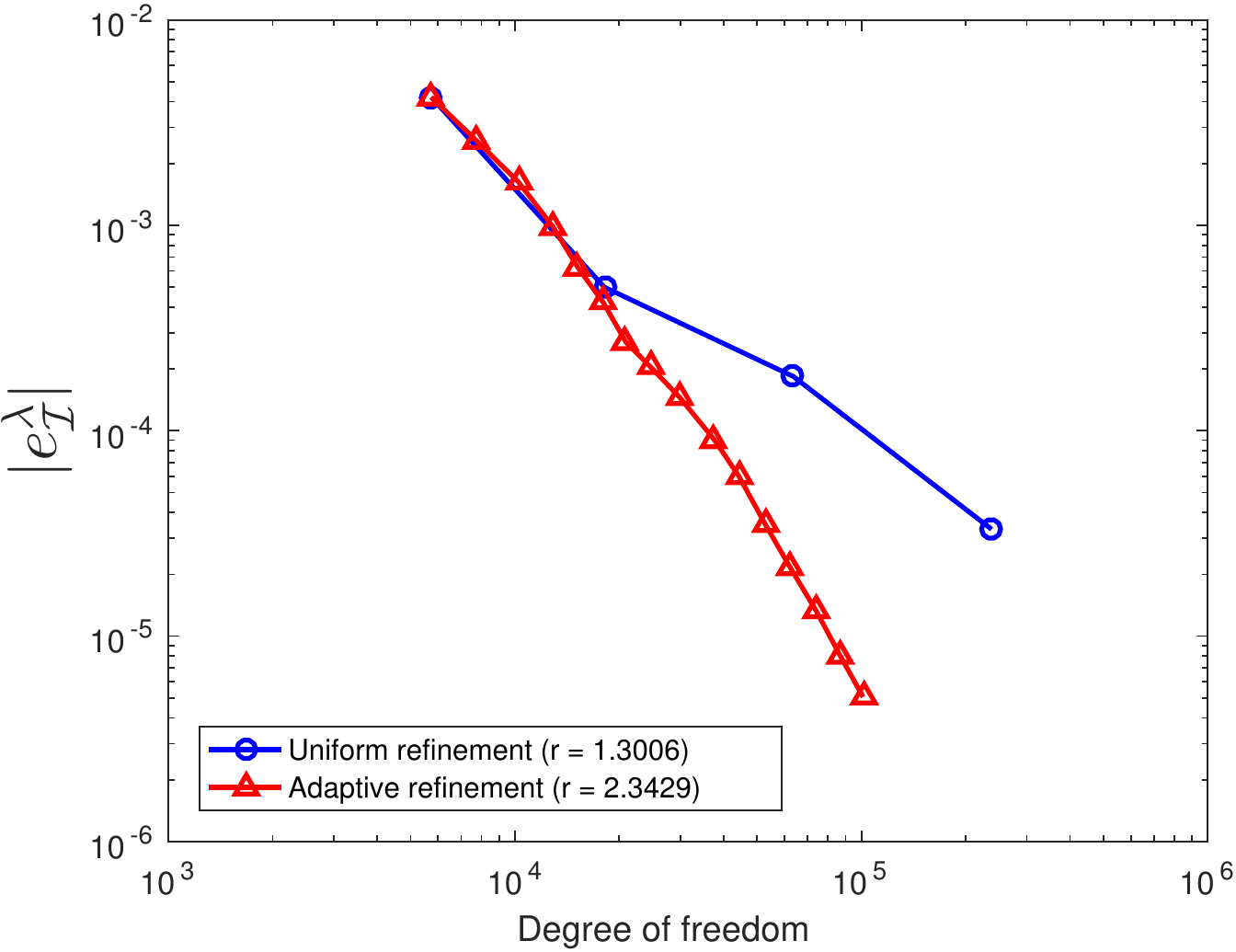}
		\caption{}
	\end{subfigure}
	\hspace{1cm}
	\begin{subfigure}[t]{0.4\textwidth}
		\centering
		\includegraphics[width=1\columnwidth]{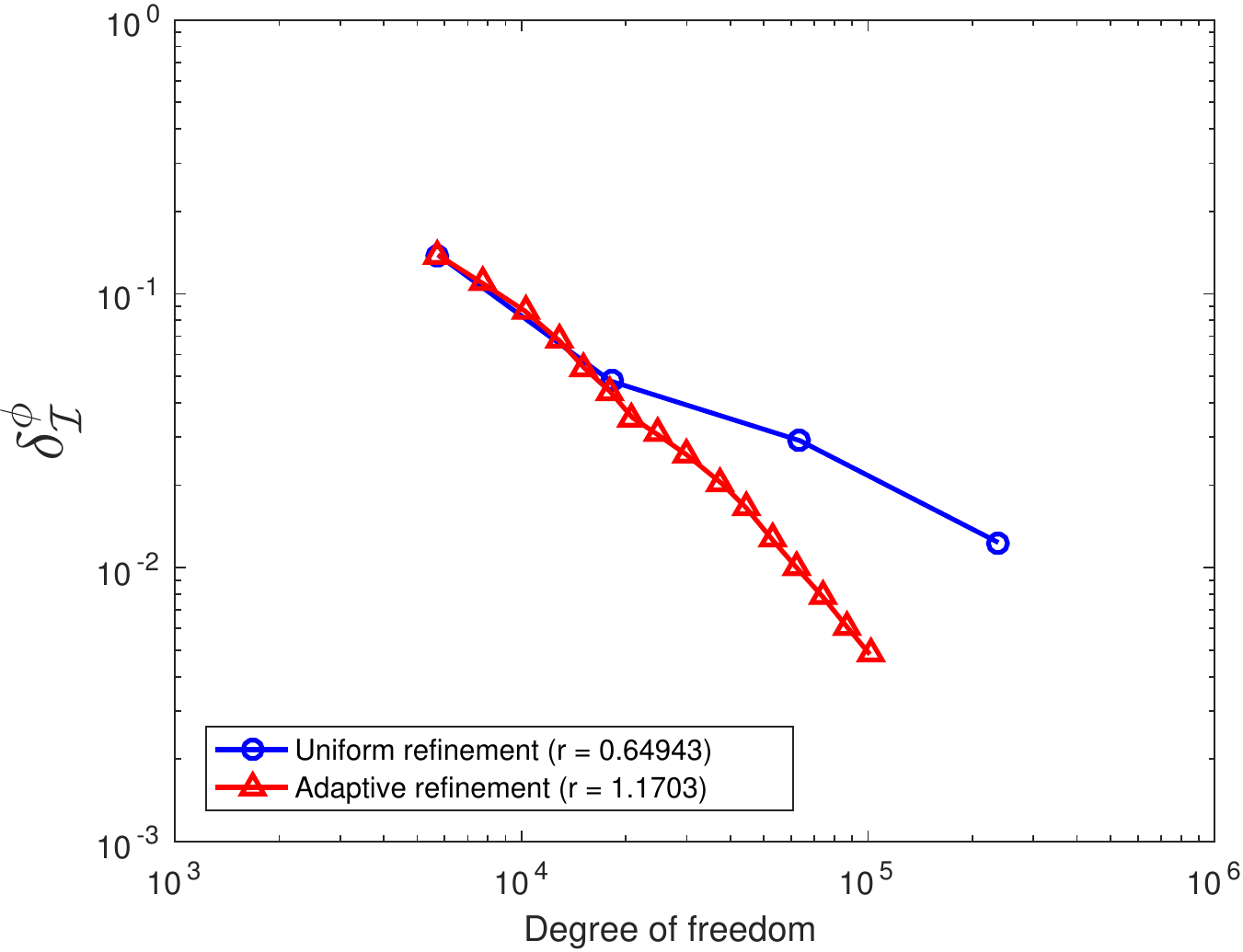}
		\caption{}
	\end{subfigure}
	\caption{The Comparison of error estimators (a)\ $|e_{\vect{\mathcal{I}}}^{\lambda}|$ and (b)\ $\delta_{\vect{\mathcal{I}}}^{\phi}$ at 1st mode obtained by adaptive refinement and uniform refinement.}
	\label{fig:conadpvsuni}	
\end{figure}

The geometry and discretization of a plate with 9 holes are presented in Fig.\ref{plateholesGeo}. The colored patches around holes are softer than black patches. The density is set as $\rho =1$ and thickness is $h=0.1$. The simply supported boundary conditions are imposed at the edges of square plate, and the edges of holes are kept free. As mentioned in Section \ref{sec:adpfreq}, the GIFT method combined with MAC is applied to guide the adaptivity for the band of frequency in Fig.\ref{fig:winplatehole}, using Algorithm \ref{Algm:adpFoI} by sweeping modes. It means that the initial mesh of adaptivity for mode(s) $(\vect{\mathcal{I}}+1)$, or $\vect{\mathcal{I}}+n$ for $n$ multiple modes, inherits the mesh at completion phase of adaptivity for mode $\vect{\mathcal{I}}$ (defined in Eq.\eqref{eq:mode_mul}). Due to the symmetric geometry and material characteristics, it is not difficult to find in Fig.\ref{fig:plateholesVmode} that double modes arise at modes (2,3), modes (7,8) and modes (10,11). Note that although mode 5 and mode 6 are very close, they are not double modes. Besides, it is observed that global modal shapes dominate from 1st to 3rd mode which leads to the nearly global refinement in adaptivity, as shown from Fig.\ref{fig:adpplateholes}(a) to Fig.\ref{fig:adpplateholes}(c). Afterwards, as it can be seen in Fig.\ref{fig:plateholesVmode}(d)-(i), local vibration gradually appears, which results in local refinements at areas around the holes. The vibrations at 10th and 11th modes distribute like an orthogonal crossing on the plate (see Fig.\ref{fig:plateholesVmode}(j,k)), and then adaptive refinement follows horizontally and vertically at conjunctions of multiple patches in Fig.\ref{fig:adpplateholes}(f). This implies that the accuracy of these $C^{0}$ coupling fields is expected to be improved. The convergence of error estimators $|e_{\vect{\mathcal{I}}}^{\lambda}|$ (defined in Eq.\eqref{eq:errvect_mul}) and $\delta_{\vect{\mathcal{I}}}^{\phi}$ (defined in Eq.\eqref{eq:errfreq_mul}) is investigated, and they both have excellent convergence from low to high modes, as seen in Fig.\ref{fig:conmodesltoh}. Note that, the plot of convergence rate for some modes (such as 4th, 6th and 9th modes) is a point, since the $|e_{\vect{\mathcal{I}}}^{\lambda}|$ and $\delta_{\vect{\mathcal{I}}}^{\phi}$ at these modes have satisfied the tolerances in Eq.\eqref{eq:tol} at the beginning of the adaptivity. Hence, no refinement is required for these modes. As illustrated in Fig.\ref{fig:conlevel}, the convergence obtained by different level of $h-$refinements is almost consistent, indicating the precision of the selected level of refinement $L_{e} = 1 $ for the computation of error estimators is sufficient. Eventually, Fig.\ref{fig:conadpvsuni} depicts an improved convergent rate achieved by local adaptive refinement, as compared with the global uniform PHT refinement generated in the GIFT framework.

\section{Conclusion}\label{sec:conclusion}
In this article, we presented a strategy of local adaptivity for the structural vibrations of Reissner-Mindlin plate. Within the context of the GIFT framework, we may make use of the geometrical descriptors given by CAD directly, and independently apply PHT splines for analysis, which allows for local refinement to be performed without the limitation occuring when using tensor-product-based shape functions. The adaptivity algorithm is fully automatised, and relies on a hierarchical posteriori error estimation strategy that makes the best of the PHT-spline element subdivision capabilities. In the frequency domain, adaptivity is performed in a mode-by-mode manner, sweeping from lower to higher frequencies, and  identifying the correspondence between coarse and fine mesh solutions using a MAC-type approach. As shown in the numerical section of the paper, super convergent solutions are obtained (i.e. faster convergence than that observed when using a uniform $h-$refinement), in particular when the problem exhibits local features that require a local refinement. We are currently working on extending the findings of this paper to the adaptivity of elasto-dynamic solutions in the time domain.
\section*{Acknowledgments}
Peng Yu thanks the support from China Scholarship Council (201406150089). The research leading to these results was partly supported by the European Research Council in the project ERC-COMBAT, grant agreement 615132. Satyendra Tomar would like to thank financial support from: FP7-PEOPLE-2011-ITN (289361), RealTCut(279578), INTER/FWO/15/10318764, and St\'ephane Bordas thanks partial funding for his time provided by the European Research Council Starting Independent Research Grant (ERC Stg grant agreement No.279578) RealTCut (Towards real time multiscale simulation of cutting in non-linear materials with applications to surgical simulation and computer guided surgery). We also thank the funding from the Luxembourg National Research Fund INTER/MOBILITY/14/8813215/CBM/Bordas, and INTER/FWO/15/10318764. Pierre Kerfriden acknowledges the financial support of the Engineering Research Network Wales. 

\appendix
\numberwithin{equation}{section}
\numberwithin{figure}{section}
\numberwithin{remark}{section}
\begin{appendices}
\section{NURBS basis functions}\label{NURBS app}
The univariate B-spline functions with $p$ order are defined in a recursive form over a knot vector $\xi$ using the following formula \cite{piegl2012nurbs}\\
for $p=0$
\begin{equation}
B_{i,0}(\xi) =
\begin{cases}
1 \hspace{1cm}\text{if} \ \xi_{i}\leq\xi<\xi_{i+1}\\
0 \hspace{1cm}\text{otherwise}
\end{cases},
\end{equation}
and $p\geqslant1$
\begin{align}
B_{i,p}(\xi) = \dfrac{\xi-\xi_{i}}{\xi_{i+p}-\xi_{i}}B_{i,p-1}(\xi)+\dfrac{\xi_{i+p+1}-\xi}{\xi_{i+p+1}-\xi_{i+1}}B_{i+1,p-1}(\xi).
\end{align}
Furthermore, the first derivative of B-spline is also given in a recursive way
\begin{align}
\dfrac{d}{d\xi}B_{i,p}(\xi) = \dfrac{p}{\xi_{i+p}-\xi_{i}}B_{i,p-1}(\xi)-\dfrac{p}{\xi_{i+p+1}-\xi_{i+1}}B_{i+1,p-1}(\xi).
\end{align}
Thus, univariate NURBS functions are defined as
\begin{align}
N_{i,p}(\xi) = \dfrac{B_{i,p}(\xi){w}_{i}}{W(\xi)} =  \dfrac{B_{i,p}(\xi)w_{i}}{\sum_{\hat{i}=1}^{n}B_{\hat{i},p}w_{\hat{i}}},
\end{align}
where $n$ is the number of basis functions and $w_{i}$ are a set of weights. Thereby, the first derivative of NURBS is defined by
\begin{align}
\dfrac{d}{d(\xi)}N_{i,p}(\xi) = w_{i}\dfrac{B_{i,p}'(\xi)W(\xi)-B_{i,p}(\xi)W'(\xi)}{W^{2}(\xi)},
\end{align}\\
where
\begin{align}
B'_{i,p}(\xi) = \dfrac{d}{d\xi}N_{i,p}(\xi), \ W'(\xi) = \sum_{\hat{i}=1}^{n}B'_{\hat{i},p}w_{\hat{i}}.
\end{align}
Given a 2D parametric space $ [0,1] \times [0,1]$, a tensor-product NURBS surface $\vect{S}_{\text{NURBS}}$ can be defined by \cite{piegl2012nurbs}
\begin{align}\label{NURBSsurface}
\vect{S}_{\text{NURBS}} = \sum_{i=1}^{n}\sum_{j=1}^{m}N_{i,j}^{p,q}(\xi,\eta)\vect{P}_{i,j} = \sum_{i=1}^{n}\sum_{j=1}^{m}\dfrac{B_{i,p}(\xi)B_{j,q}(\eta)w_{i,j}\vect{P}_{i,j}}{\sum_{\hat{i}=1}^{n}\sum_{\hat{j}=1}^{m}B_{\hat{i},p}(\xi)B_{\hat{j},q}(\eta)w_{\hat{i},\hat{j}}}, \hspace{0.5cm} (\xi,\eta)\in[0,1]\times[0,1].
\end{align}

\section{PHT-spline basis functions} \label{PHT app}
The PHT-spline with bi-cubic orders was proposed by Deng et al. in \cite{Deng2008}, and is developed recently to arbitrary degree by Anitescu et al. in \cite{anitescu2017recovery}. For brevity, only main properties of the bi-cubic PHT-spline basis functions applied in this paper and the refinement process are introduced.
\subsection{Construction of the PHT-spline basis function}
Given that all elements $\mathcal{T} = \bigcup\mathcal{T}_{e}$ are defined on a hierarchical T-mesh $\mathbb{T}$ on a parameterized domain $\mathcal{P}$. Then we can define a linear space for PHT-splines \cite{Deng2008}
\begin{align}
\mathscr{S}(p,q,\alpha,\beta,\mathbb{T}) = \left\lbrace T(\xi,\eta)\in C^{\alpha,\beta}(\mathcal{P})|T(\xi,\eta)\in\mathbb{P}_{p,q}, \forall\mathcal{T}_{e}\in\mathcal{T}\right\rbrace, 
\end{align}
where the space $ \mathbb{P}_{p,q} $ consists of all the bivariate polynomials with degree $p,q$, and $C^{\alpha,\beta}(\Omega)$ is the space involving all continuous bivariate spline functions with $C^{\alpha}$ in the $\xi$-direction and $C^{\beta}$ in the $\eta$-direction. The dimension equation of spline space $ \mathscr{S}(p,q,\alpha,\beta,\mathbb{T}) $ with $p\geqslant2\alpha+1$ and $q\geqslant2\beta+1$ is presented in \cite{deng2006dimensions}. Particularly, the bi-cubic PHT-spline space can be denoted that
\begin{align}
\text{dim}\mathscr{S}(3,3,1,1,\mathbb{T}) = 4(V^{b}+V^{+}).
\end{align}
Here $V^{b}$ and $V^{+}$ indicates the number of boundary vertices and interior crossing vertices separately. Supposed that the parametric domain $\mathcal{P} = [0,1]\times[0,1]$ is provided, the tensor-product PHT surface can be defined by
\begin{align}
\vect{S}_{\text{PHT}}=\sum_{i=1}^{n}\sum_{j=1}^{m}T_{i,j}^{p,q}(\xi,\eta)\vect{P}_{i,j}, \hspace{0.3cm} (\xi,\eta)\in\mathcal{P},
\end{align}
where $ T_{i,j}^{p,q} $ are constructed $C^{\alpha,\beta}$ continuous PHT-spline functions and $\vect{P}_{i,j}$ are control points. According to the literature \cite{anitescu2017recovery,Deng2008}, the $ T_{i,j}^{p,q} $ is generally computed through B\'ezier representation, which will be introduced subsequently. Let that $\hat{\vect{F}}$ represents linear mapping from a reference domain $\hat{\mathcal{P}}$ to parametric domain $\mathcal{P}$ that
\begin{align}
	\hat{F}: \hat{\mathcal{P}} \rightarrow \mathcal{P}, \hspace{0.3cm}   \hat{\vect{F}}(\hat{\xi},\hat{\eta}) = (\xi,\eta), \hspace{0.3cm} \hat{\mathcal{P}} = [-1,1]\times[-1,1],
\end{align} 
the basis function $ T_{i,j}^{p,q} $ can be rewritten in the form of a linear combination of Bernstein polynomials that
\begin{align}
	T_{i,j}^{p,q}(\xi,\eta) = \sum_{\hat{i}=1}^{p+1}\sum_{\hat{j}=1}^{q+1}b_{\hat{i},\hat{j}}\hat{B}_{\hat{i},\hat{j}}\circ\hat{\vect{F}}^{-1}(\xi,\eta),
\end{align}
where $ \hat{B}_{\hat{i},\hat{j}}(\hat{\xi},\hat{\eta}) = \hat{B}_{\hat{i}}(\hat{\xi})\hat{B}_{\hat{j}}(\hat{\eta})$ are the tensor product of univariate Bernstein functions, which are defined on $\hat{\mathcal{P}}$ as follows
\begin{align}
\hat{B}_{\hat{i}}(\hat{\xi}) = \dfrac{1}{2^{p}}
\begin{pmatrix}
p \\ i-1
\end{pmatrix}
(1-\xi)^{p-i+1}(1+\xi)^{i-1}, \hspace{0.5cm} i = 1,\ 2,\ldots,\ p+1.
\end{align}
The $ b_{\hat{i},\hat{j}} $ are B\'ezier ordinates obtained by a recursive method called De Casteljau’s algorithm, and readers can find the details in \cite{anitescu2017recovery,Nguyen-Thanh2011}.

\subsection{Tree structure and local refinement}\label{app:PHTrefine}
\begin{figure}
	\centering
	\includegraphics[width=1\columnwidth]{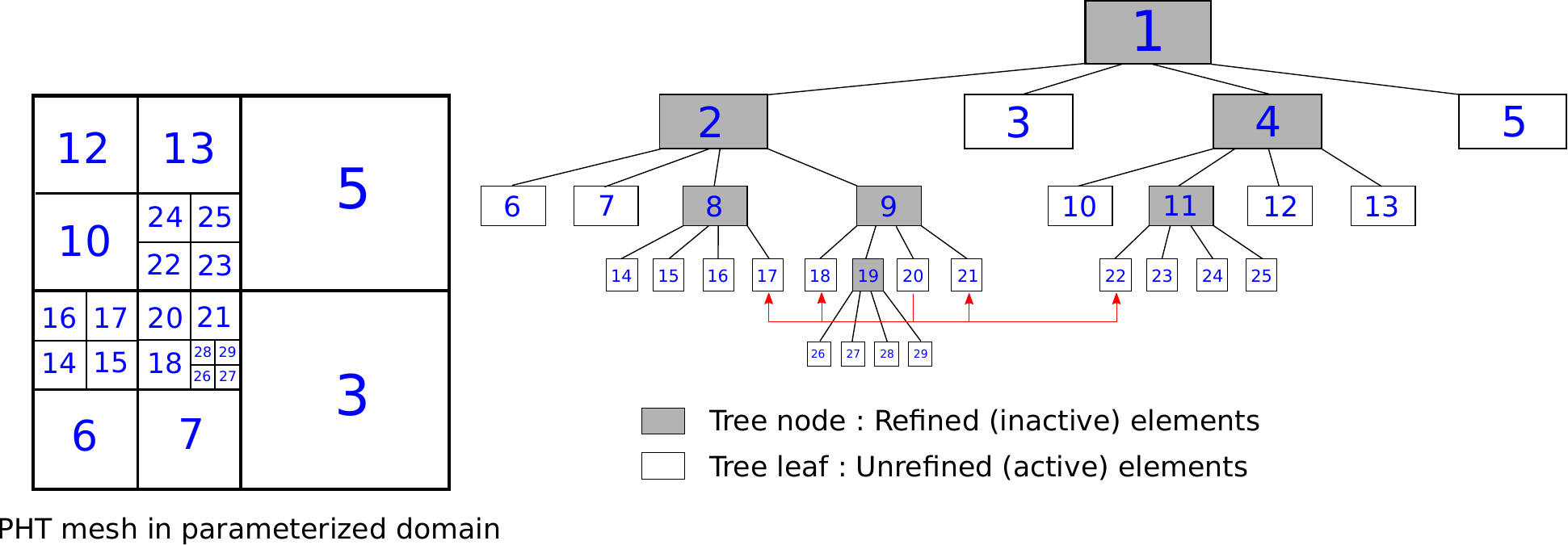}
	\caption{A typical quadtree system to represent data structure of PHT adaptive mesh. The numbering of elements is executed from coarse to fine level. The numbers in parametric space (on the left) just denote the elements without any children. With the help of quadtree structure, the relationship of all elements is readily observed. For instance, element 1 is the parent of cell 4, and element 4 is inherited by children cells 10, 11, 12, 13. The adjoint cells of element 20 are elements 17,18,21,22.}
	\label{octree}
\end{figure}

\begin{figure}
	\centering
	\begin{subfigure}[t]{0.25\textwidth}
		\centering
		\includegraphics[width=1\columnwidth]{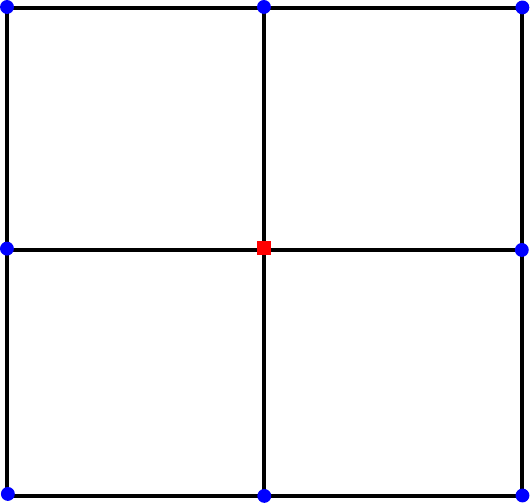}
		\caption{}
	\end{subfigure}
	\hspace{0.5cm}
	\begin{subfigure}[t]{0.25\textwidth}
		\centering
		\includegraphics[width=1\columnwidth]{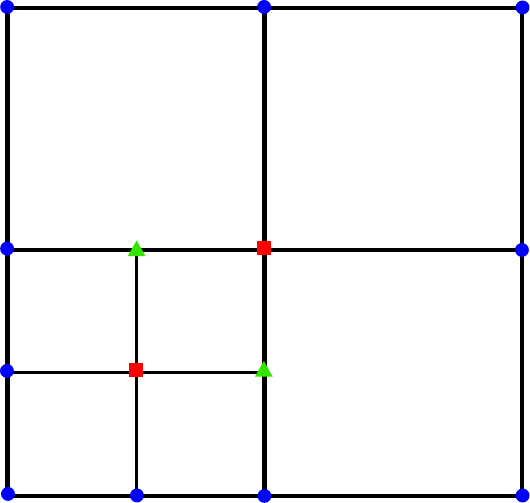}
		\caption{}
	\end{subfigure}
	\hspace{0.5cm}
	\begin{subfigure}[t]{0.25\textwidth}
		\centering
		\includegraphics[width=1\columnwidth]{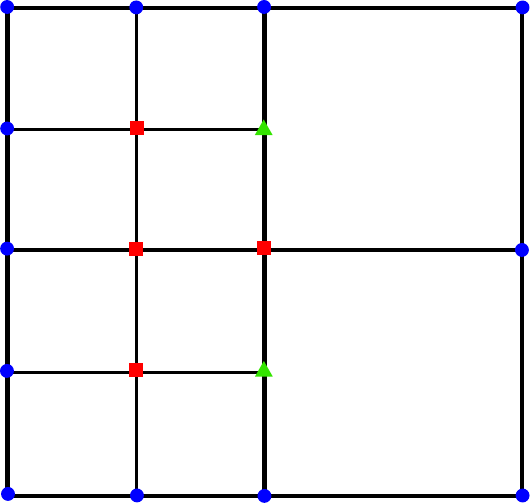}
		\caption{}
	\end{subfigure}
	\caption{An example to present vertices insertion at local refinement on a PHT mesh. Particularly, the blue dots denote boundary vertices, and the green triangles express T-junctions, which are generated when refinements are created between adjacent elements with different levels. It should be noted that T-junctions do not change basis functions until they are transfered to crossing vertices, which are indicated by red squares. A crossing vertex will lead to the truncation over $(\alpha+1)(\beta+1)$ B\'ezier ordinates around, which are set as zeros firstly, and then reset by new values on the updated spline space. (a) The initial mesh, (b) the local refinement in one cell, (c) the local refinement in the adjoint cell.}
	\label{vertexinsert}
\end{figure}

The data structure of 2D hierarchical T-mesh $\mathbb{T}$ for PHT-splines is stored within the quadtree framework, as shown in Fig.\ref{octree}. Every leaf or node of the tree represents one element $\mathcal{T}_{e}$ at different refinement levels, which reserves all the mandatory information with respect to PHT-spline basis functions, i.e., B\'ezier ordinates, numbering system of nodes and elements, refinement level, etc.. Furthermore, each leaf or node also preserves hierarchical connectivities applied to trace parent and children elements during the refinement, and adjacent connectivities which combine neighboring elements by pointers. The process of typical vertices insertions during refinements is illustrated in Fig.\ref{vertexinsert}. More details regarding the principle of the algorithm and implementation can be seen in \cite{schillinger2012isogeometric}.

\begin{figure}
	\centering
	\includegraphics[width=1\columnwidth]{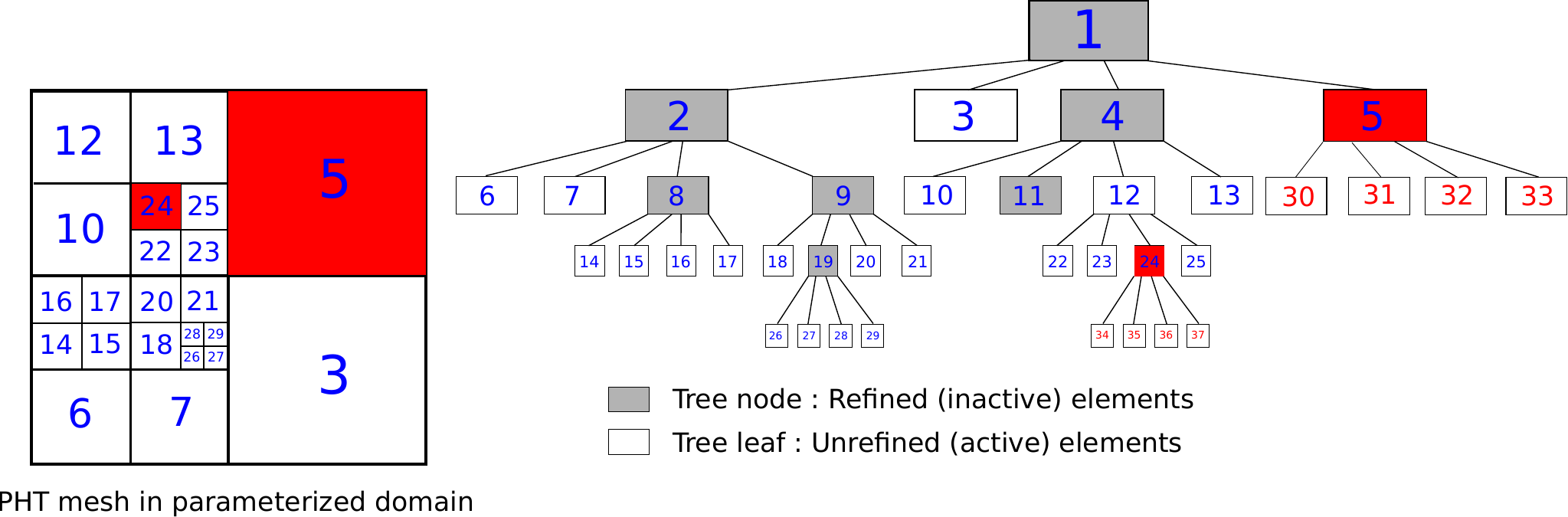}
	\caption{The direct refinement for PHT mesh in case of level $\mathcal{K} \geqslant \mathcal{K'}-1$. The cells 5 and 24 are marked and then refined. Note that refinements are always proceeded from coarse level to fine level in implementation regardless of the order of marking. So the numbering of children elements of cell 5 is smaller than those of cell 24.}
	\label{octreerefine}
\end{figure}

\begin{figure}
	\centering
	\includegraphics[width=1\columnwidth]{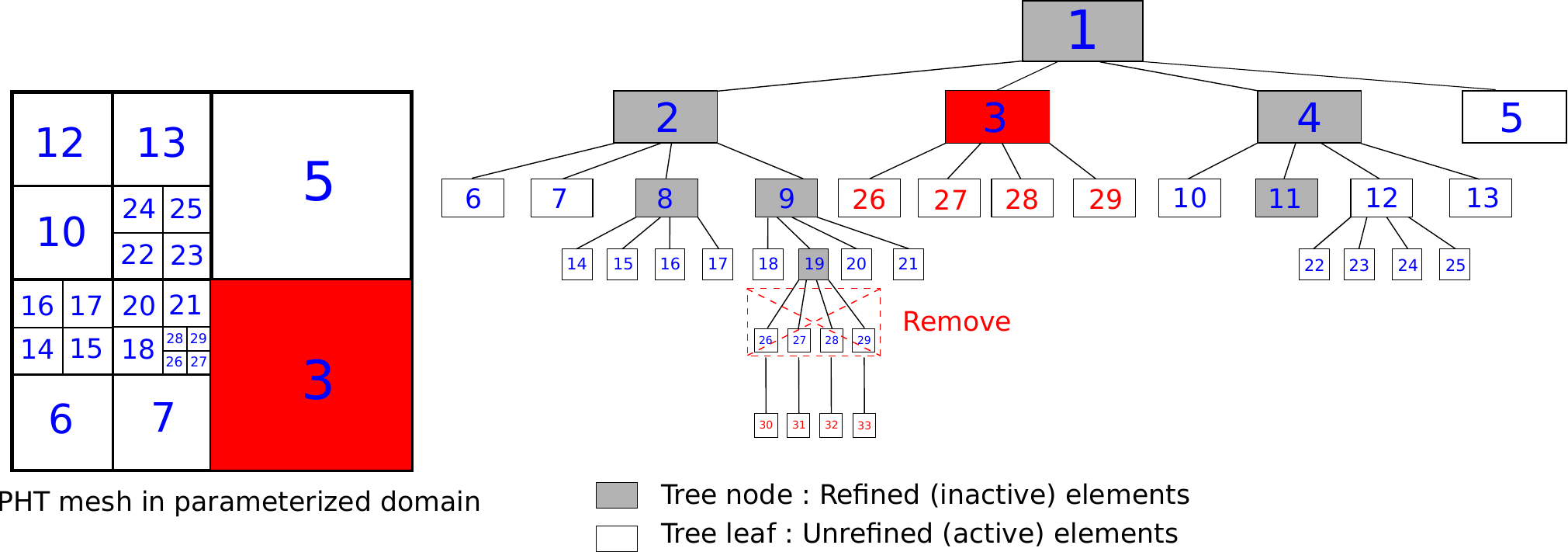}
	\caption{The refinement rule for PHT mesh in terms of level $\mathcal{K} <  \mathcal{K'}-1$. Assuming that cell 3 is considered to be refined, we have to remove the adjacent elements (26, 27, 28, 29) with level $\mathcal{K'}$ temporarily to ensure the updated maximum level $\mathcal{K''}$ of neighboring elements satisfies that $\mathcal{K} \geqslant \mathcal{K''}-1$. Afterwards, refine the cell 3 to obtain the 4 children elements with numbering (26,27,28,29). Ultimately, refine cell 19 again to acquire new children elements (30,31,32,33) as the replacement of removal elements (26,27,28,29) before.}
	\label{octreerefineremove}
\end{figure}

Before this tree system is exploited for adaptive refinement, we recall the issue of refinement procedure on PHT mesh affected by the level between the target element and its adjoint elements. As stated in \cite{Deng2008}, assumed that the level of the element $\mathcal{T}_{e}$ is $\mathcal{K}$ and maximum level of the neighboring elements is  $\mathcal{K'}$. If $\mathcal{K} \geqslant \mathcal{K'}-1$, the refinement would be straightforward as illustrated in Fig.\ref{octreerefine}. Whilst when $\mathcal{K} < \mathcal{K'}-1$, the situation would be a bit more complex but still in a good control with the application of the quadtree configuration, exhibited in Fig.\ref{octreerefineremove}.

\section{RHT-spline basis functions} \label{app:RHT}
The RHT-spline basis functions defined over 2D hierarchical T-mesh $\mathbb{T}$ are computed by \cite{Nguyen-Thanh2011}
\begin{align}
R_{i,j}^{p,q} = \dfrac{T_{i,j}^{p,q}(\xi,\eta)\hat{w}_{i,j}}{\sum_{\hat{i}=1}^{n}\sum_{\hat{j}=1}^{m}T_{\hat{i},\hat{j}}^{p,q}(\xi,\eta)\hat{w}_{\hat{i},\hat{j}}}, \hspace{0.5cm} (\xi,\eta)\in\mathcal{P}=[0,1]\times[0,1],
\end{align}
where $T_{i,j}^{p,q}(\xi,\eta)$ are PHT-spline basis functions introduced as in Appendix \ref{PHT app}, and $w_{i,j}$ are weights. Accordingly, a RHT-spline surface at level $k$ mesh is given by
\begin{align}
\vect{S}_{\text{RHT}}^{k} =\sum_{\vect{I}}\vect{R}_{I}^{k}(\xi,\eta)\vect{P}_{I}^{k},
\end{align}
where $\vect{I}$ is the multi-index, $\vect{R}_{I}^{k}(\xi,\eta)$ and $\vect{P}_{I}^{k}$ are RHT-spline basis functions and control points at level $k$ mesh. As mentioned previously, PHT-splines are utilized in the framework of GIFT where there is no need to compute control points as the geometry is characterized by NURBS. The mesh of geometry will stay at the initial stage during the refinement. However, when RHT-splines are employed in IGA scheme, it is necessary to renew the control points, as well as weights, at each refinement step. 
\subsection{Update of control points and weights with h-refinement}
Following the approach proposed by Deng et al.\cite{Deng2008} to get the new control points after refinement for the cubic PHT-spline surface, we are going to develop it for a cubic RHT-spline surface. Assuming that the control points $\vect{P}_{\vect{I}}^{k} = (x_{\vect{I}},y_{\vect{I}},z_{\vect{I}})$ are depicted over a 3D Euclidean space, we introduce the \emph{homogeneous coordinates} \cite{piegl2012nurbs} to denote $\vect{P}_{\vect{I}}^{k}$ in a four-dimensional space as follows 
\begin{align} \label{HCPtsk}
\vect{P}_{\vect{I}}^{k} = \mathbb{H}(\vect{P}_{\vect{I}}^{k(w)}) = \mathbb{H}[(\hat{w}_{\vect{I}}x_{\vect{I}}, \hat{w}_{\vect{I}}y_{\vect{I}} , \hat{w}_{\vect{I}}z_{\vect{I}}, \hat{w}_{\vect{I}})] =
\begin{cases}
(\dfrac{\hat{w}_{\vect{I}}x_{\vect{I}}}{\hat{w}_{\vect{I}}}, \dfrac{\hat{w}_{\vect{I}}y_{\vect{I}}}{\hat{w}_{\vect{I}}}, \dfrac{\hat{w}_{\vect{I}}z_{\vect{I}}}{\hat{w}_{\vect{I}}}), \hspace{0.3cm} &\text{with} \  \hat{w}_{\vect{I}} \neq 0 ,\\
x_{\vect{I}}, y_{\vect{I}}, z_{\vect{I}}, \hspace{0.3cm} &\text{with} \  \hat{w}_{\vect{I}} = 0, 
\end{cases} 
\end{align}
where $\vect{P}_{\vect{I}}^{k(w)}$ are the weighted control points on a 4D space and $\mathbb{H}$ is the mapping function. Also, we can create this analogous relationship for the surface such that $\vect{S}_{\text{RHT}}^{k} = \mathbb{H}(\vect{S}_{\text{PHT}}^{k(w)}) = \mathbb{H}(W\vect{S}_{\text{PHT}}^{k}, W)$ with
\begin{align}
\vect{S}_{\text{PHT}}^{k(w)} = \sum_{\vect{I}}\vect{T}_{I}^{k}(\xi,\eta)\vect{P}_{I}^{k(w)}, \ W = \sum_{\vect{I}}\vect{T}_{I}^{k}(\xi,\eta)\hat{w}_{I}^{k},
\end{align}
where $ \vect{S}_{\text{PHT}}^{k(w)} $ is the weighted PHT-spline surface 
at level $k$, and $W$ are weights for the surface. Note that the basis functions to represent $ \vect{S}_{\text{PHT}}^{k(w)} $ are PHT-spline basis functions so that we can directly call the algorithm \cite{Deng2008} to generate new control points and weights. To be specific, for instance, when a new vertex is inserted into an element belonging to cells $\mathcal{T}_{k}$ at level $k$, there will be $(\alpha+1)(\beta+1)$ new added basis functions $\vect{T}_{\vect{J}}^{k+1}$, weighted control points $\Delta\vect{P}_{\vect{J}}^{k+1(w)}$ created at $\mathcal{T}_{k+1}$ such that
\begin{align}
\vect{S}_{\text{PHT}}^{k+1(w)}(\xi,\eta) = \sum_{\vect{I}}^{N}\vect{\tilde{T}}_{\vect{I}}^{k+1}(\xi,\eta)\underline{\underline{\vect{P}}}_{\vect{I}}^{k+1(w)} + \sum_{\vect{J} = N +1 }^{(\alpha+1)(\beta+1)}\vect{T}_{\vect{J}}^{k+1}(\xi,\eta)\Delta\vect{P}_{\vect{J}}^{k+1(w)},
\end{align}
where $\vect{\tilde{T}}_{I}^{k+1}$ are the basis functions $\vect{T}_{I}^{k}$ represented on the $\mathcal{T}_{k+1}$ so that the relevant control points are kept as $\underline{\underline{\vect{P}}}_{\vect{I}}^{k+1(w)} = \vect{P}_{I}^{k(w)}$. Now we only consider the parametric space (by setting it as $(\xi^{*}, \eta^{*})$), dominated by the new basis. Owing to the truncation property \cite{Deng2008,anitescu2017recovery}, the basis functions $\vect{\tilde{T}}_{I}^{k+1}$ and the derivatives will vanish in this domain. Simultaneously, due to the geometry preservation during the refinement, it leads to
\begin{align}\label{WSpht}
\vect{S}_{\text{PHT}}^{k+1(w)}(\xi^{*},\eta^{*}) = \sum_{\vect{J} = N +1 }^{(\alpha+1)(\beta+1)}\vect{T}_{\vect{J}}^{k+1}(\xi^{*},\eta^{*})\Delta\vect{P}_{\vect{J}}^{k+1(w)} = \vect{S}_{\text{PHT}}^{k(w)}(\xi^{*},\eta^{*}).
\end{align} 
In order to compute $ \Delta\vect{P}_{\vect{J}}^{k+1(w)} $, we define a linear operator involving the \emph{geometric information} for the surface $\vect{S}_{\text{PHT}}^{k(w)}(\xi^{*},\eta^{*})$ such that 
\begin{align}
\mathscr G\vect{S}_{\text{PHT}}^{k(w)}(\xi^{*},\eta^{*}) = \left( \vect{S}_{\text{PHT}}^{k(w)}(\xi^{*},\eta^{*}),\  \dfrac{\partial\vect{S}_{\text{PHT}}^{k(w)}}{\partial\xi^{*}}, \ \dfrac{\partial\vect{S}_{\text{PHT}}^{k(w)}}{\partial\eta^{*}}, \ \dfrac{\partial^{2}\vect{S}_{\text{PHT}}^{k(w)}}{\partial\xi^{*}\partial\eta^{*}}\right), 
\end{align}
and then rewrite Eq.\eqref{WSpht} as
\begin{align}
\mathscr G\vect{S}_{\text{PHT}}^{k(w)}(\xi^{*},\eta^{*}) = \sum_{\vect{J} }\mathscr G\vect{T}_{\vect{J}}^{k+1}(\xi^{*},\eta^{*})\Delta\vect{P}_{\vect{J}}^{k+1(w)} = \underline{\underline{\mathbb{T}}}\cdot\Delta\vect{P}^{k+1(w)}.
\end{align}
Here, in case of cubic basis functions, the matrix $\underline{\underline{\mathbb{T}}}$ can be simply represented via the distance between the inserted vertex and the adjacent elements as follows
\begin{align}
\underline{\underline{\mathbb{T}}}(\Delta u_{1}, \Delta u_{2}, \Delta v_{1}, \Delta v_{2}) =
\begin{pmatrix}	
(1-\bar{\lambda})(1-\mu)& \hspace{0.5cm} \bar{\lambda}(1-\mu)& \hspace{0.5cm} (1-\bar{\lambda})\mu& \hspace{0.5cm} \bar{\lambda}\mu\\
-\alpha(1-\mu)& \hspace{0.5cm} \alpha(1-\mu)& \hspace{0.5cm} -\alpha\mu& \hspace{0.5cm} \alpha\mu\\
-\beta(1-\bar{\lambda})& \hspace{0.5cm} -\beta\bar{\lambda}& \hspace{0.5cm} \beta(1-\bar{\lambda})& \hspace{0.5cm} \beta\bar{\lambda}\\
\alpha\beta& \hspace{0.5cm} -\alpha\beta& \hspace{0.5cm} -\alpha\beta& \hspace{0.5cm} \alpha\beta
\end{pmatrix}, 
\end{align}
where $\alpha = \dfrac{1}{\Delta u_{1}+\Delta u_{2}}, \ \beta = \dfrac{1}{\Delta v_{1}+\Delta v_{2}}, \ \bar{\lambda} = \alpha\Delta u_{1}, \ \mu = \beta\Delta v_{1}$ can be found in \cite{Deng2008}. Thus, $ \Delta\vect{P}_{\vect{J}}^{k+1(w)} $ are obtained by
\begin{align}
\Delta\vect{P}^{k+1(w)} = (\underline{\underline{\mathbb{T}}})^{-1}\cdot\mathscr G\vect{S}_{\text{PHT}}^{k(w)},
\end{align}
and furthermore the new weighted control points are obtained by
\begin{align}
\vect{P}^{k+1(w)} = \underline{\underline{\vect{P}}}_{\vect{I}}^{k+1(w)}+\Delta\vect{P}^{k+1(w)}.
\end{align}
Consequently, the new control points at level $k+1$ mesh are updated
\begin{align}\label{HCPtsk+1}
\vect{P}^{k+1(w)} = \mathbb{H}(\vect{P}^{k+1(w)}).
\end{align}
Note that since PHT-splines are the generation of B-splines on the hierarchical T-mesh, RHT-splines at level 0 mesh are exactly NURBS. Hence, without any refinement at initial stage, a RHT surface is identical to a NURBS surface, that is, $\vect{S}_{\text{RHT}}^{0} = \vect{S}_{\text{NURBS}}^{0}$, where $\vect{S}_{\text{NURBS}}^{0}$ is defined in Eq.\eqref{NURBSsurface}, with $\vect{R}_{\vect{I}}^{0} = \vect{N}_{\vect{I}}, \vect{P}_{\vect{I}}^{0} = \vect{P}_{\vect{I}}, \hat{w}_{\vect{I}}^{0} = w_{\vect{I}}$. Afterwards, following the deduction from Eq.\eqref{HCPtsk} to Eq.\eqref{HCPtsk+1}, the added control points and weights for RHT-spline basis functions can be obtained during the hierarchical refinement.

\section{Prolongation of control variables from a coarse mesh to a refined mesh} \label{app:prolong}
Let $\vect{\phi}^{h}$ be the solution on a \emph{coarse mesh} $\mathbb{T}$, which is discretized by PHT-spline basis functions such that $ \vect{\phi}^{h} = \vect{T}\bar{\vect{\phi}}^{h} $, where $\bar{\vect{\phi}}^{h}$ are control variables. Now we aim to find the prolongation of $\bar{\vect{\phi}}^{h}$, namely, $\mathbb{P}\bar{\vect{\phi}}^{h}$, on a \emph{refined mesh} $\tilde{\mathbb{T}}$. Two methods are presented as follows.
\subsection{Update control variables in hierarchical refinement}\label{prolongation}
As the coarse mesh $\mathbb{T}$ is nested in refined mesh $\tilde{\mathbb{T}}$, namely, $\mathbb{T}\subset\tilde{\mathbb{T}}$, by recalling the algorithm of constructing new control points as from Appendix.\ref{app:RHT}, the generation of projection can be similarly interpreted as the creation of new control points. To be specific, suppose that $\Delta\bar{\vect{\phi}}^{h}$ serves as the added control variables resulting from the increase of the degrees of freedom by refinement. According to the method discussed in appendix.\ref{app:RHT}, it yields
\begin{align}
\Delta\bar{\vect{\phi}}^{h} = (\underline{\underline{\mathbb{T}}})^{-1}\mathscr G\vect{\phi}^{h}.
\end{align}
then the prolongation reads
\begin{align}
\mathbb{P}\bar{\vect{\phi}}^{h}
 = \bar{\vect{\phi}}^{h} + \Delta\bar{\vect{\phi}}^{h}.
\end{align}
\subsection{Projection}\label{Mofnorm}
Assuming that $\tilde{\vect{\phi}}^{h} = \tilde{\vect{T}}\mathbb{P}\bar{\vect{\phi}}^{h}$ is the projection of $\vect{\phi}^{h}$ from $\mathbb{T}$ onto $\tilde{\mathbb{T}}$, then the prolongation problem can also be understood that \vspace*{0.2cm}\\
\fbox{
\parbox{\textwidth}{
Find $\mathbb{P}\bar{\vect{\phi}}^{h}$ such that
\begin{align}
\underset{\mathbb{P}\bar{\vect{\phi}}^{h}}{\text{max}}\left\|\vect{\phi}^{h} - \tilde{\vect{T}}\mathbb{P}\bar{\vect{\phi}}^{h}\right\|_{m},
\end{align}
namely,
\begin{align}\label{dnorm}
\delta_{\mathbb{P}\bar{\vect{\phi}}^{h}}\left\|\vect{\phi}^{h} - \tilde{\vect{T}}\mathbb{P}\bar{\vect{\phi}}^{h}\right\|_{m}^{2}= 0,
\end{align} 
}
}\vspace*{0.2cm}
where the mass norm $\left\|\cdot\right\|_{m}$ is defined in Eq.\eqref{massnorm}.
The term $\left\| \vect{\phi}^{h} - \tilde{\vect{T}}\mathbb{P}\bar{\vect{\phi}}^{h}\right\|_{m}^{2} $ can be extended to
\begin{align}
\left\| \vect{\phi}^{h} - \tilde{\vect{T}}\mathbb{P}\bar{\vect{\phi}}^{h}\right\|_{m}^{2} = \int_{\Omega}(\vect{T}\bar{\vect{\phi}}^{h})^{T}\matr{m}\vect{T}\bar{\vect{\phi}}^{h}d\Omega - 2\int_{\Omega}(\tilde{\vect{T}}\mathbb{P}\bar{\vect{\phi}}^{h})^{T}\matr{m}\vect{T}\bar{\vect{\phi}}^{h}d\Omega + \int_{\Omega}(\tilde{\vect{T}}\mathbb{P}\bar{\vect{\phi}}^{h})^{T}\matr{m}\tilde{\vect{T}}\mathbb{P}\bar{\vect{\phi}}^{h}d\Omega.
\end{align}
Then Eq.\eqref{dnorm} is written by
\begin{align}
\delta_{\mathbb{P}\bar{\vect{\phi}}^{h}}\left\|\vect{\phi}^{h} - \tilde{\vect{T}}\mathbb{P}\bar{\vect{\phi}}^{h}\right\|_{m}^{2}=\int_{\Omega}\tilde{\vect{T}}^{T}\matr{m}\tilde{\vect{T}}\mathbb{P}\bar{\vect{\phi}}^{h}d\Omega - \int_{\Omega}\tilde{\vect{T}}^{T}\matr{m}\vect{T}\bar{\vect{\phi}}^{h}d\Omega =0.
\end{align}
Define mass matrices $\matr{M}_{\tilde{\vect{T}},\vect{T}} = \int_{\Omega}\tilde{\vect{T}}^{T}\matr{m}\vect{T}d\Omega, ~\matr{M}_{\tilde{\vect{T}}, \tilde{\vect{T}}} = \int_{\Omega}\tilde{\vect{T}}^{T}\matr{m}\tilde{\vect{T}}d\Omega$, and we obtain the prolongation as
\begin{align}
\mathbb{P}\bar{\vect{\phi}}^{h} = \matr{M}_{\tilde{\vect{T}}, \tilde{\vect{T}}}^{-1}\matr{M}_{\tilde{\vect{T}},\vect{T}}\bar{\vect{\phi}}^{h}.
\end{align}
\begin{remark}
When computing $\matr{M}_{\tilde{\vect{T}},\vect{T}} = \int_{\Omega}\tilde{\vect{T}}^{T}(\vect{\xi})\matr{m}\vect{T}(\vect{\xi})d\Omega$, suppose that the integration is proceeded in $\tilde{\mathbb{T}}$. Since the mappings $\vect{x} = \vect{F}(\vect{\xi})$ and $\vect{x} = \tilde{\vect{F}}(\vect{\xi})$ are for $\mathbb{T}$ and $\tilde{\mathbb{T}}$ respectively, actually the term $\vect{T}(\vect{\xi})$ has to be calculated through
\begin{align}
\vect{T}(\xi) = \vect{T}\circ\tilde{\vect{F}}^{-1}[\vect{F}(\vect{\xi})].
\end{align}
Owing to the geometric preservation by isogeometric method during the refinement, it yields $\tilde{\vect{F}}^{-1}[\vect{F}(\vect{\xi})] = \vect{\xi} $, which indicates that $ \vect{T}(\xi) $ can be computed directly.
\end{remark}
\begin{remark}
It is worth noting that the presented prolongation method in Appendix \ref{prolongation} is restricted to cubic PHT-spline basis functions, though, it is computationally cheap. In contrast, because of unavoidable calculation for the inverse matrix $ \matr{M}_{\tilde{\vect{T}}, \tilde{\vect{T}}}^{-1} $, the strategy in Appendix \ref{Mofnorm} is more expensive, nevertheless, it is more widely used for any arbitrary degree of PHT-splines or other basis functions.
\end{remark}	  
\end{appendices}

\section*{References}
\bibliography{mybibfile}

\begin{thebibliography}{10}
\expandafter\ifx\csname url\endcsname\relax
  \def\url#1{\texttt{#1}}\fi
\expandafter\ifx\csname urlprefix\endcsname\relax\def\urlprefix{URL }\fi
\expandafter\ifx\csname href\endcsname\relax
  \def\href#1#2{#2} \def\path#1{#1}\fi

\bibitem{hughes2005}
T.~J. Hughes, J.~A. Cottrell, Y.~Bazilevs, Isogeometric analysis: {CAD}, finite
  elements, nurbs, exact geometry and mesh refinement, Computer Methods in
  Applied Mechanics and Engineering 194~(39) (2005) 4135--4195.

\bibitem{Nguyen201589}
V.~P. Nguyen, C.~Anitescu, S.~P. Bordas, T.~Rabczuk,
  \href{http://www.sciencedirect.com/science/article/pii/S0378475415001214}{Isogeometric
  analysis: An overview and computer implementation aspects}, Mathematics and
  Computers in Simulation 117 (2015) 89 -- 116.

\bibitem{Cottrell20065257}
J.~Cottrell, A.~Reali, Y.~Bazilevs, T.~Hughes, Isogeometric analysis of
  structural vibrations, Computer Methods in Applied Mechanics and Engineering
  195~(41–43) (2006) 5257 -- 5296, {J}ohn H. Argyris Memorial Issue. Part
  \{II\}.

\bibitem{Shojaee2012}
S.~Shojaee, E.~Izadpanah, N.~Valizadeh, J.~Kiendl, {Free vibration analysis of
  thin plates by using a NURBS-based isogeometric approach}, Finite Elements in
  Analysis and Design 61 (2012) 23--34.

\bibitem{sobota2017implicit}
P.~Sobota, W.~Dornisch, R.~M{\"u}ller, S.~Klinkel, Implicit dynamic analysis
  using an isogeometric reissner--mindlin shell formulation, International
  Journal for Numerical Methods in Engineering 110~(9) (2017) 803--825.

\bibitem{NME:NME4282}
C.~H. Thai, H.~Nguyen-Xuan, N.~Nguyen-Thanh, T.-H. Le, T.~Nguyen-Thoi,
  T.~Rabczuk, \href{http://dx.doi.org/10.1002/nme.4282}{Static, free vibration,
  and buckling analysis of laminated composite reissner–mindlin plates using
  {NURBS}-based isogeometric approach}, International Journal for Numerical
  Methods in Engineering 91~(6) (2012) 571--603.

\bibitem{GiannelliJS-12-THBSplines}
C.~Giannelli, B.~J{\"u}ttler, H.~Speleers,
  \href{http://www.sciencedirect.com/science/article/pii/S0167839612000519}{{THB-splines:
  The truncated basis for hierarchical splines}}, Computer Aided Geometric
  Design 29~(7) (2012) 485 -- 498.

\bibitem{giannelli2016thb}
C.~Giannelli, B.~J{\"u}ttler, S.~K. Kleiss, A.~Mantzaflaris, B.~Simeon,
  J.~{\v{S}}peh, Thb-splines: An effective mathematical technology for adaptive
  refinement in geometric design and isogeometric analysis, Computer Methods in
  Applied Mechanics and Engineering 299 (2016) 337--365.

\bibitem{schillinger2012isogeometric}
D.~Schillinger, L.~Dede, M.~A. Scott, J.~A. Evans, M.~J. Borden, E.~Rank, T.~J.
  Hughes, An isogeometric design-through-analysis methodology based on adaptive
  hierarchical refinement of {NURBS}, immersed boundary methods, and {T}-spline
  {CAD} surfaces, Computer Methods in Applied Mechanics and Engineering 249
  (2012) 116--150.

\bibitem{johannessen2014isogeometric}
K.~A. Johannessen, T.~Kvamsdal, T.~Dokken, Isogeometric analysis using {LR}
  {B}-splines, Computer Methods in Applied Mechanics and Engineering 269 (2014)
  471--514.

\bibitem{Sederberg2004}
T.~W. Sederberg, D.~L. Cardon, G.~T. Finnigan, N.~S. North, J.~Zheng, T.~Lyche,
  \href{http://dl.acm.org/citation.cfm?id=1015706.1015715}{{T-spline
  simplification and local refinement}}, ACM Transactions on Graphics 23~(3)
  (2004) 276.

\bibitem{Bazilevs2010}
Y.~Bazilevs, V.~Calo, J.~Cottrell, J.~Evans, T.~Hughes, S.~Lipton, M.~Scott,
  T.~Sederberg,
  \href{http://www.sciencedirect.com/science/article/pii/S0045782509000875}{{Isogeometric
  analysis using T-splines}}, Computer Methods in Applied Mechanics and
  Engineering 199~(5-8) (2010) 229--263.

\bibitem{Deng2008}
J.~Deng, F.~Chen, X.~Li, C.~Hu, W.~Tong, Z.~Yang, Y.~Feng,
  \href{http://www.sciencedirect.com/science/article/pii/S1524070308000039}{{Polynomial
  splines over hierarchical T-meshes}}, Graphical Models 70~(4) (2008) 76--86.

\bibitem{Nguyen-Thanh2011}
N.~Nguyen-Thanh, H.~Nguyen-Xuan, S.~Bordas, T.~Rabczuk,
  \href{http://www.sciencedirect.com/science/article/pii/S0045782511000338}{{Isogeometric
  analysis using polynomial splines over hierarchical T-meshes for
  two-dimensional elastic solids}}, Computer Methods in Applied Mechanics and
  Engineering 200~(21-22) (2011) 1892--1908.

\bibitem{marussig2015fast}
B.~Marussig, J.~Zechner, G.~Beer, T.-P. Fries, Fast isogeometric boundary
  element method based on independent field approximation, Computer Methods in
  Applied Mechanics and Engineering 284 (2015) 458--488.

\bibitem{toshniwal2017smooth}
D.~Toshniwal, H.~Speleers, T.~J. Hughes, Smooth cubic spline spaces on
  unstructured quadrilateral meshes with particular emphasis on extraordinary
  points: Geometric design and isogeometric analysis considerations, Computer
  Methods in Applied Mechanics and Engineering 327 (2017) 411--458.

\bibitem{anitescu2017recovery}
C.~Anitescu, M.~N. Hossain, T.~Rabczuk, Recovery-based error estimation and
  adaptivity using high-order splines over hierarchical t-meshes, Computer
  Methods in Applied Mechanics and Engineering.

\bibitem{atroshchenko2017weakening}
E.~Atroshchenko, S.~Tomar, G.~Xu, S.~Bordas, Weakening the tight coupling
  between geometry and simulation in isogeometric analysis: from sub-and
  super-geometric analysis to {G}eometry {I}ndependent {F}ield approxima{T}ion
  ({GIFT}), International Journal for Numerical Methods in Engineering~(2017),
  Accepted.

\bibitem{nguyen2014adaptive}
N.~Nguyen-Thanh, J.~Muthu, X.~Zhuang, T.~Rabczuk, An adaptive three-dimensional
  rht-splines formulation in linear elasto-statics and elasto-dynamics,
  Computational Mechanics 53~(2) (2014) 369--385.

\bibitem{stein2007finite}
E.~Stein, M.~R{\"u}ter, Finite element methods for elasticity with
  error-controlled discretization and model adaptivity, Encyclopedia of
  Computational Mechanics.

\bibitem{ladeveze2013new}
P.~Ladev{\`e}ze, F.~Pled, L.~Chamoin, New bounding techniques for goal-oriented
  error estimation applied to linear problems, International journal for
  numerical methods in engineering 93~(13) (2013) 1345--1380.

\bibitem{bangerth2010adaptive}
W.~Bangerth, M.~Geiger, R.~Rannacher, Adaptive galerkin finite element methods
  for the wave equation, Computational Methods in Applied Mathematics Comput.
  Methods Appl. Math. 10~(1) (2010) 3--48.

\bibitem{gonzalez2014mesh}
O.~A. Gonz{\'a}lez-Estrada, E.~Nadal, J.~R{\'o}denas, P.~Kerfriden, S.~P.-A.
  Bordas, F.~Fuenmayor, Mesh adaptivity driven by goal-oriented locally
  equilibrated superconvergent patch recovery, Computational Mechanics 53~(5)
  (2014) 957--976.

\bibitem{allemang2003modal}
R.~J. Allemang, The {M}odal {A}ssurance {C}riterion--twenty years of use and
  abuse, {S}ound and {V}ibration 37~(8) (2003) 14--23.

\bibitem{PASTOR2012543}
M.~Pastor, M.~Binda, T.~Harčarik,
  \href{http://www.sciencedirect.com/science/article/pii/S1877705812046140}{Modal
  {A}ssurance {C}riterion}, Procedia Engineering 48 (2012) 543 -- 548.

\bibitem{zienkiewicz2000finite}
O.~C. Zienkiewicz, R.~L. Taylor, The Finite Element Method: Solid Mechanics,
  Vol.~2, Butterworth-heinemann, 2000.

\bibitem{liu2006meshfree}
Y.~Liu, Y.~Hon, K.~Liew, A meshfree hermite-type radial point interpolation
  method for {K}irchhoff plate problems, International Journal for Numerical
  Methods in Engineering 66~(7) (2006) 1153--1178.

\bibitem{fernandez2004imposing}
S.~Fern{\'a}ndez-M{\'e}ndez, A.~Huerta, Imposing essential boundary conditions
  in mesh-free methods, Computer Methods in Applied Mechanics and Engineering
  193~(12) (2004) 1257--1275.

\bibitem{nguyen2017isogeometric}
N.~Nguyen-Thanh, K.~Zhou, X.~Zhuang, P.~Areias, H.~Nguyen-Xuan, Y.~Bazilevs,
  T.~Rabczuk, Isogeometric analysis of large-deformation thin shells using
  {RHT}-splines for multiple-patch coupling, Computer Methods in Applied
  Mechanics and Engineering 316 (2017) 1157--1178.

\bibitem{nguyen2014nitsche}
V.~P. Nguyen, P.~Kerfriden, M.~Brino, S.~P. Bordas, E.~Bonisoli, Nitsche’s
  method for two and three dimensional {NURBS} patch coupling, Computational
  Mechanics 53~(6) (2014) 1163--1182.

\bibitem{dorfler1996convergent}
W.~D{\"o}rfler, A convergent adaptive algorithm for poisson’s equation, SIAM
  Journal on Numerical Analysis 33~(3) (1996) 1106--1124.

\bibitem{piegl2012nurbs}
L.~Piegl, W.~Tiller, The NURBS book, Springer Science \& Business Media, 2012.

\bibitem{liew1993transverse}
K.~Liew, Y.~Xiang, S.~Kitipornchai, Transverse vibration of thick rectangular
  plates--{I}. comprehensive sets of boundary conditions, Computers \&
  Structures 49~(1) (1993) 1--29.

\bibitem{deng2006dimensions}
J.~Deng, F.~Chen, Y.~Feng, Dimensions of spline spaces over {T}-meshes, Journal
  of Computational and Applied Mathematics 194~(2) (2006) 267--283.

\end{thebibliography}
\end{document}